%% file: main.tex
\documentclass[a4paper, 10pt, english, reqno, dvipsnames, oneside]{amsart}
\usepackage[utf8]{inputenc}

\usepackage{amsmath}
\usepackage{amsthm}
\usepackage{amssymb}
\usepackage{mathtools}  % provide more math symbols 
\usepackage{mathrsfs}   % provide \mathscr
\usepackage{amsfonts}
\usepackage{mathdots}   % provide antidiagonal dots
\usepackage[bb=libus]{mathalpha} %provide lowercase \mathbb
\usepackage{amsbsy}
\usepackage{tikz}
\usepackage{tikz-cd}
\usepackage{quiver}

\usetikzlibrary{shapes.geometric, arrows}
\usetikzlibrary{decorations.pathmorphing,shapes}
\usepackage{quiver}

\tikzstyle{startstop} = [rectangle, rounded corners, minimum width=3cm, minimum height=1cm, text centered, draw=black, fill=brown!10!pink!20, align=center]
\tikzstyle{arrow} = [->, decorate]

\usepackage{geometry}

\usepackage{tabularx}
\usepackage{makecell} %https://tex.stackexchange.com/questions/2441/how-to-add-a-forced-line-break-inside-a-table-cell

\usepackage{lmodern}

\theoremstyle{plain}
\newtheorem{theorem}{Theorem}[section]
\newtheorem{assertion}[theorem]{Assertion}
\newtheorem{proposition}[theorem]{Proposition}

\newtheorem{corollary}[theorem]{Corollary}
\newtheorem{conjecture}[theorem]{Conjecture}
\newtheorem{definition-theorem}[theorem]{Definition-Theorem}
\newtheorem{definition-proposition}[theorem]{Definition-Proposition}
\newtheorem{maintheorem}{Theorem}

\theoremstyle{definition}

\newtheorem{definition}[theorem]{Definition}

\newtheorem{assumption}[theorem]{Assumption}

\theoremstyle{remark}
\newtheorem{remark}[theorem]{Remark}

\usepackage{enumitem}   % adapt the enumerate environment
\usepackage[figurewithin=none]{caption}
\usepackage{graphicx}
\usepackage{float}
\usepackage{xstring}
\usepackage{array}
\usepackage{longtable}

\usepackage{xcolor}

\usepackage{mycom}
\usepackage{myarrows}

\makeatletter
\def\l@subsection{\@tocline{2}{0pt}{2pc}{5pc}{}}
\makeatother

\usepackage[style=alphabetic,sorting=nyt, backend=biber, backref=true]{biblatex}%<- specify style
\usepackage{url}

\usepackage[unicode, colorlinks, linktocpage=true]{hyperref}
\hypersetup{
	linkcolor = blue,
	citecolor = Red,
	urlcolor = PineGreen,
	bookmarksnumbered = true
}
\makeatletter
\def\namedlabel#1#2{\begingroup
    #2%
    \def\@currentlabel{#2}%
    \phantomsection\label{#1}\endgroup
}
\makeatother

\usepackage{orcidlink}

\title[Iwasawa theory for unitary groups]{Iwasawa Theory for Unitary Groups via Eisenstein Congruences and its Applications}
\author{Ruichen Xu \,\orcidlink{0009-0003-4555-2116}}
\address{Academy of Mathematics and Systems Science, Chinese Academy of Sciences, No. 55 Zhongguancun East Road, Beijing, 100190, China.}
\curraddr{Faculty of Mathematics,
University of Vienna, Oskar-Morgenstern-Platz 1, 1090 Vienna, Austria.}
\email{ruichen.xu@univie.ac.at}

\date{\today}

\begin{document}

\begin{abstract}
Let $p$ be an odd prime number. In this article, we study automorphic aspects of Iwasawa theory for unitary groups, focusing on the $p$-primitivity of Klingen Eisenstein series and their $p$-adic families on $\mathrm{U}(m+1,n+1)$, which is a key component in establishing lower bounds for Selmer groups in Iwasawa theory for $\mathrm{U}(m,n)$.

First, we compute the Gan--Gross--Prasad period integral of a Klingen Eisenstein series on $\mathrm{U}(m+1,n+1)$ against a cuspidal automorphic form on $\mathrm{U}(m+1,n)$, and relate it to special Rankin--Selberg $L$-values. We then interpolate these period integrals $p$-adically as both forms vary in Hida families, thereby producing a Rankin--Selberg type $p$-adic $L$-function for $\mathrm{U}(m,n) \times \mathrm{U}(m+1,n)$ as a byproduct. This method is designed to replace the computation of Fourier--Jacobi coefficients used in the existing literature to establish the $p$-primitive property of Klingen Eisenstein series, which so far treats only cases where $m+n \leq 2$.

As an application, we consider the case $\mathrm{U}(m,n)=\mathrm{U}(1,0)$ and obtain a new proof of M.-L. Hsieh’s Iwasawa main conjecture for CM elliptic curves over a totally real field $\mathcal{F}$ with CM field $\mathcal{M}$ and good ordinary reduction at all places above~$p$, for all primes $p \ge 5$ unramified in $\mathcal{F}$. Our approach removes the assumption in Hsieh's original work that $p$ does not divide the relative class number of $\mathcal{MF}/\mathcal{F}$.
\end{abstract}

    \maketitle

    \tableofcontents
    
%    \newpage
%    \pagestyle{fancy}

\newpage

\input{01-Introduction}

\newpage

\part{Preliminaries and background} \label{part:pre}

\input{AP01-PreliminariesAutoreps}
\input{AP02-PreliminariesGeometry}

\input{AP03-MainConjecture}

\newpage 

\part{The Gan--Gross--Prasad-type period integral of Klingen Eisenstein families} \label{part:automorphic}

\input{AU01-Reductions}
\input{AU02-Archimedian}
\input{AU03-Split}
\input{AU04-Summary}
\input{AU05-padicGGP}

\newpage

\part{Iwasawa main conjecture for CM fields and CM elliptic curves} \label{part:application}
\input{AR01-MoreHida}
\input{AR02-Primitivity}
\input{AR03-Lattice}
\input{TH1-Theta}
\input{AZ01-IMCU1}

\newpage 

%\bibliographystyle{alpha}
%\bibliography{mybib}
\printbibliography%<-print bibliography

\end{document}

%% file: 01-Introduction.tex
\section{Introduction} \label{sec:introduction}

\subsection{General background} \label{sec:introbackground}
In number theory, a central topic is the mysterious connection between arithmetic and analytic objects. Classical examples include the class number formula, which relates the class number of a number field to its Dedekind zeta function, and the Birch and Swinnerton-Dyer conjecture, which relates the Mordell--Weil rank of an elliptic curve over a number field to its Hasse--Weil $L$-function. Iwasawa theory, initiated by the Japanese mathematician Kenkichi Iwasawa in the 1950s, studies such connections in $p$-adic families; these connections are formulated as various versions of \emph{Iwasawa main conjectures}. 

More precisely, let $F$ be a number field, and let $\Gal_{F}$ be its absolute Galois group. Let $V$ be a finite-dimensional $E$-vector space, where $E$ is a finite extension of $\QQ_p$. Consider a geometric Galois representation
\[
\rho: \Gal_{F} \rightarrow \Aut_{E}(V)
\]
which satisfies the \emph{Panchishkin condition} formulated by R. Greenberg in \cite{greenberg1991}. Let $F_{\infty}/F$ be a $\ZZ_p$-extension, and let $\Lambda_{F} = \calO\lrbracket{\Gal(F_{\infty}/F)}$ be its Iwasawa algebra, where $\calO$ is the ring of integers of a sufficiently large extension of $\QQ_p$. There are two invariants attached to $V$: the \emph{Greenberg-type Selmer group} $\Sel^{\Gr}(V)$, which is expected to be a cotorsion $\Lambda_{F}$-module, and a conjectural \emph{$p$-adic $L$-function} $\calL_{p}(V) \in \Lambda_{F}$ predicted by \cite{coates1989on, coates1989on2}, which interpolates certain classical special $L$-values. The \emph{Greenberg-type Iwasawa main conjecture} predicts the equality
\[
\Char_{\Lambda_{F}}(\Sel^{\Gr}(V)^{\vee}) = (\calL_{p}(V)),
\]
where $(-)^{\vee}$ denotes the Pontryagin dual $\Hom_{\cts}(-, \QQ_p/\ZZ_p)$. To prove such an equality of principal ideals in $\Lambda_{F}$, the two divisibilities are typically studied using different techniques. In this article, our focus is on the divisibility providing the \emph{lower bounds of the Selmer group}, namely the inclusion
\[
\Char_{\Lambda_{F}}(\Sel^{\Gr}(V)^{\vee}) \subseteq (\calL_{p}(V)).
\]
Such inclusions are established via the machinery of \emph{Eisenstein congruences}, which dates back to the work of J.-P. Serre and K. Ribet (see, for example, \cite{ribet1976modular}).

\subsubsection{Baby example: The Ramanujan congruence}
We consider the following famous modular form
\[
    \Delta(q = \exp(2\pi \rmi z)) := q \prod_{n=1}^{\infty} (1-q^n)^{24} = \sum_{n=1}^{\infty} \tau(n) q^n \in S_{12}(\SL_2(\ZZ)),
\]
namely the \emph{modular discriminant}, which is a weight $12$ classical cuspidal eigenform of level $\SL_2(\ZZ)$. Let $p = 691$. It was first observed by Ramanujan that there is a congruence relation
\begin{equation} \label{eq:ramanujan}
     \Delta \equiv E_{12} \pmod{p=691},   
\end{equation}
where
    \[
    E_{12}(q = \exp(2 \pi \rmi z)) = - \dfrac{B_{12}}{24} + \sum_{n \geq 1} \left(\sum_{d \mid n} d^{11} \right) q^n 
    \]
is the Eisenstein series of weight $12$. Given this congruence relation, one translates it into a modulo-$p$ congruence of their Galois representations:
\[
\rho_{\Delta} \equiv \rho_{E_{12}} \equiv \begin{bmatrix}
\omega^{11} & \ast \\ 0 & 1
\end{bmatrix} \pmod {p = 691},
\]
where $\omega$ is the mod-$p$ Teichmüller character. It turns out that the Galois cohomology class $\ast \in \rmH^{1}(\Gal_{\QQ}, \FF_{p}(\omega^{11}))$ corresponds to a $p$-torsion ideal class in the ideal class group of $\QQ(\zeta_{691})$. The construction of such an ideal class can be viewed as providing a lower bound for the class group.

The story above is not yet complete: cautious readers will notice that it is essential to verify that this ideal class $\ast$ is nontrivial, which is equivalent to showing that the congruence relation \eqref{eq:ramanujan} is nontrivial. This is not difficult in this baby case: one simply considers the $p$-th Fourier coefficient of $E_{12}$, which is $1+p^{11}$, and observes that this is nonzero modulo $p$. This property is often called the \emph{$p$-primitivity} (or $p$-primitive property) of $E_{12}$.

Inspired by the baby example of the Ramanujan congruence, the idea of constructing mod-$p$ congruences between cuspidal eigenforms and Eisenstein series—namely \emph{Eisenstein congruences}—was first used in Iwasawa theory by B. Mazur and A. Wiles \cite{mazur1984class, wiles1990iwasawa} to prove the Iwasawa main conjecture for totally real number fields in the $\GL_1$ case (i.e., for class groups).

\subsubsection{General scheme of Eisenstein congruences}
Generalizing the Ramanujan congruence, the Eisenstein congruence machinery for studying the Iwasawa theory of certain automorphic forms on a reductive group $G$ consists of three steps:
\begin{enumerate}
    \item Find a larger group $\tilG$ containing $G$ as a Levi subgroup, and construct a $p$-adic family of Eisenstein series $\bfE_{V}$ over $\tilG$ such that, by the computations of Langlands and Shahidi, the constant term of $\bfE_{V}$ involves the $p$-adic $L$-function $\calL_{p}(V)$.
    \item Prove that the Eisenstein family $\bfE_{V}$ is \emph{$p$-primitive}, i.e., it is coprime to the $p$-adic $L$-function $\calL_{p}(V)$. This can be understood as a modulo-$p$ nonvanishing property of $\bfE_{V}$. Combined with (1), this implies that modulo $\calL_{p}(V)$, this Eisenstein family is congruent to a certain family of cuspidal automorphic forms $\bff$ on $\tilG$. The issue of $p$-primitivity is \emph{the most technical and challenging step} in this machinery.
    \item Translate the above congruence of automorphic forms into a congruence of the corresponding Galois representations. The Galois representation of $\bfE_{V}$ should be reducible, containing the Galois representation of $\bff$ as a constituent that is "less reducible." This should, vaguely speaking, produce "sufficiently many" cohomology classes in the Galois cohomology group $\Sel^{\Gr}(V)$ as predicted by the Iwasawa main conjecture. This step is called the \emph{lattice construction} and was axiomatized in \cite{bellaiche2009families}.
\end{enumerate}

Let $f \in S_{k}(\Gamma_{0}(N))$ be a cuspidal eigenform of weight $k \geq 2$, and let $\chi$ be an algebraic Hecke character of an imaginary quadratic field $K$ with an attached CM modular form $\theta_{\chi}$ of weight $k_{\chi}$. Assume that $p$ splits in $K$, which implies that $\theta_{\chi}$ is $p$-ordinary. In recent years, to study the arithmetic of $f$, researchers have studied the Iwasawa theory of the Galois representation $V_{f,\chi} := (V_{f} \otimes V_{\theta_{\chi}})(i)$ for $i \in \ZZ$ such that it satisfies Deligne's criticality, where $V_{f}$ (resp. $V_{\theta_{\chi}}$) is the Galois representation attached to $f$ (resp. $\theta_{\chi}$). Suppose $2 \nmid k+k_{\chi}$; then $V_{f, \chi}$ satisfies the Panchishkin condition in the following two "complementary cases":
\begin{enumerate}
    \item $k > k_{\chi}$ and $f$ is $p$-ordinary, i.e., $a_p(f)$ is a $p$-adic unit;
    \item $k < k_{\chi}$, allowing $f$ to have an arbitrary reduction type at $p$.
\end{enumerate}
The former is considered in \cite{skinner2014iwasawa}, with the Eisenstein series $\bfE_{V_{f,\chi}}$ chosen as a Klingen Eisenstein family $\bfE^{\Kling}_{f,\chi}$ over the unitary group $\GU(2,2)$, while the latter is considered in \cite{wan2020iwasawa, castella2022iwasawa} with the choice of a Klingen Eisenstein family $\bfE^{\Kling}_{f,\chi}$ over the unitary group $\GU(3,1)$. In both cases, it is crucial to establish the $p$-primitive property of $\bfE^{\Kling}_{f,\chi}$. In previous works, this has been done by computing either their Fourier coefficients or their Fourier--Jacobi coefficients. However, \emph{it remains unclear whether these computations are feasible for Klingen Eisenstein families $\bfE^{\Kling}_{f,\chi}$ over unitary groups $\rmU(m+1, n+1)$ in general.}

\subsection{Main results}

In this article, we initiate the investigation of the $p$-primitive property of Klingen Eisenstein series over unitary groups of general signatures $\rmU(m+1, n+1)$ (as constructed explicitly in \cite{wan2015families}) by computing a Gan--Gross--Prasad type period integral of them against cuspidal automorphic forms over $\rmU(m+1, n)$. We now summarize the main results of this article, beginning with the necessary notation.

Let $\calK/\calF$ be a CM extension; that is, $\calF$ is a totally real number field and $\calK$ is an imaginary quadratic extension of $\calF$. Let $p$ be an odd prime that is unramified in $\calF$ and such that each prime of $\calF$ above $p$ splits in $\calK$. Throughout, we fix an isomorphism $\iota_{p}: \barQQ_p \simeq \CC$. Let $\calK_{\infty}/\calK$ be the composition of all $\ZZ_p$-extensions of $\calK$, with Galois group $\Gamma_{\calK}$. Let $m \geq n \geq 0$ be integers. We denote by $\rmU(m,n)$ the unitary group of signature $(m,n)$ over $\calF$ associated to the skew-Hermitian matrix
\[
\begin{bmatrix}
 & & \bfone_{n} \\ & \vartheta & \\ -\bfone_{n}
\end{bmatrix},
\]
where $\vartheta$ is an $(m-n)\times (m-n)$ diagonal matrix with coefficients in $\calK$ such that $\sfi^{-1} \vartheta$ is totally positive definite, with $\sfi \in \calK$ a fixed totally imaginary element such that under the action of complex conjugation $\sfc \in \Gal(\calK/\calF)$, we have $\sfc(\sfi)= - \sfi$. We assume $\rmU(m,n)(\calF_v)$ is quasi-split for every finite place $v$ of $\calF$ (see assumption \eqref{ass:QS}).

\subsubsection{Main result 1: automorphic computations}
Let $\sigma$ be an irreducible tempered regular algebraic unitary cuspidal representation of $H := \rmU(m,n)$ such that $\sigma_{\infty}$ is a holomorphic discrete series with constant scalar weights (Assumption \ref{ass:scalarwt}) and $\sigma$ is sufficiently ramified at places above $p$ (Assumption \ref{ass:generic}). We further suppose that:
\begin{enumerate}
    \item[(irred)] the residual Galois representation of $\sigma$ is absolutely irreducible. \label{eq:irredsigma}
\end{enumerate}
The ultimate goal is to study the Iwasawa theory of the automorphic representation $\sigma$\footnote{Assumptions \ref{ass:scalarwt} and \ref{ass:generic} in this article are technical and are expected to be removed in the future.}.

Let $\chi$ be a Hecke character over $\calK$. Let $\Phi$ be a cuspidal automorphic form in $\sigma$. We can construct a Klingen Eisenstein series by pulling back a Siegel Eisenstein series $E^{\Sieg}_{\chi, s}(-)$ from the larger unitary group $G^{\bdsuit} := \rmU(m+n+1, m+n+1)$ along the cuspidal automorphic form $\Phi$ to $G^{\hsuit} := \rmU(m+1, n+1)$. We denote the Klingen Eisenstein series over $G^{\hsuit}$ by
\[
E^{\Kling}_{\sigma, \chi, s}(-) := E^{\Kling}(F^{\hsuit}(f^{\Sieg}_{s, \chi}, \Phi; -); -),
\]
where $f^{\Sieg}_{s, \chi}(-)$ is the Siegel Eisenstein section defining the Siegel Eisenstein series. The operator $F^{\hsuit}(-)$ is the so-called \emph{pullback integral}.

Let $\pi$ be an irreducible tempered regular algebraic unitary cuspidal representation of $G := \rmU(m+1, n)$ such that its residual Galois representation is absolutely irreducible. One of the main themes of this article is to compute and $p$-adically interpolate the following \emph{Gan--Gross--Prasad (GGP) type period integral} of $E^{\Kling}(-)$ over $G^{\hsuit}$ with a cuspidal automorphic form $\Psi$ in $\pi$:
\begin{equation} \label{eq:introggp}
\calP^{\Kling}(\Phi, \Psi, \chi, s) := \int_{G(\calF) \bs G(\AA_{\calF})} E^{\Kling}(F^{\hsuit}(f^{\Sieg}_{s, \chi}, \Phi; -); \jmath^{\flat}(g)) \Psi(g) \dif g,    
\end{equation}
where $\jmath^{\flat}: G \hookrightarrow G^{\hsuit}$ is a canonical embedding of unitary groups. By the cuspidality of $\Psi$, this integral converges absolutely for those values of $s$ at which the Klingen Eisenstein series is defined. 

The underlying philosophy of considering such period integrals is that to prove the modulo-$p$ nonvanishing property of $E^{\Kling}$, it suffices to choose an appropriate $\Psi$ such that \eqref{eq:introggp} is nonvanishing modulo $p$. It is widely believed that such period integrals are closely related to special values of $L$-functions, for which we have more tools to study such properties.

The first main result of this article is to turn this philosophy into genuine theorems. We unfold this period integral into a product of local integrals and show that local $L$-factors appear at the unramified places of $\calF$. For simplicity, we write $\calP^{\Kling}(\Phi, \Psi, \chi, s)$ as $\calP_{\Phi, \Psi}^{\Kling}$, and let:
\begin{enumerate}
    \item $\scV_{\calF}^{\ur}$ be the set of finite places $v$ of $\calF$ such that $\sigma_v$ and $\pi_v$ are both unramified, excluding places above $p$;
    \item $\scU_{\calF}^{\ur}$ be the set of finite places $v$ of $\calF$ such that $\sigma_v$, $\pi_v$, $\chi_v$, and the local field extension $\calK_v/\calF_v$ are all unramified, excluding places above $p$.
\end{enumerate}
We denote by $\bullet^{\bad}$ the complement of $\bullet^{\ur}$ in the set of finite places of $\calF$, for $\bullet \in \{\scV, \scU\}$. We then state our first main theorem.

\begin{maintheorem}[{Theorem \ref{thm:unrret}}] \label{mainthm:A}
We choose local Siegel Eisenstein sections $f^{\Sieg}_{v,s,\chi}$ at places $v$ of $\calF$ as in Section \ref{sec:localdoubl} for the construction of the Klingen Eisenstein series.\footnote{We remark that the choice of such local Siegel Eisenstein sections is quite canonical in applications to Iwasawa theory; they are used in \cite{skinner2014iwasawa, wan2020iwasawa, castella2022iwasawa}. In particular, we choose $f^{\Sieg}_{v,s,\chi}$ to be spherical Siegel Eisenstein sections at places in $\scU_{\calF}^{\ur}$, big-cell sections at ramified places, and certain well-chosen Siegel--Weil sections at places of $\calF$ above $p$. We closely follow \cite{wan2015families} when choosing these local Siegel Eisenstein sections.} Then
\begin{multline*}
\dfrac{\left(\calP^{\Kling}_{\Phi, \Psi}\right)^2}{\aabs{\Phi}_{\sigma, \MVW}^2 \aabs{\Psi}_{\pi, \MVW}^2} = \dfrac{1}{2^{\varkappa_{\sigma}+\varkappa_{\pi}}} \cdot \scL_{\scV_{\calF}^{\ur}}(\sigma \times \pi) L_{\scU_{\calF}^{\ur}}\left(s+\dfrac{1}{2}, \pi^{\vee}, \barchi^{\sfc} \right)^{2} \\
\times \prod_{v \in \scU_{\calF}^{\ur}} d_{N+1,v}(s, \barchi_v)^{-2} \prod_{v \in \mathscr{V}_{\calF}^{\bad}} \scI_v (\sigma_v, \pi_v) \prod_{v \in \scU_{\calF}^{\bad}} \scZ^{\dsuit, \rightt}_v(f^{\Sieg}_{v,s,\chi}, \pi_v)^2,
\end{multline*}
where the local $L$-factors are: for any set $\scS$ of places of $\calF$,
\begin{itemize}
    \item local Rankin--Selberg $L$-factors
    \[
    \scL_{\scS}(\sigma \times \pi) := \prod_{v \in \scS} \left( \dfrac{L(\frac{1}{2}, \sigma_v \times \pi_v)}{L(1, \sigma_v, \Ad) L(1, \pi_v, \Ad)} \prod_{i=1}^{m+n+1} L(i, \epsilon_{\calK_v/\calF_v}^{i}) \right),
    \]
    where $L(s, \sigma_v \times \pi_v)$ is the local $L$-factor of the Rankin--Selberg product $\BC(\sigma) \times \BC(\pi)$, with $\BC(-)$ denoting the base change to an automorphic representation of the general linear groups over $\calK$;
    \item local doubling $L$-factors
    \[
    L_{\scS}\left(s, \pi, \chi \right) := \prod_{v \in \scS} L_v(s, \BC(\pi_v) \otimes \chi_v \circ \det),
    \]
    where the right-hand side is the standard local Godement--Jacquet $L$-factor and $\BC(\pi_v)$ is the local base change from $G(\calF_v)$ to $\GL_{m+n+1}(\calK_v)$.
\end{itemize}
Here, the $\scI_v$ are local Ichino--Ikeda integrals (see \eqref{eq:definitioniiintegral}) and the $\scZ_v^{\dsuit, \rightt}$ are local doubling integrals (see \eqref{eq:defnpartialrightt}). Other unspecified factors will be defined in the text.
\end{maintheorem}

The next task is to compute the local Ichino--Ikeda integrals at bad places explicitly. We require $\sigma$ and $\pi$ to be regularly ordinary at places of $\calF$ above $p$ (see Definition \ref{defn:ord}). We impose further technical assumptions to simplify the calculation:
\begin{itemize}
    \item At archimedean places, we require $\sigma_{\infty}$ to be of constant scalar weight, with the weights of $\sigma_{\infty}$ and $\pi_{\infty}$ satisfying the ``Gan--Gross--Prasad weight interlacing property''.
    \item The local representations $\sigma_v$ and $\pi_v$ are ``disjointly ramified'' (i.e., they cannot both be ramified at any place), and the places where either $\sigma_v$ or $\pi_v$ ramifies must split in $\calK$.
    \item For $p$-adic places $v$ of $\calF$, the local representations $(\sigma_v, \chi_v)$ satisfy the generic assumption in \cite[Definition 4.42]{wan2015families}.
\end{itemize}
These are recorded more explicitly as Assumptions \ref{ass:wtinterlacing}, \ref{ass:spl}, \ref{ass:ram}, \ref{ass:scalarwt}, and \ref{ass:generic} in the text.

\begin{maintheorem}[{Theorem \ref{thm:automorphicfinal}}] \label{mainthm:B}
Under the aforementioned assumptions, we have
\begin{multline*}
\dfrac{\left(\calP^{\Kling}_{\Phi, \Psi}\right)^2}{\aabs{\Phi}_{\sigma, \MVW}^2 \aabs{\Psi}_{\pi, \MVW}^2} = \dfrac{1}{2^{\varkappa_{\sigma}+\varkappa_{\pi}}} \cdot \scL^{\scV_{\calF}^{(p)}}(\sigma \times \pi) L_{\scU_{\calF}^{\ur}}\left(s+\dfrac{1}{2}, \pi^{\vee}, \barchi^{\sfc} \right)^{2} \\
\times \prod_{v \in \mathscr{V}_{\calF}^{\ram}} \zeta_v C_{\sigma_v, \bpsi_v} C_{\sigma_v^{\vee}, \bpsi_v^{-1}} \calB_{\pi_v}^{\ess} \calB_{\sigma_v}^{\ess} \prod_{v \in \mathscr{V}_{\calF}^{(p)}} \zeta_v \Delta_\beta^2 \cdot \frG(\breve{\ulmu}, \breve{\ullam}) \frG(\breve{\check{\ulmu}}, \breve{\check{\ullam}}) \calC_{\sigma}^{\ord} \calC_{\pi}^{\ord}
\\
\times \prod_{v \in \scU_{\calF}^{\ur}}  d_{N+1,v}(s, \barchi^{\sfc}_v)^{-2} \prod_{v \in \scU_{\calF}^{\bad}} \scZ^{\dsuit}_v(f^{\Sieg}_{v,s,\chi}, \pi_v)^2.
\end{multline*}
Here $\scL^{\scV_{\calF}^{(p)}}(\sigma \times \pi)$ is the product of all corresponding local factors except for those in $\scV_{\calF}^{(p)}$, and the second line arises from the local Ichino--Ikeda integrals. These are defined explicitly in the text; we remark that the unspecified local factors involve only explicit local Gauss sums, local critical adjoint $L$-values, and the local inner products of canonically chosen vectors in local representations.
\end{maintheorem}

\subsubsection{Main result 2: $p$-adic interpolations}
The next step is to $p$-adically interpolate the Gan--Gross--Prasad type period integrals defined above, as the Klingen Eisenstein series $E^{\Kling}$ and the cuspidal automorphic form $\Psi$ vary in $p$-adic families.

Let $L/\QQ_p$ be a sufficiently large finite extension, with $\calO_L$ being its ring of integers, and let $\Lambda_{m,n}^{\circ} := \calO_{L}\lrbracket{T_{m+n}(1+p\ZZ_p)}$ be the weight algebra for Hida families over the unitary group $\GU(m,n)$, where $T_{m+n}$ is the diagonal torus in the general linear group $\GL_{m+n}$. Let $\II$ (resp. $\JJ$) be a normal domain over $\Lambda_{m,n}^{\circ}$ (resp. $\Lambda_{m+1,n}^{\circ}$), which is also a finite algebra over $\Lambda_{m,n}^{\circ}$ (resp. $\Lambda_{m+1,n}^{\circ}$). Consider $\bff$ (resp. $\bfg$), an $\II$-adic (resp. $\JJ$-adic) Hida family of tempered cuspidal ordinary eigenforms on $\GU(m,n)$ (resp. $\GU(m+1,n)$) of branch character $\eta_0$, containing $\sigma$ (resp. $\pi$). Let $\chi_0: \calK^{\times} \bs \AA_{\calK}^{\times} \rightarrow \CC^{\times}$ be a Hecke character. Let $\scU$ be a finite set of ``bad'' places of $\calF$ containing those lying above $p$ and $\infty$. We collect these inputs as a tuple
\[
\bfG := (L, \II, \JJ, \eta_0, \bff, \bfg, \chi_0, \scU),
\]
called a $p$-adic Gan--Gross--Prasad datum (see Definition \ref{defn:ggpdatum}). Given such a datum $\bfG$, we have a corresponding weight space\footnote{Here $(-)^{\ur}$ is an ``unramified extension of coefficient rings''; see Definition \ref{defn:padiceisendatum} for details.}
\[
\calX_{\bfG} := \Spec \II^{\ur}\lrbracket{\Gamma_{\calK}}(\barQQ_p) \times \Spec \JJ^{\ur}(\barQQ_p),
\]
and we $p$-adically interpolate the Gan--Gross--Prasad type period integral on a Zariski dense subset $\calX_{\bfG}^{\gen} \subseteq \calX_{\bfG}$, consisting of ``generic admissible points'' (see Definition \ref{defn:ggpweightpoint}). In particular, specializations $\sfP \times \sfQ \in \calX_{\bfG}^{\gen}$ are defined such that the automorphic representations $\sigma_{\bff_\sfP}$ and $\pi_{\bfg_{\sfQ}}$ satisfy the ``Gan--Gross--Prasad weight interlacing property''. The second main result of this article is then as follows.

\begin{maintheorem}[{Theorem \ref{thm:interpolationfinal}}] \label{mainthm:C}
With notation as above, and assuming that the residual Galois representations attached to $\bff$ and $\bfg$ are absolutely irreducible\footnote{The precise formulation of this condition is given in condition \ref{item:irredbff} in the text.}, there exists an element 
\[
\bfP_{\bfG} \in (\II^{\ur} \lrbracket{\Gamma_{\calK}} \otimes \Frac(\II^{\ur})) \otimes \Frac(\JJ^{\ur})
\]
such that for any $\sfP \times \sfQ \in \calX_{\bfG}^{\gen}$, we have
\[
(\sfP \times \sfQ)(\bfP_{\bfG}) = B(s_{\sfP}, \chi_{\sfP}) \calP^{\Kling}(\bff_{\sfP}, \bfg_{\sfQ}, \chi_{\sfP}, s_{\sfP}),
\]
where the right-hand side is the GGP type period integral of a Klingen Eisenstein series against a cuspidal automorphic form, normalized by an explicit factor $B(s_{\sfP}, \chi_{\sfP})$ defined in \eqref{eq:normalizationB}\footnote{This explicit factor $B(s_{\sfP}, \chi_{\sfP})$ is used in the renormalization of Siegel Eisenstein series (and thus Klingen Eisenstein series) to ensure algebraicity, thereby allowing for $p$-adic interpolation in families.}.
\end{maintheorem}

The main technique involves applying appropriate \emph{Hecke projectors}, introduced by Wan in \cite{wan2019iwasawatheorymathrmursblochkato}, to the family of Klingen Eisenstein series. When $n=0$, an alternative approach is to use the $p$-adic interpolation of Petersson inner products in Hida families. The method of Hecke projectors dates back to the works of Hida \cite{hida1985padic, hida1988padic, hida1991padic}, where he constructed $p$-adic $L$-functions for the Rankin--Selberg product of two modular forms and their associated Hida families.

In addition to the assumptions in Theorem \ref{mainthm:C}, if we further assume that the Hida families $\bff$ and $\bfg$ satisfy the aforementioned ramification conditions at classical points, then we obtain a more explicit interpolation formula relating $\bfP_{\bfG}^{2}$ to special $L$-values at points $\sfP \times \sfQ \in \calX_{\bfG}^{\gen}$; this is stated as Theorem \ref{thm:interpoformula} in the text.

In \cite{wan2015families}, Wan constructed $\scU$-imprimitive $p$-adic $L$-functions $\calL_{p}^{\scU}(\bfg, \chi_0)$ for Hida families $\bfg$ over unitary groups. By comparing these interpolation formulas with Theorem \ref{mainthm:B} (specifically, Theorem \ref{thm:interpoformula}), a byproduct of our computation is the construction of a primitive $p$-adic $L$-function for the Rankin--Selberg product of Hida families $\bff$ and $\bfg$ over $\rmU(m,n) \times \rmU(m+1,n)$.

\begin{maintheorem}[{Theorem \ref{thm:newpadicl}}] \label{mainthm:D}
With notation as above and under the aforementioned ramification conditions, assuming that the $p$-adic $L$-function $\calL_{p}^{\scU}(\bfg,\chi_0)$ is nonzero, there exists an element 
\[
\calL_{p}(\bff \times \bfg) \in \Frac(\II^{\ur}) \otimes \Frac(\JJ^{\ur})
\]
such that for any $\sfP \times \sfQ \in \calX_{\bfG}^{\gen}$, the value $(\sfP \times \sfQ)(\calL_{p}(\bff \times \bfg))$ equals
\begin{multline*}
\dfrac{1}{2^{\varkappa_{\sigma_{\bff_{\sfP}}}+\varkappa_{\pi_{\bfg_{\sfQ}}}}} \times \aabs{\bff_{\sfP}}_{\sigma_{\bff_{\sfP}}, \MVW}^2 \aabs{\bfg_{\sfQ}}_{\pi_{\bfg_{\sfQ}}, \MVW}^2 \cdot \prod_{v \in \scS_{\bff} \cup \scS_{\bfg}} \zeta_v C_{\sigma_{\bff_{\sfP},v}, \bpsi_v} C_{\sigma_{\bff_{\sfP}, v}^{\vee}, \bpsi_v^{-1}} \calB_{\pi_{\bfg_{\sfQ},v}}^{\ess} \calB_{\sigma_{\bff_{\sfP},v}}^{\ess} \\
\times \prod_{v \in \mathscr{V}_{\calF}^{(p)}} \zeta_v \Delta_{v,\beta_v}^2 \cdot \frG(\breve{\ulmu_{v,\sfP}}, \breve{\ullam_{v,\sfQ}}) \frG(\breve{\check{\ulmu_{v,\sfP}}}, \breve{\check{\ullam_{v,\sfQ}}}) \calC_{\sigma_{\bff_{\sfP}}}^{\ord} \calC_{\pi_{\bfg_{\sfQ}}}^{\ord} \cdot \scL^{\scV_{\calF}^{(p)}}(\sigma_{\bff_{\sfP}} \times \pi_{\bfg_{\sfQ}}),
\end{multline*}
where $\scS_{\bff}$ and $\scS_{\bfg}$ are the sets of places of $\calF$ dividing the ``tame level'' of the Hida families $\bff$ and $\bfg$, respectively. 
\end{maintheorem}

We note that there have been similar recent works on such $p$-adic $L$-functions in various generalizations, for instance \cite{liu2023anticyclotomicpadiclfunctionsrankinselberg, hsieh2023fivevariablepadiclfunctionsu3times, dimitrakopoulou2024anticyclotomicpadiclfunctionscoleman}.

\subsubsection{Main result 3: application to the Iwasawa main conjecture for CM fields}

As introduced in Section \ref{sec:introbackground}, the motivation for computing the Gan--Gross--Prasad type period integral of the Klingen Eisenstein series $E^{\Kling}$ is to prove its $p$-primitive property, which can be interpreted as a certain modulo-$p$ nonvanishing property. The significance of Theorems \ref{mainthm:A}, \ref{mainthm:B}, and \ref{mainthm:C} is, qualitatively speaking, the following transition:
\begin{figure}[H]
    \centering
\begin{tikzpicture}[node distance=2.2cm, decoration=snake]
\node (original) [startstop] {The modulo-$p$ nonvanishing property \\ of Klingen Eisenstein series $E^{\Kling}$ on $\rmU(m+1, n+1)$};
\node (GGP) [startstop, below of=original] {The modulo-$p$ nonvanishing property of the period integral $\calP^{\Kling}_{\Phi, \Psi}$ \\ for a \emph{delicately chosen} cuspidal automorphic form $\Psi$ on $\rmU(m+1, n)$};
\node (Lfunction) [startstop, below of=GGP] {The \emph{joint} modulo-$p$ nonvanishing property of \\ $L$-functions $L\left(\dfrac{1}{2}, \sigma \times \pi\right)$ and $L\left(s+\dfrac{1}{2}, \pi \times \chi \right)$ for a \\ \emph{delicately chosen} cuspidal automorphic representation $\pi$ on $\rmU(m+1, n)$};
\draw [arrow] (original) -- node[anchor=east] {reduces to} (GGP);
\draw [arrow] (GGP) -- node[anchor=east] {reduces to} (Lfunction);
\end{tikzpicture}
    \caption{The reduction process}
    \label{fig:reduction}
\end{figure}
Through this two-step reduction, we translate the problem into the modulo-$p$ nonvanishing of such special $L$-values; the inclusion of the auxiliary cuspidal automorphic representation $\pi$ on $\rmU(m+1, n)$ provides additional flexibility for the proof.

However, finding suitable candidates for such a representation $\pi$ remains challenging in general. The explicit construction of automorphic representations is a difficult problem, particularly for unitary groups of high rank. We seek a construction of $\pi$ that simultaneously satisfies the following requirements:
\begin{enumerate}
    \item The cuspidal automorphic representation $\pi$ is a \emph{transfer} of a certain cuspidal automorphic representation $\pi^{\flat}$ from a simpler group, such as a lower-rank unitary group.
    \item The local $L$-factors of $\pi^{\flat}$ (at least at the unramified finite places) capture those of $\pi$, ensuring that the $L$-functions in the final stage of Figure \ref{fig:reduction} do not become overly complicated.
\end{enumerate}
It turns out that the \emph{theta correspondence} provides a manageable way to ``transfer'' between groups as required in (1), while providing satisfying results for (2).

In this article, we implement this philosophy in the simplest non-trivial case: the Iwasawa main conjecture for cuspidal automorphic representations $\sigma$ of $\rmU(1,0)$. Equivalently, this is the Iwasawa main conjecture for Hecke characters $\chi$ of CM fields $\calK/\calF$. This problem was studied systematically in Hsieh's work \cite{hsieh2014eisenstein}, which utilizes Eisenstein congruences with Klingen Eisenstein series $E^{\Kling}_{1,\chi,s}$ on $\rmU(2,1)$ and imposes certain technical assumptions to ensure that $E^{\Kling}_{1,\chi,s}$ is $p$-primitive. These assumptions were necessary to apply the modulo-$p$ nonvanishing results for Fourier--Jacobi coefficients computed in \textit{loc.~cit.}

With the new machinery of Gan--Gross--Prasad type period integrals for $E^{\Kling}_{1,\chi,s}$, we are able to bypass the use of Fourier--Jacobi coefficients and thereby remove some of these technical assumptions. As mentioned previously, this requires the introduction of an auxiliary Hida family $\bfg$ of cuspidal representations on $\rmU(2,0)$. For this, we rely on the recent work of Y.--S.~Lee \cite{lee2025hidafamilythetalift}, who constructs a Hida family via the theta lift from $\rmU(1,0)$ to $\rmU(2,0)$\footnote{A similar construction already appears in the work of Wan \cite{wan2020iwasawa}.}.

To introduce our result on the Iwasawa main conjecture for CM fields, we begin with the following notation. Let $\calK_{\infty}/\calK$ be the compositum of the cyclotomic $\ZZ_p$-extension $\calK_{\infty}^{\cyc}/\calK$ and the anticyclotomic $\ZZ_p^{[\calF:\QQ]}$-extension $\calK_{\infty}^{\ac}$. We put $\Gamma_{\calK}^{\bullet} := \Gal(\calK_{\infty}^{\bullet}/\calK)$ for $\bullet \in \{ \quad, \cyc, \ac \}$. Let $\calK^{\prime}$ be a finite extension of $\calK$ containing $\calK(\zeta_p)$ that is disjoint from $\calK_{\infty}$. Put $\Delta := \Gal(\calK^{\prime}/\calK)$ and $\calK_{\infty}^{\prime} = \calK_{\infty} \calK^{\prime}$. Let $\psi: \Delta \rightarrow \barZZ_{p}^{\times}$ be a finite-order Hecke character over $\calK$, which is necessarily a unitary Hecke character. This is called a \emph{branch character}. We put $\calR := \barZZ_{p}^{\ur}[\psi]$ and $\Lambda := \calR\lrbracket{\Gamma_{\calK}}$.

Consider $\calM_{\Sigma}$ to be the maximal $p$-abelian $\Sigma_p$-ramified extension of $\calK_{\infty}^{\prime}$. The fundamental object is 
\[
X_{\Sigma} := \Gal(\calM_{\Sigma}/\calK_{\infty}^{\prime}) \otimes_{\ZZ_p[\Delta]\lrbracket{\Gamma_{\calK}}} \calR[\Delta]\lrbracket{\Gamma_{\calK}}.
\]
Denote by $X(\psi, \Sigma)$ the maximal $\psi$-isotypic quotient of $X_{\Sigma}$, which is a finitely generated torsion $\Lambda$-module, and let $F(\psi, \Sigma) \in \Lambda$ be its characteristic power series.

On the analytic side, we have the associated Katz $p$-adic $L$-function $\calL_{\Katz}(\psi, \Sigma)$ as an element of $\Lambda$, as recalled in Section \ref{sec:twistbypanti}. We then have the following theorem as another main result of this article.

\begin{maintheorem}[{Theorem \ref{thm:myimccm}}, Iwasawa main conjecture for CM field $\calK/\calF$] \label{mainthm:E}
Suppose that $p \nmid 6 D_{\calF} \abs{\Delta}$. Then, as ideals of $\Lambda$, the following divisibility holds:
\[
(F(\psi, \Sigma)) \subseteq (\calL_{\Katz}(\psi, \Sigma)).
\]
\end{maintheorem}
Compared with \cite[Theorem 2]{hsieh2014eisenstein}, we have been able to remove the assumption that $p \nmid h_{\calK}^{-}$ (where $h_{\calK}^{-} = h_{\calK}/h_{\calF}$ is the relative class number of $\calK/\calF$) and certain unramified properties of $\psi$ (see \ref{item:hsi} and \ref{item:hsii} in the text). By applying Hida's result \cite{hida2009quadratic} or Rubin's result \cite{rubin1991main} for the reverse inclusion, we obtain the equality under further assumptions. These are presented as Theorem \ref{thm:imchida} and Theorem \ref{thm:imcrubin}.

Furthermore, we can apply the Iwasawa main conjecture for CM fields $\calK/\calF$ to study the Iwasawa main conjecture for elliptic curves over $\calF$ with complex multiplication. Let $E$ be such an elliptic curve with complex multiplication given by the ring of integers of an imaginary quadratic field $\calM$, and assume that $E$ has good ordinary reduction at all places of $\calF$ above $p$. This implies that $p$ splits in $\calM$. We write $p \calO_{\calM} = \frp \overline{\frp}$ and let $\calK = \calF \calM$. Let $\pzK/\calK$ be the unique $\ZZ_p^{2}$-extension of $\calK$ in $\calK(E[p^{\infty}])$ and let $\Lambda_{E} := \calO_{L}\lrbracket{\Gal(\pzK/\calK)}$, where $\calO_{L}$ is a sufficiently large finite extension of $\ZZ_p$. Denote $\Lambda_{E, \calO_{\frp}} = \Lambda_{E} \otimes_{\ZZ_p} \calO_{\frp}$, with $\calO_{\frp}$ the ring defined in \cite[Section 2 before Section 2.1]{hsieh2014eisenstein}. On the algebraic side, one defines the $\frp^{\infty}$-Selmer group $\Sel_{\frp^{\infty}}(E/\pzK)$, which is a cofinitely generated discrete $\Lambda_{E}$-module. On the analytic side, one considers the two-variable $p$-adic $L$-function $\calL_{p}(E) \in \Lambda_{E, \calO_{\frp}}$ attached to $(E, \frp)$ following Perrin-Riou \cite{perrinriou}. We shall prove the following main theorem.

\begin{maintheorem}[{Theorem \ref{thm:imccmtwo}}] \label{mainthm:F}
Suppose that $p \nmid 6 D_{\calF}$ and that $E$ has good ordinary reduction at all places above $p$. Then we have the inclusion of ideals in $\Lambda_{E, \calO_{\frp}}$:
\[
\Char_{\Lambda_{E}}(\Sel_{\frp^{\infty}}(E/\pzK)^{\vee}) \subseteq (\calL_{p}(E)).
\]
Suppose further that 
\begin{enumerate}[label = \rm (Ru)]
    \item $p$ does not divide $[\calF:\QQ]$ and $\calK(E[\frp])/\calM$ is abelian; \label{item:ruintro}
\end{enumerate}
then we have the equality $\Char_{\Lambda_{E}}(\Sel_{\frp^{\infty}}(E/\pzK)^{\vee}) = (\calL_{p}(E))$ of ideals in $\Lambda_{E, \calO_{\frp}}$.
\end{maintheorem}

We project the Iwasawa main conjectures obtained in Theorem \ref{mainthm:F} to the cyclotomic direction $\calK^{\cyc}/\calK$ to obtain results on cyclotomic Iwasawa main conjectures for $E$ over $\calF$, presented as Theorem \ref{thm:imccmcyc}. As a standard application, we obtain the rank-zero $p$-converse theorem and the rank-zero refined $p$-part Birch--Swinnerton-Dyer formula, as follows, proved via Mazur's control theorem and Greenberg's refined control theorem (see, for example, \cite[Theorem 1.1, Theorem 4.1]{greenberg1999iwasawa}).

\begin{maintheorem}[{Theorem \ref{thm:pconverse}, Theorem \ref{thm:pbsd}}] \label{mainthm:G}
Suppose that $p \nmid 6 D_{\calF}$ and that $E$ has good ordinary reduction at all places above $p$. Then we have the following implication:
\[
\corank_{\ZZ_p} \Sel_{p^{\infty}}(E/\calF) = 0 \Rightarrow L(E/\calF, 1) \neq 0.
\]
Suppose further that \ref{item:ruintro} holds, then $L(E/\calF, 1) \neq 0$ implies that $\corank_{\ZZ_p} \Sel_{p^{\infty}}(E/\calF) = 0$ and
\[
v_p \left( \dfrac{L(E/\calF, 1)}{\Omega_{E}} \right) = v_p \left( \dfrac{\abs{\Sha(E/\calF)[p^{\infty}]}}{\abs{E(\calF)_{\tor}}^2} \prod_{v \nmid \infty} \Tam_{v}(E/\calF) \right).
\]
Here $\Tam_{v}(E/\calF)$ is the Tamagawa number at $v$ for a finite place $v$ of $\calF$, $\Sha(E/\calF)[p^{\infty}]$ is the $p$-primary part of the Tate--Shafarevich group of $E$ over $\calF$, and $E(\calF)_{\tor}$ is the finite group of torsion points of $E$ over $\calF$.
\end{maintheorem}

\subsection{Future works}
This article represents only the initial step toward a comprehensive Iwasawa theory for unitary groups, particularly regarding the machinery of Eisenstein congruences. A great deal of work remains to be done.

\subsubsection{Technical assumptions in the automorphic computation}
Firstly, several assumptions in this article are technical in nature and are expected to be removed or weakened in future research.
\begin{itemize}
    \item For the computation of local Ichino--Ikeda integrals, Assumptions \ref{ass:ram} and \ref{ass:spl} are expected to be removed.
    \item While it is believed that Klingen Eisenstein families over unitary groups can be constructed in full generality, such a construction is currently missing from the literature. The Klingen Eisenstein family utilized here is the one from \cite{wan2015families}, which necessitates the imposition of Assumptions \ref{ass:scalarwt} and \ref{ass:generic}.
\end{itemize}

\subsubsection{Application to higher-rank unitary groups}
Motivated by the success of our method in the lowest-rank case $\rmU(1,0)$, it is natural to pursue the Iwasawa theory for unitary groups of rank two, such as $\rmU(1,1)$, $\rmU(2,0)$, or, more generally, rank-two unitary groups whose signatures vary among the infinite places of $\calF$ (the so-called \emph{mixed-signature} case\footnote{
Let $V$ be a two-dimensional skew-Hermitian space over $\calK$ such that:
\begin{itemize}
\item for a distinguished infinite place $\infty_{0}:\calF\hookrightarrow\mathbb{R}$, the space $V_{\infty_{0}}$ has signature $(2,0)$, and
\item for the remaining infinite places $\infty_{i}$ ($i=1,\ldots,d-1$), the spaces $V_{\infty_{i}}$ have signature $(1,1)$.
\end{itemize}
The mixed-signature Iwasawa main conjecture for $\rmU(V)$ has applications to the \emph{rank-one} $p$-converse theorem and to refined $p$-parts of the Birch--Swinnerton--Dyer formula for elliptic curves over totally real number fields $\calF$, as shown in the unpublished thesis of H.~Li \cite{li2025bsd}. However, this conjecture has not yet been established in the mixed-signature setting. This mixed-signature nature is not apparent when $\calF = \QQ$, as in \cite{jetchev2017birch, skinner2020converse}.}).

More ambitiously, one may hope to treat general unitary groups $\rmU(m,n)$ with $m+n>2$. In these higher-rank cases, lower bounds for Selmer groups in the Iwasawa main conjecture remain entirely inaccessible. Progress in any of these directions is hindered by two major missing ingredients:

\begin{enumerate}[label=(\arabic*)]
\item \emph{Explicit $p$-adic families on higher-rank unitary groups.}
Even in the rank-two setting, applications require the choice of an auxiliary automorphic representation $\pi$ on a unitary group of rank three. Explicit $p$-adic families of such automorphic representations are currently largely unknown.
\item \emph{Modulo-$p$ nonvanishing of special $L$-values in higher rank.}  
Nonvanishing modulo $p$ for special values of cuspidal representations and their Rankin--Selberg products is widely regarded as a problem more delicate than the associated Iwasawa main conjectures, and essentially nothing is known in genuinely higher-rank situations.
\end{enumerate}

Although item (2) is completely out of reach at present, our results for $\rmU(1,0)$ offer partial guidance toward item (1): namely, the use of theta correspondences as a source of $p$-adic families. While the automorphic theory of the theta correspondence is well-developed, its $p$-adic aspects remain largely open. For applications to Iwasawa theory, we aim to address questions such as the following for reductive dual pairs $(\rmU(V), \rmU(W))$ and potentially other dual pairs:
\begin{enumerate}
    \item Given an algebraic (resp.\ $p$-integral, $p$-ordinary) cuspidal automorphic form $\Psi^{\flat}$ on $\rmU(V)$, when is its theta lift $\Theta(\Psi^{\flat})$ on $\rmU(W)$ algebraic (resp.\ $p$-integral, $p$-ordinary)?
    \item Given a Hida family $\bff^{\flat}$ on $\rmU(V)$ passing through $\Psi^{\flat}$, how and when can we construct a Hida (or $\Lambda$-adic) family $\bTheta(\bff^{\flat})$ interpolating the classical theta lifts $\Theta(\Psi^{\flat})$ of $\Psi^{\flat}$?
    \item How can we relate the $p$-adic $L$-function of $\bff^{\flat}$ to that of $\bTheta(\bff^{\flat})$?
\end{enumerate}

These questions have been studied sporadically in special low-rank cases, for example:
\begin{itemize}
    \item \cite{prasanna2008arithmetic} for $(\GL_2, \GO(B))$, where $B$ is a quaternion algebra;
    \item \cite{zhang2025non} for $(\rmO_{2n+1}, \Mp_{4n})$;
    \item \cite{lee2025hidafamilythetalift} for $(\GU(1,0), \GU(2,0))$;
    \item \cite{iudica2024lambdaadickudlalift} for $(\GU(1,1), \GU(2,1))$; and
    \item ongoing work of M.-L.~Hsieh and Z.~Liu on Yoshida lifts $(\mathrm{GSO}(B), \GSp(4))$.
\end{itemize}
In the near future, we hope to address the applications to rank-two unitary groups.

\subsection{Organization}
This article consists of four parts.

The first part (Part \ref{part:pre}) comprises Sections \ref{sec:preliminaries} through \ref{sec:mainconjectures}, which provide the preliminaries for this article.
\begin{itemize}
    \item In Section \ref{sec:preliminaries}, we introduce the background on unitary groups and their embeddings, followed by a discussion of Eisenstein series over unitary groups.
    \item In Section \ref{sec:formsunitary}, we recall the geometric background of modular forms, $p$-adic modular forms, and their $p$-adic families over unitary groups.
    \item In Section \ref{sec:mainconjectures}, we introduce the Greenberg-type Iwasawa main conjectures for unitary groups.
\end{itemize}

The second part (Part \ref{part:automorphic}) consists of Sections \ref{sec:ggpperiodintegral} through \ref{sec:padicggp}. In this part, we prove Theorems \ref{mainthm:A}, \ref{mainthm:B}, \ref{mainthm:C}, and \ref{mainthm:D}. This is one of the main parts of the article.
\begin{itemize}
    \item In Section \ref{sec:ggpperiodintegral}, we define the Gan--Gross--Prasad type period integral of Klingen Eisenstein series against cuspidal automorphic forms, unfold them, and relate them to local $L$-factors. A more detailed sketch of our approach is provided at the beginning of Section \ref{sec:ggpperiodintegral}. In this section, we begin by recalling the doubling methods, following both the approach of Piatetski-Shapiro and Rallis and that of Garrett. While these facts are well-known, we offer proofs for certain results where standard references are lacking.
    \item From Section \ref{sec:iiarchi} to Section \ref{sec:localdoubl}, we explicitly calculate the local Ichino--Ikeda integrals at bad places. The local doubling integrals utilized here were previously computed in \cite{wan2015families}.
    \item In Section \ref{sec:padicggp}, we first recall the construction of Siegel Eisenstein families and Klingen Eisenstein families over unitary groups. We then $p$-adically interpolate the GGP type period integral, concluding the proofs of Theorems \ref{mainthm:C} and \ref{mainthm:D}.
\end{itemize}

The third part (Part \ref{part:application}) first focuses on other inputs of Eisenstein congruences: the non-cuspidal Hida theory for unitary groups, reviewed in Section \ref{sec:noncuspidalhida} and applied in Section \ref{sec:genprimitivity}, and the lattice construction, reviewed in Section \ref{sec:lattice}. After these preparations, in Sections \ref{sec:pretheta} and \ref{sec:imcforcm} we prove Theorems \ref{mainthm:E}, \ref{mainthm:F}, and \ref{mainthm:G}.
\begin{itemize}
\item In Section \ref{sec:pretheta}, we introduce background material on local and global theta correspondences, with particular emphasis on the almost-equal-rank case.
\item In Section \ref{sec:imcforcm}, we apply our automorphic computations to the lowest-rank case $\rmU(1,0)$ and obtain results on the Iwasawa main conjectures for CM fields, including Theorem \ref{mainthm:E}. We then deduce applications to the Iwasawa main conjecture for CM elliptic curves over $\calF$ and to the rank-zero case of the Birch--Swinnerton-Dyer conjecture, including Theorems \ref{mainthm:F} and \ref{mainthm:G}.
\end{itemize}

\subsection{Notation and conventions}
Next, we introduce the basic notation and conventions used throughout this article.

\subsubsection{Fields and CM-types}
\begin{itemize}
    \item Let $\QQ, \ZZ$ (resp.\ $\QQ_p, \ZZ_p$) denote the field of (resp.\ $p$-adic) rational numbers and the ring of (resp.\ $p$-adic) integers, respectively.
    \item For a finite extension $L/\QQ_p$, we denote by $\calO_L$ its ring of integers and by $\calO_L^{\ur}$ the completion of the maximal unramified extension of $\calO_L$.
    \item Let $\calM$ be a number field; that is, a finite extension of $\QQ$. We denote by $\AA_{\calM}, \AA_{\calM, \rmf}, \AA_{\calM, \rmf}^{p}$ the ring of adèles of $\calM$, the ring of finite adèles, and the ring of finite adèles with $p$-coordinates removed, respectively. 
\end{itemize}

We fix once and for all an embedding $\iota_{\infty}: \barQQ \hookrightarrow \CC$ and an isomorphism $\iota: \CC \xrightarrow{\sim} \CC_p$, where $\CC_p$ is the completion of an algebraic closure of $\QQ_p$. Let $\iota_p := \iota \circ \iota_{\infty}: \barQQ \hookrightarrow \CC_p$ be their composition. Let $\barZZ_p$ be the $p$-adic completion of the ring $\barZZ$ of algebraic integers in $\CC_p$, let $\frm_{p}$ be the maximal ideal of $\barZZ_p$, and let $\frm := \iota^{-1}(\frm_{p})$. Every number field $\calM$ is regarded as a subfield of $\CC$ (resp.\ $\CC_{p}$) via $\iota_{\infty}$ (resp.\ $\iota_p$), and we identify $\Hom(L, \barQQ) = \Hom(L, \CC_p)$.

Throughout this article, $\calK/\calF$ denotes a CM extension; that is, $\calF$ is a totally real number field of degree $d$ over $\QQ$, and $\calK$ is an imaginary quadratic extension of $\calF$. We denote by $\KK$ (resp.\ $\AA$) the ring of adèles of $\calK$ (resp.\ $\calF$). We fix an odd rational prime $p$ that is unramified in $\calF$ and assume that every prime of $\calF$ above $p$ splits in $\calK$.
\begin{itemize}
    \item Let $\sfc \in \Gal(\calK/\calF)$ be the unique nontrivial automorphism of $\calK$, called the \emph{complex conjugation}.
    \item Let $D_{\calF}$ be the discriminant of $\calF$.
    \item Let $\calK_{\infty}/\calK$ be the composition of all $\ZZ_p$-extensions of $\calK$, with Galois group $\Gamma_{\calK}$.
    \item We fix an element $\sfi \in \calK$ such that $\sfc(\sfi) = -\sfi$.
\end{itemize}
Let $\rmI_{\calK} := \Hom(\calK, \barQQ)$. For a place $w$ of $\calK$ above $p$, we set
\[
\rmI_{w} := \{\sigma_p \in \Hom(\calK, \CC_p): \sigma_p \text{ induces } w \} = \Hom(\calK_{w}, \CC_p).
\]
We identify $\rmI_{\calK}$ with $\sqcup_{w} \rmI_{w}$ via $\sigma \mapsto \sigma_p:= \iota_p \circ \sigma$. Let $\Sigma$ be a \emph{CM-type} of $\calK$, which is a subset of $\rmI_{\calK}$ such that
\[
\Sigma \cap \Sigma \sfc = \emptyset \quad \text{and} \quad \Sigma \sqcup \Sigma \sfc = \rmI_{\calK}.
\]
Denote by $\Sigma_{p}$ the set of places of $\calK$ above $p$ induced by $\iota_p \circ \sigma$ for $\sigma \in \Sigma$. We assume that $\Sigma$ is a \emph{$p$-ordinary CM-type}; i.e.,
\[
\Sigma_p \cap \Sigma_p \sfc = \emptyset \quad \text{and} \quad \Sigma_p \sqcup \Sigma_p \sfc = \scV_{\calK}^{(p)},
\]
where $\scV_{\calK}^{(p)}$ is the set of places of $\calK$ lying above $p$. The existence of such a type is guaranteed by the assumption that every prime of $\calF$ above $p$ splits in $\calK$. Henceforth, we also denote $\Sigma^{\sfc} := \Sigma \sfc$ and $\Sigma_{p}^{\sfc} := \Sigma_{p} \sfc$, and identify the CM-type $\Sigma$ (resp.\ the $p$-adic CM-type $\Sigma_p$) with the set of archimedean places of $\calF$ (resp.\ the set of places of $\calF$ above $p$). 

\subsubsection{Places}
For a number field $\calM$, we denote the set of places by $\scV_{\calM}$, $\scS_{\calM}$, $\scT_{\calM}$, or $\scU_{\calM}$. We denote by:
\begin{itemize}
    \item $\scV_{\calF}^{\square}$ (resp.\ $\scV_{\calK}^{\square}$) the set of places of $\calF$ (resp.\ $\calK$) above a finite set $\square$ of places of $\QQ$ (resp.\ $\calF$);\footnote{When $\square = \{w\}$ is a singleton, we write $\scV_{\calF}^{(w)} = \scV_{\calF}^{\{w\}}$.}
    \item $\scV_{\calF}^{\infty}$ (resp.\ $\scV_{\calF}^{0}$) the set of archimedean (resp.\ nonarchimedean) places of $\calF$;
    \item $\scV_{\calF}^{\spl}$ (resp.\ $\scV_{\calF}^{\ns}$) the subset of $\scV_{\calF}^{0}$ consisting of those that are split (resp.\ nonsplit) in $\calK$.
\end{itemize}

We introduce three partitions of $\scV_{\calF}$. The first is used for the computation of local Ichino--Ikeda integrals:
\begin{itemize}
    \item $\scV_{\calF}^{\ur}$ is the set of finite places of $\calF$ where both $\sigma_v$ and $\pi_v$ are unramified, excluding places above $p$;
    \item $\scV_{\calF}^{\bad}$ is the complement of $\scV_{\calF}^{\ur}$ in $\scV_{\calF}$;
    \item $\scV_{\calF}^{\ram}$ is the set of finite places of $\calF$ where either $\sigma_{v}$ or $\pi_v$ is ramified, excluding places above $p$.
\end{itemize}
The second is for the doubling integrals and the $p$-adic $L$-function of $\pi$:
\begin{itemize}
    \item $\scT_{\calF}^{\ur}$ is the set of finite places where $\pi_v$, $\chi_v$, and $\calK_v/\calF_v$ are unramified;
    \item $\scT_{\calF}^{\bad}$ is the complement of $\scT_{\calF}^{\ur}$ in $\scT_{\calF}$;
    \item $\scT_{\calF}^{\ram}$ is the complement of $\scV_{\calF}^{\{p, \infty\}}$ in $\scT_{\calF}^{\bad}$.
\end{itemize}
The third is for the normalization of the Klingen Eisenstein series $E^{\Kling}_{\sigma, \chi, s}$:
\begin{itemize}
    \item $\scS_{\calF}^{\ur}$ is the set of finite places where $\sigma_v$, $\chi_v$, and $\calK_v/\calF_v$ are unramified;
    \item $\scS_{\calF}^{\bad}$ is the complement of $\scS_{\calF}^{\ur}$ in $\scS_{\calF}$;
    \item $\scS_{\calF}^{\ram}$ is the complement of $\scV_{\calF}^{\{p, \infty\}}$ in $\scS_{\calF}^{\bad}$.
\end{itemize}
We denote $\scU_{\calF}^{\ur} := \scS_{\calF}^{\ur} \cap \scT_{\calF}^{\ur}$, and let $\scU_{\calF}^{\bad}$ be the complement of $\scU_{\calF}^{\ur}$ in $\scV_{\calF}$. Then
\[
\scU_{\calF}^{\bad} = (\scS_{\calF}^{\ur} \smallsetminus \scT_{\calF}^{\ur}) \cup (\scT_{\calF}^{\ur} \smallsetminus \scS_{\calF}^{\ur}) \cup (\scS_{\calF}^{\ram} \cap \scT_{\calF}^{\ram}) \cup \scU_{\calF}^{(p)} \cup \scU_{\calF}^{(\infty)}.
\]
Under Assumption \ref{ass:ram}, $\scS_{\calF}^{\ram} \cap \scT_{\calF}^{\ram} = \emptyset$. Note that $\bullet^{\ur}$ contains all but finitely many places of $\calF$ for $\bullet \in \{\scV, \scS, \scT, \scU\}$. These partitions are illustrated in Figures \ref{fig:VF-decomposition} and \ref{fig:VF-two-axes-shifted2}.

\begin{figure}[H]
    \centering
\begin{tikzpicture}[>=stealth, thick]
\draw[dashed] (-1,0) -- (11,0);
\fill (-0.5,0) circle (4pt);
\node[below, align=center] at (-0.5,-0.3) {archimedean\\places};
\node[above] at (-0.5, 0.3) {$\scV_{\calF}^{\infty}$}; 
\fill (10.5,0) circle (4pt);
\node[below, align=center] at (10.5,-0.3) {places \\ above $p$};
\node[above] at (10.5,0.3) {$\scV_{\calF}^{(p)}$};
\node[right] at (11.5,0) {$\scV_{\calF}$};
\draw[line width=2pt, red] (1,0) -- (5,0);
\node[above] at (3,0.2) {$\scV_{\calF}^{\ur}$};
\node[below, text width=4cm, align=center] at (3,-0.3) {both $\pi_v$ and $\sigma_v$\\are unramified};
\draw[thin] (5,0.25) -- (5,-0.25);
\draw[line width=2pt, blue] (5,0) -- (9,0);
\node[above] at (7,0.2) {$\scV_{\calF}^{\ram}$};
\node[below, text width=4cm, align=center] at (7,-0.3) {one of $\pi_v$ and $\sigma_v$\\is ramified};
\end{tikzpicture}
\caption{Decomposition of $\scV_{\calF}$}
\label{fig:VF-decomposition}
\end{figure}

\begin{figure}[H]
    \centering
\begin{tikzpicture}[>=stealth, thick]
\coordinate (I) at (1,1);
\coordinate (Vbottom) at (1,0.2);
\coordinate (Vtop) at (1,11);
\coordinate (Hleft)  at (0.2,1);
\coordinate (Hright) at (12,1);
\draw[dashed] (Hleft) -- (Hright);
\draw[dashed] (Vbottom) -- (Vtop);
\fill (I) circle (3.5pt);
\node[below, align=center] at ($(I)+(-1.32,-0.32)$) {archimedean\\places};
\fill (11,1) circle (3.5pt);
\node[below, align=center] at (11,0.68) {places\\above $p$};
\node[above] at (11,1.18) {$\scS_{\calF}^{(p)}$};
\node[right] at (12.2,1) {$\scS_{\calF}$};
\draw[line width=2pt, red] (2,1) -- (5,1);
\node[below] at (3.5,0.68) {$\scS_{\calF}^{\ur}$};
\draw[thin] (5,1.25) -- (5,0.75);
\draw[line width=2pt, blue] (5,1) -- (9,1);
\node[below] at (7,0.68) {$\scS_{\calF}^{\ram}$};
\draw[line width=2pt, red] (1,2) -- (1,5);
\node[left] at (0.68, 3.5) {$\scT_{\calF}^{\ur}$};
\draw[thin] (0.75,5) -- (1.25,5);
\draw[line width=2pt, blue] (1,5) -- (1,9);
\node[left] at (0.68, 7) {$\scT_{\calF}^{\ram}$};
\fill (1,10) circle (3.5pt);
\node[left, align=center] at (0.68, 10) {places\\above $p$};
\node[right] at (1.18, 10) {$\scT_{\calF}^{(p)}$};
\node[above] at (1, 11.2) {$\scT_{\calF}$};
\fill[red!25]   (2,2) rectangle (5,5);
\fill[purple!15](2,5) rectangle (5,9);
\fill[purple!15](5,2) rectangle (9,5);
\fill[blue!15]  (5,5) rectangle (9,9);
\draw[thin] (2,2) rectangle (5,5);
\draw[thin] (2,5) rectangle (5,9);
\draw[thin] (5,2) rectangle (9,5);
\draw[thin] (5,5) rectangle (9,9);
\node at (3.5,3.5) {$\scU_{\calF}^{\ur} := {\scS^{\ur}_{\calF}} \cap {\scT^{\ur}_{\calF}}$};
\node at (3.5,7.0) {${\scS^{\ur}_{\calF}} \smallsetminus {\scT^{\ur}_{\calF}}$};
\node at (7.0,3.5) {${\scT^{\ur}_{\calF}} \smallsetminus \scS^{\ur}_{\calF}$};
\node[align=center] at (7.0,7.0) {${\scS^{\ram}_{\calF}} \cap {\scT^{\ram}_{\calF}}$, \\ empty by \\ Assumption \ref{ass:ram}};
\end{tikzpicture}
\caption{Decomposition of $\scS_{\calF}$, $\scT_{\calF}$ and $\scU_{\calF}$}
\label{fig:VF-two-axes-shifted2}
\end{figure}

For an $L$-function with Euler product $L(\cdots) = \prod_{v \in \scV_{\calF}} L_v(\cdots)$ and a subset $\scV \subseteq \scV_{\calF}$, we write
\[
L_{\scV}(\cdots) := \prod_{v \in \scV} L_v(\cdots) \quad \text{and} \quad L^{\scV}(\cdots) := \prod_{v \in \scV_{\calF} \smallsetminus \scV} L_v(\cdots).
\]

\subsubsection{Characters and automorphic representations} \label{sec:notationcharacters}
By a Hecke character $\chi$ of a number field $\calM$, we automatically mean a continuous homomorphism $\chi: \calM^{\times} \bs \AA^{\times}_{\calM} \rightarrow \CC^{\times}$.

Let $\chi: \calK^{\times} \bs \AA_{\calK}^{\times} \rightarrow \CC^{\times}$ be an algebraic Hecke character of $\calK$. We say that $\chi$ has infinity type $\kappa = \sum_{\sigma \in \rmI_{\calK}} \kappa_{\sigma} \sigma \in \ZZ[\rmI_{\calK}]$ if
    \[
    \chi_{\infty}(\alpha) = \iota_{\infty}(\alpha^{\kappa}) := \prod_{\sigma \in \rmI_{\calK}} \sigma(\alpha)^{\kappa_{\sigma}}, \quad \text{for all } \alpha \in \calK^{\times}.
    \]
We associate to $\chi$ a unique $p$-adic Hecke character $\widehat{\chi}: \calK^{\times} \bs \AA_{\calK}^{\times} \rightarrow \barZZ_p^{\times}$, its \emph{$p$-adic avatar}, defined by
\[
\widehat{\chi}(z) := \iota_p \left( \chi(z) \prod_{\sigma \in \rmI_{\calK}} z_{\sigma}^{-\kappa_{\sigma}} \right) \prod_{w \mid p} \prod_{\sigma_p \in \rmI_{w}} \sigma_{p}(z_w)^{\kappa_{\sigma}}.
\]
Conversely, $\chi$ is called the \emph{complex avatar} of $\widehat{\chi}$. For any Hecke character $\chi$ of $\calK$, we write $\chi^{\calF}$ for its restriction to $\AA^{\times}$. We define $\chi^{\sfc}$ as the conjugation on the source, $x \mapsto \chi(x^{\sfc})$, and $\barchi$ as the conjugation on the target, i.e. $x \mapsto \overline{\chi(x)}$.

For any Hecke character $\chi$ of $\calK$, we define the \emph{universal character} attached to $\chi$ as the product $\bchi := \chi \cdot \bvep_{\calK}$, where $\bvep_{\calK}: \calK^{\times} \bs \AA_{\calK}^{\times} \rightarrow \Gamma_{\calK} \hookrightarrow \Lambda_{\calK}^{\times}$ is the tautological character induced by the reciprocity law of class field theory under geometric normalization. We denote by:
\begin{itemize}
    \item $\epsilon_{\calK/\calF}: \AA^{\times} \rightarrow \CC^{\times}$ the quadratic character attached to the extension $\calK/\calF$, which decomposes into local characters $\epsilon_{\calK_w/\calF_v}: \calF_v^{\times} \rightarrow \CC^{\times}$;
    \item $\epsilon_{\cyc}: \Gal_{\QQ} \rightarrow \ZZ_p^{\times}$ the cyclotomic character;
    \item $\omega$ the composition of the cyclotomic character $\Gal_{\QQ} \rightarrow \ZZ_p^{\times}$ with $(\ZZ_p/p\ZZ_p)^{\times} \rightarrow \ZZ_p^{\times}$ where the first arrow is reduction modulo $p$ and the second is the Teichmüller lift.
\end{itemize}

For a \emph{unitary automorphic representation} of $\rmU(m,n)$ over $\calF$, we realize it as an $L^2$-subspace of the space of automorphic forms on $\rmU(m,n)$, with the inner product given by the $L^2$ inner product. We denote by $\calA(\rmU(m,n)(\calF) \bs \rmU(m,n)(\AA))$ (resp., $\calA_{\cusp}(\rmU(m,n)(\calF) \bs \rmU(m,n)(\AA))$) the space of automorphic forms (resp., the space of cuspidal automorphic forms) on $\rmU(m,n)$. For an automorphic representation $\pi$ of $\rmU(m,n)$ over $\calF$, we denote by $\chi_{\pi}: \calK^{\times} \bs \AA_{\calK}^{\times} \rightarrow \CC^{\times}$ its central character. In fact, the central character $\chi_{\pi}$ of an automorphic representation $\pi$ should be a character on the center
\[
Z(\rmU(m,n))(\AA_{\calF}) = \{ \diag[a, \ldots, a] \in \rmU(m,n)(\AA_{\calK}): a a^{\sfc} = 1, a \in \AA_{\calK}^{\times} \}.
\]
Noting that $Z(\rmU(m,n))(\AA_{\calF})$ is isomorphic to the group $\AA_{\calK}^{\times, 1}$ consisting of elements in $\AA_{\calK}^{\times}$ of norm $1$, we identify it with a character on $\AA_{\calK}^{\times, 1}$. Throughout this paper, the central character $\chi_{\pi}$ always refers to a fixed extension of the central character $\chi_{\pi}$ from $\AA_{\calK}^{\times, 1}$ to $\AA_{\calK}^{\times}$ under the aforementioned identification.

\subsubsection{General linear groups and operation on matrices}
Let $m,n$ be positive integers. For any square matrix $A$ of size $n \times n$, we put (whenever the operations makes sense):
\begin{itemize}
    \item $A^{\rmt}$ as the transpose of $A$, and $A^{-\rmt}$ as the inverse of $A^{\rmt}$,
    \item $\barA$ as the conjugation of $A$, whenever the conjugation makes sense, and
    \item $A^{\star} := \barA^{\rmt}$ and $A^{-\star} := (A^{\star})^{-1}$.
    \item We write $\Herm_{n}$ for the set of matrices in $M_{n}$ such that $A^{\star}=A$.
\end{itemize}
We write $\bfone_{n}$ as the identity square matrix of size $n \times n$. We put $M_{m,n}$ to be the set of matrices with $m$ rows and $n$ columns. We denote by:
\begin{itemize}
    \item $\GL_n$ the general linear group of matrices of size $n \times n$,
    \item $T_n$ the diagonal torus of $\GL_n$,
    \item $B_n$ the upper triangular matrices in $\GL_n$ (that is, the Borel subgroup of $\GL_n$) and $B_n^{-}$ the lower triangular matrices in $\GL_n$,
    \item $U_n$ the strictly upper triangular matrices in $\GL_n$ (that is, upper triangular matrices with $1$'s on the diagonal). It is the unipotent subgroup of $B_n$ and $B_n = T_n U_n$ is the Levi decomposition of $B_n$. In the same way we have $B_n^{-} = T_n U_n^{-}$,
    \item More generally, let $n = n_1 + \cdots + n_r$ be a partition of $n$. We put $P[n_1, \ldots, n_r]$ as the parabolic subgroup of $\GL_n$ with respect to this partition, that is,
    $$
    P[n_1, \ldots, n_r] = \begin{bmatrix}
    \GL_{n_1} & \ast & \ast & \ast \\
              & \GL_{n_2} & \ast & \ast \\
              &    & \ddots & \vdots \\
              &            &         & \GL_{n_r}
    \end{bmatrix}.
    $$
    It has a Levi decomposition $P[n_1, \ldots, n_r] = M[n_1, \ldots, n_r] U[n_1, \ldots, n_r]$, where
    \[
    M[n_1, \ldots, n_r] = \diag[\GL_{n_1}, \GL_{n_2}, \ldots, \GL_{n_r}]
    \]
    is its Levi component and
    $$
    U[n_1, \ldots, n_r] = \begin{bmatrix}
    \bfone_{n_1} & \ast & \ast & \ast \\
              & \bfone_{n_2} & \ast & \ast \\
              &    & \ddots & \vdots \\
              &            &         & \bfone_{n_r}
    \end{bmatrix}
    $$
    is its unipotent radical.
\end{itemize}
Given a square matrix of size $n \times n$, we can group its entries into blocks by a partition $n = n_1 + \cdots + n_{r}$. We denote such a partition as $n = [n_1 \mid n_2 \mid \cdots \mid n_r]$.

\subsubsection{Measures and pairings} \label{sec:measures}
We fix a Haar measure $\dif g$ on the adèle group of a reductive group $G$ over $\calF$. We take $\dif g$ to be the Tamagawa measure for definiteness. We write $\dif g = \prod_{v \in \scV_{\calF}} \dif g_v$, with $\dif g_v$ a Haar measure on $G(\calF_v)$ under the following hypothesis (following \cite[Hypothesis 1.4.4]{eischen2020padic}):
\begin{itemize}
    \item At all finite places $v$ of $\calF$ at which the group $G$ is unramified, $\dif g_v$ is the measure that gives volume $1$ to a hyperspecial maximal compact subgroup.
    \item At all finite places $v$ at which $G(\calF_v)$ is isomorphic to $\prod_{i} \GL_{n_i}(\calF_{i, w_i})$, where $\calF_{i, w_i}$ is a finite extension of $\calF_v$ with integer ring $\calO_{i}$, the measure $\dif g_v$ is the measure that gives volume $1$ to the group $\prod_{i} \GL_{n_i}(\calO_{i})$.
    \item At all finite places $v$ of $\calF$, the values of $\dif g_v$ on open compact subgroups are rational numbers.
    \item At archimedean places $v$ of $\calF$, we choose measures such that $\prod_{v \in \scV_{\calF}^{\infty}} \dif g_v$ is the Tamagawa measure.
\end{itemize}
For the reductive group $G$ over $\calF$, we let $[G]$ denote the quotient $G(\calF) \bs G(\AA)$. Endowing $G(\AA)$ with the Tamagawa measure, $[G]$ is endowed with the quotient measure by the counting measure on $G(\calF)$. Given two automorphic forms $\Phi, \Phi^{\prime}$ over $G$, we define their Petersson inner product by
\[
\lrangle{\Phi, \Phi^{\prime}}_{\Pet} := \int_{[G]} \Phi(g) \Phi^{\prime}(g) \dif g.
\]
Let $\pi$ be an irreducible unitary cuspidal automorphic representation of $G$. We denote by $\pi^{\vee}$ the contragredient representation of $\pi$. Then the Petersson inner product is a canonically defined pairing $\pi \times \pi^{\vee} \rightarrow \CC$. We have factorizations
\[
\fac_{\pi}: \pi \xrightarrow{\sim} \otimes^{\prime} \pi_v, \quad \fac_{\pi^{\vee}}: \pi^{\vee} \xrightarrow{\sim} \otimes^{\prime} \pi_v^{\vee},
\]
where $\pi_v$ is an irreducible representation of $G(\calF_v)$. We have non-degenerate canonical $G(\calF_v)$-pairings 
\[
\lrangle{-,-}_v: \pi_v \times \pi_v^{\vee} \rightarrow \CC
\]
for all $v$. We renormalize the local pairings such that for $\fac_{\pi}(\Phi) = \otimes_v \Phi_v$ and $\fac_{\pi^{\vee}}(\Phi^{\vee}) = \otimes_v \Phi_v^{\vee}$, we have
\begin{equation} \label{eq:localglobalpeter}
\lrangle{\Phi, \Phi^{\vee}}_{\Pet} = \prod_{v \in \scV_{\calF}} \lrangle{\Phi_v, \Phi_v^{\vee}}_v.    
\end{equation}

\subsection*{Acknowledgments}
This article is the (unofficial and more revised) English version of the author’s doctoral thesis at the University of Chinese Academy of Sciences \footnote{Though this acknowledgment is not a strict translation of the Chinese version.}, written under the supervision of Professor Xin Wan. The author expresses his most sincere gratitude to Professor Wan for suggesting this project and for his patient and insightful guidance during his five-year doctoral study. 

He also thanks Ashay Burungale, Yangyu Fan, Olivier Fouquet, Zhibin Geng, Ming-Lun Hsieh, Francesco Iudica, Haijun Jia, Shilin Lai, Yu-Sheng Lee, Haidong Li, Wen-Wei Li, Yifeng Liu, Peikai Qi, Loren Spice, Ye Tian, Liang Xiao, Mulun Yin, Xiaoyu Zhang, and Luochen Zhao for helpful discussions and insightful comments.

In the process of the doctoral degree application, the author's institute invited an anonymous referee group of five experts, as well as a defense committee consisting of Professors Yongquan Hu, Wen-Wei Li, Xu Shen and Yichao Tian, chaired by Professor Liang Xiao. The author sincerely thanks these referees and experts for their suggestions and for the helpful discussions during the defense, which greatly improved the quality of this article.

%% file: AP01-PreliminariesAutoreps.tex
\section{Eisenstein series over unitary groups} \label{sec:preliminaries}

In this section, we introduce the background for our automorphic computations, including the definitions of unitary groups and the Eisenstein series defined over them.

\subsection{Unitary groups}
Let $\calK/\calF$ be a CM extension. Let $V$ be an $N$-dimensional vector space over $\calK$. Let $\phi: V \times V \rightarrow \calK$ be a non-degenerate skew-Hermitian form on $V$. We linearly extend $\phi$ to any $\calF$-algebra $R$ and $R$-module $V \otimes_{\calF} R$.

\begin{definition} \label{defn:unitarygrp}
The \emph{general unitary group} attached to $(V, \phi)$ is the algebraic group $\GU(V, \phi)$ over $\calF$, whose $R$-points, for each $\calF$-algebra $R$, are given by
\[
\GU(V, \phi)(R) := \{ g \in \Aut_{\calK \otimes_{\calF} R}(V \otimes_{\calF} R): \phi(gv, gw) = \nu(g) \phi(v,w), \, v, w \in V \otimes_{\calF} R \},
\]
where $\nu(g) \in R^{\times}$ is the similitude factor. The \emph{unitary group} attached to $(V, \phi)$ is the algebraic group $\rmU(V, \phi)$ over $\calF$ defined by the subgroup where $\nu(g) = 1$.
\end{definition}

Let $\ttB$ be an ordered $\calK$-basis of $V$. Then $\phi$ can be expressed as a matrix $[\phi]_{\ttB}$ under this basis, and the unitary group $\rmU(V, \phi)$ can be identified with the matrix group
\[
\rmU(V, \phi, \ttB)(R) := \{g \in \GL_{N}(\calK \otimes_{\calF} R): g [\phi]_{\ttB} g^{\star} = [\phi]_{\ttB} \}
\]
for any $\calF$-algebra $R$.

Let $n$ be the maximum dimension of totally $\phi$-isotropic subspaces of $V$ and write $m=N-n$. By \cite[Lemma 1.5, Lemma 1.6]{shimura1997euler}, we can find an ordered $\calK$-basis of $V$:
\begin{equation} \label{eq:Wittbasis}
    \mathtt{Witt}_{m,n}: \mathtt{y}^1, \ldots, \mathtt{y}^n, \mathtt{w}^1, \ldots, \mathtt{w}^{m-n}, \mathtt{x}^1, \ldots, \mathtt{x}^n,
\end{equation}
called a \emph{Witt basis}, under which $\phi$ has the matrix representation
\[
J_{m,n} := [\phi]_{\tWitt_{m,n}} := \begin{bmatrix}
& & \bfone_{n} \\ & \vartheta & \\ -\bfone_{n} & &
\end{bmatrix}.
\]
Throughout this article, we fix such a basis of $V$ and denote $\rmU(V, \phi, \tWitt_{m,n})$ as $\rmU(m,n)$ for simplicity, though it surely depends on $\vartheta$, not only on $m$ and $n$.

\subsubsection{Unitary groups at local places} \label{sec:unitarylocal}
Let $v$ be a place of $\calF$.
\begin{itemize}
    \item Suppose $v$ is an archimedean place corresponding to a real embedding $\imath_v: \calF \hookrightarrow \RR$. We require that $\vartheta$ is an $(m-n) \times (m-n)$ diagonal matrix such that 
    \begin{equation} \tag{sgn} \label{ass:sgn}
        \imath_v(\sfi^{-1} \vartheta) \text{ is positive definite for all } v \in \scV_{\calF}^{\infty}.
    \end{equation}
    This implies that $\rmU(m,n)(\calF_v)$ is isomorphic to the real Lie group
    \begin{equation}\label{eq:realumn}
    \rmU_{m,n} := \left\{g \in \GL_{m+n}(\CC): g^{\star} \begin{bmatrix} \bfone_m & \\ & -\bfone_n \end{bmatrix} g = \begin{bmatrix}
    \bfone_m & \\ & -\bfone_{n} \end{bmatrix} \right\}        
    \end{equation}
    for each archimedean place $v$ of $\calF$.
    \item Suppose $v$ is a finite place of $\calF$ that splits as $v = w \barw$ in $\calK$. Then the isomorphism $\calK \otimes \calF_v \cong \calF_v \times \calF_v$ induces $\varrho_{v,N}: \GL_{N}(\calK \otimes \calF_v) \cong \GL_{N}(\calF_v) \times \GL_{N}(\calF_v)$. If $g \mapsto (g_1,g_2)$, then $\barg \mapsto (g_2, g_1)$, and
    \[
    \rmU(m,n)(\calF_v) \simeq \{(g_1, g_2) \in \GL_{N}(\calF_v) \times \GL_{N}(\calF_v): g_2 = J_{m,n}^{\rmt} g_1^{-1} J_{m,n}^{-1}\}.
    \]
    This group is isomorphic to $\GL_N(\calF_v)$ via projection to the first factor \footnote{This involves a choice of place $w$ of $\calK$ above $v$.}. We denote this isomorphism by $\varrho_{w, m, n}$.
    \item Suppose $v$ is a finite nonsplit place. Naturally, the $\calK$-vector space extends to a $\calK_v$-vector space $V_v$ over $\calK_v$ and so does the skew-Hermitian form extended to $\phi_v$. By \cite[p.~76]{shimura1997euler}, there exists a local basis $\tWitt_{v}$ such that
    \[
    [\phi_v]_{\tWitt_v} = \begin{bmatrix}
    & & \bfone_{n_v} \\ & \vartheta_v & \\ -\bfone_{n_v}
    \end{bmatrix}, 
    \]
    where $\vartheta_v$ is an anisotropic matrix of size $r_v \times r_v$ with $r_v \in \{0,1,2\}$. We assume the quasi-split condition:
    \begin{equation} \tag{QS} \label{ass:QS}
        r_v \leq 1 \text{ for all finite nonsplit places } v \text{ of } \calF.
    \end{equation}
\end{itemize}

\subsubsection{Hermitian basis}
It is sometimes more convenient to use another basis. One checks immediately that $\sfi \phi$ is a Hermitian form on $V$, so there exists an ordered $\calK$-basis
$$
\mathtt{Herm}_{m,n}: \tta^{1}, \ldots, \tta^{m}, \ttb^{1}, \ldots, \ttb^{n}
$$
of $V$ such that $-\sfi \phi$ has a matrix representation $[-\sfi \phi]_{\mathtt{Herm}_{m,n}} := \begin{bmatrix}
 \bfone_{m} & 0 \\ 0 & -\bfone_{n} 
\end{bmatrix}$. Then under this basis, our original skew-Hermitian form can be written as $[\phi]_{\mathtt{Herm}_{m,n}} := \begin{bmatrix}
 \sfi \cdot \bfone_{m} & 0 \\ 0 & -\sfi \cdot \bfone_{n} 
\end{bmatrix}$. In fact, we can make the transition between the Hermitian basis and Witt basis explicit. We choose appropriate Witt basis such that there exists matrix $\vartheta_{0}$ satisfying $-\sfi \theta = \vartheta_0 \vartheta_0^{\star}$. Then it is readily verified that \footnote{In fact, one computes
\[
\begin{bmatrix}
\frac{1}{2} \cdot \bfone_{n} & & -\sfi \cdot \bfone_{n} \\
 & \vartheta_0^{-1} & \\
\frac{1}{2} \cdot \bfone_{n} & & \sfi \cdot \bfone_{n}
\end{bmatrix}
\begin{bmatrix}
& & \bfone_{n} \\ & \vartheta & \\ -\bfone_{n} & &
\end{bmatrix}
\begin{bmatrix}
\frac{1}{2} \cdot \bfone_{n} & & \frac{1}{2} \cdot \bfone_{n} \\
 & \vartheta_0^{-\star} & \\
\sfi \cdot \bfone_{n} & & - \sfi \cdot \bfone_{n}
\end{bmatrix}
= 
\begin{bmatrix}
\sfi \cdot \bfone_{n} & &  \\ & \sfi \cdot \bfone_{m-n} & \\ & & -\sfi \cdot \bfone_{n}
\end{bmatrix}.
\]
}
$$
[\tty^{1}, \ldots, \tty^{n}, \ttw^{1}, \ldots, \ttw^{m-n}, \ttx^{1}, \ldots, \ttx^{n}] =
[\tta^{1}, \ldots, \tta^{m}, \ttb^{1}, \ldots, \ttb^{n}]
\begin{bmatrix}
\frac{1}{2} \cdot \bfone_{n} & & -\sfi \cdot \bfone_{n} \\
 & \vartheta_0^{-1} & \\
\frac{1}{2} \cdot \bfone_{n} & & \sfi \cdot \bfone_{n}
\end{bmatrix}.
$$
Denote the matrix on the right-hand side by $\varrho_{m,n, \vartheta_0}$; it induces an isomorphism of unitary groups
$$
\varrho_{m,n, \vartheta_0}: \rmU(V, \phi, \tWitt_{m,n}) \xrightarrow{\sim} \rmU(V, \phi, \tHerm_{m,n}), \quad h \mapsto \varrho_{m,n, \vartheta_0} h \varrho_{m,n, \vartheta_0}^{-1}.
$$
When $m=n$, the matrix $\vartheta_0$ disappears and we omit it from the notation.

\subsubsection{The $\sharp$-space}
Starting from the vector space $V$ under the basis $\tHerm_{m,n}$, we consider the $\calK$-vector space 
$$
V^{\sharp} :=  \calK \tte \oplus \calK \tta^{1} \oplus \cdots \oplus \calK \tta^{m} \oplus \calK \ttb^{1} \oplus \cdots \oplus \calK \ttb^{n},
$$
with the ordered $\calK$-basis $\tHerm^{\sharp}$ as defined above. The space $V^{\sharp}$ is equipped with a skew-Hermitian form $\phi^{\sharp}$ such that $[-\sfi \phi^{\sharp}]_{\tHerm^{\sharp}} = \diag[1, \bfone_{m}, -\bfone_{n}]$. Accordingly, we have an embedding of unitary groups
$$
\jmath^{\sharp}: \rmU(V, \phi, \tHerm_{m,n}) \hookrightarrow \rmU(V^{\sharp}, \phi^{\sharp}, \tHerm^{\sharp}), \quad h \mapsto \diag[1, h].
$$
Transitioning back to the Witt basis $\tWitt^{\sharp}$ via the transition matrix $\varrho_{m+1, n, \vartheta_0^{\sharp}}$ where $\vartheta_0^{\sharp} := \diag[1, \vartheta_0]$, the corresponding embedding $\jmath^{\sharp}$ with respect to the Witt basis $\tWitt^{\sharp}$ is defined such that the diagram

$$
\begin{tikzcd}
	{\rmU(V, \phi, \tHerm_{m,n})} && {\rmU(V^{\sharp}, \phi^{\sharp}, \tHerm^{\sharp})} \\
	{\rmU(V, \phi, \tWitt_{m,n})} && {\rmU(V, \phi, \tWitt^{\sharp})}
	\arrow["{\jmath^{\sharp}}", hook, from=1-1, to=1-3]
	\arrow["{\varrho_{m,n,\vartheta_0}}"', from=1-1, to=2-1]
	\arrow["{\varrho_{m+1,n,\vartheta_0^{\sharp}}}", from=1-3, to=2-3]
	\arrow["{\jmath^{\sharp}}"', hook, from=2-1, to=2-3]
\end{tikzcd}
$$
commutes.

\subsubsection{The $\hsuit$-space}
Next, we consider adding another line to $V^{\sharp}$. Starting from the $\calK$-ordered basis $\tHerm^{\sharp}$, we consider the $\calK$-vector space
$$
V^{\hsuit} := \calK \tte \oplus \calK \tta^{1} \oplus \cdots \oplus \calK \tta^{m} \oplus \calK \ttb^{1} \oplus \cdots \oplus \calK \ttb^{n} \oplus \calK \ttf,
$$
with the ordered $\calK$-basis $\tHerm^{\hsuit}$ as defined above. The space $V^{\hsuit}$ is equipped with a skew-Hermitian form $\phi^{\hsuit}$ such that $[-\sfi \phi^{\hsuit}]_{\tHerm^{\hsuit}} = \diag[1, \bfone_{m}, -\bfone_{n}, -1]$.
Accordingly, we have an embedding of unitary groups
$$
\jmath^{\flat}: \rmU(V^{\sharp}, \phi^{\sharp}, \tHerm^{\sharp}) \hookrightarrow \rmU(V^{\hsuit}, \phi^{\hsuit}, \tHerm^{\hsuit}), \quad g \mapsto \diag[g, 1].
$$
Returning to the Witt basis via the transition matrix $\varrho_{m+1, n+1, \vartheta_0}$, the corresponding $\jmath^{\flat}$ with respect to the Witt basis $\tWitt^{\hsuit}$ is defined such that the diagram
$$
\begin{tikzcd}
	{\rmU(V^{\sharp}, \phi^{\sharp}, \tHerm^{\sharp})} && {\rmU(V^{\hsuit}, \phi^{\hsuit}, \tHerm^{\hsuit})} \\
	{\rmU(V^{\sharp}, \phi^{\sharp}, \tWitt^{\sharp})} && {\rmU(V^{\hsuit}, \phi^{\hsuit}, \tWitt^{\hsuit})}
	\arrow["{\jmath^{\flat}}", hook, from=1-1, to=1-3]
	\arrow["{\varrho_{m+1,n,\vartheta_0^{\sharp}}}"', from=1-1, to=2-1]
	\arrow["{\varrho_{m+1, n+1, \vartheta_0}}", from=1-3, to=2-3]
	\arrow["{\jmath^{\flat}}"', hook, from=2-1, to=2-3]
\end{tikzcd}
$$
commutes.

To conclude, we have successive embeddings of unitary groups
$$
\rmU(m,n) \xrightarrow{\jmath^{\sharp}} \rmU(m+1,n) \xrightarrow{\jmath^{\flat}} \rmU(m+1, n+1).
$$

\subsubsection{Doubling space}
Let $V^{\bdsuit}$ be the direct sum $V^{\hsuit} \oplus V$ with the ordered $\calK$-basis
$$
\mathtt{Herm}^{\bdsuit, \oplus}: \tte, \mathtt{a}^1, \ldots, \mathtt{a}^{m}, \mathtt{b}^1, \ldots, \mathtt{b}^{n}, \ttf, \bar{\tta}^{1}, \ldots, \bar{\tta}^{m}, \bar{\ttb}^{1}, \ldots, \bar{\ttb}^{n} 
$$
obtained by concatenating the Hermitian basis $\tHerm^{\hsuit}$ of $V^{\hsuit}$ and the Hermitian basis $\tHerm_{m,n}$ of $V$ (where bars are added to the basis elements of $V$ for distinction). We equip it with the skew-Hermitian form $\phi^{\bdsuit}$ such that
$$
[-\sfi \phi^{\bdsuit}]_{\tHerm^{\bdsuit, \oplus}} := \diag[[-\sfi \phi^{\hsuit}]_{\tHerm^{\hsuit}}, -[-\sfi \phi]_{\tHerm}] = \diag[1, \bfone_{m}, -\bfone_{n}, -1, -\bfone_{m}, \bfone_{n}].
$$
Accordingly, we have an embedding of unitary groups
$$
\imath^{\hsuit}: \rmU(V^{\hsuit}, \phi^{\hsuit}, \tHerm^{\hsuit}) \times \rmU(V, -\phi, \tHerm_{m,n}) \hookrightarrow \rmU(V^{\bdsuit}, \phi^{\bdsuit}, \tHerm^{\bdsuit, \oplus}), \quad (g,h) \mapsto \diag[g,h].
$$
Moreover, we have another embedding of unitary groups \footnote{These two embeddings arise from two viewpoints on the ``direct sum'' $\phi^{\bdsuit}$: 
$$
\diag[\diag[1, \bfone_{m}, -\bfone_{n}, -1], \diag[-\bfone_{m}, \bfone_{n}]] = \diag[\diag[1, \bfone_{m}, -\bfone_{n}], \diag[-1, -\bfone_{m}, \bfone_{n}]].
$$
The left-hand side partition gives the embedding $\imath^{\hsuit}$ and the right-hand side gives the embedding $\imath^{\dsuit}$.}
$$
\imath^{\dsuit}: \rmU(V^{\sharp}, \phi^{\sharp}, \tHerm^{\sharp}) \times \rmU(V^{\sharp}, -\phi^{\sharp}, \tHerm^{\sharp}) \hookrightarrow \rmU(V^{\bdsuit}, \phi^{\bdsuit}, \tHerm^{\bdsuit, \oplus}), \quad (g_1, g_2) \mapsto \diag[g_1, g_2].
$$
From the definition of the embeddings above, we see immediately that
\begin{equation} \label{eq:basiscomp}
\imath^{\hsuit}(\jmath^{\flat}(g), h) = \imath^{\dsuit}(g, \jmath^{\sharp}(h))
\end{equation}
for any $h \in \rmU(V, \phi, \tHerm_{m,n})$ and $g \in \rmU(V^{\sharp}, \phi^{\sharp}, \tHerm^{\sharp})$.

We rearrange the basis to represent the corresponding metric in standard Hermitian form:
$$
\mathtt{Herm}^{\bdsuit}: \tte, \mathtt{a}^1, \ldots, \mathtt{a}^{m}, \bar{\ttb}^{1}, \ldots, \bar{\ttb}^{n}, \ttf, \bar{\tta}^{1}, \ldots, \bar{\tta}^{m}, \mathtt{b}^1, \ldots, \mathtt{b}^{n}.
$$
Then under this basis,
$$
[-\sfi \phi^{\bdsuit}]_{\mathtt{Herm}^{\bdsuit}} := \varsigma^{\bdsuit} \diag[\phi^{\hsuit}, -\phi]  (\varsigma^{\bdsuit})^{-1} = \diag[1, \bfone_{m}, \bfone_{n}, -1, -\bfone_{m}, -\bfone_{n}],
$$
where
$$
\varsigma^{\bdsuit} = \begin{bmatrix}
\bfone_{1+m} &  &  &  \\ &  & & \bfone_{n} \\  &  & \bfone_{1+m} &  \\ & \bfone_{n}  & & 
\end{bmatrix}.
$$
This induces an isomorphism of unitary groups
$$
\varsigma^{\bdsuit}: \rmU(V^{\bdsuit}, \phi^{\bdsuit}, \tHerm^{\bdsuit}) \xrightarrow{\sim} \rmU(V^{\bdsuit}, \phi^{\bdsuit}, \tHerm^{\bdsuit, \oplus}), \quad g \mapsto \varsigma^{\bdsuit} g (\varsigma^{\bdsuit})^{-1}.
$$
We then return to the Witt basis via $\varrho_{m+n+1, m+n+1}$. The embeddings $\imath^{\hsuit}$ and $\imath^{\dsuit}$ are defined such that the diagrams
$$
\begin{tikzcd}
	{\rmU(V^{\hsuit}, \phi^{\hsuit}, \tHerm^{\hsuit}) \times \rmU(V, -\phi, \tHerm_{m,n})} & {\rmU(V^{\bdsuit}, \phi^{\bdsuit}, \tHerm^{\bdsuit, \oplus})} & {\rmU(V^{\bdsuit}, \phi^{\bdsuit}, \tHerm^{\bdsuit})} \\
	{\rmU(V^{\hsuit}, \phi^{\hsuit}, \tWitt^{\hsuit}) \times \rmU(V, -\phi, \tWitt_{m,n})} && {\rmU(V^{\bdsuit}, \phi^{\bdsuit}, \tWitt^{\bdsuit})}
	\arrow["{\imath^{\hsuit}}", from=1-1, to=1-2]
	\arrow["\varrho_{m+1, n+1, \vartheta_0} \times \varrho_{m,n, \vartheta_0}"', from=1-1, to=2-1]
	\arrow["{\varsigma^{\bdsuit}}", from=1-2, to=1-3]
	\arrow["{\varrho_{m+n+1, m+n+1}}", from=1-3, to=2-3]
	\arrow["{\imath^{\hsuit}}"', from=2-1, to=2-3]
\end{tikzcd}
$$
and
$$
\begin{tikzcd}
	{\rmU(V^{\sharp}, \phi^{\sharp}, \tHerm^{\sharp}) \times \rmU(V^{\sharp}, -\phi^{\sharp}, \tHerm^{\sharp})} & {\rmU(V^{\bdsuit}, \phi^{\bdsuit}, \tHerm^{\bdsuit, \oplus})} & {\rmU(V^{\bdsuit}, \phi^{\bdsuit}, \tHerm^{\bdsuit})} \\
	{\rmU(V^{\sharp}, \phi^{\sharp}, \tWitt^{\sharp}) \times \rmU(V^{\sharp}, -\phi^{\sharp}, \tWitt^{\sharp})} && {\rmU(V^{\bdsuit}, \phi^{\bdsuit}, \tWitt^{\bdsuit})}
	\arrow["{\imath^{\dsuit}}", from=1-1, to=1-2]
	\arrow["{\varrho_{m+1, n, \vartheta_0^{\sharp}} \times \varrho_{m+1, n, \vartheta_0^{\sharp}}}"', from=1-1, to=2-1]
	\arrow["{\varsigma^{\bdsuit}}", from=1-2, to=1-3]
	\arrow["{\varrho_{m+n+1, m+n+1}}", from=1-3, to=2-3]
	\arrow["{\imath^{\dsuit}}"', from=2-1, to=2-3]
\end{tikzcd}
$$
commute. It then follows that \eqref{eq:basiscomp} holds for the embeddings under the Witt basis. To sum up, we have two embeddings of unitary groups:
\begin{align*}
\imath^{\hsuit}&: \rmU(m+1,n+1) \times \rmU(m,n) \rightarrow \rmU(m+n+1, m+n+1), \\
\imath^{\dsuit}&: \rmU(m+1,n) \times \rmU(m+1,n) \rightarrow \rmU(m+n+1, m+n+1).
\end{align*}
To simplify the notation, we set
$$
H := \rmU(m,n), \quad G := \rmU(m+1,n), \quad G^{\hsuit} := \rmU(m+1,n+1), \quad G^{\bdsuit} := \rmU(m+n+1, m+n+1)
$$
and $\bfG = H \times G$, with $\bfH$ being the image of $H$ via the diagonal embedding $\Delta^{\sharp}: H \rightarrow \bfG$ defined by $h \mapsto (h, \jmath^{\sharp}(h))$. We always write $N := m+n$.

\subsection{Eisenstein series on unitary groups}
In this section, we define Klingen Eisenstein series and Siegel Eisenstein series over unitary groups.

\subsubsection{Klingen Eisenstein series} \label{sec:Klingen}
We first define the \emph{Klingen parabolic subgroup} of $G^{\hsuit}$ as the algebraic group $P$ over $\calF$ such that, for any $\calF$-algebra $R$:
$$
P(R) = \left\{ g = 
\begin{bmatrix}
a &  & b & c & \ast \\
\ast & x^{-\star} & \ast  & \ast & \ast \\
d &  & e & f & \ast \\
h &  & l & k & \ast  \\
 &  &  &  & x
\end{bmatrix}
 \in G^{\hsuit}(R) : g_0 := \begin{bmatrix}
 a & b & c \\ d & e & f \\ h & l & k
 \end{bmatrix} \in H(R), x \in (\calK \otimes_{\calF} R)^{\times} 
\right\}.
$$
The block matrix above is written with respect to the partition $[n \mid 1 \mid m-n \mid n \mid 1]$. The subgroup $P$ has a Levi decomposition $P = MN$, where the Levi subgroup $M$ is given by
$$
M(R) = \left\{ \bfm(g_0, x) := \begin{bmatrix}
a &  & b & c &  \\
 & x^{-\star} &  &  &  \\
d &  & e & f &  \\
h &  & l & k &  \\
 &  &  &  & x
\end{bmatrix} \in P(R) \right\} \xleftarrow{\sim} H(R) \times \Res_{\calK/\calF} \GG_{\rmm}(R),
$$
with the isomorphism provided by the map $\bfm(g_0, x)$.

To construct an Eisenstein series on $G^{\hsuit}$ with respect to the Klingen parabolic subgroup $P$, we define the following Eisenstein datum.

\begin{definition} \label{defn:eisensteindatum}
An \emph{Eisenstein datum} over $H$ is a tuple $\calD = (\sigma, \psi, \chi, \scS)$, where:
\begin{itemize}
    \item $(\sigma, V_{\sigma})$ is an irreducible tempered unitary cuspidal automorphic representation of $H$ over $\calF$, whose archimedean part $\sigma_{\infty}$ is a holomorphic discrete series representation. More explicitly, for any $v \in \scV_{\calF}^{(\infty)}$, the \emph{weight} $\ulk_{v}$ of $\sigma$ is an $(m+n)$-tuple
\[
\ulk_{v} = (t_{1,v}^{+}, \ldots, t_{m,v}^{+}; t_{1, v}^{-}, \ldots ,t^{-}_{n, v}) \in \ZZ^{m+n}
\]
such that $t_{1,v}^{+} \geq \cdots \geq t_{m,v}^{+}$ and $t_{1, v}^{-} \geq \cdots \geq t^{-}_{n, v}$. These weights classify the $L$-packet of discrete series of $H(\RR) := H(\calF_{v})$, denoted by $\sigma_{v, \ulk_{v}}$. The \emph{holomorphic discrete series} corresponds to weights satisfying $t_{m,v}^{+} - t_{1,v}^{-} \geq m+n$. We identify $V_{\sigma}$ as a subspace of the space of cuspidal automorphic forms $\calA_{\cusp}(H(\calF) \bs H(\AA))$.
\item $\psi$ is a unitary Hecke character of $\calK$ such that $\psi|_{\AA^{\times}}$ is the central character of $\sigma$.
\item $\chi: \KK^{\times} \rightarrow \CC^{\times}$ is a unitary Hecke character of infinity type $\frac{\kappa}{2}(1-\sfc)$ for $\kappa \in \ZZ[\Sigma]$. We also regard $\chi$ as an automorphic form over the reductive group $\Res_{\calF}^{\calK}(\GG_{\rmm})$.
\item $\scS$ is a finite set of places $v$ of $\calF$ where at least one of $\sigma_v$, $\chi_v$, or $\calK_v/\calF_v$ is ramified, and it contains $\scV_{\calF}^{\{p, \infty\}}$.
\end{itemize}
Attached to an Eisenstein datum $\calD$ over $H$, we define the character $\xi := \psi/\chi = \psi \barchi$.
\end{definition}

Let $s$ be an arbitrary complex number. Given an Eisenstein datum $\calD$ over $H$, the exterior tensor product $\sigma \boxtimes \chi$ defines an automorphic representation of $M(\AA)$ via the identification $\bfm$, which we extend trivially to $P(\AA)$. Let $\delta_P$ be the modulus character of $P$, and let $\delta$ be a character such that $\delta^{m+n+1} = \delta_{P}$. We consider the induced representation \footnote{Henceforth, $\Ind_{P(\AA)}^{G^{\hsuit}(\AA)}$ denotes the \emph{unnormalized} smooth parabolic induction. We have added the factor of the modulus character to effectively treat it as a normalized induction. We follow the notation of \cite{skinner2014iwasawa, wan2015families} by using the modified character $\delta$ instead of $\delta_{P}$, where $\delta_{P}(\bfm(g_0,x)) = \abs{x \barx}_{\AA}^{-(m+n+1)}$. Thus, the complex variable $s$ used here relates to the variable $s_{\mathrm{ICM}}$ in \cite{skinner2006icm} by $s = (m+n+1) s_{\mathrm{ICM}}$. Our conventions are then consistent with \cite{skinner2006icm}.}
$$
I^{\Kling}(\sigma, \chi, s) := \Ind_{P(\AA)}^{G^{\hsuit}(\AA)}(\delta^{\frac{m+n+1}{2}+s} \cdot \sigma \boxtimes \chi).
$$
The representation space $I^{\Kling}(\sigma, \chi, s)$ consists of smooth functions $f_{\sigma, \chi, s}^{\Kling}: G^{\hsuit}(\AA) \rightarrow V_{\sigma}$ satisfying:
\begin{enumerate}
    \item For any $\bfm(g_0, x) \in M$, $n \in N$, and $g \in G^{\hsuit}(\AA)$, 
    $$
    f_{\sigma, \chi, s}^{\Kling}(\bfm(g_0, x) ng) = \delta(\bfm(g_0, x))^{\frac{m+n+1}{2}+s} \chi(x)\sigma(g_0)f_{\sigma, \chi, s}^{\Kling}(g).
    $$
    \item $f_{\sigma, \chi, s}^{\Kling}$ is right $K$-finite for some maximal open compact subgroup $K \subseteq G^{\hsuit}(\AA)$.
\end{enumerate}
Elements of this space are called \emph{Klingen Eisenstein sections} with respect to the datum $\calD$.

We have a natural evaluation map
$$
\ev_g : \calA_{\cusp}(H(\calF) \bs H(\AA)) \rightarrow \CC, \quad \Phi \mapsto \Phi(g)
$$
for any $g \in H(\AA)$. By our assumption, $V_{\sigma}$ lies in the space of cuspidal automorphic forms, so for every Klingen section $f_{\sigma, \chi, s}^{\Kling}$, we attach it with a scalar-valued section $f_{s, \chi, \sigma}^{\Kling, g} := \ev_g(f_{\sigma, \chi, s}^{\Kling})$. Then we define the Klingen Eisenstein series as
$$
E^{\Kling}(f_{\sigma, \chi, s}^{\Kling}, g) := \sum_{\gamma \in P(\calF) \bs G^{\hsuit}(\calF)} f_{s, \chi, \sigma}^{\Kling,1}(\gamma g).
$$
This series converges absolutely and uniformly for $(s,g)$ in compact subsets of $\{s \in \CC: \Re(s) > (m+n+1)/2\} \times G^{\hsuit}(\AA)$.

In this article, we assume $\sigma_{\infty}$ has scalar weights and $\sigma_{p}$ is sufficiently ramified following \cite[Definition 4.42]{wan2015families}. These assumptions are needed for the construction of an explicit Hida family of Klingen Eisenstein series, as in \cite{wan2015families}.

\begin{assumption}[Scalar weight assumption] \label{ass:scalarwt}
We fix $\kappa \in \ZZ \subseteq \ZZ[\Sigma]$. For any $v \in \scV_{\calF}^{(\infty)}$, we assume $\sigma_{v}$ is a holomorphic discrete series representation with scalar weight $(0, \ldots, 0; -\kappa, \ldots, -\kappa)$ ($m$ zeros and $n$ kappas). Here $\kappa \geq m+n$, and we are particularly interested in $s_{\kappa} := (\kappa - (m + n + 1))/2 \in \CC$. \footnote{We use this notation to ensure $\kappa > 0$, matching the conventions in \cite{wan2015families}.}
\end{assumption}

\begin{assumption}[Sufficiently ramified assumption] \label{ass:generic}
We assume $\sigma_v$ is $p$-ordinary and $(\sigma_{v}, \chi_v)$ is ``generic'' in the sense of \cite[Definition 4.42]{wan2015families} for all $v \in \scS_{\calF}^{(p)}$. This requires the ramification of $\pi$ at primes dividing $p$ to be sufficiently large. \footnote{This requirement informs our decision to go beyond the ``semi-stably ordinary'' case treated in \cite{liu2023anticyclotomicpadiclfunctionsrankinselberg} during our computations in Section \ref{sec:iiintegralatp}.}
\end{assumption}

\subsubsection{Siegel Eisenstein series}
As we shall see in Section \ref{sec:garrettdou}, Klingen Eisenstein series are often constructed via the pullbacks of certain Siegel Eisenstein series on the larger quasi-split unitary groups $G^{\bdsuit}$. We define the \emph{Siegel parabolic subgroup} of $G^{\bdsuit}$ as the algebraic group $Q$ over $\calF$ such that, for any $\calF$-algebra $R$:
$$
Q(R) = \left\{ g = \begin{bmatrix}
 A_g & B_g \\ 0 & D_g
\end{bmatrix} \in \GL_{2N+2}(R \otimes_{\calF} \calK) : D_g = A_g^{-\star}, A_g^{-1} B_g \in \Herm_{N+1}(R \otimes_{\calF} \calK) \right\}.
$$
The block matrix above is written with respect to the partition $[N+1 \mid N+1]$. It has a Levi decomposition $Q = M_{Q} N_{Q}$, where 
\begin{align*}
M_{Q}(R) &:= \left\{ \begin{bmatrix}
 A_g & 0 \\ 0 & A_g^{-\star}
\end{bmatrix} : A_g \in \GL_{N+1}(R \otimes_{\calF} \calK)  \right\} \cong \GL_{N+1}(R \otimes_{\calF} \calK), \\
N_{Q}(R) &:= \left\{ \begin{bmatrix}
 \bfone_{N+1} & X_g \\ 0 & \bfone_{N+1}
\end{bmatrix} : X_g \in \Herm_{N+1}(R \otimes_{\calF} \calK)  \right\} \cong \Herm_{N+1}(R \otimes_{\calF} \calK).
\end{align*}
Given any $g \in Q(R)$, we can decompose it as
$$
g = \begin{bmatrix}
 A_g & B_g \\ 0 & D_g
\end{bmatrix} = \begin{bmatrix}
 A_g & 0 \\ 0 & A_g^{-\star}
\end{bmatrix} \begin{bmatrix}
 \bfone_{N+1} & X_g \\ 0 & \bfone_{N+1}
\end{bmatrix}
$$
with $A_g \in \GL_{N+1}(R \otimes_{\calF} \calK)$ and $X_g = A_g^{-1} B_g \in \Herm_{N+1}(R \otimes_{\calF} \calK)$. For any character $\chi$ of $(R \otimes_{\calF} \calK)^{\times}$, we view it as a character on the Levi subgroup $M_{Q}$ via
$$
\chi : \begin{bmatrix}
 A_g & 0 \\ 0 & A_g^{-\star}
\end{bmatrix} \mapsto \chi(\det D_g),
$$
and extend it trivially to $N_{Q}$ to obtain a character of $Q$.

Let $v$ be a place of $\calF$. For any character $\chi_v : \calK_v^{\times} \rightarrow \CC^{\times}$ and $s \in \CC$, we define the local induction space
$$
I^{\Sieg}_{v,Q}(\chi_v, s) := \Ind_{Q(\calF_v)}^{G^{\bdsuit}(\calF_v)} (\delta_{Q}^{\frac{1}{2}+s} \cdot \chi_v),
$$
consisting of smooth functions $f_{s, \chi_v}^{\Sieg} : G^{\bdsuit}(\calF_v) \rightarrow \CC$ such that:
\begin{enumerate}
    \item For any $g \in Q(\calF_v)$ and $h \in G^{\bdsuit}(\calF_v)$, 
    $$
    f_{s, \chi_v}^{\Sieg}\left( \begin{bmatrix}
 A_g & B_g \\ 0 & D_g
\end{bmatrix} h \right) = \delta_{Q}(g)^{\frac{1}{2}+s} \chi_v(\det D_g)f_{s, \chi_v}^{\Sieg}(h) = \abs{\det A_g D_g^{-1}}_{v}^{s+\frac{n}{2}} \chi_v(\det D_g)f_{s, \chi_v}^{\Sieg}(h),
    $$
    \item $f_{s, \chi_v}^{\Sieg}$ is right $K_v$-finite, where $K_v := G^{\bdsuit}(\calO_{\calF_v})$ when $v \in \scV_{\calF}^{(0)}$, and $K_v$ is defined as the maximal compact subgroup of $G^{\bdsuit}(\calF_v)$ made explicit in \cite[Section 5.1.3]{skinner2014iwasawa} when $v \in \scV_{\calF}^{(\infty)}$.
\end{enumerate}
Such an element is called a \emph{Siegel Eisenstein section} at $v$ with respect to the datum $(\chi_v,s)$. If $v$ is a finite place and $\chi_v$ is unramified, we define the \emph{spherical Siegel Eisenstein section} $f^{\sph}_{s, \chi_v} \in I^{\Sieg}_v(\chi_v, s)$ as the unique section satisfying $f^{\sph}_{s, \chi_v}(K_v) = 1$ for $K_v = G^{\bdsuit}(\calO_{\calF_v})$.

Let $\chi : \KK^{\times} \rightarrow \CC^{\times}$ be a unitary Hecke character with the decomposition $\chi = \otimes_{v}^{\prime} \chi_v$. Consider the restricted tensor product $I^{\Sieg}(\chi, s) = \otimes_v^{\prime} I^{\Sieg}_v(\chi_v, s)$ with respect to the spherical sections $f^{\sph}_{s, \chi_v}$ at almost all finite places where $\chi_v$ is unramified. For $f_{\chi, s}^{\Sieg} \in I^{\Sieg}(\chi, s)$, the \emph{Siegel Eisenstein series} is defined as
$$
E^{\Sieg}(f_{s, \chi}^{\Sieg}, g) := \sum_{\gamma \in Q(\calF) \bs G^{\bdsuit}(\calF)} f_{s, \chi}^{\Sieg}(\gamma g).
$$
This series converges absolutely and uniformly for $(s,g)$ in compact subsets of $\{s \in \CC : \Re(s) > (N+1)/2\} \times G^{\bdsuit}(\AA)$.

%\begin{remark}
%We note that the notion of Siegel Eisenstein sections here slightly differs from \cite{skinner2014iwasawa, wan2015families}. As explained in \cite[Remark on page 170]{skinner2014iwasawa}, their representation $\chi$ on the Levi subgroup $M_Q$ is defined by
%$$
%\chi: \begin{bmatrix}
% D_g^{-\star} & 0 \\ 0 & D_g
%\end{bmatrix} \mapsto \chi(\det D_g).
%$$
%So the induced representation $I(\chi)$ in \textit{loc.cit} is actually $I^{\Sieg}_{Q}((\chi^{\rmc})^{-1},s)$ in our setup. Their convention is convenience when dealing with the functional equations of Eisenstein series.
%\end{remark}

\subsection{On multiplicity one conditions} \label{sec:multiplicityone}
We record a result regarding the multiplicity one theorem for automorphic representations of unitary groups. Let $\pi$ be an automorphic representation of a unitary group $\rmU(m,n)$ over $\calF$, and suppose that:
\begin{enumerate}[label = (BC)] 
    \item the functorial base change $\BC(\pi)$ of $\pi$ to $\GL_{m+n, \calK}$ is a cuspidal automorphic representation. \label{ass:multiplicityone}
\end{enumerate}
Then $\pi$ occurs in the space of cuspidal automorphic forms over $\rmU(m,n)$ with multiplicity one, as a consequence of \cite{mok2015endoscopic, kaletha2014endoscopicclassificationrepresentationsinner}. Consequently, up to a complex scalar, there exists a unique $\rmU(m,n)(\AA)$-invariant pairing between $\pi$ and its contragredient $\pi^{\vee}$.

In particular, let $\pi$ be any regular algebraic cuspidal automorphic representation of $\rmU(m,n)$ over $\calF$ whose residual Galois representation is absolutely irreducible; then property \ref{ass:multiplicityone} holds for $\pi$. Recall that we have assumed from the outset that the automorphic representations $\sigma$ and $\pi$ considered in this article satisfy this condition. Therefore, $\sigma$ and $\pi$ satisfy property \ref{ass:multiplicityone} and hence the multiplicity one theorem.

%% file: AP02-PreliminariesGeometry.tex
\section{Modular forms over unitary groups} \label{sec:formsunitary}
In this section, we recall the geometric background of modular forms over unitary groups, their $p$-adic analogues, and $p$-adic families. The materials are largely taken from \cite{eischen2020padic, hsieh2014eisenstein}.

\subsection{Generalities on PEL datums and unitary Shimura datums}
\subsubsection{PEL-type Shimura datums}
\begin{definition}
By a \emph{PEL-type Shimura datum}, we mean a tuple $\scP = (B, \ast, V, \psi_{V}, h)$, where
\begin{itemize}
    \item $B$ is a finite semisimple $\QQ$-algebra with a positive involution $\ast$,
    \item $V$ is a \emph{symplectic $(B, \ast)$-module}, that is, a $B$-module $V$ with a skew-symmetric nondegenerate $\QQ$-bilinear form $\psi_{V}: V \times V \rightarrow \QQ$ such that $\psi_{V}(b^{\ast}u,v) = \psi_{V}(u, bv)$ for any $u, v \in V$ and $b \in B$. 
    \item Let $\calC := \End_{B}(V)$. It carries an \emph{adjoint involution} defined by, for $\alpha \in \End_{B}(V)$, the $\alpha^{\ast} \in \calC$ such that 
    $$
    \psi_{V}(\alpha^{\ast} v, w) = \psi_{V}(v, \alpha w), \quad \text{for all } v, w \in V.
    $$
    Then $h$ is defined to be an $\RR$-algebra map $h: \CC \rightarrow \calC_{\RR}$ such that
    \begin{itemize}
        \item $h(\barz) = h(z)^{\ast}$, where the $\ast$ on the right-hand side is the adjoint involution on $\End_{B}(V)_{\RR}$,
        \item $(u,v) \mapsto \psi_V(u, h(\rmi)v)$ is positive-definite and symmetric.
    \end{itemize}
\end{itemize}
We say $\scP$ is a \emph{simple PEL-type Shimura datum} if $B$ is a simple $\QQ$-algebra. Given a PEL-type Shimura datum $\scP$, we associate with it:
\begin{itemize}
    \item the field $F$, defined to be the center of $B$,
    \item the field $F_0 := \{b \in F: b^{\ast} = b\}$, i.e., the subalgebra of $\ast$-invariants in $F$, and
    \item two algebraic groups over $\QQ$, defined as
    $$
    G(R) := \{g \in \GL_{B \otimes_{\QQ} R}(V \otimes_{\QQ} R) : \psi_{V}(gu,gv) = \nu(g) \psi_{V}(u,v), \, u,v \in V \otimes_{\QQ} R, \, \nu(g) \in R^{\times} \}
    $$
    and
    $$
    G_1(R) := \{g \in \GL_{B \otimes_{\QQ} R}(V \otimes_{\QQ} R) : \psi_{V}(gu,gv) = \nu(g) \psi_{V}(u,v), \, u,v \in V \otimes_{\QQ} R, \, \nu(g) \in (F_0 \otimes_{\QQ} R)^{\times} \}        
    $$
    for any $\QQ$-algebra $R$. Clearly, $G$ is a subgroup of $G_1$, identifying similitude factors via the natural inclusion $R^{\times} \hookrightarrow (F_0 \otimes_{\QQ} R)^{\times}$.
\end{itemize}
\end{definition}

Let $X$ be the $G(\RR)$-conjugacy class of $h^{-1}: \CC^{\times} \rightarrow G_{\RR}$. Then $(G, X)$ is a Shimura datum \`a la Deligne. For each \emph{neat} open compact subgroup $K$ of $G(\AA_{\QQ, \rmf})$, by the result of Deligne, there is an algebraic variety $\Sh_{K}(G, X)$ over the reflex field $E$ of the PEL Shimura datum $\scP$ such that 
$$
\Sh_{K}(G,X)(\CC) = G(\QQ) \bs X \times G(\AA_{\QQ, \rmf}) / K.
$$
Under certain conditions, this model is unique. We call it the \emph{Shimura variety} of the Shimura datum $(G,X)$, or the Shimura variety of the PEL Shimura datum $\scP$.

Let $p$ be a prime number. We hope to construct an integral model of the Shimura variety $\Sh_{K}(G,X)$ at $p$, that is, a smooth model over the ring $\calO_{E} \otimes_{\ZZ} \ZZ_{(p)}$. For this purpose, we need some extra data and assumptions.

\begin{definition}
By an \emph{integral PEL-type Shimura datum}, we mean a tuple $(\scP, \calO_{B}, L)$, where
\begin{itemize}
    \item $\scP = (B, \ast, V, \psi_{V}, h)$ is a PEL-type Shimura datum, 
    \item $\calO_{B}$ is a $\ZZ_{(p)}$-order in $B$ which is stable under the involution $\ast$ on $B$, and $\calO_{B} \otimes_{\ZZ_{(p)}} \ZZ_p$ is a maximal order in $B$. We require $B$ to be unramified at $p$, which means that $B_{\QQ_p}$ is isomorphic to a product of matrix algebras over unramified extensions of $\QQ_p$.
    \item $L$ is a $\ZZ_p$-lattice in $V_{\QQ_p}$ such that $L$ is stable under $\calO_B$ and $L$ is self-dual with respect to the pairing $\psi_V$.
\end{itemize}
\end{definition}

Given an integral PEL-type Shimura datum with the extra unramified condition of $B$ at $p$, $G_{\QQ_p}$ is unramified. Indeed, let
$$
K_p^{0} := \{g \in G(\QQ_p): g \cdot L = L \},
$$
i.e., the subgroup of $G(\QQ_p)$ that stabilizes the lattice $L$. Then $K_p^{0}$ is the hyperspecial subgroup of the $\ZZ_p$-points of the smooth reductive model $\calG$ of $G_{\QQ_p}$ over $\ZZ_p$.

\subsubsection{Unitary Shimura datums}
\begin{definition}
We consider a special integral PEL-type Shimura datum $(\scP = (B, \ast, V, \psi_{V}, h), \calO_B, L)$, where
\begin{itemize}
    \item $B = \calK^{m}$, the product of $m$ copies of $\calK$.
    \item $\ast$ is the complex conjugation on each factor $\calK$ of $B$.
    \item For $1 \leq i \leq m$, let $V_i$ be a finite-dimensional $\calK$-vector space of dimension $n_i$, equipped with a Hermitian form $\lrangle{-,-}_{V_i}: V_i \times V_i \rightarrow \calK$ relative to $\calK/\calF$. Let $\delta \in \calO_{\calK}$ be a totally imaginary element that is prime to $p$. We then put $\psi_{V_i} := \Tr_{\calK/\QQ} (\delta \lrangle{-,-}_{V_i})$. Then $V$ is taken to be $V_1 \times \cdots \times V_{m}$, and $\psi_V := \psi_{V_1} \times \cdots \times \psi_{V_m}$.
    \item For each $\sigma \in \scV_{\calK}^{\infty}$, $V_{i, \sigma} := V_i \otimes_{\calK, \sigma} \CC$ has a $\CC$-basis with respect to which $\psi_{V_i}$ is given by a matrix of the form $\diag[\bfone_{r_{i, \sigma}}, -\bfone_{s_{i, \sigma}}]$. Fixing such a basis, let $h_{i, \sigma}: \CC \rightarrow \End_{\RR}(V_{i,\sigma})$ be $z \mapsto \diag[z \bfone_{r_{i, \sigma}}, \barz \bfone_{s_{i, \sigma}}]$. Let $h_i := \prod_{\sigma \in \Sigma_{\calK}} h_{i, \sigma}$ and $h := \prod_{i=1}^{m} h_i$. 
    \item $\calO_{B} = \calO_{\calK}^{m}$.
    \item Let $L_i \subset V_i$ be a free $\calO_{\calK}$-module of rank $\dim_{\calK} V_{i}$, such that $\lrangle{L_i, L_i}_{\QQ} \subset \ZZ$ and $\psi_{V_i}$ is a perfect pairing on $L_i \otimes \ZZ_p$.
\end{itemize}
We call this PEL-type Shimura datum a \emph{unitary Shimura datum}. It is called a \emph{simple unitary Shimura datum} if $m = 1$.
\end{definition}

Given such a unitary Shimura datum, the corresponding objects are as follows:
\begin{itemize}
    \item the field $F = \calK^m$,
    \item the field $F_0 = \calF^{m}$,
    \item the group $G_1 = \Res^{\calF}_{\QQ} \left( \prod_{i=1}^{m} \GU(V_i, \lrangle{-,-}_{V_i}) \right)$, where for $1 \leq i \leq m$, $\GU(V_i)$ is the general unitary group over $\calF$ attached to the Hermitian space $V_i$, defined in Definition \ref{defn:unitarygrp}.
    \item the group $G$ is the subgroup of $G_1$ consisting of rational similitudes.
\end{itemize}
In particular, when $\calF = \QQ$ and $m=1$, then $G_1 = G$ is the general unitary group defined in Definition \ref{defn:unitarygrp}.

\subsubsection{Hodge structures, lattices and level subgroups at $p$} \label{sec:levelgroups}
Now we concentrate on the case of a unitary Shimura datum. For each $1 \leq i \leq m$ and $\sigma \in \scV_{\calK}^{\infty}$, $h_{i, \sigma}$ determines a pure Hodge structure of weight $-1$ on $V_{i, \sigma} := L_{i} \otimes_{\calO_{\calK}, \sigma} \CC$. Let $V_{i, \sigma}^{0}$ be the degree zero piece of the Hodge filtration. This is an $\calO_{\calK} \otimes_{\sigma} \CC$-submodule of $V$. For each $\sigma \in \scV_{\calK}^{\infty}$, let $a_{\sigma, i} = \dim_{\CC}(V_{i, \sigma}^{0} \otimes_{\calO_{\calK} \otimes_{\sigma} \CC} \CC)$ and $b_{i, \sigma} := n_i - a_{\sigma, i}$. We note that for $\sigma \in \scV_{\calK}^{\infty}$, 
$$
(a_{\sigma, i}, b_{\sigma, i}) = (b_{\sfc \sigma, i}, a_{\sfc \sigma, i}) = (r_{i, \sigma}, s_{i, \sigma}).
$$
We assume throughout the following fundamental hypothesis, called the \emph{ordinary hypothesis}:
\begin{equation}\tag{ord}
    v_{\sigma} = v_{\sigma^{\prime}} \Rightarrow r_{\sigma, i} = r_{\sigma^{\prime}, i}, \quad i = 1, \ldots, m
\end{equation}
for any $\sigma, \sigma^{\prime} \in \scV_{\calK}^{\infty}$, where $v_{\sigma}$ and $v_{\sigma^{\prime}}$ are their corresponding $p$-adic places given by the embedding $\iota_p$. This hypothesis holds in our setup since we have assumed \eqref{ass:sgn} throughout.

Then for each place $w$ of $\calK$ above $p$, we can define $(a_{w,i}, b_{w,i}) = (a_{\sigma, i}, b_{\sigma, i})$ for any $\sigma \in \Sigma_{\calK}$ corresponding to $w$. Let $L_{i, w} := L_{i} \otimes_{\calO_{\calK}} \calO_{\calK, w}$. For each $1 \leq i \leq m$, we fix an $\calO_{\calK} \otimes \ZZ_p$-decomposition $L_i \otimes \ZZ_p = L_i^{+} \oplus L_i^{-}$ such that:
\begin{itemize}
    \item $L_{i,w}^{+}$ is an $\calO_{\calK, w}$-module with $\rank_{\calO_{\calK, w}} L_{i,w}^{+} = a_{w, i}$. So $L_{i,w}^{-}$ is an $\calO_{\calK, w}$-module with $\rank_{\calO_{\calK, w}} L_{i,w}^{-} = b_{w, i}$, and $L_{i, w} = L_{i,w}^{+} \oplus L_{i,w}^{-}$.
    \item $L_{i, w}^{\pm}$ is the annihilator of $L_{i, \barw}^{\pm}$ with respect to the perfect pairing $\psi_V: L_{i,w} \times L_{i, \barw} \rightarrow \ZZ_p(1)$.
\end{itemize}
For each $i$ and $w$, we fix a basis of $L_{i,w}^{+}$ to regard it as a direct sum of copies of $\calO_{\calK, w}$. \footnote{In \cite[Section 1.8]{hsieh2014eisenstein}, Hsieh very carefully chose the precise basis of $L_{i,w}^{+}$. Here we only roughly speak of such a basis without explicitly specifying it.} Taking $\ZZ_p$-duals via $\psi_V$ yields a decomposition of $L_{i,w}^{-}$ as a direct sum of copies of $\calO_{\calK, \barw}$. The choices of these decompositions determine isomorphisms:
\[
\GL_{\calO_{\calK, w}}(L_{i,w}^{+}) \simeq \GL_{a_{w,i}}(\calO_{\calK, w}),
\]
\[
\GL_{\calO_{\calK, w}}(L_{i,w}^{-}) \simeq \GL_{b_{w,i}}(\calO_{\calK, w}),
\]
\[
\GL_{\calO_{\calK, w}}(L_{i,w}) \simeq \GL_{n_{i}}(\calO_{\calK, w}),
\]
and the embedding
$$
\GL_{\calO_{\calK, w}}(L_{i,w}^{+}) \times \GL_{\calO_{\calK, w}}(L_{i,w}^{-}) \hookrightarrow \GL_{\calO_{\calK, w}}(L_{i,w})
$$
is given by $(A,B) \mapsto \diag[A,B]$.

To define appropriate level subgroups at $p$, we start by defining:
$$
H_{\scP} = \GL_{\calO_{B} \otimes \ZZ_p} (L^{+}) \xrightarrow{\sim} \prod_{w \mid p} \prod_{i=1}^{m} \GL_{a_{w,i}}(\calO_{\calK,w}). 
$$
Let $B_{\scP} \subset H_{\scP}$ be the Borel subgroup corresponding via this isomorphism to the product of the upper-triangular Borel subgroups of the general linear groups. Let $N_{\scP}$ be its unipotent radical. Let $T_{\scP} := B_{\scP}/N_{\scP}$; this is identified by the isomorphism with diagonal matrices. Let $B^{+}_{\scP} \subset G_{\ZZ_p}$ be the Borel subgroup that stabilizes $L^{+}$ and such that:
\begin{equation} \label{eq:Bplusproj}
    B^{+}_{\scP} \twoheadrightarrow \GG_{\rmm} \times B_{\scP} \subset \GG_{\rmm} \times H_{\scP},
\end{equation}
where the map to the first factor is the similitude character $\nu$ and the second projection is the projection to $H_{\scP}$. Let $N_{\scP}^{+}$ be its unipotent radical. Under the identification:
\begin{equation} \label{eq:BpscP}
B^{+}_{\scP} \xrightarrow{\sim} \GG_{\rmm} \times \prod_{w \in \scV_{\calK}^{(p)}} \prod_{i=1}^{m} \left\{ \begin{bmatrix}
A & B \\ 0 & D
\end{bmatrix} \in \GL_{n_i}(\calO_{\calK,w}): A \in B_{a_{w, i}}(\calO_{\calK, w}), D \in B^{-}_{b_{w, i}}(\calO_{\calK, w}) \right\}.    
\end{equation}

Then we define the following level subgroups at $p$. They are all defined under the chosen basis of $L^{\bullet}$ above:
\begin{itemize}
    \item We have the \emph{hyperspecial subgroup} $K_{p}^{0} = \GG_{\rmm} \times H_{\scP}$, where $\GG_{\rmm}$ is the similitude factor part.
    \item We have the \emph{congruence level subgroup} $K_{p, r}^{0} \subset G(\ZZ_p)$ consisting of those $g$ such that $g \mod p^{r} \in B_{\scP}^{+}(\ZZ/p^r)$.
    \item We define the \emph{principal level subgroup} $K_{p,r}^{1}$ as the subgroup of $K_{p, r}^{0}$ consisting of those $g$ projecting under the surjection \eqref{eq:Bplusproj} to an element in $(\ZZ/p^r)^{\times} \times N_{\scP}^{+}(\ZZ/p^r)$.
\end{itemize}
Then it follows that $K_{p,r}^{0}/K_{p,r}^{1} \simeq T_{\scP}(\ZZ/p^r)$.

\subsection{Moduli problems}
In this subsection, we start with an integral PEL-type Shimura datum $(\scP = (B, \ast, V, \psi_{V}, h), \calO_B, L)$ and omit it from the following notations.

Let $\square$ be a finite set of primes, which is usually taken to be the empty set or $\{p\}$. Let $S$ be a locally Noetherian connected $\calO_E \otimes \ZZ_{(\square)}$-scheme.

\begin{definition}
We say that a tuple $\ulA = (A, \lambda, \iota)$ is a \emph{$\ZZ_{(\square)}$-polarized abelian scheme with an action of $\calO_{B}$} if:
\begin{itemize}
    \item $A$ is an abelian scheme over $S$.
    \item $\lambda$ is a prime-to-$\square$ polarization of $A$ over $S$.
    \item $\iota: \calO_{B} \hookrightarrow \End_{S} A \otimes_{\ZZ} \ZZ_{(\square)}$ is a homomorphism which respects the involutions on both sides: the involution $\ast$ on the left and the Rosati involution coming from $\lambda$ on the right.
\end{itemize}

Let $K$ be a compact open subgroup of $G(\AA_{\QQ, \rmf})$. We define a \emph{$K^{\square}$-level structure} $\bareta^{(\square)} = \eta^{(\square)}K$ of $\ulA$ to be a $\pi_1(S, \bars)$-invariant $K$-orbit of the isomorphism of $\calO_{B} \otimes \AA_{\QQ, \rmf}^{\square}$-modules
    $$
    \eta^{(\square)}: L \otimes \AA_{\QQ, \rmf}^{\square} \xrightarrow{\sim} \rmH_1(A_{\bars}, \AA_{\QQ, \rmf}^{\square}),
    $$
    which identifies $\psi_V$ with an $\AA_{\QQ, \rmf}^{\square, \times}$-multiple of the symplectic pairing on the Tate module $\rmH_1(A_{\bars}, \AA_{\QQ, \rmf}^{\square})$ defined by $\lambda$ and the Weil pairing.
\end{definition}

\begin{definition}
Let $\frM_{K}^{\square}$ be the following category fibered in groupoids over the category of $\calO_E \otimes_{\ZZ} \ZZ_{(\square)}$-schemes:
\begin{itemize}
    \item The objects over a scheme $S$ are quadruples $(\ulA, \bareta^{\square})$, where $\ulA$ is a $\ZZ_{(\square)}$-polarized abelian scheme with an action of $\calO_{B}$, and $\bareta^{(\square)}$ is a $K^{\square}$-level structure, such that $\Lie_{S} A$ satisfies the Kottwitz determinant condition defined by $(L \otimes \RR, \psi_V, h)$ \footnote{The requirements on the dimension of the abelian scheme $A$ and the information on the signature of the Hermitian spaces $V_i$ in the case of unitary Shimura datums are encoded in this Kottwitz determinant condition. So here the definition coincides with our usual definition, for example, in \cite[Section 2.1]{hsieh2014eisenstein}.}.
    \item The morphisms from $(\ulA, \bareta^{\square})$ to $(\ulA^{\prime}, \bareta^{\square, \prime})$ are given by a $\ZZ_{(\square)}$-isogeny $f: A \rightarrow A^{\prime}$ that is compatible with the action of $\calO_{B}$ and the level structures.
\end{itemize}
\end{definition}

We have the following well-known representability theorems in the case when $\square$ is the empty set or $\{p\}$.

\begin{theorem}[Deligne--Shimura]
When $\square$ is the empty set and $K$ is neat, $\frM_{K}^{\emptyset}$ is representable by a scheme $\calS_{G}(K)$. Moreover, in the case when $\scP$ is a unitary Shimura datum, we have
$$
\calS_{G}(K) = \bigsqcup_{G^{\prime} \in \ker^{1}(\QQ, G)} \Sh_{K}(G, X),
$$
where $\ker^{1}(\QQ, G)$ is the set of locally trivial elements of $\rmH^{1}(\QQ, G)$.
\end{theorem}

More precisely, the elements of $\ker^{1}(\QQ,G)$ classify isomorphism classes of Hermitian tuples $(V_{j}^{\prime}, \lrangle{-,-}_{V_j^{\prime}})_{1 \leq j \leq m}$ that are locally isomorphic to $(V_{j}, \lrangle{-,-}_{V_j})_{1 \leq j \leq m}$. Then $\calS_{G}(K)$ is a disjoint union of copies of $\Sh_{K}(G,X)$. This will not cause much difficulty for applications to automorphic forms, as we only need one copy of them. See \cite[Section 2.3]{eischen2020padic} for details.

\begin{theorem}[Kottwitz]
When $\square = \{p\}$, the category $\frM_{K}^{p}$ is a smooth Deligne--Mumford stack. When $K^{p}$ is neat, it is representable by a quasi-projective scheme $\calS_{G}(K^p)$. Moreover, in the case when $\scP$ is a unitary Shimura datum, we have
$$
\calS_{G}(K^p) \times_{\calO_{E} \otimes_{\ZZ} \ZZ_{(p)}} \Spec E \xrightarrow{\sim} \calS_G(K^p K_p^{0}).
$$
\end{theorem}

In the case when $\scP$ is a simple unitary Shimura datum given by a Hermitian $\calK$-vector space $V$, by abuse of notation, we call $\calS_{G}(K)$ or $\calS_{G}(K^p)$ the \emph{unitary Shimura variety attached to $\GU(V)$}. The geometric theory of modular forms and their $p$-adic theory over $\GU(V)$ are built upon $\calS_{G}(K)$ and $\calS_{G}(K^p)$.

In the following sections, we concentrate on the case of simple unitary PEL datums.

\subsection{Compactifications}
The theory of toroidal compactification of unitary Shimura varieties is developed in \cite{lan2013arithmetic}. Fixing a certain smooth projective polyhedral cone decomposition (which we do not make precise here), one can attach the \emph{toroidal compactification} $\calS_{G}^{\tor}(K^{\square})$ to $\calS_{G}(K^{\square})$. We focus on the case $\square = \{p\}$. Then we know:
\begin{itemize}
    \item The toroidal compactification $\calS_{G}^{\tor}(K^p)$ contains $\calS_{G}(K^p)$ as an open dense subscheme. The complement of $\calS_{G}(K^p)$ is a relative Cartier divisor with normal crossings. We denote by $\scI_{\calS_{G}^{\tor}(K^p)}$ the ideal sheaf of the boundary of $\calS_{G}^{\tor}(K^p)$.
    \item There is a quadruple $\underline{\scG} := (\scG, \lambda, \iota, \eta)$ over $\calS_{G}^{\tor}(K^p)$, where $\scG$ is a semi-abelian scheme with an $\calO_{\calK}$-action $\iota$ and a homomorphism $\lambda: \scG \rightarrow \scG^{\vee}$, such that $\underline{\scG}|_{\calS_{G}(K^p)} = \underline{\scA}$, the universal quadruple over $\calS_{G}(K^p)$, and $\eta$ is the level structure in the quadruple $\underline{\scA}$.
\end{itemize}

Let $\ulome := e^{\ast}\Omega_{\scG/\calS_{G}^{\tor}(K^p)}^{1}$, where $e: \calS_{G}^{\tor}(K^p) \rightarrow \scG$ is the zero section of the semi-abelian scheme $\scG$ over $\calS_{G}^{\tor}(K^p)$. Then $\ulome$ is a locally free coherent $\scO_{\calS_{G}^{\tor}(K^p)}$-module. The \emph{minimal compactification} of $\calS_{G}(K^p)$ is defined to be
$$
\calS^{\min}_{G}(K^p) := \Proj \left( \bigoplus_{k=0}^{\infty} \rmH^{0} (\calS_{G}^{\tor}(K^p), (\det \ulome)^{k}) \right).
$$
Let $\pi: \calS_{G}^{\tor}(K^p) \rightarrow \calS^{\min}_{G}(K^p)$ be the natural projection.

\subsection{Igusa schemes}
In the following, by abuse of notation, we also denote $\calS_{G}^{\tor}(K^p)$ and $\calS_{G}^{\min}(K^p)$ by their base change to $\ZZ_p$ via the map $\calO_{\calK, (p)} \rightarrow \ZZ_p$ induced by our fixed embedding $\iota_p$, and let $\calS_{G /\FF_p}^{\tor}(K^p)$ and $\calS_{G /\FF_p}^{\min}(K^p)$ be their corresponding special fibers.

\subsubsection{Hasse invariants}
Let $\Ha \in \rmH^{0}(\calS_{G, /\FF_p}^{\tor}(K^p), (\det\ulome)^{p-1})$ be the Hasse invariant defined in \cite[Section 6.3.1]{lan2013arithmetic}. In particular, for each geometric point $\bars$ of $\calS_{G /\FF_p}^{\tor}(K^p)$, the Hasse invariant of the corresponding semi-abelian scheme $\scG_{\bars}$ is non-zero if and only if the abelian part of $\scG_{\bars}$ is ordinary. Because $\pi_{\ast} \ulome$ is ample, for some $t_E > 0$, there exists an element in $\rmH^{0}(\calS_{G}^{\min}(K^p), (\pi_{\ast} \ulome)^{t_E(p-1)})$ lifting the $t_E$-th power of the push-forward of $\Ha$. We denote by $E$ the pullback under $\pi$ of any such lift, which (because $\pi^{\ast} \pi_{\ast} \ulome \simeq \ulome$) defines an element $E \in \rmH^{0}(\calS_{G}^{\tor}(K^p), (\det\ulome)^{t_E (p-1)})$.

\subsubsection{Ordinary locus and the Igusa tower}
In the Kottwitz model $\calS_{G}(K^p)$ of the unitary Shimura variety attached to $G$, the level structure at $p$ is missing. We consider the following moduli problem to add level structures at $p$ and formulate the Igusa tower attached to $G$.
\begin{definition}
Let $\frI_{K^p K_{p,n}^{1}}$ be the following category fibered in groupoids over the category of $(\calO_E \otimes_{\ZZ} \ZZ_{(p)})$-schemes:
\begin{itemize}
    \item The objects over a scheme $S$ are quadruples $(\ulA, \bareta^{p}, j_n)$, where $(\ulA, \bareta^{p}) \in \frM^{p}(S)$ and $j_n$ is a \emph{level $K_{p,n}^{1}$-structure}, defined as a $K_{p,n}^{1}$-orbit of monomorphisms of $\calO_{\calK}$-schemes over $S$:
    $$
    j_{p^n} : \bmu_{p^n} \otimes_{\ZZ} L^{+} \hookrightarrow A[p^n].
    $$
    \item The morphisms from $(\ulA, \bareta^{p}, j_n)$ to $(\ulA^{\prime}, \bareta^{p, \prime}, j_n^{\prime})$ are given by an element 
    $$
    f \in \Hom_{\frM^{p}}((\ulA, \bareta^{p}), (\ulA^{\prime}, \bareta^{p, \prime}))
    $$ 
    that is compatible with the level-$p^n$ structures of $A$ and $A^{\prime}$.
\end{itemize}
\end{definition}

\begin{theorem}
The moduli problem $\frI_{K^p K_{p,n}^{1}}$ is relatively representable over $\calS_{G}(K^p)$, and thus it is represented by a scheme $\calI_G(K^p K_{p,n}^{1})$ over the scheme $\calS_{G}(K^p)$. The scheme $\calI_G(K^p K_{p,n}^{1})$ is called the \emph{ordinary locus of level $K_{p,n}^{1}$}.
\end{theorem}

We sometimes write $\calI_G(K^p K_{p,n}^{1})$ simply as $\calI_G(K_{p,n}^{1})$ or even $\calI_n$ for short if it does not cause any confusion.

Let $\calI_{G}^{\tor}(K^p K_{p,n}^{1})$ be the partial toroidal compactification of the ordinary locus $\calI_G(K^p K_{p,n}^{1})$ (\cite[Theorem 5.2.1.1]{lan2018compactification}). It is obtained by gluing the toroidal boundary charts parameterizing degenerating families, defined in \cite[Definition 3.4.2.0]{lan2018compactification}, to $\calI_{G}(K^p K_{p,n}^{1})$.

Let $\calI_{G, m}^{\tor}(K^p K_{p,n}^{1})$ be the base change of $\calI_{G}^{\tor}(K^p K_{p,n}^{1})$ to $\ZZ/p^m \ZZ$. By \cite[Lemma 6.3.2.7]{lan2018compactification}, $\calS_{G, m}^{\tor}(K^p)[1/E]$ agrees with the ordinary locus in \cite[Theorem 5.2.1.1]{lan2018compactification} for the hyperspecial level at $p$, and by \cite[Corollary 5.2.2.3]{lan2018compactification}, the map 
$$
\pi_{m,n}: \calI_{G, m}^{\tor}(K^p K_{p,n}^{1}) \rightarrow \calS_{G, m}^{\tor}(K^p)[1/E]
$$
forgetting the $K_{p,n}^{1}$-level structure is finite étale.

We gather the above objects into the following diagram, referred to as the diagram of \emph{Igusa towers}.
\begin{equation} \label{tikz:igusa}
  \begin{tikzcd}
	& \vdots && \vdots \\
	\cdots & {\calI_{G, m+1}^{\tor}(K^pK_{p,n+1}^{1})} && {\calI_{G, m}^{\tor}(K^pK_{p,n+1}^{1})} \\
	\cdots & {\calI_{G, m+1}^{\tor}(K^pK_{p,n}^{1})} && {\calI_{G, m}^{\tor}(K^pK_{p,n}^{1})} \\
	\cdots & {\calS_{G, m+1}^{\tor}[1/E]} && {\calS_{G, m}^{\tor}[1/E]} && {\calS_{G}^{\tor}[1/E]} \\
	\cdots & {\Spec \ZZ/p^{m+1}} && {\Spec \ZZ/p^m} && {\Spec \ZZ_p}
	\arrow[from=1-2, to=2-2]
	\arrow[from=1-4, to=2-4]
	\arrow[from=2-2, to=2-1]
	\arrow[from=2-2, to=3-2]
	\arrow[from=2-4, to=2-2]
	\arrow[from=2-4, to=3-4]
	\arrow[from=3-2, to=3-1]
	\arrow["{\pi_{m+1,n}}"', from=3-2, to=4-2]
	\arrow[from=3-4, to=3-2]
	\arrow["{\pi_{m,n}}"', from=3-4, to=4-4]
	\arrow[from=4-2, to=4-1]
	\arrow[from=4-2, to=5-2]
	\arrow[from=4-4, to=4-2]
	\arrow[from=4-4, to=4-6]
	\arrow[from=4-4, to=5-4]
	\arrow[from=4-6, to=5-6]
	\arrow[from=5-2, to=5-1]
	\arrow["{\mod p}"', from=5-4, to=5-2]
	\arrow[from=5-4, to=5-6]
\end{tikzcd} .
\end{equation}
To simplify the notation, when the level group $K^p$ away from $p$ is clear or fixed in the context, we write $\calI_{m,n}^{\tor}$ for $\calI_{G, m}^{\tor}(K^pK_{p,n}^{1})$.

\subsection{Modular forms on unitary groups}
\subsubsection{$p$-adic modular forms} \label{sec:padicmodularforms}
We define the space of \emph{mod $p^m$ automorphic forms} on $G$ of level $n$ by
$$
V_{G, m,n} := \rmH^{0}(\calI_{m,n}^{\tor}, \scO_{\calI_{m,n}^{\tor}}).
$$
We let $\scI_{m,n}^{\tor} := \pi_{m,n}^{\ast} \scI_{\calS_m^{\tor}}$, and similarly define the space of \emph{mod $p^m$ cuspidal automorphic forms} on $G$ of level $n$ by
$$
V^{0}_{G, m,n} := \rmH^{0}(\calI_{m,n}^{\tor}, \scI_{m,n}^{\tor}).
$$
Then we define $p$-adic automorphic forms by passing to the limit.
\begin{definition} \label{defn:padicautomorphicform}
We define the space of \emph{$p$-adic automorphic forms (resp. cuspidal automorphic forms) with torsion coefficients} as
$$
\calV_G := \varinjlim_{m} \varinjlim_{n} V_{G, m,n} \quad (\text{resp. } \calV^{0}_G := \varinjlim_{m} \varinjlim_{n} V_{G, m,n}^0).
$$
We define the space of \emph{$p$-adic automorphic forms (resp. cuspidal automorphic forms) with integral coefficients} as
$$
V_G := \varprojlim_{m} \varinjlim_{n} V_{G, m,n} \quad (\text{resp. } V^{0}_G := \varprojlim_{m} \varinjlim_{n} V_{G, m,n}^0).
$$
For any $p$-adic ring $R$, i.e., $R$ satisfies $R \simeq \varprojlim_{n} R/p^n R$, we have the \emph{$p$-adic automorphic forms (resp. cuspidal automorphic forms)} with coefficient ring $R$ defined by base change to $R$, denoted by $V_G(R)$ and $V_G^{0}(R)$ respectively.
\end{definition}

Recall that we previously defined several algebraic groups $H_{\scP}$, $B_{\scP}$, and $T_{\scP}$ attached to a unitary Shimura datum.

\begin{definition}
For any sufficiently large finite extension $L/\QQ_p$, we define the \emph{weight algebra} with coefficient field $L$ of the unitary Shimura datum $\scP$ as the completed group algebra $\Lambda_{\scP} := \calO_{L} \lrbracket{T_{\scP}(\ZZ_p)}$.
\end{definition}

The space of $p$-adic automorphic forms carries several actions:
\begin{itemize}
    \item Recall that $T_{\scP}(\ZZ/p^n)$ naturally identifies with $K_{p,n}^{0}/K_{p,n}^{1}$. Hence, the group $T_{\scP}(\ZZ_p)$ acts on $\calV_G, \calV^0_G, V_{G}(R)$, and $V_{G}^0(R)$, making these spaces into $\Lambda_{\scP}$-modules.
    \item The action of $G(\AA_{\QQ, \rmf}^{p})$ on the Igusa tower gives an action of $G(\AA_{\QQ, \rmf}^{p})$ on the space of $p$-adic automorphic forms. Let $K^{p}$ be an open compact subgroup of $G(\AA_{\QQ,\rmf}^{p})$. We denote the submodules fixed by $K^{p}$ by $V_G(K^{p}, R)$, $V_G^{0}(K^{p}, R)$, and so on. These $p$-adic automorphic forms are said to be of \emph{tame level} $K^p$.
\end{itemize}

\begin{definition} \label{defn:padicweight}
We define a \emph{$p$-adic weight} $\tau$ as a $\barQQ_p$-valued character of $\Lambda_{\scP}$. In other words, under the isomorphism $T_{\scP}(\ZZ_p) \simeq \prod_{v \in \scV_{\calF}^{(p)}} T_{r_v + s_v}(\calO_{\calF,v})$, a $p$-adic weight $\tau$ is a collection $\{\tau_v\}_{v \in \scV_{\calF}^{(p)}}$ with each $\tau_v$ given by
\begin{align*}
\tau_v: T_{r_v + s_v}(\calO_{\calF, v}) &\rightarrow \barQQ_p^{\times}, \\
\diag[a_1, \ldots, a_{r_v+s_v}] &\mapsto \tau_{1,v}^{+}(a_1) \cdots \tau_{r_v, v}^{+}(a_{r_v}) \tau_{1,v}^{-}(a_{r_v + 1})\cdots\tau_{s_v, v}^{-}(a_{r_v+s_v}),    
\end{align*}
where $\tau_{i,v}^{\pm}: \calO_{\calF,v}^{\times} \rightarrow \barQQ_p^{\times}$ are continuous characters. We say a $p$-adic weight $\tau$ is \emph{arithmetic} if, for each $v \in \scV_{\calF}^{(p)}$, $\tau_{i,v}^{\pm}: \calO_{\calF, v}^{\times} \rightarrow \barQQ_p^{\times}$ is given by a product of an algebraic character and a finite order character, i.e., $\tau_{i,v}^{\pm}(a) = \epsilon^{\pm}_{i, \tau, v}(a) a^{t_{i,\tau,v}^{\pm}}$ for finite order characters $\epsilon_{i,\tau,v}^{\pm}: \calO_{\calF, v}^{\times} \rightarrow \barQQ_p^{\times}$ and $t_{i,\tau, v}^{\pm} \in \ZZ$.
\end{definition}

Given a $p$-adic weight $\tau$, we denote by $\QQ_p(\tau)$ the finite field extension of $\QQ_p$ generated by the image of $\tau$. We put $V_{m,n}[\tau]$ as the subspace of $V_{m,n} \otimes_{\ZZ_p} \calO_{\QQ_p(\tau)}$ on which $\Lambda_{\scP}$ acts by the inverse of the character $\tau$. Similarly, we define the spaces $V_{m,n}^{0}[\tau]$, $\calV[\tau]$, and $\calV^{0}[\tau]$. These $p$-adic automorphic forms are said to have \emph{weight} $\tau$.

\subsubsection{Classical automorphic forms}
Though we are free to use the space of $p$-adic modular forms, which is larger than the space of classical automorphic forms, automorphic computations are more frequently performed in the classical way. In this subsection, we consider the classical automorphic forms over unitary groups.

Let $\ult := (t_1^{+}, \ldots, t_r^{+}, t_{1}^{-}, \ldots, t_{s}^{-}) \in \ZZ^{r+s}$ be any $(r+s)$-tuple. We define, for any algebra $R$, the space
\[
W_{\ult}(R) := \{f \in R[\GL_r \times \GL_s] : f(tn_{+}g) = k^{-1}(t) f(g), t \in T_r \times T_s, n_{+} \in U_r \times U_{s}^{-} \},
\]
where $R[\GL_r \times \GL_s]$ denotes the polynomial functions on $\GL_r \times \GL_s$ with coefficients in $R$, and $k$ is regarded as an algebraic character on $T_r \times T_s$ defined by
$$
k (\diag[a_1, \ldots, a_r], \diag[a_{r+1}, \ldots, a_{r+s}]) = a_1^{t_{1}^{+}} \cdots a_r^{t_{r}^{+}} a_{r+1}^{t_{1}^{-}} \cdots a_{r+s}^{t_{s}^{-}}.
$$
Then $W_{\ult}$ is a free $R$-module and is the algebraic representation of $\GL_{r}(R) \times \GL_{s}(R)$ with minimal weight $-k$ with respect to $U_r \times U_{s}^{-}$. We regard it as
$$
W_{\ult} = W_{\ult}^{+} \boxtimes W_{\ult}^{-}
$$
as an algebraic representation of $\GL_r \times \GL_s$.

Recall that $\ulome := e^{\ast}\Omega_{\scG/\calS_{G}^{\tor}(K^p)}^{1}$ has the decomposition $\ulome = \ulome^{+} \oplus \ulome^{-}$, where $\ulome^{+}$ (resp. $\ulome^{-}$) is the subsheaf of $\ulome$ on which $\iota(b)$ acts by $b$ (resp. $\barb$) for all $b \in \calO_{\calK}$. Because $p$ is unramified in $\calK$, $\ulome^{+}$ (resp. $\ulome^{-}$) is locally free of rank $r$ (resp. $s$). Set
\[
\omega_{\ult}^{+} = \uIsom_{\calS^{\tor}_{G}}(\scO_{\calS^{\tor}_{G}}^{\oplus r}, \ulome^{+}) \times^{\GL_{r}} W_{\ult}^{+}, \quad \omega_{\ult}^{-} = \uIsom_{\calS^{\tor}_{G}}(\scO_{\calS^{\tor}_{G}}^{\oplus s}, \ulome^{-}) \times^{\GL_{s}} W_{\ult}^{-},
\]
and put $\omega_{\ult} = \omega_{\ult}^{+} \otimes \omega_{\ult}^{-}$. 

Suppose for each $v \in \scV_{\calF}^{(p)}$, we are given an $(r_v + s_v)$-tuple $\ult_{v}$ as above and, by abuse of notation, let $\ult = \{\ult_v\}_{v \in \scV_{\calF}^{(p)}}$. Then, following the above process, we have a well-defined $\omega_{\ult}$ provided $\ult$ is a parallel weight.

Let $v \in \scV_{\calF}^{(p)}$ and let $\uleps_v := (\epsilon_{1,v}^{+}, \ldots, \epsilon_{r,v}^{+}, \epsilon_{1,v}^{-}, \ldots, \epsilon_{s,v}^{-})$ be an $(r+s)$-tuple of finite order characters $\calO_{\calF_v}^{\times} \rightarrow \barQQ_p^{\times}$, and let $\uleps = \{\uleps_{v}\}_{v \in \scV_{\calF}^{(p)}}$. Let $F/\QQ_p$ be a finite extension containing the values of all $\epsilon_{i}^{\pm}$.

\begin{definition}
We define the space of \emph{classical automorphic forms} on $G$ of \emph{weight} $\ult$, \emph{level} $K^{p}K_{p,n}^{1}$, and \emph{nebentypus} $\uleps$ as elements in the space
\[
M_{\ult}(K^p K_{p,n}^{1}; \uleps; F) := (\rmH^{0}(\calS_{G}^{\tor}(K^p K_{p,n}^{1}), \omega_{\ult}) \otimes_{\calO_{\calK}} F)[\uleps].
\]
Similarly, we have the space of \emph{classical cuspidal automorphic forms}
\[
M_{\ult}^{0}(K^p K_{p,n}^{1}; \uleps; F) := (\rmH^{0}(\calS_{G}^{\tor}(K^p K_{p,n}^{1}), \omega_{\ult} \otimes \scI_{\calS_{G}^{\tor}(K^p)}) \otimes_{\calO_{\calK}} F)[\uleps].
\]
\end{definition}
One has embeddings (see, for example, \cite[Section 3.3]{hsieh2014eisenstein}):
\begin{subequations} \label{eq:embedclassical}
\begin{equation}
    M_{\ult}(K^p K_{p,n}^{1}; \uleps; F) \hookrightarrow V_{G}(K, F)[\tau_{\ult, \uleps}]
\end{equation}
and
\begin{equation} 
M_{\ult}^{0}(K^p K_{p,n}^{1}; \uleps; F) \hookrightarrow V_{G}^{0}(K, F)[\tau_{\ult, \uleps}],    
\end{equation}
\end{subequations}
where $\tau_{\ult, \uleps}$ is the $p$-adic weight given by $\ult$ and $\uleps$ following the notation in Definition \ref{defn:padicweight}.

\subsection{Hecke operators on modular forms} \label{sec:Heckeoperators}
One can refer to \cite[Section 3.7]{hsieh2014eisenstein} for carefully presented details; here we provide only a brief sketch. Let $K_{j}^{p}$ be open compact subgroups of $G(\AA_{\QQ, \rmf}^{p})$ for $j = 1, 2$, such that $K_j^{p} \cdot G(\ZZ_p)$ are neat.

\subsubsection{Hecke operators away from $p$}
Let $g \in G(\AA_{\QQ, \rmf}^{p})$. We define the double coset operator
\[
[K_2^p g K_1^p]: V_{G}(K_1^p)[\tau] \rightarrow  V_{G}(K_2^p)[\tau] 
\]
through the action of $G(\AA_{\QQ, \rmf}^{p})$ on the space of modular forms:
\[
[K_2^p g K_1^p] f = \sum_{g_j} [g_j]^{\ast} f, \quad \text{where } K_2^p g K_1^p = \bigsqcup_{g_j} g_j K_1^p.
\]
Similarly, for neat open compact subgroups $K_j^{p} \cdot K_{p,r}^{1}$, we can define the corresponding double coset operators in the same way and denote them by $[K_{2, r}^p g K_{1,r}^p]$. When $K_1^p = K_2^p$ is understood, we write $T(g)$ instead of $[K_1^p g K_1^p]$ and $T_r(g)$ instead of $[K_{1,r}^p g K_{1,r}^p]$.

\subsubsection{Hida's operators at $p$}
We define certain elements $t_{j}$ in $B_{\scP}^{+}(\QQ_p)$ as $t_j := (1, (t_{w,j}))$ under the isomorphism \eqref{eq:BpscP}, for $1 \leq j \leq r+s$, by:
\[
t_{w,j} := \begin{cases}
\diag[p \bfone_{j}, \bfone_{r+s-j}], &\quad j \leq a_w \\
\diag[p \bfone_{a_w}, \bfone_{r+s-j}, p \bfone_{j-a_w}], &\quad j > a_w
\end{cases}.
\]
In \cite[Section 8.3.1]{hida2004padic}, Hida defined an action of the double cosets $u_{w,j} = B_{\scP}^{-}(\ZZ_p) t_{w,j} B_{\scP}(\ZZ_p)$ on the modules of $p$-adic modular forms and cusp forms via correspondences on the Igusa tower. For any $v \in \scV_{\calF}^{(p)}$, we put:
$$
\UU_v^{\prime} = \prod_{w \in \scV_{\calK}^{(v)}} \prod_{j=1}^{r+s} u_{w,j}, \quad \text{and } \UU_p^{\prime} = \prod_{v \in \scV_{\calF}^{(p)}} \UU_v^{\prime}.
$$
\begin{definition}
We define a projector $\ee := \lim_{n \to \infty} (\UU_p^{\prime})^{n!}$, called \emph{Hida's ordinary projector}. For any $p$-adic ring $R$, we define the submodule of \emph{ordinary $p$-adic automorphic forms} (resp. \emph{ordinary $p$-adic cuspidal automorphic forms}) over $R$ as
\[
V_{G}^{\ord}(K^p, R) := \ee V_{G}(K^p, R) \quad (\text{resp. } V_{G}^{0,\ord}(K^p, R) := \ee V_{G}^{0}(K^p, R) ).
\]
\end{definition}

We can further define ordinary $p$-adic automorphic forms of weight $\tau$ for a $p$-adic weight $\tau$ by adding decorations to the notation, which we shall not list here.

\subsection{Interlude: review on $p$-adic measures}
Next, we shall introduce $p$-adic families of ($p$-adic) automorphic forms. We prefer the setup of $p$-adic measures among various different but essentially equivalent perspectives. Here we briefly review the basic notions of $p$-adic measures.
\begin{itemize}
    \item Let $R$ be a $p$-adic ring. This will act as the base coefficient ring.
    \item Let $M$ be a $p$-adically complete $R$-module. 
    \item Let $Y$ be a compact abelian group with a totally disconnected topology; i.e., $Y$ is a profinite abelian group.
\end{itemize}
Then we denote by $\scC(Y, R)$ the $R$-algebra of continuous $R$-valued functions on $Y$, equipped with the topology of uniform convergence. An $M$-valued \emph{$p$-adic measure} on $Y$ is a continuous $R$-linear map
$$
\mu: \scC(Y,R) \rightarrow M, \quad f \mapsto \mu(f):= \int_{Y} f(y) \,\mathrm{d}\mu(y).
$$
The set of $M$-valued $p$-adic measures on $Y$ is a $p$-adically complete $R$-module and is denoted by $\Meas_{R}(Y,M)$. We often omit the coefficient ring $R$ in the notation. 

Here are some operations on $p$-adic measures:
\begin{itemize}
    \item \emph{Base change}: Let $R^{\prime}$ be an $R$-algebra which is also $p$-adically complete. Since $\scC(Y, R^{\prime}) \simeq \scC(Y, R) \hotimes_R R^{\prime}$, there is a natural map $\Meas_{R}(Y,M) \rightarrow \Meas_{R^{\prime}}(Y, M \hotimes_R R^{\prime})$. If the structure map $R \rightarrow R^{\prime}$ is injective, we may view $\Meas_{R}(Y,M)$ as a subset of $\Meas_{R^{\prime}}(Y, M \hotimes_R R^{\prime})$.
    \item \emph{Dirac measures}: Given $y \in Y$, we define $\delta_{y}: f \mapsto f(y)$ as a $p$-adic measure in $\Meas(Y,R)$.
    \item \emph{Action of continuous functions}: Let $h \in \scC(Y, R)$. For $\mu \in \Meas(Y,M)$, we define 
    $$
    h \ast \mu: f \mapsto \int_{Y} f(y)h(y) \,\mathrm{d}\mu(y).
    $$
    \item \emph{Convolution}: If $Y$ is equipped with the structure of an abelian group (written multiplicatively), then we define the convolution on $\Meas(Y,R)$ as
    $$
    \mu_1 \ast \mu_2: f \mapsto \int_Y \int_Y f(yz) \,\mathrm{d}\mu_1(y) \,\mathrm{d}\mu_2(z).
    $$
    If $f: Y \rightarrow R^{\times}$ is a continuous multiplicative character, then we have
    $$
    \int_{Y} f \,\mathrm{d}(\mu_1 \ast \mu_2) = \left( \int_Y f \,\mathrm{d}\mu_1 \right) \left( \int_Y f \,\mathrm{d}\mu_2 \right).
    $$
    We can define the convolution of measures in $\Meas(Y,M)$ in the same way whenever there is an appropriate "product" on $M$.
    \item \emph{Product of the test space}: If $Y = Y_1 \times Y_2$ is a product of profinite abelian groups, then there is a natural isomorphism
    \begin{equation} \label{eq:productmeasure}
    \Meas(Y, M) \simeq \Meas(Y_1, \Meas(Y_2, M)).
    \end{equation}
\end{itemize}
Moreover, we note that $\Meas(Y, R)$ can be identified with the completed group algebra $R\lrbracket{Y}$. In practice, we often encounter the following cases:
\begin{itemize}
    \item $Y$ is a product of a finite abelian group $\Delta$ with several copies of $\ZZ_p$. For example, $Y$ could be the Galois group $\Gamma_{\calK}$, the group of units $U_{\calK, p}^{(r)} := 1 + p^{r} \calO_{\calK,p}$ for $r \geq 1$, or the torus $T_{n}(\ZZ_p)$ for some positive integer $n$.
    \item $R$ is often the ring of integers $\calO_L$ for a finite extension $L/\QQ_p$.
    \item $M$ is often taken to be $R$ or the space of $p$-adic modular forms over certain reductive groups.
\end{itemize}

\subsection{Hida families on unitary groups}
Recall that $\scP$ is a simple unitary PEL datum given by a Hermitian space $V$. Let $R$ be a $p$-adic ring.
\begin{definition} \label{defn:Hidafam}
We define a \emph{$p$-adic family of automorphic forms} on $G$ of level $K^{p}$ over $R$ to be a $p$-adic measure $\bff \in \Meas_R(T_{\scP}(\ZZ_p), V_{G}(K^p, R))$ such that 
\begin{equation} \label{eq:Hidafamilycompact}
\int_{T_{\scP}(\ZZ_p)} t \cdot \phi(y) \,\mathrm{d}\bff(y) = t \cdot \int_{T_{\scP}(\ZZ_p)} \phi(y) \,\mathrm{d}\bff(y)    
\end{equation}
for any $\phi \in \scC(T_{\scP}(\ZZ_p), R)$ and $t \in T_{\scP}(\ZZ_p)$. A \emph{$p$-adic family of cuspidal automorphic forms} (resp. \emph{Hida family}, \emph{cuspidal Hida family}) on $G$ of level $K^p$ over $R$ is a $p$-adic family of automorphic forms taking values in 
\[
V_{G}^{0}(K^p, R) \quad (\text{resp. } V_{G}^{\ord}(K^p, R), \quad V_{G}^{0, \ord}(K^p, R)).
\]
We denote the space of such families by $\calM^{\bullet}_{G}(K^{p}, R)$ with decorations $\bullet \in \{\emptyset, 0, \ord\}$.
\end{definition}

In particular, let $\tau$ be a $p$-adic weight such that $R$ contains all values of $\tau$; then for $\bff \in \calM_{G}(K^p, R)$, the value
\[
\bff_{\tau} := \int_{T_{\scP}(\ZZ_p)} \tau \,\mathrm{d}\bff \in V_{G}(K^p, R)[\tau]
\]
is called the \emph{specialization} of $\bff$ at $\tau$. We remark that $\calM^{\bullet}_{G}(K^{p}, R)$ becomes a $\Lambda_{\scP}$-module by the action of $T_{\scP}(\ZZ_p)$ induced by its action on $V_G(K^p, R)$. This is well-defined by the requirement \eqref{eq:Hidafamilycompact}.

Recall that we have a decomposition $\ZZ_p^{\times} = \bmu(\ZZ_p) \times (1+p\ZZ_p)$ (assuming $p \neq 2$), where $\bmu(\ZZ_p)$ is the group of roots of unity in $\ZZ_p$. This induces a decomposition of $T_{\scP}(\ZZ_p)$ as
$$
T_{\scP}(\ZZ_p) = \Delta \times T^{\circ}_{\scP},
$$
where $\Delta = T_{\scP}(\bmu(\ZZ_p))$ is the torsion subgroup and $T^{\circ}_{\scP} := T_{\scP}(1+p\ZZ_p)$ is the identity component. With respect to this decomposition, we rewrite the weight algebra as $\Lambda_{\scP} \simeq \calO_{L}[\Delta]\lrbracket{T_{\scP}^{\circ}}$. Since $|\Delta|$ is prime to $p$, the space $\calM^{\bullet}_{G}(K^{p}, R)$ decomposes into a direct sum of isotypic pieces for the $\calO_L$-characters $\eta$ of $\Delta$:
\[
\calM^{\bullet}_{G}(K^{p}, R) = \bigoplus_{\eta \in \widehat{\Delta}} \calM^{\bullet}_{G}(K^{p}, R, \eta).
\]
The characters $\eta$ are called \emph{branch characters}. Each $\calM^{\bullet}_{G}(K^{p}, R, \eta)$ is a $\Lambda_{\scP}^{\circ} := \calO_{L}\lrbracket{T_{\scP}^{\circ}}$-module.

Let $\II$ be a normal domain which is also a finite algebra over $\Lambda_{\scP}$ (resp. $\Lambda_{\scP}^{\circ}$). We define the space of \emph{$\II$-adic Hida families} by the base change of $\calM^{\bullet}_{G}(K^{p}, \calO_{L})$ (resp. $\calM^{\bullet}_{G}(K^{p}, \calO_{L}, \eta)$) to $\II$, forming the spaces $\calM^{\bullet}_{G}(K^{p}, \II)$ and $\calM^{\bullet}_{G}(K^{p}, \eta, \II)$.

Hecke operators on the space of $p$-adic families of modular forms are induced by those defined in Section \ref{sec:Heckeoperators}. Let $\scS$ be a finite set of places of $\calF$ such that $K^{p}$ is maximal outside $\scS$. For any finite torsion-free $\Lambda_{\scP}^{\circ}$-algebra $A$, we define the \emph{unramified ordinary cuspidal Hecke algebra} $\bfh^{0, \ord, \scS}(A)$ as the $A$-subalgebra of $\End_A(\calM^{0, \ord}_{G}(K^{p}, A))$ generated by Hecke operators away from $\scV_{\calF}^{(p)} \cup \scS$ and the $\bfU^{\prime}_v$-operators for $v \in \scV_{\calF}^{(p)}$.

%\begin{remark}[Other $p$-adic families of automorphic forms]
%Hida families defined in Definition \ref{defn:Hidafam} are possibly the easiest classes of $p$-adic families of automorphic forms. There are ways of generalizing it.
%\begin{itemize}
%    \item \emph{Allowing general parabolics}. Here Hida's ordinary projector $\ee$ is constructed via Borel double cosets. In \cite{MR1734126}, Hida introduced the notion of $P$-ordinary modular forms on a reductive group $F$ with parabolic subgroup $P$ of $G$, generalizing the case $P=B$ in this article. A $p$-adic automorphic form can be $P$-ordinary without being $B$-ordinary. There are works on this aspects, for example \cite{MR4201778}. The deep work \cite{MR4522696} built up the “noncuspidal Hida theory for semiordinary modular forms” and proved the Iwasawa main conjecture of such $p$-adic families.
%    \item \emph{Allowing general $p$-slopes}. Here the ordinary $p$-adic automorphic forms are defined to be of $p$-slope zero. There are generalizations to the $p$-adic family of automorphic forms of finite $p$-slope, as \emph{Coleman families}.
%\end{itemize}
%\end{remark}

%% file: AP03-MainConjecture.tex
\section{Greenberg-type Iwasawa main conjectures for unitary groups} \label{sec:mainconjectures}

In this section, we formulate the Greenberg-type Iwasawa main conjectures for unitary groups.

\subsection{Galois representations} \label{sec:Galoisreps}
Dating back to \cite{skinner2012galois}, which relies on the works \cite{shin2011galois, morel2010cohomology}, it is possible to attach a Galois representation to a cuspidal automorphic representation $\sigma$ as part of the Eisenstein datum $\calD$.

\begin{theorem} \label{thm:galoisrep}
For some finite extension $L/\QQ_p$, there exists a Galois representation
\[
R_{p}(\sigma): \Gal_{\calK} \rightarrow \GL_{m+n}(L),
\]
such that:
\begin{enumerate}[label = \rm (\arabic*)]
    \item $R_p(\sigma)^{\sfc} \cong R_{p}(\sigma)^{\vee} \otimes \sigma_{p}(\chi_{\pi}^{1+\sfc}) \epsilon_{\cyc}^{1-m-n}$, where $\sigma_{p}(-)$ denotes the Galois character associated via class field theory, and $\chi_{\pi}$ denote the central character of $\sigma$.
    \item $R_{p}(\sigma)$ is unramified at all places away from $p$ where $\sigma$ is unramified. For such a place $w$, we have
    \[
    \det(1- R_{p}(\sigma)(\Frob_{w}) q_w^{-s}) = L \left(\BC(\sigma)_{w} \otimes \chi_{\pi, w}^{\sfc}, s+ \frac{1-m-n}{2} \right)^{-1},
    \]
    where $\Frob_{w}$ is the geometric Frobenius automorphism.
\end{enumerate}
\end{theorem}
We write $V_{\sigma}$ or simply $V$ for the representation space underlying $R_{p}(\sigma)$ \footnote{Note that we also used $V_{\sigma}$ for the representation space of $\sigma$ as an automorphic representation over $\rmU(m,n)$; we expect this will not cause confusion.}. It is well-known that there exists an $\calO_L$-lattice $T \subset V$ which is stable under the action of $\Gal_{\calK}$. We always assume the following:
\begin{itemize}
    \item[\namedlabel{item:irred}{(irred)}] There is a $\Gal_{\calK}$-stable lattice $T$ such that the residue Galois representation $\barR_{p}(\sigma)$ is absolutely irreducible. 
    \item[\namedlabel{item:dist}{(dist)}] The characters on the diagonal in \eqref{eq:Rpupper} are distinct modulo the maximal ideal of $\calO_{L}$. 
\end{itemize}
Under \ref{item:irred}, the $\Gal_{\calK}$-stable lattice $T$ is unique up to homothety.

We further suppose that $\sigma$ is unramified and ordinary at all places above $p$. Then $R_{p}(\sigma)$ satisfies the following property for $v = v_0 \barv_0 \in \scV_{\calF}^{(p)}$ with $v_0$ in our $p$-adic CM-type $\Sigma_p$:
\begin{equation} \label{eq:Rpupper}
R_p(\sigma)|_{\Gal_{\calK_{v_0}}} \simeq \begin{bmatrix}
\xi_{1, v_0} \epsilon_{\cyc}^{-\kappa_{m+n,v_0}} & \ast & \ast \\
 & \ddots & \ast \\
 & & \xi_{m+n, v_0} \epsilon_{\cyc}^{-\kappa_{1,v_0}}
\end{bmatrix},
\end{equation}
where the $\xi_{i, v_0}$ are unramified characters, and
\begin{itemize}
    \item for $i = 1, \ldots, n$, $\kappa_{i, v_0}$ is defined as $\frac{m+n}{2} + n - i - t_{n+1-i, v}^{-}$;
    \item for $i = 1, \ldots, m$, $\kappa_{n+i, v_0}$ is defined as $\frac{m+n}{2} + m+n - i - t_{m+1-i, v}^{+}$.
\end{itemize}
Using Theorem \ref{thm:galoisrep}(1), we deduce that $R_{p}(\sigma)|_{\Gal_{\calK_{\barv_0}}}$ is also equivalent to an upper triangular representation as above. By \eqref{eq:Rpupper}, there is a cotorsion-free $\Gal_{\calK_{v_0}}$-stable submodule $T_{v_0}^{+} \subset T_{v_0}$ of rank $m$ corresponding to the first $m$ rows. Similarly, there is a cotorsion-free $\Gal_{\calK_{\barv_0}}$-stable submodule $T_{\barv_0}^{+} \subset T_{\barv_0}$ of rank $n$ corresponding to the first $n$ rows.

Suppose $\sigma$ is a specialization of an $\II$-adic cuspidal Hida family $\bff$, where $\II$ is a normal domain and a finite $\Lambda_{m,n}^{\circ}$-algebra, where $\Lambda_{m,n}^{\circ} := \calO_{L} \lrbracket{T_{m,n}(1+p\ZZ_p)}$ is the weight algebra. We assume that \footnote{See Definition \ref{defn:specializations} for the definitions of specializations and arithmetic specializations.}:
\begin{itemize}
    \item[\namedlabel{item:irredbff}{(IRRED)}] The specialization $\sigma_{\bff_{\sfP}}$ satisfies \ref{item:irred} for some (and hence all) arithmetic specializations $\sfP$.
    \item[\namedlabel{item:distbff}{(DIST)}] The specialization $\sigma_{\bff_{\sfP}}$ satisfies \ref{item:dist} for some (and hence all) arithmetic specializations $\sfP$.
\end{itemize}
From the theory of pseudo-representations, it follows from \ref{item:irredbff} that there exists a continuous $\II$-linear Galois representation $(R_{p}(\bff), T_{\bff})$ with 
\[
R_{p}(\bff): \Gal_{\calK} \rightarrow \GL_{\II}(T_{\bff})
\]
and $T_{\bff}$ a free $\II$-module of rank $m+n$, interpolating the representations $R_{p}(\sigma_{\bff_{\sfP}})$ for arithmetic specializations $\sfP$ (possibly after extension of scalars). Let $V_{\bff} := T_{\bff} \otimes_{\II} \Frac(\II)$. We use the properties of each $\bff_{\sfP}$ to deduce the existence of a cotorsion-free $\Gal_{\calF_{v_0}}$-stable $\Frac(\II)$-submodule $V_{\bff, v_0}^{+}$ of $V_{\bff, v_0}$ of rank $m$. If \ref{item:distbff} holds, then $T_{\bff, v_0}^{+} := T_{\bff, v_0} \cap V_{\bff, v_0}^{+}$ is a free $\II$-summand of $T_{\bff, v_0}$ of rank $m$. Similarly, we have a free $\II$-summand $T_{\bff, \barv_0}^{+}$ of $T_{\bff, \barv_0}$ of rank $n$. \footnote{Regarding how to obtain the direct summand $T_{\bff, v_0}^{+}$ satisfying the above conditions by assuming \ref{item:distbff}, one may refer to the parallel discussion in the proof for cuspidal ordinary Hida families for $\GL_{2,\QQ}$ in \cite[Theorem 4.3.2]{hida2012geometric}, where the condition ``$\varepsilon \not\equiv \delta \pmod{\frm_{\II}}$'' is parallel to Assumption \ref{item:distbff} here.}

\subsection{Eisenstein datum and $p$-adic Eisenstein datum}
Recall that we defined the Eisenstein datum over the group $H = \rmU(m,n)$ in Definition \ref{defn:eisensteindatum}. We also have a $p$-adic variation of Eisenstein data, defined as follows.

\begin{definition} \label{defn:padiceisendatum}
A \emph{$p$-adic Eisenstein datum} over $H$ is a tuple $\bfD = (L, \II, \eta_0, \bff, \psi_0, \chi_0, \scS)$, where:
\begin{itemize}
    \item $L/\QQ_p$ is a finite extension.
    \item $\II$ is a normal domain over $\Lambda_{m,n}^{\circ}$, which is also a finite $\Lambda_{m,n}^{\circ}$-algebra.
    \item $\eta_0$ is a branch character, i.e., a finite order character of $\Delta$, the torsion part of $T_{m,n}(\ZZ_p)$.
    \item $\bff$ is an $\II$-adic Hida family of tempered cuspidal ordinary eigenforms on $H$, of tame level group $K_H^p \leq H(\AA_{\rmf}^{p})$ and branch character $\eta_0$ \footnote{In \cite{wan2015families}, the tempered condition is included in the assumption (TEMPERED) in \cite[Section 5B2]{wan2015families}. Here we include this assumption in the definition of the $p$-adic Eisenstein datum.}, satisfying \ref{item:irredbff} and \ref{item:distbff}.
    \item $\psi_0$ is a unitary Hecke character $\calK^{\times} \bs \KK^{\times} \rightarrow \CC^{\times}$ such that $\psi_0|_{\AA}$ is the central character of $\bff$.
    \item $\chi_0$ is a unitary Hecke character $\calK^{\times} \bs \KK^{\times} \rightarrow \CC^{\times}$.
    \item $\scS := \scS_{\bff} \cup \scS_{\chi_0} \cup \scS_{\calF}^{(p)}$, where $\scS_{\bff}$ (resp. $\scS_{\chi_0}$) is the set of finite places $v$ of $\calF$ away from $p$ such that $K_{H,v} \neq H(\calF_v)$ (resp. $\chi_0$ is ramified).
\end{itemize}
Let $\II^{\ur}$ be the normalization of an irreducible component of $\II \hotimes_{\calO_{L}} \calO_{L}^{\ur}$\footnote{More precisely, there is a bijection between the minimal primes of $\II \hotimes_{\calO_{L}} \calO_{L}^{\ur}$ and the irreducible components of $\Spec (\II \hotimes_{\calO_{L}} \calO_{L}^{\ur})$. Let $\frP$ be a minimal prime ideal of $\II \hotimes_{\calO_{L}} \calO_{L}^{\ur}$; then the corresponding irreducible component is isomorphic to $\Spec (\II \hotimes_{\calO_{L}} \calO_{L}^{\ur}/\frP)$. The latter quotient ring is an integral domain, and we take its integral closure in its fraction field; this is the $\II^{\ur}$ in the text. As remarked in \cite[p.~1957]{wan2015families}, for each such irreducible component, we can perform the subsequent construction. \label{fn:iurdomain}}. Given a $p$-adic Eisenstein datum $\bfD$, we define the corresponding Iwasawa algebra as $\Lambda_{\bfD} := \II^{\ur}\lrbracket{\Gamma_{\calK}}$. It is a $\Lambda_{m,n}^{\circ}$-algebra. Attached to $\bfD$, we define $\xi_{0} := \psi_0/\chi_0$, and $\bxi_0 = \bpsi_0/\bchi_0$ is the universal character attached to $\xi_0$ with target in $\Lambda_{\bfD}^{\times}$.
\end{definition}

The Iwasawa algebra $\Lambda_{\bfD}$ will be used as the weight space for $p$-adic interpolations. We define
\begin{equation} \label{eq:weightspace}
 \calX_{\bfD} := \Hom_{\cts}(\Lambda_{\bfD}, \barQQ_p) = \Spec \Lambda_{\bfD}(\barQQ_p)
\end{equation}
as the set consisting of continuous $\barQQ_p$-valued characters of $\Lambda_{\bfD}$.

\begin{definition} \label{defn:specializations}
Let $\sfP \in \calX_{\bfD}$ be a \emph{specialization} of $\bfD$.
\begin{enumerate}[label = \rm (\arabic*)]
    \item We call $\sfP$ an \emph{arithmetic point} if there exists an integer $\kappa_{\sfP} \geq m+n+2$ satisfying the following two conditions:
\begin{itemize}
    \item Let $\tau_{\sfP}$ be the pullback of $\sfP$ to $\Lambda_{m,n}^{\circ}$ along the canonical maps $\Lambda_{m,n}^{\circ} \hookrightarrow \II \rightarrow \II \hotimes \calO_L^{\ur} \rightarrow \II^{\ur}$; then it is a $p$-adic weight in the sense of Definition \ref{defn:padicweight}, with
    $$
    (t_{1, \tau_{\sfP}, v}^{+}, \ldots, t_{m, \tau_{\sfP}, v}^{+}; t_{1, \tau_{\sfP}, v}^{-}, \ldots, t_{n, \tau_{\sfP}, v}^{-}) = (0, \ldots, 0; -\kappa_{\sfP}, \ldots, -\kappa_{\sfP})
    $$
    for any $v \in \scV_{\calF}^{(p)}$.
    \item Let $\chi_{\sfP} := \sfP \circ \bchi_0$ be the composition of $\sfP$ with the universal character $\bchi_0$, i.e.,
    \begin{equation} \label{eq:definechiP}
    \calK^{\times} \bs \KK^{\times} \xrightarrow{\bchi_0} \Lambda_{\bfD} \xrightarrow{\sfP} \barQQ_{p}.
    \end{equation}
    Then $\chi_{\sfP}$ is a Hecke character of infinite type $\frac{\kappa_{\sfP}}{2}(\Sigma-\Sigma^\sfc)$.
\end{itemize}
The set of such points is denoted by $\calX_{\bfD}^{\ari}$.
\item We call $\sfP$ a \emph{classical arithmetic point} if $\sfP$ is an arithmetic point and the specialization $\bff_{\sfP} := \bff_{\tau_{\sfP}} \in V_{H}^{0, \ord}[\tau_{\sfP}]$ is a classical automorphic form under \eqref{eq:embedclassical} such that it generates an irreducible tempered unitary cuspidal automorphic representation $\sigma_{\bff_{\sfP}}$ of $H$. \footnote{By the requirement $\kappa_{\sfP} \geq m+n+2$, the archimedean part $\sigma_{\bff_{\sfP}, \infty}$ is indeed a holomorphic discrete series.} The set of such points is denoted by $\calX_{\bfD}^{\cls}$.
\item We call $\sfP$ a \emph{generic classical arithmetic point} if $\sfP$ is a classical arithmetic point such that $(\sigma_{\bff_{\sfP}}, \chi_\sfP)$ satisfies the generic condition \cite[Definition 4.42]{wan2015families}. The set of such points is denoted by $\calX_{\bfD}^{\gen}$.
\end{enumerate}
Given a point $\sfP \in \calX^{\cls}_{\bfD}$, we have an Eisenstein datum $\bfD_{\sfP} := (\sigma_{\bff_{\sfP}}, \psi_{\sfP}, \chi_{\sfP}, \scS)$ in the sense of Definition \ref{defn:eisensteindatum} as the specialization of $\bfD$ at $\sfP$, where $\psi_{\sfP}$ and $\chi_{\sfP}$ are defined following \eqref{eq:definechiP}. We write $s_{\sfP} := (\kappa_{\sfP} - (m+n+1))/2$.
\end{definition}

We know that $\calX_{\bfD}^{\ari}$, $\calX_{\bfD}^{\cls}$, and $\calX_{\bfD}^{\gen}$ are all Zariski dense in $\calX_{\bfD}$ \footnote{This is not quite trivial; it depends on the classicality theorem in Hida theory. See \cite[Remark 5.2]{wan2015families} for the justifications.}. Since we have assumed that $\bff$ satisfies \ref{item:irredbff}, for every $\sfP \in \calX_{\bfD}^{\cls}$, the cuspidal automorphic representation $\sigma_{\bff_\sfP}$ satisfies the multiplicity one theorem.

Let $\calD = (\sigma, \psi, \chi, \scS)$ be an Eisenstein datum over $H$ as defined in Definition \ref{defn:eisensteindatum}, and let $\bfD$ be a $p$-adic Eisenstein datum as defined in Definition \ref{defn:padiceisendatum}. In this section, we shall introduce the formulation of the Greenberg-type Iwasawa main conjecture for $\calD$ and $\bfD$.

\subsection{Greenberg-type Iwasawa main conjecture for unitary groups}
We formulate the Greenberg-type Iwasawa main conjecture for $\calD$ and $\bfD$.

\subsubsection{Algebraic side: Greenberg-type Selmer groups}
We are now ready to define the Greenberg-type Selmer group, following the axiomatic framework in \cite[Chapter 3]{skinner2014iwasawa}. Let $A$ be a profinite $\ZZ_p$-algebra, and assume that $T$ is a free $A$-module of finite rank, equipped with a continuous action of $\Gal_{\calF}$. We write $T_{v}$ for $T$ equipped with the restricted action of $\Gal_{\calF_{v}}$. Suppose that for each place $v \in \scV_{\calF}^{(p)}$, we are given a $\Gal_{\calF_{v}}$-stable free $A$-direct summand $T_{v}^{+} \subseteq T_v$.

\begin{definition} \label{defn:selmergroups}
Let $\scS$ be a finite set of places of $\calF$. We define the \emph{$\scS$-primitive Greenberg-type Selmer group} $\Sel_{\Gr}^{\scS}(T, \{T_{v}^{+}\}_{v \in \scV_{\calF}^{(p)}}, \calF, A)$ to be the kernel of the restriction map:
\[
\rmH^{1}(\calF, T \otimes_{A} A^{\ast}) \rightarrow \prod_{v \not\in \scS, v \nmid p} \rmH^{1}(I_{v}, T \otimes_{A} A^{\ast}) \times \prod_{v \mid p} \rmH^{1}(I_{v}, (T_v/T_{v}^{+}) \otimes_{A} A^{\ast}).
\]
We also define the Pontryagin dual:
\[
X_{\Gr}^{\scS}(T, \{T_{v}^{+}\}_{v \in \scV_{\calF}^{(p)}}, \calF, A) := \Hom_{A}(\Sel_{\Gr}^{\scS}(T, \{T_{v}^{+}\}_{v \in \scV_{\calF}^{(p)}}, \calF, A), A^{\ast}).
\]
\end{definition}

\begin{remark}
In more memorable terms, when defining the Greenberg-type Selmer group, we impose the relaxed Selmer condition at places in $\scS \setminus \scV_{\calF}^{(p)}$, the unramified Selmer condition at places outside $\scS \cup \scV_{\calF}^{(p)}$, and the Greenberg Selmer condition at places above $p$.
\end{remark}

For a profinite $\ZZ_p$-algebra $A$, we define the \emph{Iwasawa algebra with $A$-coefficients} as $\Lambda_{\calK, A} := \ZZ_p \lrbracket{\Gamma_{\calK}} \otimes_{\ZZ_p} A$. Recall that in Section \ref{sec:notationcharacters}, we defined the character $\bvep_{\calK}: \calK^{\times} \bs \AA_{\calK}^{\times} \rightarrow \Gamma_{\calK} \hookrightarrow \Lambda_{\calK}^{\times}$ given by class field theory. Furthermore, we define the \emph{Iwasawa deformation} of an $A$-module $T$ as
\[
T_{\Iw} := T \otimes_{A} \Lambda_{\calK, A}(\bvep_{\calK}^{-1}),
\]
which is a free $\Lambda_{\calK, A}$-module of finite rank \footnote{Note that, although not explicitly indicated in the notation, $T_{\Iw}$ clearly depends on the choice of the coefficient ring $A$, the number field $\calK$, and the (multiple) $\ZZ_p$-extension over $\calK$. As we shall see, when defining the Greenberg-type Selmer group for a cuspidal automorphic representation $\sigma$, the coefficient ring $A$ is taken to be the coefficient ring $\calO$ of its Galois representation; for an $\II$-adic cuspidal Hida family $\bff$, the coefficient ring $A$ is taken to be $\II$.}.

We define the \emph{Greenberg-type Selmer group} of a cuspidal automorphic representation $\sigma$ as
\[
\Sel_{\Gr}^{\scS}(\sigma/\calK_{\infty}) := \Sel_{\Gr}^{\scS}(T_{\Iw}, \{T_{v_0, \Iw}^{+}, T_{\barv_0, \Iw}^{+}\}_{v \in \scV_{\calF}^{(p)}}, \calK, \Lambda_{\calK, \calO}).
\]
If $\bff$ is an $\II$-adic cuspidal Hida family on $\rmU(m,n)$, then we can define the \emph{Greenberg-type Selmer group} of the Hida family $\bff$ as
\[
\Sel_{\Gr}^{\scS}(\bff/\calK_{\infty}) := \Sel_{\Gr}^{\scS}(T_{\bff, \Iw}, \{T_{\bff, v_0, \Iw}^{+}, T_{\bff, \barv_0, \Iw}^{+}\}_{v \in \scV_{\calF}^{(p)}}, \calK, \II\lrbracket{\Gamma_{\calK}}).
\]

Suppose $\chi: \calK^{\times} \bs \AA_{\calK}^{\times} \to \CC^{\times}$ is an algebraic Hecke character; we define the \emph{$\chi$-twisted Greenberg-type Selmer groups}
\[
\Sel_{\Gr}^{\scS}(\sigma/\calK_{\infty}, \chi) \quad \text{and} \quad \Sel_{\Gr}^{\scS}(\bff/\calK_{\infty}, \chi)
\]
by replacing the Galois-stable lattices in their definitions with their $\chi$-twists. We also have 
\[
X_{\Gr}^{\scS}(\sigma/\calK_{\infty}, \chi) \quad \text{and} \quad X_{\Gr}^{\scS}(\bff/\calK_{\infty}, \chi)
\]
defined following Definition \ref{defn:selmergroups}. We remark that \cite[Proposition 3.4]{skinner2014iwasawa} also specifies other equivalent ways of defining these Selmer groups.

\subsubsection{Analytic side: $p$-adic $L$-functions}
There are various works regarding $p$-adic $L$-functions over unitary groups. In \cite{eischen2020padic}, the authors defined an $\scS$-imprimitive $p$-adic $L$-function $\calL^{\scS}_{\mathtt{EHLS}}(\sigma, \calK, \chi) \in \calO^{\ur}_{L}\lrbracket{\Gamma_{\calK}}$ interpolating certain special values of $L^{\scS}(\sigma^{\vee}, \chi_{0}^{-\sfc}, s)$. Putting $\sigma$ into a Hida family $\bff$, there is a recent generalization of \cite{eischen2020padic} by Marcil \cite{marcil2024PadicLfunctionsPordinary}. In this article, we prefer to work with the $p$-adic $L$-function of Hida families $\bff$ constructed by Wan in \cite{wan2015families}, denoted by 
\[
\calL^{\scS}_{\mathtt{Wan}}(\bff, \calK, \chi_{0}) \in \II^{\ur}\lrbracket{\Gamma_{\calK}} \otimes_{\II^{\ur}} \Frac(\II^{\ur}),
\]
which was proved to be an element in $\II^{\ur}\lrbracket{\Gamma_{\calK}}$ for $n = 0$ in \cite{wan2015families}. It satisfies the following interpolation property:
\begin{equation} \label{eq:wanpadicl}
\sfP(\calL^{\scS}_{\mathtt{Wan}}(\bff, \calK, \chi_{0})) = \left(\dfrac{(-2)^{-n}(2 \pi \rmi)^{n \kappa_{\sfP}} (2/\pi)^{n(n-1)/2}}{\prod_{j=0}^{n-1} \Gamma(\kappa_{\sfP}-j)} \right)^{-d} C_{\bff_{\sfP}}^{(p)} \Omega_{\sfP}^{-1} L^{\scS} \left( \sigma_{\bff_{\sfP}}^{\vee} \times \bchi_{0, \sfP}^{-\sfc}, \dfrac{\kappa_{\sfP}-(m+n)}{2} \right)    
\end{equation}
for $\sfP \in \calX_{\bfD}^{\gen}$ of weight $\kappa_{\sfP}$, where $C_{\bff_{\sfP}}^{(p)}$ is an explicit constant factor specified in \cite[Theorem 1.1]{wan2015families}, and 
\begin{equation} \label{eq:defineperiod}
\Omega_{\sfP} = \begin{cases}
(\Omega_{p}^{\infty}/\Omega_{\infty}^{\infty})^{-m \kappa_{\sfP}}, &\quad n = 0, \\
\lrangle{\bff_{\sfP}^{\vee}, \bff_{\sfP}}_{\Pet}, &\quad n \neq 0,
\end{cases}    
\end{equation}
with $\Omega_p \in (\ZZ_p^{\ur})^{\scV_{\calF}^{\infty}}$ the \emph{$p$-adic period} and $\Omega_{\infty} \in \CC^{\scV_{\calF}^{\infty}}$ the \emph{CM period} \footnote{We refer to \cite{hida2004non} for precise definitions.}. We denote by $\Omega_{\infty}^{\infty}$ the product of $d$ elements of $\Omega_{\infty}$ and define $\Omega_p^{\infty}$ similarly. Note that when $n = 0$, we normalize using the CM period, and when $n \neq 0$, we normalize using the Petersson period.

\subsubsection{Greenberg-type Iwasawa main conjecture for Eisenstein data}
We expect the following conjectures to hold, which are known as the Greenberg-type Iwasawa main conjectures. 

\begin{conjecture}[Greenberg-type Iwasawa main conjecture for $\calD$] \label{conj:twovariableimc}
With notations as above and assuming \ref{item:irred} and \ref{item:dist}:
\begin{enumerate}[label = \rm (\arabic*)]
    \item The dual Greenberg-type Selmer group $X_{\Gr}^{\scS}(\sigma/\calK_{\infty}, \xi_{0}^{-\sfc} \epsilon_{\cyc})$ is a finitely generated torsion $\calO_{L}\lrbracket{\Gamma_{\calK}}$-module.
    \item The $p$-adic $L$-function $\calL_{\mathtt{EHLS}}^{\scS}(\sigma, \calK, \chi_{0})$ belongs to $\calO_{L}^{\ur}\lrbracket{\Gamma_{\calK}}$.
    \item There is an equality of principal ideals of $\calO_{L}^{\ur}\lrbracket{\Gamma_{\calK}}$:
    \[
    \Char_{\Lambda_{\calK, \calO}}(X_{\Gr}^{\scS}(\sigma/\calK_{\infty}, \xi_{0}^{-\sfc} \epsilon_{\cyc})) = (\calL^{\scS}_{\mathtt{EHLS}}(\sigma, \calK, \chi_0)).
    \]
    The left-hand side is a principal ideal in $\calO_{L}\lrbracket{\Gamma_{\calK}}$, regarded here as a principal ideal in $\calO_{L}^{\ur}\lrbracket{\Gamma_{\calK}}$.
\end{enumerate}
\end{conjecture}

\begin{conjecture}[Greenberg-type Iwasawa main conjecture for $\bfD$] \label{conj:threevariableimc}
With notations as above and assuming \ref{item:irredbff} and \ref{item:distbff}:
\begin{enumerate}[label = \rm (\arabic*)]
    \item The dual Greenberg-type Selmer group $X_{\Gr}^{\scS}(\bff/\calK_{\infty}, \xi_{0}^{-\sfc} \epsilon_{\cyc})$ is a finitely generated torsion $\II\lrbracket{\Gamma_{\calK}}$-module.
    \item The $p$-adic $L$-function $\calL_{\mathtt{Wan}}^{\scS}(\bff, \calK, \chi_{0})$ belongs to $\II^{\ur}\lrbracket{\Gamma_{\calK}}$.
    \item There is an equality of principal ideals of $\II^{\ur}\lrbracket{\Gamma_{\calK}}$:
    \[
    \Char_{\II^{\ur}\lrbracket{\Gamma_{\calK}}}(X_{\Gr}^{\scS}(\bff/\calK_{\infty}, \xi_0^{-\sfc} \epsilon_{\cyc})) = (\calL^{\scS}_{\mathtt{Wan}}(\bff, \calK, \chi_0)).
    \]
    The left-hand side is a principal ideal in $\II\lrbracket{\Gamma_{\calK}}$, regarded here as a principal ideal in $\II^{\ur}\lrbracket{\Gamma_{\calK}}$.
\end{enumerate}
\end{conjecture}

\begin{remark}
In informal terms, Conjecture \ref{conj:twovariableimc} is often called the ``two-variable Iwasawa main conjecture,'' since when $\calK$ is an imaginary quadratic extension of $\QQ$, $\Lambda_{\calK}$ is isomorphic to a two-variable power series ring over $\ZZ_p$. Conjecture \ref{conj:threevariableimc} is then called the ``three-variable Iwasawa main conjecture'' because $\sigma$ has been deformed into a Hida family. This terminology can be confusing, as $\II$ is an algebra over the weight algebra $\Lambda_{m,n}^{\circ}$, which is itself a power series ring in $m+n$ variables over $\calO$.
\end{remark} 

Currently, the Iwasawa main conjectures are proved by establishing two divisibilities individually, each with completely different machinery:
\begin{itemize}
    \item \emph{The upper bound of Selmer groups}: This is the inclusion ``$\supseteq$'' in Conjecture \ref{conj:twovariableimc} and \ref{conj:threevariableimc}. It is often proved by the machinery of Euler systems.
    \item \emph{The lower bound of Selmer groups}: This is the inclusion ``$\subseteq$'' in Conjecture \ref{conj:twovariableimc} and \ref{conj:threevariableimc}. It is often proved by the machinery of Eisenstein congruences, which is the topic of our article.
\end{itemize}

We summarize the main breakthroughs of Greenberg-type Iwasawa main conjectures in recent years, yet ignoring technical hypotheses. We also include the case of Iwasawa theory over the symplectic group $\Sp(4)$, in which case the Greenberg-type Iwasawa main conjecture can be formulated in a parallel way.

\begin{table}[H]
\centering
\begin{tabular}{|c|c|c|}
\hline
\multicolumn{1}{|c|}{}       & Lower bounds & Upper bounds \\ \hline
$\rmU(1,0)$            & \cite{hsieh2014eisenstein}     &  \makecell{\cite{rubin1991main} \\ (and \cite{hida1994anticyclotomic} for anticyclotomic characters)}   \\ \hline
$\rmU(1,1)$         &   \makecell{\cite{skinner2014iwasawa, burungale2025base} for $\calF=\QQ$   \\ and \cite{wan2015hilbert} in general}          &         \makecell{\cite{kato2004padic} for $\calF = \QQ$ \\ (and \cite{loeffler2025iwasawa} for $\calF$ being real quadratic fields) }     \\ \hline
$\rmU(2,0)$      &   \makecell{ \cite{wan2020iwasawa, castella2022iwasawa} for $\calF = \QQ$   \\ (also \cite{wan2019iwasawatheorymathrmursblochkato, fouquet2022iwasawa}) }         &   \makecell{ Same as above for Jacquet-Langlands \\ transfer of elliptic newforms } \\ \hline \hline
$\Sp(4)$ & \cite{urban2001selmer, urban2006groupes} & \cite{vincentelli2025euler} \\ \hline
\end{tabular}
\end{table}
Recently, there have been positive results for upper bounds of Selmer groups in higher rank cases \cite{liu2022beilinson, liu2024iwasawa}, but the lower bounds of Selmer groups are still out of reach in higher rank cases.

\subsection{Eisenstein congruences} \label{sec:eisensteincongruences}
We explain the machinery of Eisenstein congruences more explicitly here. To prove the lower bound of Selmer groups, Eisenstein congruences break the inclusion into two steps:
\begin{equation} \label{eq:inserteisenstein}
    \Char_{\Lambda_{\calK}}(X_{\Gr}^{\scS}) \subseteq \calE_{\bfD} \subseteq (\calL^{\scS}_{p})
\end{equation}
where we write the Greenberg-type Selmer group $X_{\Gr}^{\scS}$ and $p$-adic $L$-functions $\calL^{\scS}_{p}$ (of \texttt{EHLS} and \texttt{Wan}) in short. Here $\calE_{\bfD}$ is a suitably defined ``\emph{Eisenstein ideal}'' attached to $\bfD$, which measures the congruence relations between the ordinary Klingen Eisenstein family $\bfE_{\bfD}^{\Kling}$ over $\rmU(m+1, n+1)$ and cuspidal Hida families over $\rmU(m+1, n+1)$. The Klingen Eisenstein family has the advantage that its constant term is essentially the $p$-adic $L$-function $\calL^{\scS}_{p}$.

These two successive inclusions in \eqref{eq:inserteisenstein} are expected to be proved in the following two steps:
\begin{enumerate}[label = (\arabic*)]
    \item The inclusion $\calE_{\bfD} \subseteq (\calL^{\scS}_{p})$ can be regarded as ``purely automorphic'' and is the most crucial part of the machinery of Eisenstein congruences. It has two key inputs:
    \begin{enumerate}[label = (1.\alph*)]
        \item A fundamental exact sequence of Eisenstein congruences that arises from Hida theory or appropriate theory of $p$-adic families of modular forms. \label{item:1a}
        \item Certain $p$-primitivity results of Eisenstein families. \label{item:1b}
    \end{enumerate}
    \item The inclusion $\Char_{\Lambda_{\calK}}(X_{\Gr}^{\scS}) \subseteq \calE_{\bfD}$ is proved by the so-called ``\emph{lattice construction},'' which has been built up axiomatically in \cite[Chapter 1]{bellaiche2009families}. Morally speaking, this is done by translating the congruences of Eisenstein series with cusp forms into congruences of the corresponding Galois representations, which provide sufficiently many Galois cohomology classes in $X_{\Gr}^{\scS}$. \label{item:2}
\end{enumerate}
In this article, we shall introduce \ref{item:1a} in Section \ref{sec:noncuspidalhida} and \ref{item:2} in Section \ref{sec:lattice}. The main theme of this article is \ref{item:1b}, which will be introduced in Part \ref{part:automorphic}.

%% file: AU01-Reductions.tex
\section{The Gan--Gross--Prasad-type period integral} \label{sec:ggpperiodintegral}
We inherit all notations and conventions from previous sections. Let
\begin{itemize}
    \item $(\pi, V_{\pi})$ be an irreducible unitary cuspidal automorphic representation of $G$ over $\calF$ and identify $V_{\pi}$ with a subspace of $\calA_{\cusp}(G(\calF) \bs G(\AA))$, and
    \item $E^{\Kling}(f^{\Kling}_{\sigma, \chi, s}, -)$ be a Klingen Eisenstein series on $G^{\hsuit}$, defined in Section \ref{sec:Klingen}. 
\end{itemize}
Let $\Psi \in V_{\pi}$. We define the \emph{Gan--Gross--Prasad-type period integral} (GGP period integral) of $E^{\Kling}$ with the cusp form $\Psi$ as the integral
\begin{equation} \label{eq:GGPdefn}
\calP(\Psi, E^{\Kling}(f^{\Kling}_{\sigma, \chi, s}, -)) := \int_{G(\calF) \bs G(\AA)} E^{\Kling}(f^{\Kling}_{\sigma, \chi, s}, \jmath^{\flat}(g)) \Psi(g) \dif g.    
\end{equation}
By the cuspidality of $\Psi$, this converges absolutely for those values of $s$ at which $E^{\Kling}(f^{\Kling}_{\sigma, \chi, s}, -)$ is defined. \footnote{We thank Wen-Wei Li for his helpful hint on this issue.} In fact, we can also regard this integral as integrating the automorphic form 
\[
E^{\Kling}(f_{\sigma,\chi,s}^{\Kling}, ?_{1}) \boxtimes \Psi(?_{2}) \quad \text{ of the group } \bfG 
\]
over the subgroup $\bfH$ of $\bfG$, so it is a \emph{Bessel-type period integral}.

In this section, we consider in particular the Klingen Eisenstein series constructed by pulling back from Siegel Eisenstein series (see Section \ref{sec:garrettdou}). Our plan is as follows:
\begin{enumerate}
    \item We first unfold the GGP period integral \eqref{eq:GGPdefn} to reduce it to a certain ``cuspidal GGP period integral''. See Proposition \ref{prop:redcusp}.
    \item Then we square the integral and use the Ichino--Ikeda formula (Theorem \ref{thm:ikedaichino}) to decompose it into products of local integrals, involving local Ichino--Ikeda integrals and local doubling integrals. See Theorem \ref{thm:breaklocal}.
    \item We further invoke the unramified computations of such local integrals to see that the Rankin--Selberg local $L$-factors and standard local $L$-factors arise at ``good'' places (which covers all but finitely many places) of $\calF$. See Theorem \ref{thm:unrret}.
\end{enumerate}
There are finitely many ``bad'' places that remain, which we shall deal with in the sections to follow. Before carrying out this three-step scheme, we review several basic facts concerning the doubling method. Although these facts are well known to experts, the specific formulation required in our setting is not readily available in the literature. For this reason, and for the sake of completeness, we include proofs of certain statements.

\subsection{Doubling method, à la Piatetski-Shapiro and Rallis}
Here we briefly recall the doubling method à la Piatetski-Shapiro and Rallis, first introduced in \cite{gelbart1987explicit}.

Let $(\pi, V_{\pi})$ be an irreducible tempered unitary cuspidal automorphic representation of $G$ and $(\pi^{\vee}, V_{\pi^{\vee}})$ be its contragredient representation. We thus identify $V_{\pi}$ with a subspace of $\calA_{\cusp}(G(\calF) \bs G(\AA))$. Recall that we assumed that $\pi$ satisfies \ref{item:irred} and hence the property \ref{ass:multiplicityone}. We define the global \emph{doubling integral à la Piatetski-Shapiro and Rallis} as
$$
Z^{\dsuit}(f^{\Sieg}_{s,\chi}, \Psi, \Psi^{\vee}) := \int_{[G^{\dsuit}]} E^{\Sieg}(f^{\Sieg}_{s,\chi}, \imath^{\dsuit}(g,h)) \Psi(g) \Psi^{\vee}(h) \chi^{-1}(\det h) \dif g \dif h.
$$
where $\Psi \in V_{\pi}$ and $\Psi^{\vee} \in V_{\pi^{\vee}}$. It converges wherever the Eisenstein series is defined, by the cuspidality of $\Psi$ and $\Psi^{\vee}$.

\subsubsection{Basic identity of Piatetski-Shapiro and Rallis}
We define $G^{\dsuit}$ as the image of $G \times G$ via the canonical doubling embedding $\imath^{\dsuit}: G \times G \hookrightarrow G^{\bdsuit}$ and $G^{\Delta}$ as the image of the diagonal embedding $G \rightarrow G \times G$ composed with $\imath^{\dsuit}$.

The following fundamental result is well-known; it appears in the proof of the ``Basic Identity'' in \cite[page 3]{gelbart1987explicit}, though the character $\chi$ was absent there. For the reader's convenience, we provide a sketch of the proof here, adapted from the proof of \cite[Theorem 4.3.4]{eischen2024automorphic} in the setup of unitary groups.

\begin{theorem}[Basic identity of Piatetski-Shapiro and Rallis] \label{thm:basicid}
\[
Z^{\dsuit}(f^{\Sieg}_{s,\chi}, \Psi, \Psi^{\vee}) = \int_{G^{\Delta}(\calF) \bs G^{\dsuit}(\AA)} f^{\Sieg}_{s,\chi}(\imath^{\dsuit}(g,h)) \Psi(g) \Psi^{\vee}(h) \chi^{-1}(\det h) \dif g \dif h.
\]
\end{theorem}

\begin{proof}
The theorem follows from an analysis of the orbits of $G^{\dsuit}$ acting on $\calX := Q \bs G^{\bdsuit}$ by multiplication on the right. We write $(G^{\dsuit})^{\gamma}$ for the stabilizer of a point $\gamma \in \calX$. We can rewrite the Siegel Eisenstein series $E^{\Sieg}(f_{s,\chi}^{\Sieg},-)$ by grouping the summands by orbits:
\[
E^{\Sieg}(f_{s,\chi}^{\Sieg}, \eta) = \sum_{[\gamma] \in Q(\calF) \bs G^{\bdsuit}(\calF)/G^{\dsuit}(\calF)} \left( \sum_{[\gamma_0] \in (G^{\dsuit})^{\gamma}(\calF) \bs G^{\dsuit}(\calF)} f_{s,\chi}^{\Sieg}(\gamma \gamma_0 \eta) \right). 
\]
Here $[\gamma]$ denotes the orbit of $Q(\calF) \gamma \in \calX(\calF)$ under the right action of $G^{\dsuit}(\calF)$. Inserting this expression into the doubling integral, we see that $Z^{\dsuit}(f^{\Sieg}_{s,\chi}, \Psi, \Psi^{\vee})$ equals
\begin{align*}
& \sum_{[\gamma] \in Q(\calF) \bs G^{\bdsuit}(\calF)/G^{\dsuit}(\calF)} \left( \sum_{[\gamma_0] \in (G^{\dsuit})^{\gamma}(\calF) \bs G^{\dsuit}(\calF)} \int_{[G^{\dsuit}]} f^{\Sieg}_{s,\chi}(\gamma \gamma_0 \imath^{\dsuit}(g,h)) \Psi(g) \Psi^{\vee}(h) \chi^{-1}(\det h) \dif g \dif h \right) \\
&= \sum_{[\gamma] \in Q(\calF) \bs G^{\bdsuit}(\calF)/G^{\dsuit}(\calF)} \left( \int_{(G^{\dsuit})^{\gamma}(\calF) \bs G^{\dsuit}(\AA)} f^{\Sieg}_{s,\chi}(\gamma \imath^{\dsuit}(g,h)) \Psi(g) \Psi^{\vee}(h) \chi^{-1}(\det h) \dif g \dif h \right).
\end{align*}
Temporarily denote the integral in the summand as $I(\gamma)$. Note that for each $\gamma \in G^{\bdsuit}(\calF)$,
\begin{align*}
(G^{\dsuit})^{\gamma}(\calF) &= \{\imath^{\dsuit}(g,h) \in G^{\dsuit}(\calF) : Q(\calF) \gamma \imath^{\dsuit}(g,h) = Q(\calF) \gamma \} \\
&= \{\imath^{\dsuit}(g,h) \in G^{\dsuit}(\calF): \gamma \imath^{\dsuit}(g,h) : \gamma \imath^{\dsuit}(g,h) \gamma^{-1} \in Q(\calF)  \}.
\end{align*}
We first consider the orbit $\gamma = 1$. In this case, the stabilizer is 
\[
(G^{\dsuit})^{1}(\calF) = Q(\calF) \cap G^{\dsuit}(\calF) = \{\imath^{\dsuit}(g,g): g \in G(\calF) \} = G^{\Delta}(\calF),
\]
and
\[
I(1) = \int_{G^{\Delta}(\calF) \bs G^{\dsuit}(\AA)} f^{\Sieg}_{s,\chi}(\imath^{\dsuit}(g,h)) \Psi(g) \Psi^{\vee}(h) \chi^{-1}(\det h) \dif g \dif h.
\]
The remaining orbits (i.e., orbits with $\gamma \neq 1$) are \emph{negligible} in the sense of \cite[Part A, Chapter I]{gelbart1987explicit}, and $I(\gamma) = 0$ for these cases. This vanishing result follows from the cuspidality of $\Psi$ and $\Psi^{\vee}$, and the fact that $(G^{\dsuit})^{\gamma}$ contains the unipotent radical of a proper parabolic subgroup of $G^{\dsuit}(\calF)$ as a normal subgroup. Details can be found in \cite{gelbart1987explicit} or in the proof of \cite[Theorem 4.3.4]{eischen2024automorphic}. We remark that, as seen in the proof of \cite[Theorem 4.3.4]{eischen2024automorphic}, although we have additionally introduced the character $\chi$, $\chi \circ \det$ is trivial on unipotent radicals; thus, the application of the proof in \cite[Part A, Chapter I]{gelbart1987explicit} to negligible orbits remains unobstructed.
\end{proof}

\subsubsection{Partial doubling integrals}
Following \cite[page 174]{skinner2014iwasawa}, we define the partial doubling integrals:
\[
Z^{\dsuit, \leftt}(f^{\Sieg}_{s,\chi}, \Psi^{\vee}; g) := \int_{[G]} E^{\Sieg}(f^{\Sieg}_{s,\chi}, \imath^{\dsuit}(g,h)) \Psi^{\vee}(h) \chi^{-1}(\det h) \dif h,
\]
and
\[
Z^{\dsuit, \rightt}(f^{\Sieg}_{s,\chi}, \Psi, h) := \int_{[G]} E^{\Sieg}(f^{\Sieg}_{s,\chi}, \imath^{\dsuit}(g,h)) \Psi(g) \dif g.
\]
These converge by the cuspidality of $\Psi$ and $\Psi^{\vee}$. We then have:
\begin{equation} \label{eq:pairSU}
 Z^{\dsuit}(f^{\Sieg}_{s,\chi}, \Psi, \Psi^{\vee}) = \lrangle{\Psi, Z^{\dsuit, \leftt}(f^{\Sieg}_{s,\chi}, \Psi^{\vee}; -)}_{\Pet} = \lrangle{Z^{\dsuit, \rightt}(f^{\Sieg}_{s,\chi}, \Psi, -), \Psi^{\vee} \cdot \chi^{-1}(\det -)}_{\Pet}.   
\end{equation}

We shall use the following corollary of the basic identity.
\begin{corollary} \label{coro:partialdou}
With notations as above:
\begin{enumerate}[label= \rm (\arabic*)]
    \item Let $g, g_0 \in G(\AA)$, then
    \begin{align*}
        Z^{\dsuit, \leftt}(f^{\Sieg}_{s,\chi}, \Psi^{\vee}; g) &= \int_{G(\AA)} f^{\Sieg}_{s,\chi}(\imath^{\diamondsuit}(g, h)) \chi^{-1}(\det h) \Psi^{\vee}(h) \dif h \\
        &= \int_{G(\AA)} f^{\Sieg}_{s,\chi}(\imath^{\diamondsuit}(g_0, h)) \chi^{-1}(\det h) \Psi^{\vee}(g g_0^{-1} h) \dif h.
    \end{align*}
    \item Let $h, h_0 \in G(\AA)$, then
    \begin{align*}
    Z^{\dsuit, \rightt}(f^{\Sieg}_{s,\chi}, \Psi; h) &= \int_{G(\AA)} f^{\Sieg}_{s,\chi}(\imath^{\diamondsuit}(g, h)) \Psi(g) \dif g \\
    &= \chi(\det h h_0^{-1}) \int_{G(\AA)} f^{\Sieg}_{s,\chi}(\imath^{\diamondsuit}(g, h_0)) \Psi(h h_0^{-1} g) \dif g.  
    \end{align*}
    \item As a result, for $g, g_0, h, h_0 \in G(\AA)$,
    \begin{align*}
    Z^{\dsuit}(f^{\Sieg}_{s,\chi}, \Psi, \Psi^{\vee}) &= \int_{G(\AA)} f^{\Sieg}_{s,\chi}(\imath^{\diamondsuit}(g_0, h)) \lrangle{\Psi, \pi^{\vee}(g_0^{-1} h)\Psi^{\vee}}_{\Pet} \dif h \\
    &= \chi(\det h_0^{-1}) \int_{G(\AA)} f^{\Sieg}_{s,\chi}(\imath^{\diamondsuit}(g, h_0)) \lrangle{\pi(h_0^{-1} g)\Psi, \Psi^{\vee}}_{\Pet} \dif g.
    \end{align*}
\end{enumerate}
\end{corollary}
\begin{proof}
We note that there are two ways of identifying $G^{\Delta}(\calF) \bs (G \times G)(\AA)$ and $G(\AA) \times (G(\calF) \bs G(\AA))$ by sending $(g,h)$ to either $(g,h)$ or $(h,g)$. Theorem \ref{thm:basicid}, together with \eqref{eq:pairSU}, gives the first equalities of (1) and (2) accordingly, noting that the Petersson pairing $\lrangle{-,-}_{\Pet}$ is perfect.

To deduce the second equality of (1), we observe:
\[
f^{\Sieg}_{s,\chi}(\imath^{\diamondsuit}(g, h)) = f^{\Sieg}_{s,\chi}(\imath^{\diamondsuit}(g g_0^{-1}, g g_0^{-1})\imath^{\diamondsuit}(g_0, g_0 g^{-1} h)) = \chi(\det gg_0^{-1}) f^{\Sieg}_{s,\chi}(\imath^{\diamondsuit}(g_0 , g_0 g^{-1} h)).
\]
Substituting this into the integral from (1), we have:
\begin{align*}
    Z^{\dsuit, \leftt}(f^{\Sieg}_{s,\chi}, \Psi^{\vee}; g) &= \chi(\det g g_0^{-1}) \int_{G(\AA)} f^{\Sieg}_{s,\chi}(\imath^{\diamondsuit}(g_0 , g_0 g^{-1} h)) \chi^{-1}(\det h) \Psi^{\vee}(h) \dif h \\
&= \chi(\det gg_0^{-1}) \int_{G(\AA)} f^{\Sieg}_{s,\chi}(\imath^{\diamondsuit}(g_0, h^{\prime})) \chi^{-1}(\det gg_0^{-1} h^{\prime}) \Psi^{\vee}(gg_0^{-1} h^{\prime}) \dif h^{\prime} \\
&= \int_{G(\AA)} f^{\Sieg}_{s,\chi}(\imath^{\diamondsuit}(g_0, h^{\prime})) \chi^{-1}(\det h^{\prime}) \Psi^{\vee}(gg_0^{-1} h^{\prime}) \dif h^{\prime}
\end{align*}
via the change of variables $h^{\prime} = g_0 g^{-1} h$. This proves (1). The second equality of (2) is deduced similarly using
\[
f^{\Sieg}_{s,\chi}(\imath^{\diamondsuit}(g, h)) = \chi(\det hh_0^{-1}) f^{\Sieg}_{s,\chi}(\imath^{\diamondsuit}(h_0 h^{-1} g,h_0)).
\]
Putting into the integral on the right hand side of the first equality of (2), we see that
\begin{align*}
    Z^{\dsuit, \rightt}(f^{\Sieg}_{s,\chi}, \Psi; h) &= \chi(\det h h_0^{-1}) \int_{G(\AA)} f^{\Sieg}_{s,\chi}(\imath^{\diamondsuit}(h_0 h^{-1} g,h_0))
\Psi(g) \dif g \\
&= \chi(\det hh_0^{-1}) \int_{G(\AA)} f^{\Sieg}_{s,\chi}(\imath^{\diamondsuit}(g^{\prime}, h_0)) \Psi(h h_0^{-1} g^{\prime}) \dif g^{\prime}
\end{align*}
by change of variables $g^{\prime} = h_0 h^{-1} g$. This shows the second equality in (2). The equalities in (3) follow from (1) and (2) by applying \eqref{eq:pairSU}.
\end{proof}

By the uniqueness of $G$-invariant pairings between $\pi$ and $\pi^{\vee}$ provided by \ref{ass:multiplicityone}, and part (3) of Corollary \ref{coro:partialdou}, we obtain the following result.

\begin{corollary} \label{cor:doubintlocal}
Under the assumptions above, if $\Psi = \otimes^{\prime}_v \Psi_v \in \pi$ and $\Psi^{\vee} = \otimes^{\prime}_v \Psi^{\vee}_v \in \pi^{\vee}$ are factorizable, then
\[
Z^{\dsuit}(f^{\Sieg}_{s,\chi}, \Psi, \Psi^{\vee}) = \lrangle{\Psi, \Psi^{\vee}}_{\Pet} \prod_{v \in \scV_{\calF}} \frac{Z^{\dsuit}_v(f^{\Sieg}_{s,\chi,v}, \Psi_v, \Psi_{v}^{\vee})}{\lrangle{\Psi_v, \Psi^{\vee}_v}},
\]
where
\begin{align*}
    Z^{\dsuit}_v(f^{\Sieg}_{v, s,\chi}, \Psi_v, \Psi_v^{\vee}) &:= \int_{G(\calF_v)} f^{\Sieg}_{v,s,\chi}(\imath^{\diamondsuit}(g_0, h_v)) \lrangle{\Psi_v, \pi_v^{\vee}(g_0^{-1} h_v) \Psi_v^{\vee}} \dif h_v \\
    &= \chi_v^{-1}(\det h_0) \int_{G(\calF_v)} f^{\Sieg}_{v,s,\chi}(\imath^{\diamondsuit}(g_v, h_0)) \lrangle{\pi_v(h_0^{-1} g_v) \Psi_v, \Psi_v^{\vee}} \dif g_v.   
\end{align*}
\end{corollary}

\subsubsection{Local partial doubling integrals}
We also consider the local counterparts of the doubling integrals as vector-valued integrals (namely the Gelfand--Pettis integral; see, for example, \cite[Chapter 14]{garrett2018modern} for a precise definition).

Let $g_0, h_0 \in G(\calF_v)$. In the spirit of Corollary \ref{coro:partialdou}, we define the \emph{local partial doubling integrals} as follows:
\begin{align} \label{eq:defnpartialleftt}
Z_v^{\dsuit, \leftt}(f^{\Sieg}_{s,\chi,v}, \Psi^{\vee}_v; g_v) &:= \int_{G(\calF_v)} f^{\Sieg}_{s,\chi,v}(\imath^{\diamondsuit}(g_v, h_v)) \chi_v^{-1}(\det h_v) \pi_v^{\vee}(h_v) \Psi^{\vee}_v \dif h_v \nonumber \\
&= \int_{G(\calF_v)} f^{\Sieg}_{s,\chi,v}(\imath^{\diamondsuit}(g_0, h_v)) \chi_v^{-1}(\det h_v) \pi^{\vee}_v(g_v g_0^{-1} h_v) \Psi_v^{\vee} \dif h_v \quad \in \pi_{v}^{\vee},    
\end{align}
and
\begin{align} \label{eq:defnpartialrightt}
Z_v^{\dsuit, \rightt}(f^{\Sieg}_{s,\chi,v}, \Psi_v; h_v) &:= \int_{G(\calF_{v})} f^{\Sieg}_{s,\chi,v}(\imath^{\diamondsuit}(g_v, h_v)) \pi_v(g_v) \Psi_v \dif g_v \nonumber \\
&= \chi_v^{-1}(\det h_0) \int_{G(\calF_v)} f^{\Sieg}_{s,\chi,v}(\imath^{\diamondsuit}(g_v, h_0)) \pi_v(h_v h_0^{-1} g_v) \Psi_v \dif g_v \quad \in \pi_{v}.    
\end{align}
Here, the integrals $Z_v^{\dsuit, \leftt}$ and $Z_v^{\dsuit, \rightt}$ are understood as vectors in $\pi_{v}^{\vee}$ and $\pi_v$, respectively; the second equalities in \eqref{eq:defnpartialleftt} and \eqref{eq:defnpartialrightt} follow from the same change of variables used in Corollary \ref{coro:partialdou}. \footnote{We remark that $Z_v^{\dsuit, \rightt}(f^{\Sieg}_{s,\chi,v}, \Psi_v; - )$ is precisely the integral $F_{\Psi}^{\prime}(f^{\Sieg}_{s,\chi,v}; s, -)$ in the notation of \cite[Section 3B3]{wan2015families}.} By the definition of vector-valued integrals, they satisfy:
\begin{equation} \label{eq:localdoublingleft}
\lrangle{\Psi_{v}, Z_v^{\dsuit, \leftt}(f^{\Sieg}_{s,\chi,v}, \Psi^{\vee}_v; g_v)} = \int_{G(\calF_v)} f^{\Sieg}_{s,\chi,v}(\imath^{\diamondsuit}(g_0, h_v)) \chi_v^{-1}(\det h_v) \lrangle{\Psi_{v}, \pi_v^{\vee}(g_v g_0^{-1} h_v) \Psi_v^{\vee}} \dif h_v,    
\end{equation}
and
\begin{equation} \label{eq:localdoublingright}
\lrangle{Z_v^{\dsuit, \rightt}(f^{\Sieg}_{s,\chi,v}, \Psi_v; h_v), \Psi_{v}^{\vee}} = \chi_v(\det h_v h_0^{-1}) \int_{G(\calF_v)} f^{\Sieg}_{s,\chi,v}(\imath^{\diamondsuit}(g_v, h_0)) \lrangle{\pi_v(h_v h_0^{-1} g_v) \Psi_v, \Psi_v^{\vee}} \dif g_v.  
\end{equation}

\begin{assertion} \label{assertion:localdoubling}
With our choice of local Siegel Eisenstein sections $f_{s,\chi,v}^{\Sieg}$ at each place $v$ of $\calF$ as in Section \ref{sec:localdoubl}, there exists a scalar $\scZ^{\dsuit, \rightt}_v(f^{\Sieg}_{v,s,\chi}, \pi_v) \in \CC$, independent of $\Psi_v \in \pi_v$, such that:
\[
Z_v^{\dsuit, \rightt}(f^{\Sieg}_{s,\chi,v}, \Psi_v; \bfone_{m+n+1}) = \scZ^{\dsuit, \rightt}_v(f^{\Sieg}_{v,s,\chi}, \pi_v) \cdot \Psi_{v}
\]
for all $\Psi_v \in \pi_v$.
\end{assertion}
%We remark that $g_0$ and $h_0$ can certainly be $1 \in G(\calF_v)$, and actually this is the only case we shall use. In other applications, a flexibility to choose $g_0, h_0 \in G(\calF_v)$ may be convenient for computations. By the uniqueness of $G$-invariant pairings between $\pi$ and $\pi^{\vee}$ granted by \eqref{ass:multiplicityone}, we see the local quotients is independent of the choice of local vectors $\Psi_v$ and $\Psi_v^{\vee}$, but only depend on the local Siegel Eisenstein section $f^{\Sieg}_{v,s,\chi}$ and the local representation $\pi_v$. We therefore denote
%$$
%\scZ^{\dsuit}_v(f^{\Sieg}_{v,s,\chi}, \pi_v) := \frac{Z^{\dsuit}_v(f^{\Sieg}_{s,\chi,v}, \Psi_v, \Psi_{v}^{\vee})}{\lrangle{\Psi_v, \Psi^{\vee}_v}}
%$$
%for simplicity.

\subsection{Doubling method, à la Garrett} \label{sec:garrettdou}
The primary method for explicitly constructing Klingen Eisenstein series is via pullbacks of Siegel Eisenstein series. This approach generalizes the doubling method à la Piatetski-Shapiro and Rallis; it was initiated by Garrett in \cite{garrett1984pullbacks, garrett1989integral} and was later further developed by Shimura \cite{shimura1997euler}.

We define the \emph{pullback integral} formally as the vector-valued integral
\begin{equation} \label{eq:scalarpb}
F^{\heartsuit}(f^{\Sieg}_{s,\chi}, \Phi; g_0) := \int_{H(\AA)} f^{\Sieg}_{s,\chi}(\imath^{\hsuit}(g_0, h)) \chi^{-1}(\det h) \pi(h)\Phi \dif h,
\end{equation}
with values in $V_{\pi}$, for any Siegel section $f^{\Sieg}_{s,\chi} \in I^{\Sieg}(\chi, s)$ and any $g_0 \in G^{\hsuit}(\AA)$. Composing this with $\ev_1$ (the evaluation map at $\bfone_{m+n+1}$ from $V_{\pi}$ to $\CC$), we formally have
\[
F^{\heartsuit, 1}(f^{\Sieg}_{s,\chi}, \Phi; g_0) := \int_{H(\AA)} f^{\Sieg}_{s,\chi}(\imath^{\hsuit}(g_0, h)) \chi^{-1}(\det h) \Phi(h) \dif h.
\]

\begin{proposition}\label{prop:pbkling}
With the above notations:
\begin{enumerate}[label = \rm (\arabic*)]
    \item The function $g_0 \mapsto F^{\heartsuit}(f^{\Sieg}_{s,\chi}, \Phi; g_0)$ defines a Klingen Eisenstein section in $I^{\Kling}(\sigma, \chi, s)$ whenever the integral exists. The integral $F^{\heartsuit, 1}(f^{\Sieg}_{s,\chi}, \Phi; g_0)$ converges for $(s,g_0)$ in compact subsets of $\{\Re(s) > (N+1)/2 \} \times G^{\hsuit}(\AA)$. \footnote{There appears to be a typographical error in the convergence range stated in \cite[Proposition 3.5 (ii)]{wan2015families}. Specifically, $\{\Re(z) > r+s+1/2\}$ should be $\{\Re(z) > (r+s+1)/2\}$ in the notation of \cite{wan2015families}.}
    \item Moreover, we have the identity:
\[
\int_{H(\calF) \bs H(\AA)} E^{\Sieg}(f^{\Sieg}_{s,\chi}; \imath^{\hsuit}(g_0, h)) \barchi(\det h) \Phi(h) \dif h = E^{\Kling}(F^{\heartsuit}(f^{\Sieg}_{s,\chi}, \Phi; -); g_0).
\]
\end{enumerate}
\end{proposition}

\begin{proof}
This is \cite[Proposition 3.5]{wan2015families}. The convergence issue is discussed in the proof therein. The fact that $F^{\heartsuit, 1}$ defines a Klingen Eisenstein section is well-known and can be directly verified. The reader may refer to the proof of \cite[Theorem 2.6]{zhang2019fourier} for the case $H = \rmU(2,0)$; the computation there can be generalized to the broader case. Note that as we have assumed $\chi$ is a unitary Hecke character, we have $\chi^{-1} = \barchi$.
\end{proof}

\subsection{Reduction to the cuspidal GGP period integral}
We write
\begin{equation} \label{eq:periodKLE}
    \calP^{\Kling}_{\Phi, \Psi} := \calP^{\Kling}(\Phi, \Psi, \chi, s) := 
    \calP(\Psi, E^{\Kling}(F^{\heartsuit}(f^{\Sieg}_{s,\chi}, \Phi; -)))
\end{equation}
for simplicity. By Proposition \ref{prop:pbkling}, 
\[
\calP^{\Kling}_{\Phi, \Psi} = \int_{[G]} \int_{[H]} \Psi(g) E^{\Sieg}(f^{\Sieg}_{s,\chi}; \imath^{\hsuit}(\jmath^{\flat}(g), h)) \chi^{-1}(\det h) \Phi(h) \dif h \dif g .
\]
Interchanging the two integrals, we isolate the inner part:
\[
\calP^{\Kling}_{\Phi, \Psi} = \int_{[H]} \left( \int_{[G]} \Psi(g) E^{\Sieg}(f^{\Sieg}_{s,\chi}; \imath^{\heartsuit}(\jmath^{\flat}(g), h)) \dif g \right) \chi^{-1}(\det h) \Phi(h) \dif h .
\]
By the observation in \eqref{eq:basiscomp}, we have
\[
\calP^{\Kling}_{\Phi, \Psi} = \int_{[H]} \left( \int_{[G]} \Psi(g) E^{\Sieg}(f^{\Sieg}_{s,\chi}; \imath^{\dsuit}(g,\jmath^{\sharp}(h))) \dif g \right) \chi^{-1}(\det h) \Phi(h) \dif h .
\]
The inner integral is precisely $Z^{\dsuit, \rightt}(f^{\Sieg}_{s,\chi}, \Psi, \jmath^{\sharp}(h))$. By Corollary \ref{coro:partialdou} (2), setting $h_0 = 1$ there, we obtain
\begin{align*}
 \calP^{\Kling}_{\Phi, \Psi} &= \int_{[H]} \left( \chi(\det \jmath^{\sharp}(h)) \int_{G(\AA)} f^{\Sieg}_{s,\chi}( \imath^{\diamondsuit}(g, 1)) \Psi(\jmath^{\sharp}(h) g) \dif g \right) \chi^{-1}(\det h) \Phi(h) \dif h   \\
 &= \int_{[H]} \int_{G(\AA)} f^{\Sieg}_{s,\chi}( \imath^{\diamondsuit}(g, 1)) \Psi(\jmath^{\sharp}(h) g) \Phi(h) \dif g \dif h.
\end{align*}
Interchanging the two integrals once more yields
\begin{align*}
    \calP^{\Kling}_{\Phi, \Psi} = \int_{G(\AA)} f^{\Sieg}_{s,\chi}(\imath^{\diamondsuit}(g,1))
\left( \int_{[H]} \Phi(h) (\pi(g)\Psi)(\jmath^{\sharp}(h)) \dif h \right) \dif g.
\end{align*}
The inner integral is a \emph{cuspidal GGP period integral}. In general, for $\Phi \in \sigma$ and $\Psi \in \pi$, we define the cuspidal GGP period integral of $\Psi$ with $\Phi$ as
\[
\calP^{\sharp}(\Phi, \Psi) := \int_{[H]} \Phi(h) \Psi(\jmath^{\sharp}(h)) \dif h,
\]
which converges absolutely by the cuspidality of $\Phi$ and $\Psi$. To sum up, we have proved the following:

\begin{proposition} \label{prop:redcusp}
Let $\calP^{\Kling}_{\Phi, \Psi}$ be the GGP period integral defined in \eqref{eq:periodKLE}. Then
\[
\calP^{\Kling}_{\Phi, \Psi} = \int_{G(\AA)} f^{\Sieg}_{s,\chi}(\imath^{\dsuit}(g,1)) \calP^{\sharp}(\Phi, \pi(g)\Psi) \dif g.
\]
\end{proposition}

\subsection{Break into local integrals}
Before conducting concrete computations, we make the following preparations.

\subsubsection{Contragredient, conjugation, and MVW involutions}
Given an irreducible cuspidal automorphic representation $\pi \subseteq \calA_{0}(G(\calF) \bs G(\AA))$ and its complex conjugation $\bar{\pi} \subseteq \calA_{0}(G(\calF) \bs G(\AA))$, which is isomorphic to the contragredient representation $\pi^{\vee}$ \footnote{See, for example, \cite[Proposition 8.9.6]{goldfeld2011automorphici} for the case of $\GL_2(\AA_{\QQ})$. The proof generalizes to any reductive group.}, we fix factorizations:
\[
\fac_{\pi}: \pi \cong \otimes^{\prime}_{v \in \scV_{\calF}} \pi_v, \quad \fac_{\bar{\pi}}: \bar{\pi} \cong \otimes^{\prime}_{v \in \scV_{\calF}} \pi^{\vee}_v,
\]
where the restricted tensor product is taken with respect to spherical elements $\Psi_{v}^{\circ}$ and $\Psi_{v}^{\vee, \circ}$ at the places $v$ where $\pi_v$ is unramified.

Besides complex conjugation, another model for $\pi^{\vee}$ is established by the MVW involution \cite{moeglin1987correspondances}. Consider a unitary group $\rmU(V, \phi)$. There exists a unique element $\delta \in \Aut_{\calK}(V)$ such that $\phi(\delta v, \delta w) = \phi(w, v)$ for all $v,w \in V$. Conjugation by $\delta$ gives an automorphism of the group $\rmU(V, \phi)$:
\[
(-)^{\dagger}: g \mapsto \delta g \delta^{-1}.
\]
We fix an MVW involution $\dagger$ on $\bfG$ that stabilizes $H$. For every $v \in \scV_{\calF}^{\spl}$, we fix a standard isomorphism $\bfG(\calF_v) \cong \GL_{n}(\calF_v) \times \GL_{n+1}(\calF_v)$, under which $\dagger$ coincides with the transpose-inverse. For $\Psi = \otimes_{v}^{\prime} \Psi_v \in \pi$, the function $\Psi^{\dagger}(g) := \Psi(g^{\dagger})$ belongs to $\pi^{\vee}$. By a change of variables, we have:
\begin{equation} \label{eq:MVWperiod}
\calP^{\sharp}(\Phi^{\dagger}, \Psi^{\dagger}) = \calP^{\sharp}(\Phi, \Psi).
\end{equation}
We also define $g^{\ddagger} := \delta^{-1} g \delta$ and $\Psi^{\ddagger}(g) := \Psi(g^{\ddagger})$. Since $(\Psi^{\ddagger})^{\dagger} = \Psi$, it follows that $\Psi^{\ddagger} \in \pi^{\vee}$.

\subsubsection{Ichino--Ikeda formula}
Recall $\bfG = H \times G$. Let $\Pi = \sigma \boxtimes \pi$ be the cuspidal representation of $\bfG(\AA)$.

The product $L$-series associated with $\sigma$ and $\pi$ is defined as
\[
L(s, \sigma \times \pi) := L^{\JPSS}(s, \BC(\sigma) \times \BC(\pi)),
\]
where $\BC(\sigma)$ (resp. $\BC(\pi)$) is the functorial lift of $\sigma$ (resp. $\pi$) to an automorphic representation of $\GL_{m+n}(\KK)$ (resp. $\GL_{m+n+1}(\KK)$). The right-hand side is the $L$-factor defined by Jacquet, Piatetski-Shapiro, and Shalika in \cite{jacquet1983rankin}. Let $L(s, \sigma, \Ad)$ denote the adjoint $L$-series for $\sigma$.

Assuming that both $\sigma$ and $\pi$ are tempered, we set
\[
\scL(\sigma \times \pi) = \dfrac{L(\frac{1}{2}, \sigma \times \pi)}{L(1, \sigma, \Ad) L(1, \pi, \Ad)} \prod_{i=1}^{m+n+1} L(i, \epsilon_{\calK/\calF}^{i})
\]
and
\[
\scL(\sigma_v \times \pi_v) = \dfrac{L(\frac{1}{2}, \sigma_v \times \pi_v)}{L(1, \sigma_v, \Ad) L(1, \pi_v, \Ad)} \prod_{i=1}^{m+n+1} L(i, \epsilon_{\calK_v/\calF_v}^{i}) 
\]
for any $v \in \scV_{\calF}$. 

Let $\Xi$ be a cusp form on $\bfG$; we define the integral
\[
\calP(\Xi) := \int_{[\bfH]} \Xi(h) \dif h.
\]
If $\Xi = \otimes^{\prime}_v \Xi_v \in \Pi$ and $\Xi^{\prime} = \otimes^{\prime}_v \Xi^{\prime}_v \in \Pi^{\vee}$ are factorizable, we define the \emph{local Ichino--Ikeda integral}
\begin{equation} \label{eq:definitioniiintegral}
 I(\Xi_v, \Xi_v^{\prime}) := \int_{H(F_v)} \lrangle{\Pi_v(\Delta^{\sharp}(h_v))\Xi_v, \Xi^{\prime}_v}_{\Pi_v} \dif h_v.   
\end{equation}
This integral is convergent if $\Pi_v$ is tempered. Note that $\calP(\Psi \boxtimes \Phi) = \calP^{\sharp}(\Phi, \Psi)$.

Recently, there has been significant progress on the global Gan--Gross--Prasad conjecture for unitary groups (see, e.g., \cite{ichino2010period, harris2014refined, beuzart2021isolation, beuzart2022global}). We adopt the version reinterpreted in \cite[Theorem 4.2]{hsieh2023fivevariablepadiclfunctionsu3times} \footnote{Following our convention in Section \ref{sec:measures}, the constant $C_H$ in \cite[Theorem 4.2]{hsieh2023fivevariablepadiclfunctionsu3times} is taken to be $1$ here.}.

\begin{theorem}[Ichino--Ikeda formula] \label{thm:ikedaichino}
Let $\Pi$ be an irreducible tempered cuspidal automorphic representation of $\bfG(\AA)$. If $\Xi = \otimes^{\prime}_v \Xi_v \in \Pi$ and $\Xi^{\prime} = \otimes^{\prime}_v \Xi^{\prime}_v \in \Pi^{\vee}$ are factorizable, then
\[
\dfrac{\calP(\Xi)\calP(\Xi^{\prime})}{(\Xi, \Xi^{\prime})_{\Pet}} = \dfrac{\scL(\sigma \times \pi)}{2^{\varkappa_{\sigma}+\varkappa_{\pi}}} \prod_{v \in \mathscr{V}_{\calF}} \dfrac{I(\Xi_v, \Xi_v^{\prime})}{\scL(\sigma_v \times \pi_v) \lrangle{\Xi_v, \Xi_v^{\prime}}_{\Pi_v} },
\]
where $2^{\varkappa_{\pi}}$ (resp. $2^{\varkappa_{\sigma}}$) is the order of the component group associated with the $L$-parameter of $\pi$ (resp. $\sigma$).
\end{theorem}

\subsubsection{Reduction to local Ichino--Ikeda integrals}
We square the GGP period integral, using Proposition \ref{prop:redcusp}, to obtain
\begin{align*}
(\calP^{\Kling}_{\Phi, \Psi})^2 
& = \int_{G(\AA) \times G(\AA)} f^{\Sieg}_{s,\chi}( \imath^{\diamondsuit}(g_1,1)) f^{\Sieg}_{s,\chi}( \imath^{\diamondsuit}(g_2,1)) \calP^{\sharp}(\Phi, \pi(g_1)\Psi) \calP^{\sharp}(\Phi, \pi(g_2)\Psi) \dif g_1 \dif g_2 \\
&\xlongequal{\eqref{eq:MVWperiod}} \int_{G(\AA) \times G(\AA)} f^{\Sieg}_{s,\chi}( \imath^{\diamondsuit}(g_1,1)) f^{\Sieg}_{s,\chi}( \imath^{\diamondsuit}(g_2,1)) \calP^{\sharp}(\Phi, \pi(g_1)\Psi) \calP^{\sharp}(\Phi^{\dagger}, (\pi(g_2)\Psi)^{\dagger}) \dif g_1 \dif g_2 .
\end{align*}

\begin{remark}[On complex conjugation]
Using a naive approach with complex conjugation, we have:
\begin{align*}
\abs{\calP^{\Kling}_{\Phi, \Psi}}^2 &= \calP^{\Kling}_{\Phi, \Psi} \cdot \overline{\calP^{\Kling}_{\Phi, \Psi}} \\
&= \int_{G(\AA) \times G(\AA)} f^{\Sieg}_{s,\chi}( \imath^{\diamondsuit}(g_1,1)) \overline{f^{\Sieg}_{s,\chi}( \imath^{\diamondsuit}(g_2,1))} \calP^{\sharp}(\Phi, \pi(g_1)\Psi) \overline{\calP^{\sharp}(\Phi, \pi(g_2)\Psi)} \dif g_1 \dif g_2 \\
& = \int_{G(\AA) \times G(\AA)} f^{\Sieg}_{s,\chi}( \imath^{\diamondsuit}(g_1,1)) \overline{f^{\Sieg}_{s,\chi}( \imath^{\diamondsuit}(g_2,1))} \calP^{\sharp}(\Phi, \pi(g_1)\Psi) \calP^{\sharp}(\bar{\Phi}, \pi^{\vee}(g_2)\bar{\Psi}) \dif g_1 \dif g_2.
\end{align*}
While this is acceptable for purely automorphic computations, complex conjugation behaves poorly in $p$-adic families, making it unsuitable for $p$-adic interpolation.
\end{remark}

We now focus on the product $\calP^{\sharp}(\Phi, \pi(g_1)\Psi) \calP^{\sharp}(\Phi^{\dagger}, (\pi(g_2)\Psi)^{\dagger})$. Applying Theorem \ref{thm:ikedaichino} with $\Xi := \Phi \boxtimes \pi(g_1)\Psi$ and $\Xi^{\prime} := \Phi^{\dagger} \boxtimes (\pi(g_2)\Psi)^{\dagger}$, we have:
\begin{multline} \label{eq:useii}
\dfrac{\calP^{\sharp}(\Phi, \pi(g_1)\Psi) \calP^{\sharp}(\Phi^{\dagger}, (\pi(g_2)\Psi)^{\dagger})}{\lrangle{\Phi, \Phi^{\dagger}}_{\Pet} \lrangle{\pi(g_1)\Psi, (\pi(g_2)\Psi)^{\dagger}}_{\pi, \Pet}} \\
= \dfrac{\scL(\sigma \times \pi)}{2^{\varkappa_{\sigma}+\varkappa_{\pi}}} \prod_{v \in \mathscr{V}_{\calF}} \dfrac{I(\Phi_v \boxtimes \pi(g_1)\Psi_v, \Phi_v^{\dagger} \boxtimes (\pi(g_2)\Psi)^{\dagger}_v)}{\scL(\sigma_v \times \pi_v) \lrangle{\Phi_v \boxtimes \pi(g_1)\Psi_v, \Phi^{\dagger}_v \boxtimes (\pi_v(g_2)\Psi_v)^{\dagger}}_{\Pi_v}}.    
\end{multline}
Hereafter, we denote $\aabs{\Phi}_{\sigma, \dagger}^{2} := \lrangle{\Phi, \Phi^{\dagger}}_{\Pet}$. By the uniqueness of the $\bfG$-invariant pairings $\Pi \times \Pi^{\vee} \rightarrow \CC$ (Assumption \ref{ass:multiplicityone}), the local ratio
\begin{equation} \label{eq:iiindependent}
\scI_v (\sigma_v, \pi_v) := \dfrac{I(\Phi_v \boxtimes \pi(g_1)\Psi_v, \Phi_v^{\dagger} \boxtimes (\pi_v(g_2)\Psi_v)^{\dagger})}{\lrangle{\Phi_v \boxtimes \pi(g_1)\Psi_v, \Phi^{\dagger}_v \boxtimes (\pi_v(g_2)\Psi_v)^{\dagger}}_{\Pi_v}}
\end{equation}
depends only on the local representation $\Pi_v$ and not on the particular vectors. Thus, we have:
\begin{equation} \label{eq:useiishort}
    \dfrac{\calP^{\sharp}(\Phi, \pi(g_1)\Psi) \calP^{\sharp}(\Phi^{\dagger}, (\pi(g_2)\Psi)^{\dagger})}{\aabs{\Phi}_{\sigma, \dagger}^{2} \lrangle{\pi(g_1)\Psi, (\pi(g_2)\Psi)^{\dagger}}_{\pi, \Pet}} = \dfrac{\scL(\sigma \times \pi)}{2^{\varkappa_{\sigma}+\varkappa_{\pi}}} \prod_{v \in \mathscr{V}_{\calF}} \scL(\sigma_v \times \pi_v)^{-1} \scI_v (\sigma_v, \pi_v).
\end{equation}

\subsubsection{Reduction to the doubling integrals}
To simplify the notation in this subsection, we denote the right-hand side of \eqref{eq:useiishort} by $(\star)$. We have seen that this term is independent of $g_1$ and $g_2$. Therefore, in the expression for $(\calP^{\Kling}_{\Phi, \Psi})^2$, we may pull it out of the integrals to obtain:
\begin{align*}
(\calP^{\Kling}_{\Phi, \Psi})^2 &= (\star) \cdot \aabs{\Phi}_{\sigma, \dagger}^{2} \int_{G(\AA) \times G(\AA)} f^{\Sieg}_{s,\chi}(\imath^{\diamondsuit}(g_1,1)) f^{\Sieg}_{s,\chi}(\imath^{\diamondsuit}(g_2,1)) \lrangle{\pi(g_1)\Psi, (\pi(g_2)\Psi)^{\dagger}}_{\pi, \Pet} \dif g_1 \dif g_2 \\
&= (\star) \cdot \aabs{\Phi}_{\sigma, \dagger}^{2} \int_{G(\AA)} \left( \int_{G(\AA)} f^{\Sieg}_{s,\chi}( \imath^{\diamondsuit}(g_1,1)) \lrangle{\pi(g_1)\Psi, (\pi(g_2)\Psi)^{\dagger}}_{\pi, \Pet} \dif g_1 \right) f^{\Sieg}_{s,\chi}(\imath^{\diamondsuit}(g_2,1)) \dif g_2.
\end{align*}

We emphasize to the reader that, from this point forward, we choose the \emph{specific} local Siegel Eisenstein sections $f_{s,\chi,v}^{\Sieg}$ defined in Section \ref{sec:localdoubl} to ensure the validity of Assertion \ref{assertion:localdoubling}.

The integral in the parentheses is precisely $Z^{\dsuit}(\Psi, (\pi(g_2)\Psi)^{\dagger}, f^{\Sieg}_{s,\chi})$ by Corollary \ref{coro:partialdou} (3) (setting $h_0=1$). Decomposing this into local doubling integrals via Corollary \ref{cor:doubintlocal}, we have:
\begin{align*}
 (\calP^{\Kling}_{\Phi, \Psi})^2 &= (\star) \cdot \aabs{\Phi}_{\sigma, \dagger}^{2} \int_{G(\AA)} \left( \prod_{v \in \scV_{\calF}} \int_{G(\calF_{v})} f^{\Sieg}_{s,\chi,v}( \imath^{\diamondsuit}(g_1,1)) \lrangle{\pi_v(g_1)\Psi_v, (\pi_v(g_2)\Psi_v)^{\dagger}} \dif g_1 \right) f^{\Sieg}_{s,\chi}(\imath^{\diamondsuit}(g_2,1)) \dif g_2 \\
 &\xlongequal{\eqref{eq:localdoublingright}} (\star) \cdot \aabs{\Phi}_{\sigma, \dagger}^{2} \int_{G(\AA)} \left( \prod_{v \in \scV_{\calF}} \lrangle{Z_v^{\dsuit, \rightt}(f^{\Sieg}_{s,\chi,v}, \Psi_v; 1), (\pi_v(g_2)\Psi_v)^{\dagger}} \right) f^{\Sieg}_{s,\chi}(\imath^{\diamondsuit}(g_2,1)) \dif g_2 \\
 &\xlongequal{\text{Assertion } \ref{assertion:localdoubling}} (\star) \cdot \aabs{\Phi}_{\sigma, \dagger}^{2} \int_{G(\AA)} \left( \prod_{v \in \scV_{\calF}} \scZ_v^{\dsuit, \rightt}(f^{\Sieg}_{s,\chi,v}, \pi_v) \lrangle{\Psi_v, (\pi_v(g_2)\Psi_v)^{\dagger}} \right) f^{\Sieg}_{s,\chi}(\imath^{\diamondsuit}(g_2,1)) \dif g_2 \\
 &\xlongequal{\eqref{eq:localglobalpeter}} (\star) \cdot \aabs{\Phi}_{\sigma, \dagger}^{2} \left( \prod_{v \in \scV_{\calF}} \scZ_v^{\dsuit, \rightt}(f^{\Sieg}_{s,\chi,v}, \pi_v)\right) \int_{G(\AA)} f^{\Sieg}_{s,\chi}(\imath^{\diamondsuit}(g_2,1)) \lrangle{\Psi, (\pi(g_2)\Psi)^{\dagger}}_{\pi, \Pet} \dif g_2.
\end{align*}
By a change of variables, we observe that $\lrangle{\Psi, (\pi(g_2)\Psi)^{\dagger}}_{\pi, \Pet} = \lrangle{\pi(g_2)\Psi, \Psi^{\ddagger}}_{\pi, \Pet}$. The remaining integral is then $Z^{\dsuit}(\Psi, \Psi^{\ddagger}, f^{\Sieg}_{s,\chi})$ by Corollary \ref{coro:partialdou} (3). Decomposing this into local doubling integrals again and repeating the process above, we obtain:
\[
(\calP^{\Kling}_{\Phi, \Psi})^2 = (\star) \cdot \aabs{\Phi}_{\sigma, \dagger}^{2} \aabs{\Psi}_{\pi, \ddagger}^{2} \cdot \prod_{v \in \scV_{\calF}} \scZ^{\dsuit, \rightt}_v(f^{\Sieg}_{v,s,\chi}, \pi_v)^{2},
\]
where $\aabs{\Psi}_{\pi, \ddagger}^{2} := \lrangle{\Psi, \Psi^{\ddagger}}_{\pi, \Pet}$. Since $\aabs{\Psi}_{\pi, \ddagger}^{2} = \aabs{\Psi}_{\pi, \dagger}^{2}$, we denote both by $\aabs{\Psi}_{\pi, \MVW}^{2}$ from now on.

To conclude, we have established the following decomposition:

\begin{theorem} \label{thm:breaklocal}
With notations as above, we have:
\begin{equation}
\dfrac{(\calP^{\Kling}_{\Phi, \Psi})^2}{\aabs{\Phi}_{\sigma, \MVW}^{2} \aabs{\Psi}_{\pi, \MVW}^{2}} = \dfrac{\scL(\sigma \times \pi)}{2^{\varkappa_{\sigma}+\varkappa_{\pi}}} \prod_{v \in \mathscr{V}_{\calF}} \scZ^{\dsuit, \rightt}_v(f^{\Sieg}_{v,s,\chi}, \pi_v)^2 \cdot \dfrac{\scI_v (\sigma_v, \pi_v)}{\scL(\sigma_v \times \pi_v)}.
\end{equation}
\end{theorem}

\subsection{Unramified Computations}
We have the following standard results relating local integrals to local $L$-factors at unramified places. 

\begin{theorem} \label{thm:unrii}
Let $v \in \scV_{\calF}^{\ur}$. Then $\scI_v (\sigma_v, \pi_v) = \scL_{v}(\sigma_v \times \pi_v)$.
\end{theorem}
\begin{proof}
This is \cite[Theorem 2.12]{harris2014refined} for the case where $v$ is unramified in $\calK$. When $v$ ramifies in $\calK$, the result is proved in \cite[Appendix B]{hsieh2023fivevariablepadiclfunctionsu3times} for the case $H = \rmU(2,0)$ and $G = \rmU(3,0)$, which generalizes to unitary groups of arbitrary signature.
\end{proof}

\begin{theorem}[{\cite[Lemma 4.5]{wan2015families}}] \label{thm:urrdou}
Let $v \in \scT_{\calF}^{\ur}$. Then
\[
\scZ^{\dsuit, \rightt}_v(f^{\Sieg, \sph}_{v,s,\chi}, \pi_v) = d_{N+1,v}(s, \barchi^{\sfc}_v)^{-1} \cdot L_{v} \left(s+\dfrac{1}{2}, \pi_v^{\vee}, \barchi^{\sfc}_v \right),
\]
where:
\begin{itemize}
    \item $d_{n,v}(s,\chi_v) := \prod_{j=1}^{n} L_v(2s+j, \chi^{\calF} \cdot \epsilon_{\calK/\calF}^{n-j})$ is the product of local $L$-factors for Hecke characters over $\calF$. \footnote{The correction factor $d_{N+1,v}(s, \barchi_v) = d_{N+1,v}(s, \barchi^{\sfc}_v)$ because, by definition, it depends only on the restriction of $\barchi$ to $\calF$. We use the conjugation $\sfc$ to maintain a natural appearance in the identity.}
    \item $L_{v}(s+\dfrac{1}{2}, \pi_v, \chi_v) = L_v(s, \BC(\pi_v) \otimes \chi_v \circ \det)$, where the right-hand side is the standard local Godement--Jacquet $L$-factor and $\BC(\pi_v)$ is the local base change from $G(\calF_v)$ to $\GL_{N+1}(\calK_v)$.
\end{itemize}
\end{theorem}

Combining Theorem \ref{thm:breaklocal}, \ref{thm:unrii}, and \ref{thm:urrdou}, we obtain the following result. Here, we apply Theorem \ref{thm:urrdou} only to the primes in $\scU_{\calF}$, due to the considerations in Remark \ref{rem:modificationsections}.

\begin{theorem} \label{thm:unrret}
Assume that \ref{item:irred} holds for $\sigma$ and $\pi$. With local Siegel Eisenstein sections chosen as in Section \ref{sec:localdoubl}, we have:
\begin{multline*}
\dfrac{\left(\calP^{\Kling}_{\Phi, \Psi}\right)^2}{\aabs{\Phi}_{\sigma, \MVW}^2 \aabs{\Psi}_{\pi, \MVW}^2} = \dfrac{1}{2^{\varkappa_{\sigma}+\varkappa_{\pi}}} \cdot \scL_{\scV_{\calF}^{\ur}}(\sigma \times \pi) L_{\scU_{\calF}^{\ur}}\left(s+\dfrac{1}{2}, \pi^{\vee}, \barchi^{\sfc} \right)^{2} \\
\times \prod_{v \in \scU_{\calF}^{\ur}} d_{N+1,v}(s, \barchi_v)^{-2} \prod_{v \in \mathscr{V}_{\calF}^{\bad}} \scI_v (\sigma_v, \pi_v) \prod_{v \in \scU_{\calF}^{\bad}} \scZ^{\dsuit, \rightt}_v(f^{\Sieg}_{v,s,\chi}, \pi_v)^2,
\end{multline*}
where $\scL_{\scV_{\calF}^{\ur}}(\cdots)$ and $L_{\scU_{\calF}^{\ur}}(\cdots)$ denote the products of the local factors $\scL_{v}(\cdots)$ and $L_{v}(\cdots)$ for $v$ in $\scV_{\calF}^{\ur}$ or $\scU_{\calF}^{\ur}$, respectively.
\end{theorem}

\subsection{Further assumptions}
In the following sections, we shall:
\begin{enumerate}
    \item Compute the local Ichino--Ikeda integrals at bad places $v \in \scV_{\calF}^{\bad}$ (Sections \ref{sec:iiarchi} through \ref{sec:iispl}).
    \item Select Siegel Eisenstein sections at bad places $v \in \scS_{\calF}^{\bad}$ such that the resulting Siegel Eisenstein series (and consequently the Klingen Eisenstein series via pullback) possess $p$-adic interpolatable $q$-expansions (see Section \ref{sec:padicggp}).
\end{enumerate}
To facilitate step (1), we impose additional assumptions on the representations $\sigma$ and $\pi$.

\begin{assumption}[Weight interlacing assumption] \label{ass:wtinterlacing}
We assume that for any $v \in \scV_{\calF}^{\infty}$, the representations $\sigma_v$ and $\pi_v^{\vee}$ are irreducible discrete series representations of $\rmU_{m,n}$ and $\rmU_{m,n+1}$ respectively, and satisfy the Gan--Gross--Prasad weight interlacing property (to be recalled in Section \ref{sec:iiarchi}).
\end{assumption}

\begin{assumption}[Splitting assumption] \label{ass:spl}
We assume that every $v \in \scV_{\calF}^{\ram}$ splits in $\calK$.
\end{assumption}

\begin{assumption}[Disjointly ramified assumption] \label{ass:ram}
We assume that for every $v \in \scV_{\calF}^{\ram}$, at least one of $\pi_v$ or $\sigma_v$ is unramified (i.e., they are not both ramified).
\end{assumption}

Assumption \ref{ass:spl} is a version of the "Heegner hypothesis"; it is well-known that there exist infinitely many imaginary quadratic extensions $\calK$ of $\calF$ satisfying this condition. If $\sigma$ and $\pi$ are viewed as representations associated with modular forms of tame levels $N_\sigma$ and $N_{\pi}$, Assumption \ref{ass:ram} is equivalent to requiring $N_{\sigma}$ and $N_{\pi}$ to be coprime.

For simplicity, where the context is clear, we shall adopt the following local notation:
\begin{itemize}
    \item $H := H(\calF_v)$, $G := G(\calF_v)$, and similarly for other groups.
    \item $\sigma := \sigma_v$, $\pi := \pi_v$, and similarly for local representations.
    \item The extension $\calK_v / \calF_v$ is denoted by $K/F$. We write $\fro, \frp, \varpi$ for the ring of integers, the maximal ideal, and a fixed uniformizer of $F$, respectively, and let $q$ be the cardinality of the residue field $\fro/\frp$. The absolute value $\abs{-}$ on $F$ is normalized by $\abs{\varpi} = q^{-1}$.
\end{itemize}

%% file: AU02-Archimedian.tex
\section{Local Ichino--Ikeda integrals at archimedean places} \label{sec:iiarchi}
In this section, we consider $v \in \scV_{\calF}^{\infty}$, an archimedean place of $\calF$.

\subsection{Weight interlacing property}
As promised, we first introduce the weight interlacing assumption (Assumption \ref{ass:wtinterlacing}). The primary reference for this material is \cite{he2017gan}.

Consider the real Lie group $\rmU_{m,n}$ as defined in \eqref{eq:realumn}. A maximal compact subgroup $K \cong \rmU_{m,0} \times \rmU_{n,0}$ is given by $K = \{g \in \rmU_{m,n}: g = g^{-\star}\}$. Let $\frg_0$ be the Lie algebra of $\rmU_{m,n}$ and $\frt_0$ be the Cartan subalgebra of $\frg_0$ consisting of diagonal matrices. Let $\frg := \frg_0 \otimes_{\RR} \CC$ and $\frt := \frt_0 \otimes_{\RR} \CC$ be their complexifications. We identify $\frt^{\ast}$ with $\CC^{m+n}$ via the basis $\varepsilon_1, \ldots, \varepsilon_{m+n}$ defined by:
\[
\varepsilon_i(\diag[a_1, \ldots, a_{m+n}]) = a_i,
\]
and define a bilinear form $\lrangle{-,-}: \frt^{\ast} \times \frt^{\ast} \rightarrow \CC$ by 
\[
\lrangle{x,y} := \sum_{i=1}^{m+n} x_i y_i
\]
for $x, y \in \frt^{\ast} \cong \CC^{m+n}$. Let $\Delta$ be the set of roots of $\frt$ in $\frg$, namely $\Delta = \{\pm(\varepsilon_i - \varepsilon_j) : 1 \leq i < j \leq m+n \}$. Let $\Delta_{\rmc}$ be the set of compact roots and take the positive system $\Delta_{\rmc}^{+}$ as:
\[
\Delta_{\rmc}^{+} = \{ \varepsilon_i - \varepsilon_j : 1 \leq i < j \leq m \} \cup \{ \varepsilon_i - \varepsilon_j : m+1 \leq i < j \leq m+n \}.
\]
The discrete series representations of $\rmU_{m,n}$ are parameterized by Harish-Chandra parameters that are dominant for $\Delta_{\rmc}^{+}$:
\[
\ullam = (\lambda_1^{+}, \ldots, \lambda_{m}^{+}, \lambda_{1}^{-}, \ldots, \lambda_{n}^{-}) \in \sqrt{-1} \frt_{0}^{\ast},
\]
where $\lambda_i^{\circ}, \lambda_j^{\bullet} \in \ZZ + \frac{m+n-1}{2}$, $\lambda_i^{\circ} \neq \lambda_j^{\bullet}$ for $\circ, \bullet \in \{+,-\}$ and all distinct $i, j$, with $\lambda_1^{+} > \cdots > \lambda_m^{+}$ and $\lambda_{1}^{-} > \cdots > \lambda_n^{-}$. We denote by $\DS(\ullam)$ the discrete series representation parameterized by $\ullam$. 

Let $\sigma = \DS(\ullam)$ be a discrete series of $\rmU_{m,n}$ and $\pi = \DS(\uleta)$ be a discrete series of $\rmU_{m+1,n}$. We denote the signs in $\uleta$ by $\oplus, \ominus$. We say $\sigma$ and $\pi$ satisfy the \emph{Gan--Gross--Prasad weight interlacing relation} if, when the components of $\ullam$ and $\uleta$ are lined up in descending order, the resulting sequence of signs consists only of the following adjacent pairs:
\[
(\oplus, +), (+, \oplus), (-, \ominus), (\ominus, -), (+, -), (-,+), (\oplus, \ominus), (\ominus, \oplus).
\]

The local Gan--Gross--Prasad conjecture for discrete series at archimedean places, now a theorem of He \cite{he2017gan}, is as follows.

\begin{theorem}[{\cite[Theorem 1.1]{he2017gan}}] \label{thm:heGGP}
The discrete series $\DS(\ullam)$ of $\rmU_{m,n}$ appears as a subrepresentation of the restriction $\DS(\uleta)|_{\rmU_{m,n}}$ if and only if $\DS(\ullam)$ and $\DS(\uleta)$ satisfy the Gan--Gross--Prasad weight interlacing relation.
\end{theorem}

\section{Local Ichino--Ikeda Integrals}
According to the adjointness of the corresponding functors, we have
\[
\Hom_{\bfH}(\Pi^{\vee}, \CC) = \Hom_{H}\left( \pi^{\vee}|_{H} \boxtimes \sigma^{\vee}, \CC \right) 
= \Hom_{H}(\pi^{\vee}|_{H}, \Hom_{H}(\sigma^{\vee}, \CC))
= \Hom_{H}(\pi^{\vee}|_{H}, \sigma)
\]
Note that since $\pi$ is an irreducible admissible representation of $G$, its dual representation $\pi^{\vee}$ is also irreducible and admissible; in the last equality, we identify $(\sigma^{\vee})^{\vee}$ with $\sigma$. By Schur's Lemma, we see that $\Hom_{\bfH}(\Pi^{\vee}, \CC)$ is nonzero if and only if $\sigma$ is a subrepresentation of $\pi^{\vee}|_{H}$. Therefore, by Theorem \ref{thm:heGGP}, this holds if and only if $\sigma$ and $\pi^{\vee}$ satisfy the GGP weight interlacing property.

Furthermore, the condition ``$\Hom_{\bfH}(\Pi^{\vee}, \CC) \neq \{0\}$'' implies the existence of a nontrivial $\bfH$-invariant functional $\Lambda^{\circ}: \Pi^{\vee} \rightarrow \CC$, such that for any $v \in \Pi^{\vee}$ and $h \in \bfH$, 
\[
\Lambda^{\circ}(h v) = h \Lambda^{\circ}(v) = \Lambda^{\circ}(v),
\]
noting that the $\bfH$-action on $\CC$ is trivial. Viewing $\Lambda^{\circ}$ as an element in the contragredient representation $(\Pi^{\vee})^{\vee}$, this means $\Lambda^{\circ}$ is an $\bfH$-invariant vector. We identify $(\Pi^{\vee})^{\vee}$ with $\Pi$, so under this identification, $\Lambda^{\circ}$ corresponds to an $\bfH$-invariant vector $\Xi^{\circ} \in \Pi$. Thus, directly from the definition, we have
$$
I(\Xi^{\circ}, \Xi^{\vee}) := \int_{H} \lrangle{\Pi(\Delta^{\sharp}(h))\Xi^{\circ}, \Xi^{\vee}}_{\Pi} \dif h = \lrangle{\Xi^{\circ}, \Xi^{\vee}}_{\Pi}
$$
for any $\Xi^{\vee} \in \Pi^{\vee}$. \footnote{The above discussion is an expansion of the assertion in the middle of \cite[p. 38]{hsieh2023fivevariablepadiclfunctionsu3times}. In particular, \cite[Eq. (1.1)]{hsieh2023fivevariablepadiclfunctionsu3times} is the special case of ``$\sigma$ and $\pi^{\vee}$ satisfying the GGP weight interlacing property'' for definite unitary groups $\rmU_{n,0} \times \rmU_{n-1, 0}$.} Therefore, we have the following proposition:

\begin{proposition}
Under Assumption \ref{ass:wtinterlacing}, for any $v \in \scV_{\calF}^{\infty}$, $\scI_v(\sigma_v, \pi_v) = 1$.
\end{proposition}

We point out that the weight interlacing assumption (i.e., Assumption \ref{ass:wtinterlacing}) ensures that the local Ichino--Ikeda integral is nonzero.

%% file: AU03-Split.tex
\section{Local Ichino--Ikeda integrals at split primes} \label{sec:iispl}
In this section, we treat the local Ichino--Ikeda integrals at places $v \in \scV_{\calF}^{\spl}$, i.e., finite places of $\calF$ that split in $\calK$. This includes places of $\calF$ above $p$.

In Section \ref{sec:unitarylocal}, we established isomorphisms
\[
\varrho_{w, m, n}: H \xrightarrow{\sim} \GL_{m+n}(F), \quad \varrho_{w, m, n+1}: G \xrightarrow{\sim} \GL_{m+n+1}(F).
\]
For convenience when dealing with local Ichino--Ikeda integrals, we adjust the isomorphisms $\varrho_{w, m, n}$ and $\varrho_{w, m, n+1}$ such that the following diagram commutes:
\[\begin{tikzcd}
	{H} && {G} \\
	{\GL_{m+n}(F)} && {\GL_{m+n+1}(F)}
	\arrow["{\jmath^{\sharp}}", hook, from=1-1, to=1-3]
	\arrow["{\varrho_{w, m, n}}"', from=1-1, to=2-1]
	\arrow["{\varrho_{w, m, n+1}}", from=1-3, to=2-3]
	\arrow["{\jmath^{\GL}}"', hook, from=2-1, to=2-3]
\end{tikzcd}\]
where
$$
\jmath^{\GL}: \GL_{m+n}(\calF_v) \hookrightarrow \GL_{m+n+1}(\calF_v), \quad g \mapsto \diag[g,1].
$$
Under the isomorphisms $\varrho_{w, m, n}$ and $\varrho_{w, m, n+1}$, we regard $\sigma$ and $\pi$ as representations of $\GL_{m+n}(F)$ and $\GL_{m+n+1}(F)$, respectively.

\subsection{Basic representation theory of $\GL_n(F)$}
In this subsection, we briefly review some basic representation theory of $\GL_n(F)$. For this subsection only, let $m, n, l$ denote general positive integers (unrelated to the unitary groups previously discussed).

We fix an additive character $\bpsi: F \rightarrow \CC^{\times}$ that is trivial on $\fro$ and nontrivial on $\frp^{-1}$. Let $\scS(F)$ be the space of locally constant, compactly supported functions on $F$. The Fourier transform of $\phi \in \scS(F)$ is defined by 
$$
\widehat{\phi}(y) = \int_F \phi(x) \bpsi(-xy) \dif x.
$$
The measure $\dif x$ is chosen such that $\widehat{\widehat{\phi}}(x)=\phi(-x)$. 

Let $\tau: F^{\times} \rightarrow \CC^{\times}$ be an arbitrary multiplicative character of $F^{\times}$. We let $\frg(\tau, \bpsi, y)$ denote the Gauss sum
$$
\frg(\tau, \bpsi, y) := \int_{\fro^{\times}} \tau(z) \bpsi(yz) \dif^{\times} z
$$
for $y \in F$. When $y=1$ and $\bpsi$ is clear from the context, we shall simply write $\frg(\tau)$.

\subsubsection{Whittaker models}
We extend $\bpsi$ to $U_n(F)$ by the rule 
$$
\bpsi(u) := \prod_{i=1}^{n-1} \bpsi(u_{i,i+1}), \quad \text{ for } u=(u_{i,j}) \in U_n(F).
$$
This is a generic character. Note that it is trivial on $U_n(\fro)$.

Let $\pi$ be an irreducible admissible representation of $\GL_{n}(F)$. It is called \emph{generic} if
$$
\Hom_{\GL_n(F)}(\pi, \Ind_{U_n(F)}^{\GL_n(F)} \bpsi) \neq 0.
$$
It is known that tempered representations are always generic. By Frobenius reciprocity, this means there exists a nonzero linear form $\lambda: V_{\pi} \rightarrow \CC$ such that
$$
\lambda(\pi(u)v) = \bpsi(u) \lambda(v), \quad v \in V_{\pi}, u \in U_n(F).
$$
For a generic $\pi$, the space $\Hom_{\GL_n(F)}(\pi, \Ind_{U_n(F)}^{\GL_n(F)} \bpsi)$ is one-dimensional. The \emph{Whittaker model} $\calW_{\bpsi}(\pi)$ of $\pi$ with respect to $\bpsi$ is defined as
\[
\calW_{\bpsi}(\pi) := \{W_v: \GL_n(F) \rightarrow \CC : W_v(g) = \lambda(\pi(g)v), v \in V_{\pi} \}.
\]
Then $\calW_{\bpsi}(\pi)$ is independent of the choice of $\lambda$, and for $u \in U_n(F)$, $g \in \GL_n(F)$, 
$$
W_v(ug) = \bpsi(u)W_v(g), \quad W_v(g) = W_{\pi(g)v}(\bfone_n).
$$
The map $v \mapsto W_v$ yields an isomorphism $\sfW_{\pi}: V_{\pi} \rightarrow \calW_{\bpsi}(\pi)$. One can define an invariant perfect pairing 
$$
\lrangle{- , -}: \calW_{\bpsi}(\pi) \otimes \calW_{\bpsi^{-1}}(\pi^{\vee}) \rightarrow \CC
$$
such that $\lrangle{\phi, \phi^{\vee}} = \lrangle{\sfW_{\pi}(\phi), \sfW_{\pi^{\vee}}(\phi^{\vee})}$, where the pairing on the left is the canonical pairing between $\pi$ and $\pi^{\vee}$.

Given $W \in\calW_{\bpsi}(\pi)$, we define $\tilW \in \calW_{\bpsi^{-1}}(\pi^\vee)$ by $\widetilde{W}(g)=W(w_{n} g^{-\rmt})$, where $w_n := \begin{bmatrix}
& & & 1 \\ & & 1 & \\ & \iddots &  & \\ 1 & & &
\end{bmatrix}$ is the longest Weyl element.

\subsubsection{The JPSS integrals and local JPSS $L$-factors}
Let $\pi$ be an irreducible admissible generic representation of $\GL_{m+1}(F)$. Let $n$ be a positive integer such that $n \leq m$, and set $l=m-n$. 
Let $\sigma$ be an irreducible admissible generic representation of $\GL_n(F)$ with central character $\omega_\sigma$. We associate to the Whittaker functions $W \in \calW_{\bpsi}(\pi)$ and $W' \in \calW_{\bpsi^{-1}}(\sigma)$ the local zeta integrals \footnote{There appears to be a typographical error in \cite[page 37]{hsieh2023fivevariablepadiclfunctionsu3times}, where it is written that $W^{\prime} \in\calW_{\bpsi}(\sigma)$.}
\begin{align*}
Z(s,W,W')&=\int_{U_n(F) \bs \GL_n(F)} W \left(\begin{bmatrix} h & \\ & \bfone_{l+1} \end{bmatrix}\right)W'(h) \abs{\det h}^{s-\frac{l+1}{2}} \dif h, \\
\widetilde{Z}(s,\widetilde{W},\widetilde{W'})&=\int_{U_n(F) \bs \GL_n(F)}\int_{\rmM_{l \times n}(F)} \widetilde{W}\left(\begin{bmatrix} h & & \\ x & \bfone_l & \\ & & 1 \end{bmatrix}\right)\widetilde{W'}(h) \abs{\det h}^{s-\frac{l+1}{2}} \dif x \dif h, 
\end{align*}
which converge absolutely for $\Re s \gg 0$, where $\dif h$ is the Haar measure on $\GL_n(F)$ normalized such that $\vol(\GL_n(\fro)) = 1$. 

We denote by $L^\JPSS(s,\pi\times\sigma)$, $\varepsilon^\JPSS(s,\pi\times\sigma,\bpsi)$, and $\gamma^\JPSS(s,\pi\times\sigma,\bpsi)$ the $L$, epsilon, and gamma factors associated with $\pi$ and $\sigma$. These local factors were studied extensively in \cite{jacquet1983rankin}. The gamma factor is defined as the proportionality constant of the functional equation 
\begin{equation}
    Z(1-s,\pi^\vee(w_{m+1,n})\widetilde{W},\widetilde{W'})=\omega_\sigma(-1)^m\gamma^\JPSS(s,\pi\times\sigma,\bpsi)\widetilde{Z}(s,W,W'),
\end{equation}
where 
$$
w_{m+1,n}=\begin{bmatrix} 
 \bfone_n & \\ & w_{m-n+1} 
 \end{bmatrix}.
$$ 

\begin{remark}
We are primarily interested in the case $n=m$ (so $l = 0$). By a change of variables, we have
\begin{equation} \label{eq:invarianceJPSS}
Z(s,W, \sigma(h_0)W^{\prime}) = \abs{\det h_0}^{\frac{1}{2}-s} Z(s, \pi(\jmath^{\GL}(h_0^{-1}))W, W^{\prime})
\end{equation}
for any $h_0 \in \GL_{n}(F)$. This property is referred to as the \emph{invariance of JPSS integrals}.
\end{remark}

When $\pi$ and $\sigma$ are viewed as representations of unitary groups over the split quadratic algebra $K = F \oplus F$, 
\begin{equation} \label{eq:uniJPSS}
L(s,\pi\times\sigma)=L^\JPSS(s,\pi\times\sigma)L^\JPSS(s,\pi^\vee\times\sigma^\vee).
\end{equation}
When $n=1$ and $\chi$ is a character of $F^\times$, we have 
\begin{subequations}
\begin{equation}
L^{\GJ}(s,\pi\otimes\chi)=L^\JPSS(s,\pi\times\chi),      
\end{equation}
\begin{equation}
\varepsilon^{\GJ}(s,\pi\otimes\chi,\bpsi)=\varepsilon^\JPSS(s,\pi\times\chi,\bpsi),
\end{equation}
\begin{equation}
\gamma^{\GJ}(s,\pi\otimes\chi,\bpsi)=\gamma^\JPSS(s,\pi\times\chi,\bpsi),
\end{equation}
\end{subequations}
where the local factors on the left-hand side are the Godement--Jacquet $L$-factors studied in \cite{godement1972zeta}. Furthermore, recall that
\begin{equation} \label{eq:uniad}
L(s, \pi, \Ad) = L^{\JPSS}(s, \pi \times \pi^{\vee}), \quad L(s, \sigma, \Ad) = L^{\JPSS}(s, \sigma \times \sigma^{\vee}),      
\end{equation}
where $\pi$ and $\sigma$ are regarded as representations of unitary groups on the left-hand side and as representations of general linear groups on the right-hand side. Combining \eqref{eq:uniJPSS} and \eqref{eq:uniad}, we obtain
\begin{equation} \label{eq:GLcalL}
    \scL(\pi \times \sigma) = \dfrac{L^\JPSS\left(\dfrac{1}{2},\pi \times \sigma \right) L^\JPSS\left(\dfrac{1}{2},\pi^{\vee} \times \sigma^{\vee} \right)}{L^{\JPSS}\left(1, \pi \times \pi^{\vee} \right) L^{\JPSS}\left(1, \sigma \times \sigma^{\vee} \right)}.
\end{equation}

\subsubsection{The naive local Rankin--Selberg $L$-factor}
We maintain the notation of the previous section and consider the special case $n = m$ (so $l = 0$). By definition, the local Godement--Jacquet $L$-function $L^{\GJ}(s,\sigma)$ is of the form $P_{\sigma}(q^{-s})^{-1}$, where $P_{\sigma} \in \CC[X]$ is a polynomial of degree at most $n$ satisfying $P_{\sigma}(0) = 1$. We may then find $n$ complex numbers $\{\alpha_i\}_{i=1}^{n}$ such that
$$
L^{\GJ}(s, \sigma) = \prod_{i=1}^{n}(1-\alpha_i q^{-s})^{-1}.
$$
We call the set $\{\alpha_i\}_{i=1}^{n}$ the \emph{Langlands parameters} of $\sigma$. Let $\{\gamma_j\}_{j=1}^{n+1}$ be the Langlands parameters of $\pi$. We then define
$$
L^{\RS}(s, \pi \times \sigma) := \prod_{i=1}^{n} \prod_{j=1}^{n+1}(1-\alpha_i \gamma_j q^{-s})^{-1}
$$
as the \emph{naive local Rankin--Selberg $L$-factor}. Morally speaking, $L^{\RS}(s, \pi \times \sigma)$ captures the ``unramified part'' of $L^{\JPSS}(s, \pi \times \sigma)$.

We require Bernstein--Zelevinsky's classification theorem for irreducible admissible representations of $\GL_n(F)$. For our purposes, we restrict ourselves to the classification of tempered representations \footnote{See, for example, \cite[Theorems 14.6.4--14.6.5]{goldfeld2011automorphicii} for the case $\GL_n(\QQ_p)$}. This classification consists of two parts:
\begin{itemize}
    \item \textit{Discrete series}. Let $r, d$ be two positive integers such that $rd=n$, and let $\eta$ be an irreducible supercuspidal representation of $\GL_{r}(F)$ with a unitary central character. Then
    $$
    \Ind_{P[r,\ldots, r]}^{\GL_n(F)}(\eta, \balpha_F^{(1-d)/2} \eta, \ldots, \balpha_F^{(d-1)/2} \eta)
    $$
    has a unique irreducible quotient, denoted by $\BZ(\eta, d)$. Every discrete series representation of $\GL_n(F)$ is isomorphic to some $\BZ(\eta,d)$.
    \item \textit{Tempered representations}. Let $n = r_1 + \cdots + r_s$ be a partition of $n$ and let $\tau_i$ be a discrete series representation of $\GL_{r_i}(F)$ for $i = 1, \ldots, s$. Then the induction
    $$
    \tau_1 \boxplus \cdots \boxplus \tau_s := \Ind_{P[r_1, \ldots, r_s]}^{\GL_n(F)}(\tau_1, \ldots, \tau_s)
    $$
    is irreducible and tempered. Every tempered representation is isomorphic to a representation of this form. The operation ``$\boxplus$'' is called the \emph{isobaric sum}.
\end{itemize}
Let $\sigma = \tau_1 \boxplus \cdots \boxplus \tau_s$ be a tempered representation of $\GL_{n}(F)$. Let $(\tau_{j_1}, \ldots, \tau_{j_r})$ be the tuple of unramified isobaric summands of $\sigma$ with indices $j_k$ in increasing order. We define the \emph{unramified socle} of $\sigma$ as the isobaric sum $\sigma_{\ur} := \boxplus_{k=1}^{r} \tau_{j_k}$, which is an unramified tempered representation of $\GL_{r_{j_1} + \cdots + r_{j_r}}(F)$. Clearly, if $\sigma$ itself is unramified, then $\sigma_{\ur} = \sigma$.

\begin{proposition} \label{prop:RSJPSS}
Under the conventions and notation established above:
\begin{enumerate}[label = \rm (\arabic*)]
    \item The local $L$-factors $L^{\RS}(s, - \times -)$ and $L^{\JPSS}(s, - \times -)$ are bi-additive with respect to the isobaric sum; that is,
    \begin{align*}
    L^{\bullet}(s, \pi \times (\sigma \boxplus \sigma^{\prime})) &=  L^{\bullet}(s, \pi \times \sigma) L^{\bullet}(s, \pi \times \sigma^{\prime}), \\
    L^{\bullet}(s, (\pi \boxplus \pi^{\prime}) \times \sigma) &=  L^{\bullet}(s, \pi \times \sigma) L^{\bullet}(s, \pi^{\prime} \times \sigma),
    \end{align*}
for $\bullet \in \{\RS, \JPSS\}$ and all irreducible admissible representations $\pi, \pi^{\prime}, \sigma, \sigma^{\prime}$ of general linear groups.
    \item If $\sigma$ and $\pi$ are unramified irreducible admissible representations of $\GL_n(F)$ and $\GL_{n+1}(F)$ respectively, then $L^{\RS}(s, \pi \times \sigma) = L^{\JPSS}(s, \pi \times \sigma)$.
    \item If either $\sigma$ or $\pi$ is a ramified discrete series representation of $\GL_n(F)$ or $\GL_{n+1}(F)$ respectively, then $L^{\RS}(s, \pi \times \sigma) = 1$.
\end{enumerate}
In particular, we have $L^{\RS}(s, \pi \times \sigma) = L^{\JPSS}(s, \pi \times \sigma_{\ur}) = L^{\JPSS}(s, \pi_{\ur} \times \sigma)$.
\end{proposition}

\begin{proof}
The bi-additivity of $L^{\RS}$ follows directly from its definition, while that of $L^{\JPSS}$ is proven in \cite[Section 9.5, Theorem]{jacquet1983rankin}. The result in (2) is stated in \cite[Equation (14), p. 371]{jacquet1983rankin}. For (3), it follows from the definition of $L^{\RS}$ that it suffices to show that if $\tau$ is a ramified discrete series representation of $\GL_{n}(F)$, then $L^{\GJ}(s, \tau) = 1$. (This is equivalent to the statement that all Langlands parameters of $\tau$ are zero.) To prove this, let $\tau = \BZ(\eta, d)$ by the Bernstein--Zelevinsky classification. Then we have
\begin{equation*}
    L^{\GJ}(s, \tau) = \begin{cases}
    L(s, \balpha_{F}^{1-n} \eta) , &\quad d = n, \\
    1, &\quad d \neq n,
    \end{cases}
\end{equation*}
where the $L$-factor in the first case is the local Hecke $L$-factor. If $d \neq n$, the claim is immediate. If $d = n$ and $\tau$ is ramified, then $\eta$ must be ramified, which implies $L(s, \balpha_{F}^{1-n} \eta) = 1$. The final equalities in the proposition follow from (1)--(3) and the definition of the unramified socle.
\end{proof}

\subsubsection{Essential Whittaker vectors}
We now review the theory of the essential Whittaker vector associated with an irreducible admissible generic representation $\pi$ of $\GL_{n}(F)$. Given an open compact subgroup $\Gamma \subset \GL_{n}(F)$ and a character $\calX: \Gamma \rightarrow \CC^{\times}$, we define
$$
\calW_{\bpsi}(\pi, \Gamma, \calX) = \{W \in \calW_{\psi}(\pi): \pi(\gamma)W = \calX(\gamma)W \text{ for all } \gamma \in \Gamma\}.
$$

For any integer $f \geq 0$, we consider the following two compact open subgroups of $\GL_n(F)$:
\begin{align*}
    K_1(\frp^{f}) &:= \left\{ g \in \GL_{n}(\fro): g \equiv \begin{bmatrix}
 \ast & \ast \\ 0_{1 \times (n-1)} & 1
\end{bmatrix} \pmod{\frp^{f}} \right\}, \\
    K_0(\frp^{f}) &:= \left\{ g \in \GL_{n}(\fro): g \equiv \begin{bmatrix}
 \ast & \ast \\ 0_{1 \times (n-1)} & \ast
\end{bmatrix} \pmod{\frp^{f}} \right\},
\end{align*}
where $K_1(\frp^{f})$ is a normal subgroup of $K_0(\frp^{f})$ with quotient $K_0(\frp^{f})/K_1(\frp^{f}) \cong (\fro/\frp^{f})^{\times}$.

Let $c(\pi)$ denote the exponent of the conductor of $\pi$; i.e., the epsilon factor of $\pi$ satisfies
$$
\varepsilon\left(s+\dfrac{1}{2}, \pi, \bpsi\right) = q^{-c(\pi)s} \varepsilon\left(\dfrac{1}{2}, \pi, \bpsi \right).
$$
Following \cite{jacquet1981conducteur, jacquet2012correction}, we have
$$
\dim_{\CC} \calW_{\bpsi}(\pi, K_1(\frp^{c(\pi)}), \bfone) = 1.
$$
This subspace of $\calW_{\bpsi}(\pi)$ is called the \emph{essential line} of $\calW_{\bpsi}(\pi)$. We define the \emph{normalized essential Whittaker vector} $W_{\pi, \bpsi}^{\ess}$ of $\pi$ to be the unique vector in the essential line satisfying
$$
W_{\pi, \bpsi}^{\ess}\left( \begin{bmatrix}
 gh & \\ & 1
\end{bmatrix} \right) = W_{\pi, \bpsi}^{\ess}\left( \begin{bmatrix}
 g & \\ & 1
\end{bmatrix} \right) \, \text{ for all } h \in \GL_{n-1}(\fro).
$$

\begin{remark}
\begin{enumerate}[label = \rm (\arabic*)]
    \item If $\pi$ is unramified, then $c(\pi)=0$ and $W_{\pi, \bpsi}^{\ess}$ is the normalized spherical function.
    \item The group $K_0(\frp^{c(\pi)})$ acts on the essential line via a character related to the central character $\omega_{\pi}$. Specifically, for $g = (g_{i,j}) \in K_0(\frp^{c(\pi)})$, let
    \[
    \omega_{\pi}^{\downarrow}(g) = \begin{cases}
    1, &\quad \text{ if } c(\pi) = 0, \\
    \omega_{\pi}(g_{n,n}), &\quad \text{ if } c(\pi) > 0.
    \end{cases}
    \]
    Then $\omega_{\pi}^{\downarrow}$ is a character of $K_0(\frp^{c(\pi)})$ trivial on $K_1(\frp^{c(\pi)})$, and
    \[
    \pi(g) W_{\pi}^{\ess} = \omega_{\pi}^{\downarrow}(g) W_{\pi}^{\ess} \, \text{ for all } g \in K_0(\frp^{c(\pi)}).
    \]
    \item The conductor of a multiplicative character $\chi$ of $F^{\times}$, denoted by $c(\chi)$, is defined as follows: if $\chi$ is unramified, $c(\chi)=0$; otherwise, $c(\chi)$ is the least integer $n \geq 1$ such that $\chi|_{1+\frp^n} = 1$.
\end{enumerate}
\end{remark}

\subsubsection{Test vector problem}
Let $\sigma$ and $\pi$ be irreducible admissible tempered representations of $\GL_{n}(F)$ and $\GL_{n+1}(F)$, respectively. Previously, we defined the JPSS integrals $Z(s, W, W^{\prime})$ for Whittaker functions $W \in \calW_{\bpsi}(\pi)$ and $W' \in \calW_{\bpsi^{-1}}(\sigma)$. The test vector problem is to find Whittaker functions $W^{\circ} \in \calW_{\bpsi}(\pi)$ and $W^{\prime, \circ} \in \calW_{\bpsi^{-1}}(\sigma)$ such that
$$
Z(s, W^{\circ}, W^{\prime, \circ}) = L^{\JPSS}(s, \pi \times \sigma).
$$
In this section, we introduce two partial results on this problem, which are sufficient to compute the local Ichino--Ikeda integral under Assumption \ref{ass:ram}.

\begin{theorem}[{\cite[Théorème, p. 208]{jacquet1981conducteur}}] \label{thm:testJPSSur}
Suppose $\sigma$ is unramified. Then 
$$
Z\left(\dfrac{1}{2}, W_{\pi,\bpsi}^{\ess}, W_{\sigma, \bpsi^{-1}}^{\ess}\right) = L^\JPSS\left(\dfrac{1}{2},\pi \times \sigma \right).
$$
\end{theorem}

The case where $\sigma$ is ramified is more complicated. For this case, Booker, Krishnamurthy, and Lee \cite{booker2020test} modified $W_{\pi, \psi}^{\ess}$ through a process of unipotent averaging. This method dates back to the work of Schmidt \cite{schmidt1993relative}, among others.

Let $\frn$, $\frq$, and $\frc$ denote the conductors of $\pi$, $\sigma$, and the central character $\omega_{\sigma}$, respectively. Consider $\beta = (\beta_1, \ldots, \beta_{n}) \in F^{n}$ with $\beta_i \in \frq^{-1}$ for $i=1, \ldots, n$. Let $\bfu(\beta) = (u_{i,j})$ denote the matrix
$$
\bfu(\beta) = \begin{bmatrix}
1 &  & & & \beta_1 \\
  & 1 & & &  \beta_2\\
  & & \ddots &&  \vdots  \\
  & & & 1 & \beta_{n} \\
  &  &  &  & 1
\end{bmatrix} \in U_{n+1}(F).
$$
We define the unipotent averaging operator $\Theta$ on $W \in \calW_{\bpsi}(\pi)$ by
$$
\Theta W := \dfrac{1}{[\fro:\frq]^{n-1}} \sum_{(\beta_1, \ldots, \beta_{n-1}) \in (\frq^{-1}/\fro)^{n-1}} \pi(\bfu(\beta_1, \ldots, \beta_{n-1}, \varpi^{-v(\frc)})) W \in \calW_{\bpsi}(\pi).
$$
(When $n=1$, we understand there to be a single summand, so that $\Theta W = \pi(\varpi^{-v(\frc)})W$).

Then we have the following theorem.

\begin{theorem} \label{thm:testJPSSram}
With notation and conventions as above, under Assumption \ref{ass:ram}, if $\sigma_v$ is ramified, then
$$
Z\left(\dfrac{1}{2}, \Theta W_{\pi,\bpsi}^{\ess}, W_{\sigma, \bpsi^{-1}}^{\ess}\right) = C_{\sigma, \bpsi} \cdot L^\JPSS\left(\dfrac{1}{2},\pi\times\sigma \right),
$$
where $C_{\sigma, \bpsi}$ is the explicit non-zero constant
$$
C_{\sigma, \bpsi} := \dfrac{\frg(\omega_{\sigma}, \bpsi, \varpi^{-v(\frc)})}{[\GL_{n}(\fro):K_0(\frq)]} \neq 0.
$$
\end{theorem}

\begin{proof}
The main theorem of \cite{booker2020test} states that
$$
Z\left(s, \Theta W_{\pi,\bpsi}^{\ess}, W_{\sigma, \bpsi^{-1}}^{\ess}\right) = C_{\sigma, \bpsi} \cdot L^\RS\left(s,\pi\times\sigma\right),
$$
with the constant $C_{\sigma, \bpsi}$ specified as on \cite[p. 46]{booker2020test}. Under Assumption \ref{ass:ram}, if $\sigma_v$ is ramified, then $\pi_v$ must be unramified. Therefore, by Proposition \ref{prop:RSJPSS}, we have $L^{\RS}(s, \pi \times \sigma) = L^{\JPSS}(s, \pi_{\ur} \times \sigma) = L^{\JPSS}(s, \pi \times \sigma)$. Evaluation at $s=1/2$ yields the result.
\end{proof}

\subsection{The splitting lemma}
We now return to the computation of the local Ichino--Ikeda integrals.

Let $\pi$ be an irreducible admissible tempered representation of $\GL_{m+n+1}(F)$ and let $\sigma$ be that of $\GL_{m+n}(F)$. We consider the integral
$$
J(W_1, W_2, W_1^{\prime}, W_2^{\prime}) := \int_{\GL_{m+n}(F)} \lrangle{\pi(\jmath^{\GL}(h)) W_1, W_2}_{\pi} \lrangle{\sigma(h)W_1^{\prime}, W_2^{\prime}}_{\sigma} \dif h,
$$
where 
$$
W_1 \in \calW_{\bpsi}(\pi), \, W_2 \in \calW_{\bpsi^{-1}}(\pi^{\vee}), \, W_1^{\prime} \in \calW_{\bpsi^{-1}}(\sigma), \, W_2^{\prime} \in \calW_{\bpsi}(\sigma^{\vee}).
$$
This integral is known to converge absolutely. 

An essential tool in this context is the splitting lemma of Wei Zhang \footnote{Note that in \cite{zhang2014automorphic}, Zhang uses unnormalized local Haar measures, while here we use normalized ones, following the reformulation in \cite[Lemma 5.2]{hsieh2023fivevariablepadiclfunctionsu3times}. Interested readers may refer to \cite[Remark 5.3]{hsieh2023fivevariablepadiclfunctionsu3times} for familiar special cases of this lemma.}, stated as follows:

\begin{theorem}[Splitting lemma, {\cite[Proposition 4.10]{zhang2014automorphic}}]
With notation as above, we have
\[
J(W_1, W_2, W_1^{\prime}, W_2^{\prime}) =  Z\left(\dfrac{1}{2}, W_1, W_1^{\prime}\right) Z\left(\dfrac{1}{2}, W_2, W_2^{\prime}\right) \prod_{i=1}^{m+n-1} \zeta_{F}(i). 
\]
\end{theorem}

Consequently, the computation of local Ichino--Ikeda integrals reduces to the evaluation of local JPSS integrals. Indeed, under the identifications of unitary groups with general linear groups via $\varrho_{w, m, n}$ and $\varrho_{w, m, n+1}$, and the use of local Whittaker models, we obtain the following corollary:

\begin{corollary}\label{cor:splitii}
With notation as above, we have
\[ 
\scI (\sigma, \pi) = \dfrac{J(W_1, W_2, W_1^{\prime}, W_2^{\prime})}{\lrangle{W_1, W_2} \lrangle{W_1^{\prime}, W_2^{\prime}}} = \dfrac{Z\left(1/2, W_1, W_1^{\prime}\right) Z\left(1/2, W_2, W_2^{\prime}\right)}{\lrangle{W_1, W_2} \lrangle{W_1^{\prime}, W_2^{\prime}}} \prod_{i=1}^{m+n-1} \zeta_{F}(i).     
\]
\end{corollary}

\subsection{Local Ichino--Ikeda integrals at ramified places outside $p$}
Let
$$
\calB_{\pi}(W_1, W_2) := \dfrac{\lrangle{W_1, W_2}}{L^{\JPSS}(1, \pi \times \pi^{\vee})}, \quad \calB_{\sigma}(W_1^{\prime}, W_2^{\prime}) := \dfrac{\lrangle{W_1^{\prime}, W_2^{\prime}}}{L^{\JPSS}(1, \sigma \times \sigma^{\vee})}
$$
and define \footnote{We note that this notation differs slightly from the definition of $\calB_{\pi_l}$ and $\calB_{\sigma_l}$ in \cite[Sect. 4.7]{hsieh2023fivevariablepadiclfunctionsu3times}.}
$$
\calB_{\pi}^{\ess} := \calB_{\pi}(W_{\pi,\bpsi}^{\ess}, W_{\pi^{\vee},\bpsi^{-1}}^{\ess}), \quad \calB_{\sigma}^{\ess} := \calB_{\sigma}(W_{\sigma,\bpsi^{-1}}^{\ess}, W_{\sigma^{\vee},\bpsi}^{\ess}).
$$
These can be regarded as invariants of the local representations $\pi$ and $\sigma$.

\subsubsection{The case where $\sigma$ is unramified}
We first consider the case where $\sigma$ is unramified. By taking
$$
W_1 = W_{\pi,\bpsi}^{\ess}, \quad W_1^{\prime} = W_{\sigma, \bpsi^{-1}}^{\ess}, \quad W_2 = W_{\pi^{\vee},\bpsi^{-1}}^{\ess}, \quad W_2^{\prime} = W_{\sigma^{\vee}, \bpsi}^{\ess}
$$
in Corollary \ref{cor:splitii} and applying Theorem \ref{thm:testJPSSur} along with \eqref{eq:GLcalL}, we obtain the following:

\begin{proposition}
With notation as above, let $v \in \scV_{\calF}^{\ram}$ such that $\sigma_v$ is unramified. Then
\[ 
\scI (\sigma, \pi) = \scL(\pi \times \sigma) \cdot \calB_{\pi}^{\ess} \calB_{\sigma}^{\ess} \cdot \prod_{i=1}^{m+n-1} \zeta_{F}(i).      
\]
\end{proposition}

\subsubsection{The case where $\sigma$ is ramified}
By taking
$$
W_1 = \Theta W_{\pi,\bpsi}^{\ess}, \quad W_1^{\prime} = W_{\sigma, \bpsi^{-1}}^{\ess}, \quad W_2 = \Theta W_{\pi^{\vee},\bpsi^{-1}}^{\ess}, \quad W_2^{\prime} = W_{\sigma^{\vee}, \bpsi}^{\ess}
$$
in Corollary \ref{cor:splitii} and applying Theorem \ref{thm:testJPSSram} along with \eqref{eq:GLcalL}, we obtain the following:

\begin{proposition}
With notation as above, let $v \in \scV_{\calF}^{\ram}$ such that $\sigma_v$ is ramified. Under Assumption \ref{ass:ram}, we have
\[ 
\scI (\sigma, \pi) = \scL(\pi \times \sigma) \cdot C_{\sigma, \bpsi} C_{\sigma^{\vee}, \bpsi^{-1}} \cdot \calB_{\pi}^{\ess} \calB_{\sigma}^{\ess} \cdot \prod_{i=1}^{m+n-1} \zeta_{F}(i).      
\]
\end{proposition}

\subsection{Further representation theory of $\GL_n(F)$: the ordinary condition} \label{sec:iiintegralatp} 
To handle the local Ichino--Ikeda integrals at places $v$ of $\calF$ above $p$, we introduce the background on the ordinary condition for irreducible admissible representations of $\GL_n(F)$. As before, for this subsection only, we let $m, n$ be general positive integers (unrelated to the unitary groups discussed previously). 

Recall that at the beginning of this article, we fixed an identification $\iota_p: \CC \xrightarrow{\sim} \barQQ_p$. Through $\iota_p$, we regard characters of $F$ as being valued in $\barQQ_p$, and every irreducible admissible representation $\pi$ of $\GL_n(F)$ is considered to have coefficients in $\barQQ_p$.

\subsubsection{The ordinary line}
Let $\ulmu = (\mu_1, \ldots, \mu_n)$ be a tuple of characters $F^{\times} \rightarrow \barQQ_p^{\times}$, which naturally induces a character $\mu$ of $T_n(F)$. We define a character $\mu^{\natural}$ of $T_n(F)$ by
$$
\mu^{\natural}: \diag[x_1, \ldots, x_n] \mapsto \prod_{i=1}^{n} \abs{x_i}^{n-i} \mu_i(x_i),
$$
which we extend to a character of $B_n(F)$ by inflation. We then define the algebraically induced principal series representation \footnote{The representation $\abs{\det(-)}^{\frac{1-n}{2}} \otimes I_{B_n}^{\GL_n}(\mu)$ coincides with the normalized induction of $\mu$ from $B_n(F)$ to $\GL_n(F)$.}
$$
I_{B_n}^{\GL_n}(\mu) := \{f: \GL_n(F) \rightarrow \barQQ_p^{\times} \text{ locally constant} : f(bg) = \mu^{\natural}(b)f(g) \text{ for all } b \in B_n(F), g \in \GL_n(F) \}
$$
as an admissible representation of $\GL_n(F)$ acting via right translation.

\begin{definition} \label{defn:ord}
Let $\pi$ be an irreducible admissible representation of $\GL_n(F)$.
\begin{enumerate}
    \item We say $\pi$ is \emph{ordinary} if there exists a unique tuple $\ulmu$ of admissible characters such that $\abs{x}^{i-1} \mu_{n+1-i}(x)$ are $p$-adic units for $1 \leq i \leq n$ and all $x \in F^{\times}$, and $\pi$ is isomorphic to $I_{B_n}^{\GL_n}(\mu)$.
    \item We say $\pi$ is \emph{semi-stably ordinary} if, in addition, the characters $\mu_i$ are all unramified.
    \item We say $\pi$ is \emph{regularly ordinary} if $\pi$ is ordinary and the $p$-adic valuations of $\mu_i(\varpi)$ are all distinct. We say $\pi$ is \emph{regularly semi-stably ordinary} if it is both regularly ordinary and semi-stably ordinary.
\end{enumerate}
\end{definition}

\begin{remark} \label{rem:ordcontra}
Note that if $\pi$ satisfies the properties in Definition \ref{defn:ord}, then so does its contragredient $\pi^{\vee}$ with respect to the tuple $\check{\ulmu} := (\balpha_F^{1-n} \mu_n^{-1}, \ldots, \balpha_F^{1-n} \mu_1^{-1})$.
\end{remark}

\begin{definition}
For $x \in \fro \cap F^{\times}$, let $\bfd(x) := \diag[x^{n-1}, \ldots, x, 1] \in \GL_n(F)$. We define the operator $\VV_{n}^{x}$ on the space $\pi^{U_n(\fro)}$ by
\[
\VV_{n}^{x} := \sum_{u \in U_n(\fro) /( U_n(\fro) \,\cap\, \bfd(x) U_n(\fro) \bfd(x)^{-1})} \pi(u \bfd(x)).
\]
\end{definition}

We remark that for $x = \varpi$, the operator $\VV_{n}^{\varpi}$ corresponds to the $\UU_{\frp}$-operator (cf. \cite[Section 1.4]{januszewski2024nonabelian}).

We have the following proposition, proved in \cite[Lemma 4.4]{liu2023anticyclotomicpadiclfunctionsrankinselberg}. \footnote{We note that in the statement of \cite[Lemma 4.4]{liu2023anticyclotomicpadiclfunctionsrankinselberg}, the eigenvalue is given by 
$$
\left( \prod_{m=1}^{n-1} \prod_{i=1}^{m} \abs{x}^{i-1} \mu_i(x) \right),
$$
after reordering $(\mu_1, \ldots, \mu_n)$ into $(\mu_{n}, \mu_{n-1}, \ldots, \mu_1)$. Our ordinary condition in Definition \ref{defn:ord} is adjusted accordingly. While \cite{liu2023anticyclotomicpadiclfunctionsrankinselberg} requires $\pi$ to be a subrepresentation of $I_{B_n}^{\GL_n}(\mu)$, our Definition \ref{defn:ord} utilizes the quotient representation. Since we assume $\pi$ is irreducible, the ordering of the $\mu_i$ does not affect the outcome. We thank Yifeng Liu for his guidance on this matter.}

\begin{proposition} \label{prop:ordinaryline}
Suppose that $\pi$ is regularly ordinary. Then there exists a unique (up to scalar) non-zero element $f^{\ord}_{\pi} \in \pi^{U_n(\fro)}$ such that
$$
\VV_n^{x} \cdot f^{\ord}_{\pi} = \left(\prod_{m=1}^{n-1} \prod_{i=1}^{m} \abs{x}^{i-1} \mu_{n+1-i}(x) \right) f^{\ord}_{\pi}
$$
for every $x \in \fro \cap F^{\times}$. In particular, the $\VV_{n}^{x}$-eigenvalue is a $p$-adic unit.
\end{proposition}

The one-dimensional $\barQQ_p$-subspace of $\pi^{U_n(\fro)}$ generated by $f^{\ord}_{\pi}$ is called the \emph{ordinary line} of $\pi$, denoted by $\pi^{\ord}$, and a non-zero element thereof is an \emph{ordinary vector}. Let $\calW_{\bpsi}(\pi)$ be the Whittaker model of $\pi$ with respect to $\bpsi$. For any $f \in V_{\pi}$ supported on the big cell $B_n(F) w_n B_n(F)$, we consider the integral 
$$
W_f(g) := \int_{U_n(F)} f(w_n u g) \overline{\bpsi}(u) \dif u,
$$
where $\dif u$ is the normalized Haar measure on $U_n(F)$. By a result of Rodier, this integral converges and extends uniquely to an intertwining operator 
$$
\sfW: I_{B_n}^{\GL_n} (\mu) \rightarrow \Ind_{U_n(F)}^{\GL_n(F)} \bpsi.
$$
See \cite[Corollary 1.8]{casselman1980unramified} for details. By the uniqueness of Whittaker models for $I_{B_n}^{\GL_n} (\mu)$, $\sfW$ provides the Whittaker model of $\pi$, i.e.,
$$
\sfW_{\pi}: I_{B_n}^{\GL_n} (\mu) \rightarrow \calW_{\bpsi}(\pi) \hookrightarrow \Ind_{U_n(F)}^{\GL_n(F)} \bpsi.
$$
The image of the ordinary line $\pi^{\ord}$ in $\calW_{\bpsi}(\pi)$ is denoted by $\calW^{\ord}_{\bpsi}(\pi)$.

We can explicitly construct a canonical ordinary vector. We define the big-cell section\footnote{It is called the \emph{big cell section} since it is supported on the big cell of the Bruhat decomposition, i.e. the cell of the longest Weyl element. We use the notation “$\dagger$” to denote such sections.} $f^{\dagger}_{\pi}: \GL_n(F) \rightarrow \barQQ_p^{\times}$ by
$$
f^{\dagger}_{\pi}(g) := \begin{cases}
\mu^{\natural}(b), &\quad \text{ if } g=b w_n r \text{ with } b \in B_n(F), \, r \in U_n(\fro), \\
0, &\quad \text{ otherwise}.
\end{cases}
$$

Let $W_{\pi}^{\dagger}$ be the Whittaker vector corresponding to $f^{\dagger}_{\pi}$. Here are some basic properties of $f^{\dagger}_{\pi}$:

\begin{proposition} \label{prop:bigcellord}
With notation as above, we have:
\begin{enumerate}[label = \rm (\arabic*)]
    \item $f^{\dagger}_{\pi}$ is an ordinary vector.
    \item $W^{\dagger}_{\pi}(\bfone_n) = \vol(U_n(\fro)) \neq 0$.
\end{enumerate}
\end{proposition}

\begin{proof}
It is immediate that $f_{\pi}^{\dagger} \in \pi^{U_{n}(\fro)}$. We fix $x \in \fro \cap F^{\times}$. Let $g \in \GL_n(F)$ be such that $\VV_{n}^{x} f_{\pi}^{\dagger}(g) \neq 0$. By the definition of $\VV_{n}^{x}$, there exists $u \in U_{n}(\fro)$ such that $f_{\pi}^{\dagger}(gu \bfd(x)) \neq 0$, implying $gu \bfd(x) \in B_n(F) w_n U_n(\fro)$. One verifies that for any $u = (u_{i,j}) \in U_n(\fro)$:
\begin{subequations}
\begin{equation} \label{eq:conjdx}
\bfd(x) u \bfd(x)^{-1} = \begin{bmatrix}
1 & x u_{1,2} & x^2 u_{1,3}  & \cdots & x^{n-1} u_{1,n} \\ 
  & 1 & x u_{2,3} & \cdots  & x^{n-2} u_{2,n} \\
  &  & \ddots & \vdots & \vdots \\
  &   &         &   1      &  x u_{n-1, n} \\    
  &   &         &  &  1   
\end{bmatrix} \in U_n(\fro)
\end{equation}
and
\begin{equation} \label{eq:conjdxinv}
\bfd(x)^{-1} u \bfd(x) = \begin{bmatrix}
1 & x^{-1} u_{1,2} & x^{-2} u_{1,3}  & \cdots & x^{1-n} u_{1,n} \\ 
  & 1 & x^{-1} u_{2,3} & \cdots  & x^{2-n} u_{2,n} \\
  &  & \ddots & \vdots & \vdots \\
  &   &         &   1      &  x^{-1} u_{n-1, n} \\    
  &   &         &  &  1   
\end{bmatrix} \in U_n \left(\frp^{-(n-1) v(x)} \fro \right).
\end{equation}
\end{subequations}
It follows from \eqref{eq:conjdx} that $g \in B_n(F) w_n U_n(\fro)$. Thus, it suffices to compute $\VV_{n}^{x} f^{\dagger}_{\pi}(g)$ for $g \in B_n(F) w_n U_n(\fro)$. Moreover, since $\VV_{n}^{x} f^{\dagger}_{\pi}$ is $U_n(\fro)$-invariant, we may assume $g \in B_n(F) w_n$.

By \eqref{eq:conjdx}, the quotient $U_n(\fro) / (\bfd(x) U_n(\fro) \bfd(x)^{-1})$ has a complete set of representatives given by:
\begin{equation} \label{eq:repofhecke}
U_n(\fro) / (\bfd(x) U_n(\fro) \bfd(x)^{-1}) \xleftrightarrow{1:1} \begin{bmatrix}
1 & \fro/\frp^{v(x)} & \fro/\frp^{2v(x)}  & \cdots & \fro/\frp^{(n-1)v(x)} \\ 
  & 1 & \fro/\frp^{v(x)} & \cdots  & \fro/\frp^{(n-2)v(x)} \\
  &  & \ddots & \vdots & \vdots \\
  &   &         &   1      &  \fro/\frp^{v(x)} \\    
  &   &         &  &  1   
\end{bmatrix}.
\end{equation}
Here, for any positive integer $t$, we denote
\[
\fro/\frp^{t} := \left\{ \sum_{k=0}^{t-1} s_k \varpi^{k}: s_k \in \fro/\frp \right\}.
\]
It remains to compute $\pi(u \bfd(x)) f^{\dagger}_{\pi}(g)$ for $g = b w_n \in B_n(F) w_n$ and $u$ in the set of representatives \eqref{eq:repofhecke}. We have
\begin{align*}
\pi(u \bfd(x)) f^{\dagger}_{\pi} (g) &= f^{\dagger}_{\pi}(b w_n u \bfd(x)) \\
&= f^{\dagger}_{\pi}(b (w_n \bfd(x) w_n^{-1}) w_n (\bfd(x)^{-1} u \bfd(x))).
\end{align*}
By \eqref{eq:conjdxinv} and the definition of $f_{\pi}^{\dagger}$, we see that $\bfd(x)^{-1} u \bfd(x) \in U_n(\fro)$ if and only if $u = 1$. Furthermore, 
\[
w_n \bfd(x) w_n^{-1} = \diag[1, x, x^2, \ldots, x^{n-1}].
\]
Consequently, we obtain \footnote{The third equality is obtained by representing the factors in a matrix where rows are indexed by $i = 1, \ldots, n$ and columns by $m = 0, 1, \ldots, n-1$, and taking the product of the non-zero entries:
\[
\begin{pmatrix}
1 & \mu_n(x) & \mu_n(x) & \cdots & \cdots & \mu_n(x) \\
  & 1 & \abs{x} \mu_{n-1}(x) & \cdots & \cdots  & \abs{x} \mu_{n-1}(x) \\
  & & \ddots & \vdots & \vdots  & \vdots \\
  & & & 1 & \abs{x}^{n-3} \mu_3(x)&  \abs{x}^{n-3} \mu_3(x) \\
  & & & & 1 & \abs{x}^{n-2} \mu_2(x) \\
  & & & & & 1
\end{pmatrix}.
\]
}
\begin{align*}
\VV_{n}^{x} f^{\dagger}_{\pi} (g) &= \mu^{\natural}(\diag[1, x, x^2, \ldots, x^{n-1}]) f^{\dagger}_{\pi}(g) \\ 
&= \prod_{i=1}^{n} \left(\abs{x}^{n-i} \mu_i(x) \right)^{i-1} f^{\dagger}_{\pi}(g) = \left(\prod_{m=1}^{n-1} \prod_{i=1}^{m} \abs{x}^{i-1} \mu_{n+1-i}(x) \right) f^{\dagger}_{\pi}(g).
\end{align*}
This verifies (1). For (2), 
\[
W_{\pi}^{\dagger}(\bfone_n) = \int_{U_n(F)} f_{\pi}^{\dagger}(w_n u) \overline{\bpsi}(u) \dif u = \int_{U_n(\fro)} \dif u = \vol(U_n(\fro)) \neq 0,
\]
since the integrand vanishes for $u \notin U_n(\fro)$ and $f_{\pi}^{\dagger}(w_n u) = 1$ for $u \in U_n(\fro)$.
\end{proof}

\subsubsection{Iwahori type}
For any integer $\alpha \geq 0$, we define the \emph{Iwahori subgroup} of level $\alpha$, denoted by $\Iw_{\alpha}$, as the subgroup of matrices in $\GL_n(\fro)$ that are upper triangular modulo $\varpi^{\alpha}$.

The subgroup $\Iw_{\alpha}$ acts naturally on $\calW_{\bpsi}(\pi)^{U_n(\fro)}$. Since $\calW^{\ord}_{\bpsi}(\pi)$ is one-dimensional, there exists a character $\vartheta: \Iw_{\alpha} \rightarrow \CC^{\times}$ such that
$$
W(gr) = \vartheta(r) W(g), \quad \forall g \in \GL_n(F), r \in \Iw_{\alpha}.
$$
This is called the \emph{Iwahori type} of $\pi$. Using the decomposition $\Iw_{\alpha} = [\Iw_{\alpha}, \Iw_{\alpha}] \cdot T_n(\fro)$, where $[\Iw_{\alpha}, \Iw_{\alpha}]$ is the derived subgroup (consisting of matrices that are strictly upper triangular modulo $\varpi^{\alpha}$), the character $\vartheta$ is determined by its values on $T_n(\fro)$. Thus, there exist characters $\theta_1, \ldots, \theta_n: F^{\times} \rightarrow \CC^{\times}$ such that
$$
\vartheta((r_{ij})_{i,j}) = \prod_{i=1}^{n} \theta_i(r_{ii}).
$$

The following proposition characterizes the Iwahori type of $\pi$. \footnote{We thank Loren Spice for helpful guidance regarding this section.}

\begin{proposition} \label{prop:iwahoriact}
Let $\pi = I_{B_n}^{\GL_n}(\mu)$ be a regularly ordinary representation of $\GL_n(F)$. Then the Iwahori type of $\pi$ is $(\breve{\mu})^{\natural}$, where for $\ulmu = (\mu_1, \ldots, \mu_n)$, we define $\breve{\ulmu} := (\mu_n, \mu_{n-1}, \ldots, \mu_1)$.
\end{proposition}
\begin{proof}
Since $W_{\pi}^{\dagger}(\bfone_n) = \vol(U_n(\fro))$, it follows that for all $r \in \Iw_\alpha \cap T_n(F) = T_n(\fro)$, $W_{\pi}^{\dagger}(r)$ is given by
\begin{align*}
\int_{U_n(F)} f_\pi^{\dagger}(w_n u r) \overline{\bpsi}(u) \dif u &= \int_{U_n(F)} f_\pi^{\dagger}((w_n r w_n^{-1}) w_n (r^{-1} u r)) \overline{\bpsi}(u) \dif u \\
&= \mu^{\natural}(w_n r w_n^{-1}) \int_{U_n(F)} f_\pi^{\dagger} (w_n u^{\prime}) \overline{\bpsi}(u^{\prime}) \dif u^{\prime} \\
&= (\breve{\mu})^{\natural}(r) \vol(U_n(\fro)) = (\breve{\mu})^{\natural}(r) W^{\dagger}_{\pi}(\bfone_n).
\end{align*}
The second equality holds because $r$ normalizes $U_n(\fro)$, and the third follows from the identity $w_n \diag[r_1, \ldots, r_n] w_n^{-1} = \diag[r_n, r_{n-1}, \ldots, r_1]$. Thus, $r$ acts by multiplication by $(\breve{\mu})^{\natural}(r)$.
\end{proof}

\begin{remark}
In particular, if $\pi$ is regularly semi-stably ordinary (i.e., all $\mu_i$ are unramified), then $\mu(\fro^{\times}) = 1$. Consequently, $\Iw_{\alpha}$ acts trivially on the ordinary line, recovering \cite[Lemma 4.8]{liu2023anticyclotomicpadiclfunctionsrankinselberg}.
\end{remark}

\subsection{Test vector problems: Januszewski's generalization of the local Birch lemma}
We first introduce Fabian Januszewski's generalization of the local Birch Lemma (cf. \cite{januszewski2024nonabelian}). Only in this subsection, let $\sigma = I_{B_{n}}^{\GL_{n}}(\ullam)$ and $\pi = I_{B_{n+1}}^{\GL_{n+1}}(\ulmu)$ be regular ordinary representations of $\GL_{n}(F)$ and $\GL_{n+1}(F)$, respectively. To apply Januszewski's results, we impose the following assumption on the local representations $\sigma$ and $\pi$.

\begin{assumption}[Constant conductor assumption] \label{ass:conductor}
    The principal series representations $\sigma = I_{B_{n}}^{\GL_{n}}(\ullam)$ and $\pi = I_{B_{n+1}}^{\GL_{n+1}}(\ulmu)$ satisfy the following: for all $1 \leq i \leq n+1$ and all $1 \leq j \leq n$, the conductors of $\mu_i \lambda_j$ are nontrivial and equal.
\end{assumption}

This can be viewed as a ``ramification interlacing condition'' for the local representations $\sigma$ and $\pi$. In fact, fixing the local representation $\sigma$, if we choose local characters $\mu_{1}, \ldots, \mu_{n+1}: F^{\times} \rightarrow \barQQ_{p}^{\times}$ with sufficiently large conductors, then the pair of local representations $\pi := I_{B_{n+1}}^{\GL_{n+1}}(\ulmu)$ and $\sigma$ satisfies the constant conductor assumption. The converse also holds.

Next is the generalized local Birch lemma by F. Januszewski \cite[Theorem 2.8]{januszewski2024nonabelian}. Here we state a simplified version. For $x \in F^{\times}$, let
$$
h_n := \begin{bmatrix}
 w_n & 1 \\  & 1
\end{bmatrix}, \quad \bft(x) := \diag[x^{n}, x^{n-1}, \ldots, x] \in \GL_n(F).
$$
We denote $h_n^{(x)} := (\jmath^{\GL}(\bft(x)))^{-1} h_n \jmath^{\GL}(\bft(x))$.

\begin{theorem}[Local Birch Lemma] \label{thm:janbirch}
With notation as above, and requiring the local representations $\sigma$ and $\pi$ to satisfy Assumption \ref{ass:conductor}. Let $W_{\pi} \in \calW_{\bpsi}(\pi)$ and $W_{\sigma} \in \calW_{\bpsi^{-1}}(\sigma)$ be Whittaker vectors of Iwahori types $\mu$ and $\lambda$, respectively. Then:
$$
    Z \left(\frac{1}{2},  \pi(h_n^{\left(-\varpi^{\alpha}\right)})W_{\pi}, W_{\sigma} \right) =  \Delta_\beta \frG(\ulmu, \ullam) W_{\pi}(\jmath^{\GL}(\bft(\varpi^{\alpha-\beta}))) W_{\sigma}(\bft(\varpi^{\alpha-\beta})),
$$
where $\alpha$ is the level of the Iwahori subgroup $\Iw_{\alpha}$ we have fixed, $\alpha, \beta > 0$ are sufficiently large integers (to be determined in the proof), and
$$
\Delta_\beta := \prod_{j=1}^{n}(1-q^{-j}) q^{-\frac{(n+2)(n+1)n \beta}{6}}, \quad \frG(\ulmu, \ullam) := \prod_{j=1}^{n}\prod_{i=1}^{j} \mu_{j+1-i}\lambda_{i}(\varpi^{\beta}) \frg(\mu_{j+1-i}\lambda_{i}).
$$
\end{theorem}

\begin{proof}
This is an application of Januszewski's generalization of the local Birch lemma (i.e., \cite[Theorem 2.8]{januszewski2024nonabelian}). We fix an auxiliary local character $\kappa: F^{\times} \rightarrow \barQQ_p^{\times}$ satisfying $\kappa(\varpi) = 1$, whose conductor $\varpi^{\gamma}$ satisfies $\gamma \leq \alpha$; for determinacy, we require $\gamma = \alpha$. Note that the Whittaker function $W_{\pi}$ has Iwahori type $\mu$, and the Whittaker function $W_{\sigma} \cdot (\kappa^{-1} \circ \det)$, i.e. the Whittaker function $W_{\sigma}$ twisted by the character $\kappa^{-1}$, has Iwahori type $\lambda \kappa^{-1}$. \footnote{This is not quite direct. To see this, we only need to verify that when the conductor exponent $\gamma$ of the character $\kappa$ satisfies $\gamma \leq \alpha$, for $g = (g_{i,j}) \in \Iw_{\alpha}$, we always have $\kappa( \det g ) = \kappa(g_{1,1} \cdots g_{n,n})$. In fact, we write
\[
\det g = \sum_{\tau \in \frS_{n}} \prod_{i=1}^{n} g_{i, \tau(i)} = \prod_{i=1}^{n} g_{i,i} + \sum_{\tau \neq \id{}} \prod_{i=1}^{n} g_{i, \tau(i)}.
\]
For a non-identity permutation $\tau \in \frS_{n}$, there must exist $j \in \{1, \ldots, n\}$ such that $\tau(j) < j$. In this case, $g_{j, \tau(j)}$ is located to the bottom-left of the main diagonal of $g$, so by the definition of the Iwahori subgroup, $g_{j, \tau(j)} \in \varpi^{\alpha} \fro$. Since all $g_{ij} \in \fro$ and all $g_{i,i} \in \fro^{\times}$ (because $g \pmod{\varpi}$ is also an upper triangular matrix, and its determinant is non-zero since $g \in \GL_n(\fro)$), we obtain
\[
\det g \in  \left( \prod_{i=1}^{n} g_{i,i} \right) + \varpi^{\alpha} \fro = \left( \prod_{i=1}^{n} g_{i,i} \right) (1+ \varpi^{\alpha} \fro).
\]
Recall that we required the conductor exponent $\gamma$ of $\kappa$ to satisfy $\gamma \leq \alpha$, hence $\kappa(1+ \varpi^{\alpha} \fro) = 1$, so $\kappa( \det g ) = \kappa(g_{1,1} \cdots g_{n,n})$, as desired. We point out that in the previous version of the preprint \cite[Section 5.5]{xu2025} written by the author in 2024, the Whittaker functions $W_{\pi}$ and $W_{\sigma}$ were directly substituted into Januszewski's result with an attempt to cancel the highly ramified local character $\kappa$ using a correction matrix $\bfs(l)$—this attempt was incorrect and futile. In this paper, when $\sigma$ is a regular ordinary representation and $W_{\sigma}$ is taken as the big cell Whittaker function $W_{\sigma}^{\dagger}$, its twist $W_{\sigma}^{\dagger} \cdot (\kappa^{-1} \circ \det)$ generally no longer corresponds to the big cell Whittaker function of the principal series representation $I_{B_{n}}^{\GL_{n}}(\mu \kappa^{-1})$. However, following the above discussion, when $\gamma \leq \alpha$, it remains a Whittaker function of that representation and has Iwahori type of level $\alpha$, so \cite[Theorem 2.8]{januszewski2024nonabelian} can still be applied. We also note that if restricted to the special case where both $W_{\pi}$ and $W_{\sigma}$ are the corresponding big cell Whittaker functions, \cite[Proposition 5.10]{hsieh2023fivevariablepadiclfunctionsu3times} provides a direct calculation for this local JPSS integral in the $\GL_2 \times \GL_3$ case, but generalizing that calculation to higher rank groups does not appear to be straightforward.
}
By Assumption \ref{ass:conductor}, for all $1 \leq i \leq n+1$ and all $1 \leq j \leq n$, the conductors of $\kappa \mu_i (\kappa^{-1} \lambda_j) = \mu_i \lambda_j$ are nontrivial and equal; let $\varpi^{\beta}$ denote a generator of this conductor. This allows us to apply \cite[Theorem 2.8]{januszewski2024nonabelian} to these two Whittaker functions and take the character $\chi$ therein to be the auxiliary character $\kappa$ mentioned above. We obtain \footnote{Note that we have already assumed $\kappa(\varpi) = 1$, thus $\kappa(\det \bft(-\varpi^{\alpha})) = 1$.}
\begin{multline*}
 \int_{U_n(F) \bs \GL_n(F)} \pi(h_n \jmath^{\GL}(\bft(-\varpi^{\alpha}))) W_{\pi}(\jmath^{\GL}(g)) \cdot \sigma(\bft(-\varpi^{\alpha}))W_{\sigma}(g) \dif g \\ =  \Delta_\beta \frG(\ulmu, \ullam) W_{\pi}(\jmath^{\GL}(\bft(\varpi^{\alpha-\beta}))) W_{\sigma}(\bft(\varpi^{\alpha-\beta})).    
\end{multline*}
The integral on the left side is:
\begin{align*}
    &Z \left(\frac{1}{2}, \pi( h_n \jmath^{\GL}(\bft(-\varpi^{\alpha})))W_{\pi}, \sigma(\bft(-\varpi^{\alpha}))W_{\sigma} \right) \\
    &\xlongequal{\eqref{eq:invarianceJPSS}} Z \left(\frac{1}{2}, \pi(\jmath^{\GL}(\bft(-\varpi^{\alpha}))^{-1})  h_n \jmath^{\GL}(\bft_{-\varpi^{\alpha}}))W_{\pi}, W_{\sigma} \right) \\
    &= Z \left(\frac{1}{2},  \pi(\jmath^{\GL}(\bft(-\varpi^{\alpha}))^{-1} h_n \jmath^{\GL}(\bft(-\varpi^{\alpha})))W_{\pi}, W_{\sigma} \right) \\
    &= Z \left(\frac{1}{2},  \pi(h_n^{\left(-\varpi^{\alpha}\right)})W_{\pi}, W_{\sigma} \right),
\end{align*}
as desired.
\end{proof}

\subsection{Local Ichino--Ikeda integrals at places above $p$}
We now return to the context of local Ichino--Ikeda integrals. Let $\sigma = I_{B_{m+n}}^{\GL_{m+n}}(\ullam)$ and $\pi = I_{B_{m+n+1}}^{\GL_{m+n+1}}(\ulmu)$ be regularly ordinary representations of $\GL_{m+n}(F)$ and $\GL_{m+n+1}(F)$, respectively.

We set
\[
W_1 = \pi(h_{m+n}^{(-\varpi^{\alpha})})W_{\pi}^{\dagger}, 
\quad W_1^{\prime} = W_{\sigma}^{\dagger}, \quad W_2 = \pi^{\vee}(h_{m+n}^{(-\varpi^{\alpha})})W_{\pi^{\vee}}^{\dagger}, \quad W_2^{\prime} = W_{\sigma^{\vee}}^{\dagger}
\]
in Corollary \ref{cor:splitii}. By Theorem \ref{thm:janbirch}, together with Remark \ref{rem:ordcontra} and Proposition \ref{prop:iwahoriact}, we obtain the following:

\begin{proposition}
With notation and conventions as above, let $\kappa: F^{\times} \rightarrow \barQQ_p^{\times}$ be a sufficiently ramified character with $\kappa(\varpi)=1$. Then
\[
\scI (\sigma, \pi) = 
\Delta_\beta^2 \cdot \frG(\breve{\ulmu}, \breve{\ullam}) \frG(\breve{\check{\ulmu}}, \breve{\check{\ullam}}) \cdot \calC_{\sigma}^{\ord} \calC_{\pi}^{\ord} \cdot \prod_{i=1}^{m+n-1} \zeta_{F}(i),
\]
where
\[
\calC_{\pi}^{\ord} := \dfrac{W_{\pi}^{\dagger}(\jmath^{\GL}(\bft(\varpi^{\alpha-\beta}))) W_{\pi^{\vee}}^{\dagger}(\jmath^{\GL}(\bft(\varpi^{\alpha-\beta})))}{\lrangle{W_{\pi}^{\dagger}, W_{\pi^{\vee}}^{\dagger}}} \text{ and } 
\calC_{\sigma}^{\ord} := \dfrac{W_{\sigma}^{\dagger}(\bft(\varpi^{\alpha-\beta}))  W_{\sigma^{\vee}}^{\dagger}(\bft(\varpi^{\alpha-\beta}))}{\lrangle{W_{\sigma}^{\dagger}, W_{\sigma^{\vee}}^{\dagger}}}.
\]
\end{proposition}

%% file: AU04-Summary.tex
\section{A summary of automorphic computations} \label{sec:localdoubl}
In this section, we describe the choice of the appropriate local Siegel sections $f_{s,\chi, v}^{\Sieg}$ from \cite{wan2015families} and summarize the results obtained thus far.

\subsection{The choice of Siegel Eisenstein sections}
Recall that we have previously established Assumption \ref{ass:scalarwt} and Assumption \ref{ass:generic}, which allow us to directly utilize the Hida family of Klingen Eisenstein series and the results of the local doubling integral computations in \cite{wan2015families}.

\subsubsection{Archimedean places}
Let $v \in \scV_{\calF}^{\infty}$ be an archimedean place of $\calF$. Following \cite[Section 4A2]{wan2015families}, let 
\[
\bfi := \diag \left[\dfrac{1}{2} \rmi \bfone_{n}, \rmi, \dfrac{1}{2}\vartheta, \dfrac{1}{2} \rmi \bfone_{n} \right]
\]
be the distinguished point in the unbounded realization of the symmetric domain for $\GU(m+n+1, m+n+1)$. We define the Siegel Eisenstein section as
\[
f^{\Sieg}_{v, s, \chi}\left(g = \begin{bmatrix}
A_g & B_g \\ C_g & D_g
\end{bmatrix} \right) = \det(C_g \bfi + D_g)^{-\kappa} \abs{\det(C_g \bfi + D_g)}^{\kappa-2s-(m+n)}, \quad g \in G^{\bdsuit}(\calF_v),
\]
which depends on the weight $\kappa$.

\subsubsection{Unramified places}
For $v \in \scU_{\calF}^{\ur}$, we choose $f_{v,s,\chi}^{\Sieg}$ to be the spherical section $f_{v,s,\chi}^{\Sieg, \sph}$.

\subsubsection{Ramified places away from $p$}
At finite places $v \in \scU_{\calF}^{\bad}$ not dividing $p$, let $f_{s, \chi,v}^{\dagger}$ be the big-cell section defined as the Siegel section supported on the big cell $Q(\calF_v) w_{N+1} N_{Q}(\calO_{\calF_v})$ such that $f_{s, \chi,v}^{\dagger}(w_{N+1} N_{Q}(\calO_{\calF_v})) = 1$. We set
    \[
    \gamma_v := \begin{bmatrix}
    \bfone_{N+1} & \gamma_v^{\UR} \\ 0 & \bfone_{N+1} 
    \end{bmatrix}, \quad \gamma_v^{\UR} := \begin{bmatrix}
    0 & 0 & 0 & x^{-1} \bfone_{n} \\ 0 & 0 & 0 & 0 \\ 0 & 0 & (y \bary)^{-1} \bfone_{m-n} & 0 \\ \barx^{-1} \bfone_{n} & 0 & 0 & 0
    \end{bmatrix},
    \]
where $x$ and $y$ are fixed constants in $\calK$ divisible by a sufficiently high power of $\varpi_{\calK_v}$. We then define the Siegel section as $f_{s,\chi,v}^{\Sieg}(g) = f_{v, s,\chi}^{\dagger}(g \gamma_v)$.

\subsubsection{Places above $p$}
Let $v \in \scV_{\calF}^{(p)}$. We note that in \cite[Section 4D]{wan2015families}, $p$ is assumed to split completely in $\calK$, whereas we only assume that $p$ is unramified in $\calF$ and that every place of $\calF$ above $p$ splits in $\calK$. The more restrictive assumption in \textit{loc.\ cit.\ } is primarily to simplify notation and can be generalized to our setting through careful bookkeeping.

Recall that we began with a Hecke character $\chi$ of $\calK$, which decomposes as $\chi = \otimes_{w \in \scS_{\calK}} \chi_w$. For $v \in \scS_{\calF}^{(p)}$ splitting as $v = w \barw$ in $\calK$, we set
\[
\chi_1 := \chi_w, \quad \chi_2 := \chi_{\barw}^{-1},
\]
identifying $\calK_{w} = \calK_{\barw} = \calF_{v}$. These local characters determine a character
\[
\chi_{w,s}: B_{2(N+1)}(\calK_w) \rightarrow \CC^{\times}, \quad \begin{bmatrix}
A & B \\ 0 & D
\end{bmatrix} \mapsto \chi_1(\det D) \chi_2(\det A) \abs{\det AD^{-1}}_w^{s},
\]
where the matrix is partitioned as $(N+1) \times (N+1)$. One verifies that for each 
\begin{equation} \label{eq:GLSiegelsection}
f_{w,s} \in \Ind_{B_{2(N+1)}(\calK_w)}^{\GL_{2(N+1)}(\calK_w)}(\chi_{w,s}),
\end{equation}
we have $f_{v,s}(g_v) := f_{w,s}(g_w) \in I_v^{\Sieg}((\chi_1, \chi_2), s)$, where $g_w$ is the projection of $g_v$ under $\varrho_{w, N+1, N+1}$. This identification reduces the construction to sections in \eqref{eq:GLSiegelsection} over $\GL_{2(N+1)}(\calK_w)$. Furthermore, we have the local representation $\sigma_v = I_{B_{N}}^{\GL_{N}}(\ullam)$, where $\ullam = (\lambda_{1}, \ldots, \lambda_{m+n})$, and the conductor of each $\lambda_{i}$ is $\varpi_v^{t_i}$. We define
\[
\xi_i := \begin{cases}
\lambda_i \chi_1^{-1}, &\quad 1 \leq i \leq m, \\
1, &\quad i=m+1, \\
\lambda_{i-1}^{-1} \chi_2, &\quad m+2 \leq i \leq m+n+1.
\end{cases}
\]

Let $\tilde{f}_{v,s,\chi}^{\dagger}$ be the big-cell Siegel Eisenstein section such that \footnote{This corresponds to $\tilde{f}_t$ in the notation of \cite{wan2015families}.}
\begin{itemize}
    \item it is supported on the big cell $Q(\calF_v) w_{N+1} K_{Q}(p^{t})$, where $K_{Q}(p^t)$ is the subgroup of $\GL_{2(N+1)}(\calO_{\calF_v})$ consisting of matrices that are blockwise upper triangular modulo $p^t$ under the partition $(N+1) \times (N+1)$, and $t$ is such that the conductor of $(\chi^{\calF})_v$ is $\varpi_v^{t}$;
    \item $\tilde{f}_{v,s,\chi}^{\dagger}\left(w_{N+1} \begin{bmatrix}
A & B \\ C & D
\end{bmatrix} \right) = \chi_v(\det D)$ for $ \begin{bmatrix}
A & B \\ C & D
\end{bmatrix} \in K_{Q}(p^t)$.
\end{itemize}
Following \cite[Section 4D4]{wan2015families}, we define the Siegel Eisenstein section \footnote{There are subtle differences between the Siegel Eisenstein series defined here and those in \cite[Section 4D4]{wan2015families}. Following \cite[Section 4.7]{wan2019iwasawatheorymathrmursblochkato}, we decide to add the factor $\left( \frg(\barchi_{\sfP, v})^{m+n+1} c_{m+n+1}(\chi_{\sfP, v}, -s_{\sfP}) \right)^{-1}$ directly to the Siegel Eisenstein section, rather than including it in the normalization factor $B_{\scD}$ at the beginning of \cite[Section 5C1]{wan2015families}.}
\begin{multline*}
f_{v,s,\chi}^{\Sieg} := p^{- \sum_{i=1}^{m}i t_i - \sum_{i=1}^{n}i t_{m+i}} \prod_{i=1}^{m} \frg(\xi_i) \xi_{i}(-1) \prod_{i=1}^{n} \frg(\xi_{m+1+i}) \xi_{m+1+i}(-1) \\
\times \sum_{A, B, C, D, E} \prod_{i=1}^{m-n} \barxi_i \left(\dfrac{\det A_i}{\det A_{i-1}} p^{t_i} \right) \prod_{i=1}^{n} \barxi_{m-n+i} \left(\dfrac{\det D_i}{\det D_{i-1}} p^{t_{m-n+i}} \right) \prod_{i=1}^{n} \barxi_{m+1+i} \left(\dfrac{\det E_i}{\det E_{i-1}} p^{t_{m+1+i}} \right)   \\
\times \left( \frg((\barchi^{\calF})_{v})^{m+n+1} c_{m+n+1}((\chi^{\calF})_{v}, -s) \right)^{-1} \tilde{f}^{\dagger}_{v,s,\chi}\left( g w_{\Borel}^{-1} \begin{bmatrix}
\bfone_{N+1} & \circledast \\ 0 & \bfone_{N+1}
\end{bmatrix} w_{\Borel} \right),
\end{multline*}
where 
\[
\circledast := \begin{bmatrix}
 & & C & D \\ & & & \\ & & A & B \\ E & & &
\end{bmatrix}
\]
is partitioned as $[n \mid 1 \mid m-n \mid n]$ and $A, B, C, D, E$ range over the set defined in \cite[Lemma 4.29]{wan2015families}. Here $A_i$ is the $i$-th upper-left minor of $A$, $D_i$ is the $(m-n+i)$-th upper-left minor of $\begin{bmatrix} A & B \\ C & D \end{bmatrix}$, and $E_i$ is the $i$-th upper-left minor of $E$. The element $w_{\Borel} \in G^{\bdsuit}(\calF_v)$ has projection $\diag[w_{m}, 1] \in \GL_{2(N+1)}(\calK_w)$. The factor $c_{m+n+1}(-)$ is as defined in \cite[(13)]{wan2015families}.

\subsection{Local doubling integrals} \label{sec:doublingintegrals}
The local doubling integrals under these choices of local Siegel sections were computed in \cite{wan2015families}:
\begin{itemize}
    \item For $v \in \scV_{\calF}^{\infty}$, the integral is computed in \cite[Lemma 4.3 (ii)]{wan2015families}.
    \item For $v \in \scU_{\calF}^{\ur}$, it is computed in \cite[Lemma 4.5]{wan2015families}.
    \item For $v \in \scU_{\calF}^{\bad}$ not lying above $p$, it is computed in \cite[Lemma 4.10(ii)]{wan2015families}. Note that although the local representation $\pi_v$ might be unramified at $v$, the computation remains valid.
    \item For $v \in \scV_{\calF}^{(p)}$, it is computed throughout \cite[pp.~2005--2017]{wan2015families}, with the results summarized on \cite[p.~2021]{wan2015families} (see also \cite[Section 4.7]{wan2019iwasawatheorymathrmursblochkato}).
\end{itemize}
These computations show that Assertion \ref{assertion:localdoubling} holds under our choice of Siegel Eisenstein sections. To avoid overly cumbersome notation, we do not record the explicit expressions of the local doubling integrals in this article.

\begin{remark} \label{rem:modificationsections}
The only difference between our choices and those in \cite{wan2015families} is that we take spherical sections only at places in $\scU_{\calF}^{\ur}$, which is a (possibly proper) subset of the set $\scS_{\calF}^{\ur}$ used in \cite{wan2015families}. At the remaining unramified places $v \in \scS_{\calF}^{\ur} \smallsetminus \scU_{\calF}^{\ur}$, we use the big-cell section $f^{\dagger}_{s,\chi,v}$ in place of the spherical one. This modification affects neither the construction of the $p$-adic $L$-functions nor the ordinary Klingen Eisenstein family in \cite{wan2015families}. In particular, the computations of the local doubling integrals in \cite[Section~4C]{wan2015families} still apply at these places, even if $\pi_v$ is unramified.

Our reason for making this slight change is that the evaluation of the local doubling integrals $\scZ^{\dsuit, \rightt}_v(f^{\Sieg, \sph}_{v,s,\chi}, \pi_v)$ for ramified $\pi_v$ is not currently available. Moreover, this adjustment is required in Theorem~\ref{thm:newpadicl} to compare the interpolation formula for $\calL_{\Wan}^{\scU}(\bfg,\calK,\chi_0)$ with our computation of the Gan--Gross--Prasad type period integral.

One could avoid any modification by imposing the stronger assumption that $\pi$ is ramified only at places above $p$. However, such an assumption would prevent us from applying our machinery to the Iwasawa theory for CM fields in Part~\ref{part:application}, and is therefore too restrictive for our purposes \footnote{More precisely, in Lee's construction of the pullback theta lift, the ideal $\frs$ of $\calO_{\calK}$ is chosen such that $\frs = \frf \overline{\frf}$ with $\frf + \overline{\frf} = \calO_{\calK}$. This prevents $\frs$ from being the zero ideal; hence, the resulting ordinary CM family $\bTheta_{\lambda, \frs}$ must be ramified at some places of $\calF$ away from $p$.}. 
\end{remark}

\subsection{Summary of automorphic computations}
Combining all the calculations in this part, we obtain the following main theorem.
\begin{theorem} \label{thm:automorphicfinal}
With notation as above, assume \ref{item:irred} for $\sigma$ and $\pi$, as well as Assumptions \ref{ass:wtinterlacing}, \ref{ass:spl}, \ref{ass:ram}, \ref{ass:scalarwt}, and \ref{ass:generic}. We choose local Siegel Eisenstein sections $f^{\Sieg}_{v,s,\chi}$ at places $v$ of $\calF$ as in Section \ref{sec:localdoubl} for the construction of the Klingen Eisenstein series. Suppose $\sigma$ and $\pi$ are both regularly ordinary at all places $v$ above $p$. Then
\begin{multline*}
\dfrac{\left(\calP^{\Kling}_{\Phi, \Psi}\right)^2}{\aabs{\Phi}_{\sigma, \MVW}^2 \aabs{\Psi}_{\pi, \MVW}^2} = \dfrac{1}{2^{\varkappa_{\sigma}+\varkappa_{\pi}}} \cdot \scL^{\scV_{\calF}^{(p)}}(\sigma \times \pi) L_{\scU_{\calF}^{\ur}}\left(s+\dfrac{1}{2}, \pi^{\vee}, \barchi^{\sfc} \right)^{2} \\
\times \prod_{v \in \mathscr{V}_{\calF}^{\ram}} \zeta_v C_{\sigma_v, \bpsi_v} C_{\sigma_v^{\vee}, \bpsi_v^{-1}} \calB_{\pi_v}^{\ess} \calB_{\sigma_v}^{\ess} \prod_{v \in \mathscr{V}_{\calF}^{(p)}} \zeta_v \Delta_\beta^2 \cdot \frG(\breve{\ulmu}, \breve{\ullam}) \frG(\breve{\check{\ulmu}}, \breve{\check{\ullam}}) \calC_{\sigma}^{\ord} \calC_{\pi}^{\ord}
\\
\times \prod_{v \in \scU_{\calF}^{\ur}} d_{N+1,v}(s, \barchi^{\sfc}_v)^{-2} \prod_{v \in \scU_{\calF}^{\bad}} \scZ^{\dsuit}_v(f^{\Sieg}_{v,s,\chi}, \pi_v)^2.
\end{multline*}
Here, for $v \in \scV_{\calF}^{\ram}$ such that $\sigma_v$ is unramified, we set $C_{\sigma_v, \bpsi_v} = C_{\sigma_v^{\vee}, \bpsi_v^{-1}} = 1$, and we define $\zeta_v := \prod_{i=1}^{m+n-1} \zeta_{\calF_v}(i)$.
\end{theorem}

%% file: AU05-padicGGP.tex
\section{The $p$-adic interpolation of Gan--Gross--Prasad type period integrals} \label{sec:padicggp}
Let $V$ be an $N$-dimensional vector space over $\calK$ equipped with a non-degenerate skew-Hermitian form $\phi$. Then $\sfi \phi$ is an Hermitian form on $V$, and we regard $V$ as an Hermitian space over $\calK$. In Section \ref{sec:preliminaries} and Section \ref{sec:formsunitary}, we considered the spaces $V^{\sharp}$, $V^{\hsuit}$, and $V^{\bdsuit}$. Let $\scP, \scP^{\sharp}, \scP^{\hsuit}$, and $\scP^{\bdsuit}$ be the simple unitary Shimura data attached to these Hermitian spaces, and let $H, G, G^{\hsuit}$, and $G^{\bdsuit}$ be the unitary groups attached to them, respectively. We maintain the assumptions \eqref{ass:sgn} and \eqref{ass:QS} from Section \ref{sec:unitarylocal}. We denote the algebraic groups $H_{\bullet}, B_{\bullet}, B^{+}_{\bullet}$, and $T_{\bullet}$ for these unitary Shimura data by $H_{r,s}, B_{r,s}, B^{+}_{r,s}$, and $T_{r,s}$ if the corresponding unitary group has signature $(r,s)$, adopting similar notation for other cases.

\begin{remark}[On unitary groups, {\cite[Remark 2.1]{wan2015families}}] \label{rem:unitarygrp}
As recalled in Section \ref{sec:formsunitary}, to obtain Shimura varieties for studying $p$-adic automorphic forms and Galois representations, we use the general unitary group $G$ defined over $\QQ$, which is smaller than the general unitary group defined in Definition \ref{defn:unitarygrp}. However, this group is less convenient for local automorphic computations since we cannot treat the primes of $\calF$ independently. Consequently, for automorphic computations, we implicitly consider automorphic forms on the larger general unitary groups as in Definition \ref{defn:unitarygrp} and then restrict them to the smaller ones. For the algebraic construction, we only perform pullbacks for unitary groups rather than general unitary groups.
\end{remark}

\subsection{Gan--Gross--Prasad datum}
We define the following datum for the $p$-adic interpolation of the Gan--Gross--Prasad period integral.

\begin{definition} \label{defn:ggpdatum}
A \emph{$p$-adic Gan--Gross--Prasad datum} is a tuple 
\[\bfG = (L, \II, \JJ, \eta_0, \bff, \bfg, \psi_{0}, \psi_{0}^{\prime}, \chi_0, \scS, \scT, \scU), \]
where
\begin{itemize}
    \item the tuple $\bfD_{H} := (L, \II, \eta_0, \bff, \psi_0, \chi_0, \scS)$ is a $p$-adic Eisenstein datum over $H$;
    \item the tuple $\bfD_{G} := (L, \JJ, \eta_0, \bfg, \psi_0^{\prime}, \chi_0, \scT)$ is a $p$-adic Eisenstein datum over $G$;
    \item $\scU$ is a finite set of places of $\calF$ containing $\scS \cup \scT$.
\end{itemize}
Recall that we defined $\scS_{\bff}$ (resp. $\scT_{\bfg}$) in Definition \ref{defn:padiceisendatum} as the set of finite places $v$ of $\calF$ away from $p$ such that $K_{H,v} \neq H(\calF_v)$ (resp. $K_{G,v} \neq G(\calF_v)$); by definition, $\scU$ contains $\scS_{\bff}$ and $\scT_{\bfg}$.
\end{definition}

We define the \emph{weight space of the $p$-adic GGP datum} $\bfG$ as the product space $\calX_{\bfG} := \calX_{\bfD_{H}} \times \calX_{\bfD_{G}}$, where $\calX_{\bfD_{H}}$ (resp. $\calX_{\bfD_{G}}$) is the weight space for the $p$-adic Eisenstein datum $\bfD_{H}$ (resp. $\bfD_{G}$) defined in \eqref{eq:weightspace}.

\begin{definition} \label{defn:ggpweightpoint}
Let $\sfP \times \sfQ \in \calX_{\bfG}$, where $\sfP \in \calX_{\bfD_{H}}$ and $\sfQ \in \calX_{\bfD_{G}}$. Following Definition \ref{defn:specializations}:
\begin{enumerate}[label = \rm (\arabic*)]
    \item we call $\sfP \times \sfQ$ an \emph{admissible point} if $\sfP$ and $\sfQ$ are classical arithmetic points and the pair $(\sigma_{\bff_{\sfP}, \infty}, \pi_{\bfg_{\sfQ}, \infty})$ satisfies the weight interlacing assumption (Assumption \ref{ass:wtinterlacing});
    \item we call $\sfP \times \sfQ$ a \emph{generic admissible point} if it is an admissible point, $\sfP \in \calX_{\bfD_{H}}^{\gen}$, and $\sfQ \in \calX_{\bfD_{G}}^{\gen}$.
\end{enumerate}
The sets of such points are denoted by $\calX_{\bfG}^{\adm}$ and $\calX_{\bfG}^{\gen}$, respectively. Both are Zariski dense in $\calX_{\bfG}$.
\end{definition}

\subsection{The $p$-adic family of Eisenstein series}
In this subsection, we recall the construction of $p$-adic families of Siegel Eisenstein series and Klingen Eisenstein series from \cite{wan2015families}, with the slight modification discussed in Remark \ref{rem:modificationsections}.

\subsubsection{The Siegel Eisenstein family}
To perform $p$-adic interpolation, we require a renormalization of Siegel Eisenstein series. For $\sfP \in \calX_{\bfD_{H}}^{\gen}$, we define the normalization factor
\begin{equation} \label{eq:normalizationB}
B(s_{\sfP},\chi_{\sfP}) := \Omega_{\sfP}^{-1} \left(\dfrac{(-2)^{-n}(2 \pi \rmi)^{n \kappa_{\sfP}} (2/\pi)^{n(n-1)/2}}{\prod_{j=0}^{n-1} \Gamma(\kappa_{\sfP}-j)} \right)^{-d} \cdot \prod_{v \in \scU_{\calF}^{0} \smallsetminus \scU} d_{m+n+1,v}(s_{\sfP}, \barchi^{\sfc}_{\sfP}),
\end{equation}
where the period factor $\Omega_{\sfP}$ is as defined in \eqref{eq:defineperiod}. We define the \emph{normalized Siegel Eisenstein series} by
\[
E^{\Sieg}_{\chi_{\sfP}}(-) = B(s_{\sfP}, \chi_{\sfP}) \cdot E^{\Sieg}(f^{\Sieg}_{s_{\sfP}, \chi_{\sfP}}, -)
\]
and the \emph{normalized Klingen Eisenstein series} by
\[
E^{\Kling}_{\bfD_{\sfP}}(-) = B(s_{\sfP}, \chi_{\sfP}) \cdot E^{\Kling}(F^{\hsuit}(f^{\Sieg}_{s_{\sfP}, \chi_{\sfP}}, \bff_{\sfP}; -), -),
\]
where $\bfD_{\sfP}$ is the specialized Eisenstein datum over $H$ at $\sfP$.

\begin{theorem}[{\cite[Lemma 5.7]{wan2015families}}]
Attached to the $p$-adic Eisenstein datum $\bfD$, there exists a $p$-adic measure 
\[
\bfE^{\Sieg}_{\chi_0} \in \Meas(\Gamma_{\calK} \times T_{m,n}^{\circ}, V_{G^{\bdsuit}}^{\ord}(K_{G^{\bdsuit}}^{p}, \calO_{L}^{\ur}))
\]
such that for any generic classical arithmetic point $\sfP \in \calX^{\gen}_{\bfD}$, we have
\[
\int_{\Gamma_{\calK} \times T_{m,n}^{\circ}} \sfP \dif \bfE^{\Sieg}_{\chi_0} = E^{\Sieg}_{\chi_{\sfP}, s_{\sfP}}.
\]
\end{theorem}

The construction proceeds by first interpolating the Fourier coefficients of the Siegel Eisenstein series and then applying the $q$-expansion principle. The normalization factor $B(s_{\sfP},\chi_{\sfP})$ ensures that these Fourier coefficients are integral and $p$-adically interpolatable.

\subsubsection{Hecke projectors}
We recall the notion of Hecke projectors attached to a $p$-adic Eisenstein datum $\bfD$, following \cite[Section 5B]{wan2015families}.

\begin{definition}[Dual Hida family]
Let $\vee: \Lambda_{m,n}^{\circ} \rightarrow \Lambda_{m,n}^{\circ}$ be the $\calO_{L}$-involution sending $\diag[a_1, \ldots, a_n]$ to $\diag[a_n^{-1}, \ldots, a_1^{-1}]$. We define $\II^{\vee}$ to be the ring $\II$ with the $\Lambda_{m,n}^{\circ}$-algebra structure given by the composition of $\vee$ with the structure map $\Lambda_{m,n}^{\circ} \rightarrow \II$. An $\II^{\vee}$-adic cuspidal Hida family $\bff^{\vee}$ is a \emph{dual Hida family} of $\bff$ if, for all generic classical arithmetic points $\sfP \in \calX_{\bfD}^{\gen}$, we have $(\bff^{\vee})_{\sfP} \in \sigma_{\bff_{\sfP}}^{\vee}$.
\end{definition}

\begin{definition}[Hecke projectors]
Let $\bff$ be an $\II$-adic Hida family over $H$. A Hecke operator $\bbone_{\bff}$ in the Hecke algebra $\bfh^{0, \ord, \scS_{\bfD}}(\Frac(\II))$ is a \emph{Hecke projector of $\bff$} if, for any $\sfP \in \calX_{\bfD}^{\gen}$ and any $\II$-adic family of automorphic forms $\bff^{\prime}$, the specialization of $(\bbone_{\bff} \circ \ee)\bff^{\prime}$ at $\sfP$ is the projection (under the Petersson inner product) of $\bff^{\prime}_{\sfP}$ to the one-dimensional line in $V_{\sigma_{\bff}}$ spanned by $\bff_{\sfP}$. The scalar $(\Proj_{\bff})_{\sfP}(\bff^{\prime})$ is defined such that
\[
((\bbone_{\bff} \circ \bfe)\bff^{\prime})_{\sfP} = \dfrac{(\Proj_{\bff})_{\sfP}(\bff^{\prime})}{\aabs{\bff_{\sfP}}_{\pi_{\bff_{\sfP}}}} \cdot \bff_{\sfP}.
\]
This defines a map $\Proj_{\bff}: \calM_{H}(K^p_H, \II) \rightarrow \Frac(\II)$.
\end{definition}

\begin{theorem}[Wan]
Assume notation as above and suppose that $\bff$ satisfies \ref{item:irredbff}. Then the dual Hida family $\bff^{\vee}$ and the Hecke projector $\bbone_{\bff}$ exist.
\end{theorem}

This is proved in \cite[Chapter 6]{wan2019iwasawatheorymathrmursblochkato}. More specifically, the Hecke projector is constructed in \cite[Equation (28)]{wan2019iwasawatheorymathrmursblochkato}, and the dual Hida family is constructed in the proof of \cite[Theorem 6.8]{wan2019iwasawatheorymathrmursblochkato}. We remark that the existence of such projectors was well-known in the lower-rank cases of $\rmU(1,1)$ and $\rmU(2,0)$ in the classical theory prior to Wan's general construction. One may refer to, for example, \cite[Remark 5.6]{wan2015families} for further justification.

\subsubsection{The Klingen Eisenstein family}
Following the recipe of the pullback formula (Proposition \ref{prop:pbkling}) and the $p$-adic interpolation construction in \cite[Section 5C]{wan2015families}, we first restrict the Siegel Eisenstein family to $G^{\hsuit} \times H$. Then, by applying a twist of $\chi_{\sfP}^{-1}(\det g_2)$ to the $H$-factor for every specialization $\sfP$, we obtain
\[
\bfE^{\Sieg}_{\chi} |^{\star}_{G^{\hsuit} \times H} \in \Meas(\Gamma_{\calK} \times T_{m,n}^{\circ}, V_{G^{\hsuit}}(K_{G^{\hsuit}}^p, \calO_{L}^{\ur}) \otimes V_{H}(K_{H}^p, \calO_{L}^{\ur})).
\]
By our choice of $f^{\Sieg}_{s,\chi,v}$ for $v \in \scV_{\calF}^{(p)}$ and \eqref{eq:productmeasure}, this further implies
\[
\bfE^{\Sieg}_{\chi} |^{\star}_{G^{\hsuit} \times H} \in \Meas(\Gamma_{\calK}, V_{G^{\hsuit}}(K_{G^{\hsuit}}^p, \calO_{L}^{\ur}) \otimes \Meas(T_{m,n}^{\circ}, V_{H}(K_{H}^p, \calO_{L}^{\ur}))).
\]
We apply the Hecke projector $\Proj_{\bff^{\vee}}$ to the $V_{H}$-part of $\bfE^{\Sieg}_{\chi} |^{\star}_{G^{\hsuit} \times H}$, which yields
\begin{equation} \label{eq:constructklingen}
\bfE_{\bfD_{H}}^{\Kling} := \Proj_{\bff^{\vee}}(\bfE^{\Sieg}_{\chi} |^{\star}_{G^{\hsuit} \times H}) \in \Meas(\Gamma_{\calK}, V_{G^{\hsuit}}(K_{G^{\hsuit}}^p, \calO_{L}^{\ur})) \otimes \Frac(\II^{\ur}).
\end{equation}
This is the \emph{Klingen Eisenstein family} associated to the $p$-adic Eisenstein datum $\bfD_{H}$, satisfying the property that for any generic classical arithmetic point $\sfP \in \calX^{\gen}_{\bfD_{H}}$, we have
\[
\int_{\Gamma_{\calK}} \sfP \dif \bfE^{\Kling}_{\bfD_{H}} = B(s_{\sfP}, \chi_{\sfP}) E^{\Kling}_{\bff_{\sfP}, \chi_{\sfP}, s_{\sfP}}.
\]
In fact, $\bfE_{\bfD_{H}}^{\Kling}$ is an ordinary family, i.e., it lies in $\Meas(\Gamma_{\calK}, V_{G^{\hsuit}}^{\ord}(K_{G^{\hsuit}}^p, \calO_{L}^{\ur})) \otimes \Frac(\II^{\ur})$.

In the case where $H$ is a definite unitary group (i.e., $H = \rmU(m,0)$), there is an alternative approach that avoids Hecke projectors by using the $p$-adic interpolation of the Petersson inner product. This approach exploits the fact that the Shimura variety of $\rmU(m,0)$ is a finite set (see, e.g., \cite[Definition 5.1 (ii)]{wan2015families}). This technique was introduced in \cite{hsieh2014eisenstein} for $\rmU(2,0)$ and generalized to $\rmU(m,0)$ in \cite[Section 5C3]{wan2015families}.

\begin{proposition}[{\cite[Proposition 5.9]{wan2015families}}] \label{prop:padicpetersson}
Let $K^{p}$ be a neat tame level group in $H(\AA_{\QQ, \rmf}^{p})$. There exists a $\Lambda_{\bfD}$-linear pairing
$$
\bfB_{K^p}: \calM_H(K^p; \II) \times \calM_H(K^p; \II^{\vee}) \rightarrow \II
$$
such that for any $\sfP \in \calX_{\bfD}^{\gen}$, we have $\sfP(\bfB_{K^p}(\bff, \bff^{\prime})) = \lrangle{\bff_{\sfP}, \bff^{\prime}_{\sfP}}_{\sigma_{\bff_{\sfP}}, \Pet}$ for any $\bff \in \calM_H(K^p; \II)$ and $\bff^{\prime} \in \calM_H(K^p; \II^{\vee})$.
\end{proposition}

In this case, the Klingen Eisenstein family can be defined by replacing \eqref{eq:constructklingen} with 
\begin{equation*}
\bfE_{\bfD_{H}}^{\Kling} := \bfB_{K_H^p}(\bfE^{\Sieg}_{\chi} |^{\star}_{G^{\hsuit} \times H}, \bff^{\vee}) \in \Meas(\Gamma_{\calK}, V_{G^{\hsuit}}(K_{G^{\hsuit}}^p, \calO_{L}^{\ur})) \otimes \II.
\end{equation*}
Note that in this definite case, the Klingen Eisenstein family $\bfE_{\bfD_{H}}^{\Kling}$ is $\II^{\ur}$-integral. Furthermore, $\bfE_{\bfD_{H}}^{\Kling}$ is an ordinary family, i.e., it lies in $\Meas(\Gamma_{\calK}, V_{G^{\hsuit}}^{\ord}(K_{G^{\hsuit}}^p, \calO_{L}^{\ur})) \otimes \Frac(\II)$.

\subsection{The $p$-adic interpolation of GGP period integrals}
By applying appropriate Hecke projectors to $p$-adically interpolate the GGP period integral, we obtain the following theorem.

\begin{theorem} \label{thm:interpolationfinal}
There exists an element 
\[
\bfP_{\bfG} \in (\Lambda_{\bfD_{H}} \otimes \Frac(\II^{\ur})) \otimes \Frac(\JJ^{\ur})
\]
such that for any $\sfP \times \sfQ \in \calX_{\bfG}^{\gen}$,
\[
(\sfP \times \sfQ)(\bfP_{\bfG}) = B(s_{\sfP}, \chi_{\sfP}) \calP^{\Kling}(\bff_{\sfP}, \bfg_{\sfQ}, \chi_{\sfP}, s_{\sfP}),
\]
where the right-hand side is the GGP period integral of the Klingen Eisenstein series with a cuspidal automorphic form, as defined in \eqref{eq:periodKLE} and normalized by the factor in \eqref{eq:normalizationB}. When $n=0$, we have $\bfP_{\bfG} \in \Lambda_{\bfD_{H}} \otimes \JJ^{\ur}$.
\end{theorem}

\begin{proof}
Recall the Klingen Eisenstein family
\[
\bfE^{\Kling}_{\bfD_{H}} \in \Meas(\Gamma_{\calK}, V_{G^{\hsuit}}(K_{G^{\hsuit}}^p, \calO_{L}^{\ur})) \otimes \Frac(\II^{\ur}) \subset \Meas(\Gamma_{\calK}, \calM_{G^{\hsuit}}(K_{G^{\hsuit}}^p, \calO_{L}^{\ur})) \otimes \Frac(\II^{\ur}).
\]
Restricting it to $G$, we obtain
\[
\bfE^{\Kling}_{\bfD_{H}}|_{G} \in \Meas(\Gamma_{\calK}, \calM_{G}(K_{G}^p, \calO_{L}^{\ur})) \otimes \Frac(\II^{\ur}).
\]
We then apply the Hecke projector $\Proj_{\bfg}$ to obtain
\begin{multline} \label{eq:projEKling}
\Proj_{\bfg}(\bfE^{\Kling}_{\bfD_{H}}|_{G}) \in \Meas(\Gamma_{\calK}, \Frac(\JJ^{\ur})) \otimes \Frac(\II^{\ur}) \\
= \Frac(\JJ^{\ur})\lrbracket{\Gamma_{\calK}} \otimes \Frac(\II^{\ur}) = (\Lambda_{\bfD_{H}} \otimes \Frac(\II^{\ur})) \otimes \Frac(\JJ^{\ur}).
\end{multline}
This element is the desired $\bfP_{\bfG}$ satisfying the stated interpolation property. When $n=0$, instead of using the Hecke projector $\Proj_{\bfg}$, we apply $\bfB_{K_{G}^{p}}(-, \bfg)$ to obtain
\[
\bfB_{K^p_{G}}(\bfE^{\Kling}_{\bfD_{H}}|_{G}, \bfg) \in \Meas(\Gamma_{\calK}, \JJ^{\ur}) = \JJ^{\ur}\lrbracket{\Gamma_{\calK}} \otimes \II^{\ur} = \Lambda_{\bfD_{H}} \otimes \JJ^{\ur},
\]
as required.
\end{proof}

To obtain a precise interpolation formula, we impose the following ``$p$-adically automorphic assumptions.''
\begin{assumption}[Automorphic assumptions] \label{ass:automorphic}
We assume that $\bff$ and $\bfg$ both satisfy condition \ref{item:irredbff}, that the sets $\scS_{\bff}$ and $\scS_{\bfg}$ are disjoint, and that every place in $\scS_{\bff} \cup \scS_{\bfg}$ splits in $\calK$.
\end{assumption}

In previous sections, we related the period integral $\calP^{\Kling}(\bff_{\sfP}, \bfg_{\sfQ}, \chi_{\sfP}, s_{\sfP})$ to certain $L$-values with explicit local factors at bad places. Combining Theorem \ref{thm:interpolationfinal} and Theorem \ref{thm:automorphicfinal}, we obtain the following theorem.

\begin{theorem} \label{thm:interpoformula}
With notation as above, assume Assumption \ref{ass:automorphic}. Then for any $\sfP \times \sfQ \in \calX_{\bfG}^{\gen}$, $(\sfP \times \sfQ)(\bfP_{\bfG})^2$ equals
\begin{multline*}
\dfrac{1}{2^{\varkappa_{\sigma_{\bff_{\sfP}}}+\varkappa_{\pi_{\bfg_{\sfQ}}}}} \cdot \scL^{\scV_{\calF}^{(p)}}(\sigma_{\bff_{\sfP}} \times \pi_{\bfg_{\sfQ}}) L^{\scU}\left(s_{\sfP}+\dfrac{1}{2}, \pi_{\bfg_{\sfQ},v}^{\vee}, \barchi^{\sfc}_{\sfP,v} \right)^2 \Omega_{\sfP}^{-2} \\
\times \left(\dfrac{(-2)^{-n}(2 \pi \rmi)^{n \kappa_{\sfP}} (2/\pi)^{n(n-1)/2}}{\prod_{j=0}^{n-1} \Gamma(\kappa_{\sfP}-j)} \right)^{-2d} \prod_{v \in \scS_{\bff} \cup \scS_{\bfg}} \zeta_v C_{\sigma_{\bff_{\sfP,v}, \bpsi_v}} C_{\sigma_{\bff_{\sfP}, v}^{\vee}, \bpsi_v^{-1}} \calB_{\pi_{\bfg_{\sfQ},v}}^{\ess} \calB_{\sigma_{\bff_{\sfP},v}}^{\ess} \\
\times \prod_{v \in \mathscr{V}_{\calF}^{(p)}} \zeta_v \Delta_{v,\beta_v}^2 \cdot \frG(\breve{\ulmu_{v,\sfP}}, \breve{\ullam_{v,\sfQ}}) \frG(\breve{\check{\ulmu_{v,\sfP}}}, \breve{\check{\ullam_{v,\sfQ}}}) \calC_{\sigma_{\bff_{\sfP}}}^{\ord} \calC_{\pi_{\bfg_{\sfQ}}}^{\ord}
\\
\times \aabs{\bff_{\sfP}}_{\sigma_{\bff_{\sfP}}, \MVW}^2 \aabs{\bfg_{\sfQ}}_{\pi_{\bfg_{\sfQ}}, \MVW}^2 \prod_{v \in \scU} \scZ^{\dsuit, \rightt}_v(f^{\Sieg}_{v,s,\chi}, \pi_{\bfg_{\sfQ},v})^2. 
\end{multline*}
Here, $\ulmu_{v, \sfP}$ is the tuple appearing in the local representation $\sigma_{\bff_{\sfP}, v}$ (and similarly for other notations); $\kappa_v$ are auxiliary characters $\calF_v^{\times} \rightarrow \barQQ_p^{\times}$ that are sufficiently ramified for $v \in \scV_{\calF}^{(p)}$ (which may depend on $\sfP$ and $\sfQ$); and $\scZ^{\dsuit, \rightt}_v(f^{\Sieg}_{v,s,\chi}, \pi_{\bfg_{\sfQ},v})$ are local doubling integrals that can be written explicitly (see Section \ref{sec:doublingintegrals}).
\end{theorem}

\begin{proof}
This follows from a straightforward combination of Theorem \ref{thm:interpolationfinal} and \ref{thm:automorphicfinal}, using the definition of the normalization factor in \eqref{eq:normalizationB}. Note that the assumptions in the automorphic computations (i.e., Assumptions \ref{ass:wtinterlacing}, \ref{ass:spl}, \ref{ass:ram}, \ref{ass:scalarwt}, and \ref{ass:generic}) are satisfied by the definitions of the $p$-adic Eisenstein datum and GGP datum (i.e., Definitions \ref{defn:padiceisendatum} and \ref{defn:ggpdatum}), together with Assumption \ref{ass:automorphic}.
\end{proof}

\subsection{A $p$-adic $L$-function of the Rankin--Selberg product of Hida families}
In \cite{wan2015families}, Wan constructed not only Klingen Eisenstein families but also $p$-adic $L$-functions of Hida families over unitary groups.\footnote{Certain typographical errors in \cite[Theorem 1.1 (1)]{wan2015families} were corrected in \cite[Theorem 6.8]{wan2019iwasawatheorymathrmursblochkato}. We also remark that David Marcil \cite{marcil2024PadicLfunctionsPordinary} recently constructed the $p$-adic $L$-function for $P$-ordinary Hida families over unitary groups, which is more general than Wan's construction.} In this section, under the modification discussed in Remark \ref{rem:modificationsections}, we denote the (modified) $p$-adic $L$-function $\calL_{\Wan}^{\scU}(\bfg, \calK, \chi_0)$ attached to 
\[
\bfD_{G}^{\prime} := (L, \JJ, \eta_0, \bfg, \psi_0^{\prime}, \chi_0, \scU)
\]
simply as $\calL_{p}^{\scU}(\bfg,\chi_0)$. Since we have replaced $\scT$ in $\bfD_{G}$ with $\scU$, this $p$-adic $L$-function is $\scU$-imprimitive rather than $\scT$-imprimitive. Moreover, by the construction of $\calL_{p}^{\scU}(\bfg,\chi_0)$ in \cite{wan2015families}, for any $\sfQ \in \calX_{\bfD_{G}^{\prime}}^{\gen}$, we have
\[
C_{\bfg_\sfQ}^{(p)} = \prod_{v \in \scU} \scZ^{\dsuit, \rightt}_v(f^{\Sieg}_{v,s,\chi}, \pi_{\bfg_{\sfQ},v}),
\]
where $C_{\bfg_{\sfQ}}^{(p)}$ is the constant appearing in the interpolation formula \eqref{eq:wanpadicl}. Comparing \eqref{eq:wanpadicl} with the interpolation formula in Theorem \ref{thm:interpoformula} yields the following result, which is another main result of this article.

\begin{theorem} \label{thm:newpadicl}
With notation as above, assume Assumption \ref{ass:automorphic}. If the $p$-adic $L$-function $\calL_{p}^{\scU}(\bfg,\chi_0)$ is nonzero, then there exists an element 
\[
\calL_{p}(\bff \times \bfg) \in \Frac(\II^{\ur}) \otimes \Frac(\JJ^{\ur})
\]
such that for any $\sfP \times \sfQ \in \calX_{\bfG}^{\gen}$, $(\sfP \times \sfQ)(\calL_{p}(\bff \times \bfg))$ equals
\begin{multline*}
\dfrac{1}{2^{\varkappa_{\sigma_{\bff_{\sfP}}}+\varkappa_{\pi_{\bfg_{\sfQ}}}}} \times \aabs{\bff_{\sfP}}_{\sigma_{\bff_{\sfP}}, \MVW}^2 \aabs{\bfg_{\sfQ}}_{\pi_{\bfg_{\sfQ}}, \MVW}^2 \cdot \prod_{v \in \scS_{\bff} \cup \scS_{\bfg}} \zeta_v C_{\sigma_{\bff_{\sfP,v}, \bpsi_v}} C_{\sigma_{\bff_{\sfP}, v}^{\vee}, \bpsi_v^{-1}} \calB_{\pi_{\bfg_{\sfQ},v}}^{\ess} \calB_{\sigma_{\bff_{\sfP},v}}^{\ess} \\
\times \prod_{v \in \mathscr{V}_{\calF}^{(p)}} \zeta_v \Delta_{v,\beta_v}^2 \cdot \frG(\breve{\ulmu_{v,\sfP}}, \breve{\ullam_{v,\sfQ}}) \frG(\breve{\check{\ulmu_{v,\sfP}}}, \breve{\check{\ullam_{v,\sfQ}}}) \calC_{\sigma_{\bff_{\sfP}}}^{\ord} \calC_{\pi_{\bfg_{\sfQ}}}^{\ord} \cdot \scL^{\scV_{\calF}^{(p)}}(\sigma_{\bff_{\sfP}} \times \pi_{\bfg_{\sfQ}}).
\end{multline*}
Furthermore, when $n = 0$, we have $\calL_{p}(\bff \times \bfg) \in \II^{\ur} \otimes \Frac(\JJ^{\ur})$.
\end{theorem}

\begin{proof}
Recalling \eqref{eq:projEKling} in the construction of $\bfP_{\bfG}$, we identify 
\[
\bfP_{\bfG} \in (\Lambda_{\bfD_{H}} \otimes \Frac(\II^{\ur})) \otimes \Frac(\JJ^{\ur}) = \Frac(\II^{\ur}) \otimes (\Lambda_{\bfD^{\prime}_{G}} \otimes \Frac(\JJ^{\ur})).
\]
We then define \footnote{Note that the ring $\II^{\ur}\lrbracket{\Gamma_{\calK}} \otimes_{\JJ^{\ur}} \Frac(\JJ^{\ur})$ in which the element $\calL_{p}^{\scU}(\bfg,\chi_0)$ resides is an integral domain, since $\JJ^{\ur}$ is an integral domain as set in Footnote \ref{fn:iurdomain}.}
\begin{equation} \label{eq:constructmypadicl}
\calL_{p}(\bff \times \bfg) := \bfP_{\bfG}^2 \cdot (\calL_{p}^{\scU}(\bfg,\chi_0))^{-2} \in \Frac(\II^{\ur}) \otimes \Frac(\Lambda_{\bfD^{\prime}_{G}}),
\end{equation}
which satisfies the interpolation formula stated in the theorem. Since the specializations along the $\ZZ_p$-extension line are trivial (i.e., $\chi_{\sfP}$ does not appear in the interpolation formula), it follows that $\calL_{p}(\bff \times \bfg) \in \Frac(\II^{\ur}) \otimes \Frac(\JJ^{\ur})$. If $n = 0$, then $\bfP_{\bfG} \in \Lambda_{\bfD_{H}} \otimes \JJ^{\ur}$, which implies $\calL_{p}(\bff \times \bfg) \in \II^{\ur} \otimes \Frac(\JJ^{\ur})$ by construction \eqref{eq:constructmypadicl}.
\end{proof}

As a corollary, we obtain the following factorization:
\[
\bfP_{\bfG}^{2} = \calL_{p}(\bff \times \bfg) \cdot (\calL_{p}^{\scU}(\bfg,\chi_0))^2,
\]
which is a $p$-adically interpolated version of Theorem \ref{thm:automorphicfinal}.

The element $\calL_{p}(\bff \times \bfg)$ can be regarded as the $p$-adic $L$-function of the Rankin--Selberg product associated with the Hida family $\bff$ on $H=\rmU(m,n)$ and the Hida family $\bfg$ on $G=\rmU(m+1,n)$. Regarding this $p$-adic $L$-function, we make the following remarks:
\begin{itemize}
    \item Concerning the special $L$-values and local factors in the interpolation formula, their algebraicity can be found in \cite[Section 4.3]{hsieh2023fivevariablepadiclfunctionsu3times}, particularly the discussion in the proof of \cite[Proposition 4.4]{hsieh2023fivevariablepadiclfunctionsu3times}.
    \item We emphasize that it interpolates the special Rankin–Selberg $L$-values with only the Euler factors of the places of $\calF$ above $p$ removed (in Iwasawa theory, such $p$-adic $L$-functions are often referred to as \emph{primitive} $p$-adic $L$-functions). Strictly speaking, this relies on the disjoint ramification hypothesis for $\bff$ and $\bfg$ (Assumption \ref{ass:automorphic}, derived from Assumption \ref{ass:ram}). This assumption allows us to avoid local factors at primes where local Ichino--Ikeda integrals and local Rankin--Selberg $L$-factors cannot currently be compared.
    \item We have not yet checked whether this interpolation formula is consistent with the conjectured form of $p$-adic $L$-functions predicted by Coates and Perrin-Riou in \cite{coates1989on, coates1989on2}. Another disadvantage is that we cannot verify whether $\calL_{p}(\bff \times \bfg)$ is $\II^{\ur} \otimes \JJ^{\ur}$-integral, even in the case $n = 0$, because we divided by $\calL_{p}^{\scU}(\bfg, \chi_0)$ in the construction.
\end{itemize}
We note that there has been similar recent work in various directions, such as \cite{liu2023anticyclotomicpadiclfunctionsrankinselberg, hsieh2023fivevariablepadiclfunctionsu3times, dimitrakopoulou2024anticyclotomicpadiclfunctionscoleman, harris2025padic}.

%% file: AR01-MoreHida.tex
\section{Noncuspidal Hida theory for unitary groups} \label{sec:noncuspidalhida}
In this section, we review the fundamental exact sequence of Eisenstein congruences for unitary groups. There have been successful results for symplectic groups in \cite{liu2020noncuspidal} (specifically, \cite[Theorem 2.3.1]{liu2020noncuspidal}), and these are expected to hold, following the same proof techniques, for all groups arising from PEL-type Shimura varieties. The case of unitary groups\footnote{We remark that full results have been established for $\rmU(m,m)$ in \cite[Chapter 6]{skinner2014iwasawa} and for $\rmU(m,1)$ in \cite[Chapter 4]{hsieh2014eisenstein}.} was addressed in \cite[Section 3.2]{wan2019iwasawatheorymathrmursblochkato}, which we shall recall here.

\subsection{More on the toroidal compactification of unitary Shimura varieties}
The study of the fundamental exact sequence of Eisenstein congruences relies on a solid understanding of the geometric structure of the toroidal compactification of unitary Shimura varieties and the Igusa tower over them.

\subsubsection{Cusp labels}
We consider the noncuspidal Hida theory of a unitary group $G := \GU(m,n)$ attached to a skew-Hermitian space $(V, \phi)$ with the ordered Witt basis
\[
\tWitt_{m,n}: \tty^{1}, \ldots, \tty^{n}; \ttw^{1}, \ldots, \ttw^{m-n}; \ttx^{1}, \ldots, \ttx^{n}.
\]
We let $X^{\vee} := \calO_{\calK} \ttx^{1} \oplus \cdots \oplus \calO_{\calK}\ttx^{n}$, $ Y:= \calO_{\calK} \tty^{1} \oplus \cdots \oplus \calO_{\calK}\tty^{n}$, and $W$ be an $\calO_{\calK}$-maximal lattice such that $W_{p} := W \otimes_{\ZZ} \ZZ_p = (\calO_{\calK}\otimes_{\ZZ} \ZZ_p) \ttw^{1} \oplus \cdots \oplus (\calO_{\calK}\otimes_{\ZZ} \ZZ_p) \ttw^{m-n}$. We define an $\calO_{\calK}$-lattice of $V$ as
\[
\VV := Y \oplus W \oplus X^{\vee}.
\]
Fix a neat open compact subgroup $K := K_{G}$ of $G(\AA_{\calF}^{\infty})$, and let $\Gamma := K_{G} \cap G(\calO_{\calF})$. For simplicity, we assume $\Gamma$ is a principal congruence subgroup $\Gamma(N)$ of some level $N$. 

Let $\frC_{\VV}$ be the set of cotorsion-free isotropic $\calO_{\calK}$-submodules of $V$; the group $G(\calO_{\calF})$ acts naturally on $\frC_{\VV}$. The quotient set $\frC_{\VV}/\Gamma$ is called the set of \emph{cusp labels of level $\Gamma$}.

Let $\calS_{G}^{\tor}(K^p)$ be the toroidal compactification of the unitary Shimura variety $\calS_{G}(K^p)$. It admits a stratification indexed by cusp labels. For any cusp label $\frV \in \frC_{\VV}/\Gamma$, there is a locally closed subscheme $Z_{\frV} \subset \calS_{G}^{\tor}(K^p)$. For $0 \leq r \leq m+n$, denote
\[
Z_{r} := \bigsqcup_{\frV \in \frC_{\VV}/\Gamma, \, \rank \frV > r} Z_{\frV}.
\]
For $0 \leq r \leq m+n$, let $\scI^{r}_{\calS_{G}^{\tor}(K^p)}$ be the sheaf of ideals associated with the closed subscheme $Z_{r}$. In particular, when $r = 0$, the closed subscheme $Z_{0}$ is the boundary of $\calS_{G}^{\tor}(K^p)$, and thus $\scI^{0}_{\calS_{G}^{\tor}(K^p)}$ coincides with the sheaf $\scI_{\calS_{G}^{\tor}(K^p)}$ defined previously.

\subsubsection{Modular forms on unitary groups}
Recall the Igusa tower \eqref{tikz:igusa} introduced earlier. Similar to Section \ref{sec:padicmodularforms}, let $0 \leq r \leq m+n$ and let $\scI_{m,n}^{r, \tor} := \pi_{m,n}^{\ast} \scI_{\calS_m^{\tor}}^{r}$. We define the space of \emph{mod $p^m$ cuspidal automorphic forms vanishing along the strata of rank $> r$} on $G$ of level $n$ as
\[
V^{r}_{G, m,n} := \rmH^{0}(\calI_{m,n}^{\tor}, \scI_{m,n}^{r, \tor}).
\]
We define $p$-adic automorphic forms by passing to the limit.
\begin{definition}
The space of \emph{$p$-adic automorphic forms vanishing along the strata of rank $> r$ with torsion coefficients} is defined as
\[
\calV^{r}_G := \varinjlim_{m} \varinjlim_{n} V_{G, m,n}^r.
\]
The space of \emph{$p$-adic automorphic forms vanishing along the strata of rank $> r$ with integral coefficients} is defined as
\[
V^{r}_G := \varprojlim_{m} \varinjlim_{n} V_{G, m,n}^r.
\]
For any $p$-adic ring $R$ (i.e., $R \simeq \varprojlim_{n} R/p^n R$), the space of \emph{$p$-adic automorphic forms vanishing along the strata of rank $> r$} with coefficient ring $R$ is denoted by $V_G^{r}(R)$ and obtained via base change to $R$.
\end{definition}
When $r = 0$, we recover the notion of $p$-adic cuspidal automorphic forms from Definition \ref{defn:padicautomorphicform}.

\subsection{Main results in noncuspidal Hida theory for unitary groups}
We now state the main theorem of noncuspidal Hida theory for unitary groups, which generalizes \cite{liu2020noncuspidal} as established in \cite{wan2019iwasawatheorymathrmursblochkato}. Since $p$-ordinarity in this article is defined with respect to the Borel subgroup rather than general parabolic subgroups, we restrict our attention to the case $r = 1$, in contrast to \cite[Theorem 2.3.1]{liu2020noncuspidal}.

Furthermore, as noted in \cite[p.\ 14]{wan2019iwasawatheorymathrmursblochkato}, obtaining the Hida control theorem for noncuspidal families requires working with a weight space smaller than $\Lambda_{m,n}^{\circ}$ that consists of parallel weights. Specifically, a character
\[
\ultau = (\tau_1^{+}, \ldots, \tau_{m}^{+}; \tau_{1}^{-}, \ldots, \tau_{n}^{-})
\]
of $T_{m,n}(\calO_{\calF, p})$ is called a \emph{parallel character} (following \cite[Definition 3.3]{wan2019iwasawatheorymathrmursblochkato}) if 
\[
\tau_{1,v}^{-}/\tau_{m,v}^{+} = \tau_{\para} \circ \Norm_{\calF_v/\QQ_p}
\]
for some character $\tau_{\para}: \ZZ_p^{\times} \rightarrow \barQQ_p^{\times}$ independent of $v$, for all places $v$ of $\calF$ above $p$. We denote the subspace of modular forms with parallel nebentypus by the superscript ``$\para$.'' The following theorem is proved in \cite[Section 3.2]{wan2019iwasawatheorymathrmursblochkato}. (While the result in \cite{wan2019iwasawatheorymathrmursblochkato} assumes that $p$ splits completely in $\calK$, the proof remains valid in our setting where $p$ is unramified in $\calF$, subject to careful bookkeeping.) 

Recall that $\Lambda_{m,n} := \calO_{\calK}\lrbracket{T_{m,n}(\ZZ_p)}$ and $\Lambda_{m,n}^{\circ} := \calO_{\calK}\lrbracket{T_{m,n}(1+p\ZZ_p)}$.

\begin{definition-theorem} \label{thm:hidamain}
With notation as above:
\begin{enumerate}[label = \rm (\arabic*)]
    \item (Freeness) An ordinary projector $\ee = \ee^2$ is well-defined on $\calV_{G}^{1}$. The Pontryagin dual of its ordinary part,
    \[
    \bfV_{G}^{1, \para, \ord} := \Hom_{\ZZ_p}(\ee \calV_{G}^{1, \para}, \QQ_p/\ZZ_p),
    \]
    carries a natural $\Lambda_{m,n}$-module structure and is finite free over $\Lambda_{m,n}^{\circ}$.
    \item (Control theorem) For $r = 0, 1$, define
    \[
    \calM_{G}^{r, \para \ord} := \Hom_{\Lambda_{m,n}^{\circ}}(\bfV_{G}^{r, \para, \ord}, \Lambda_{m,n}^{\circ}),
    \]
    whose elements are called \emph{Hida families over $G$ vanishing along the strata of rank $> r$}. Let $\tau \in \Hom_{\cont}(T_{m,n}(\ZZ_p), \barQQ_p^{\times})$ be an arithmetic parallel weight with dominant algebraic part $\ult$ and finite-order part $\uleps$, and let $\frp_{\tau}$ be the corresponding prime ideal of $\Lambda_{m,n}$. Then
    \[
    \calM_{G}^{r, \para, \ord} \otimes_{\Lambda_{m,n}} \Lambda_{m,n}/\frp_{\tau} \xrightarrow{\sim} \varprojlim_{m} \varinjlim_{n} \ee V_{G, m,n}^{r, \para}[\tau],
    \]
    which is equivariant under the action of the unramified Hecke algebra away from $Np$ and the $\UU_p$-operators.
    \item (Fundamental exact sequence of Eisenstein congruences) There is an exact sequence
    \[
    0 \rightarrow \calM_{G}^{0, \para, \ord} \rightarrow \calM_{G}^{1, \para \ord} \xrightarrow{\oplus_{\frV} \Phi_{\frV}} \bigoplus_{\frV \in \frC_{\VV}/\Gamma, \rank \frV = 1} \calM_{\GU(m-1, n-1), \frV}^{0, \ord} \otimes_{\Lambda_{m-1, n-1}} \Lambda_{m,n} \rightarrow 0,
    \]
    where $\calM_{\GU(m-1, n-1), \frV}^{0, \ord}$ is the $\Lambda_{m-1, n-1}$-module of families of $p$-adic ordinary cuspidal automorphic forms over $\GU(m-1, n-1)$ defined over the strata $Y_{\frV}$. This stratum is the image of $Z_{\frV}$ under the natural projection $\pi: \calS_{G}^{\tor}(K^p_{G}) \rightarrow \calS_{G}^{\min}(K^p_G)$ and is isomorphic to the unitary Shimura variety of $\GU(m-1, n-1)$ as a scheme. The map $\Phi_{\frV}$ is the \emph{Siegel operator} obtained by restriction to the stratum indexed by $\frV$.
\end{enumerate}
\end{definition-theorem}

\begin{corollary}[Fundamental exact sequence of Eisenstein congruences] \label{thm:fundaexact}
Let $A$ be a finite torsion-free $\Lambda_{m,n}^{\circ}$-algebra. There is a short exact sequence
\[
0 \rightarrow \calM_{G}^{0, \para, \ord} \otimes_{\Lambda_{m,n}} A \rightarrow \calM_{G}^{1, \para, \ord} \otimes_{\Lambda_{m,n}} A \rightarrow \bigoplus_{\frV \in \frC_{\VV}/\Gamma, \rank \frV = 1} \calM_{\frV}^{0, \ord} \otimes_{\Lambda_{m-1, n-1}} A \rightarrow 0.
\]
\end{corollary}
\begin{proof}
This was established for $A = \Lambda_{m,n}^{\circ}$ in Theorem \ref{thm:hidamain}(3). The result for $A$ follows since the $\Lambda_{m,n}^{\circ}$-modules of modular forms are free, as shown in Theorem \ref{thm:hidamain}(1).
\end{proof}

%% file: AR02-Primitivity.tex
\section{Primitivity of Klingen Eisenstein families over unitary groups} \label{sec:genprimitivity}
In this section, we briefly review the machinery of Eisenstein congruences.

\subsection{The Klingen Eisenstein ideal}
We define the central object in the theory of Eisenstein congruences: the Klingen Eisenstein ideal. This ideal measures the congruences between the cuspidal Hecke eigensystems in $\bfh_{\bfD}$ and the Hecke eigensystem of the Klingen Eisenstein family.

\subsubsection{The Klingen Eisenstein ideal}
We have constructed a Hida family of Klingen Eisenstein series
\[
\bfE_{\bfD}^{\Kling} \in \Meas(\Gamma_{\calK}, V_{G^{\hsuit}}(K_{G^{\hsuit}}^p, \calO_{L}^{\ur})) \otimes \Frac(\II^{\ur})
\]
attached to a $p$-adic Eisenstein datum $\bfD$. When $n=0$, this family lies in $\Meas(\Gamma_{\calK}, V_{G^{\hsuit}}(K_{G^{\hsuit}}^p, \calO_{L}^{\ur})) \otimes \II$.

By the natural pairing between $V_{G^{\hsuit}}^{\ord}$ and the space $\calV_{G^{\hsuit}}^{\ord,\ast}$ from Theorem \ref{thm:hidamain}, there is a natural map
\[
\Meas(\Gamma_{\calK}, V_{G^{\hsuit}}(K_{G^{\hsuit}}^p, \calO_{L}^{\ur})) \times \calV_{G^{\hsuit}}^{\ord,\ast} \rightarrow \Meas(\Gamma_{\calK}, \calO_{L}^{\ur}) \cong \calO_{L}^{\ur}\lrbracket{\Gamma_{\calK}},
\]
which induces an $\calO_{L}^{\ur}\lrbracket{\Gamma_{\calK}}$-linear, Hecke-equivariant map
\[
\Meas(\Gamma_{\calK}, V_{G^{\hsuit}}(K_{G^{\hsuit}}^p, \calO_{L}^{\ur})) \rightarrow \Hom_{\Lambda_{m+1, n+1}}(\calV_{G^{\hsuit}}^{\ord,\ast} ,\calO_{L}^{\ur}\lrbracket{\Gamma_{\calK}}) = \calM_{G^{\hsuit}}^{\ord} \otimes_{\Lambda_{m+1, n+1}} \calO_{L}^{\ur}\lrbracket{\Gamma_{\calK}}.
\]
Here, $\calO_{L}^{\ur}\lrbracket{\Gamma_{\calK}}$ is equipped with a $\Lambda_{m+1, n+1}$-algebra structure via the action of $T_{m+1, n+1}$ on $\bfE_{\bfD}^{\Kling}$. Further tensoring with $\II^{\ur}$, we obtain
\[
\Meas(\Gamma_{\calK}, V_{G^{\hsuit}}(K_{G^{\hsuit}}^p, \calO_{L}^{\ur})) \otimes \II \rightarrow \calM_{G^{\hsuit}}^{\ord} \otimes_{\Lambda_{m+1, n+1}} \II^{\ur}\lrbracket{\Gamma_{\calK}} =: \calM_{G^{\hsuit}}^{\ord}(K_{G^{\hsuit}}^p, \II^{\ur}\lrbracket{\Gamma_{\calK}}).
\]
Consequently, the Klingen Eisenstein family can be viewed as an element 
\begin{equation} \label{eq:newklingen}
\bfE_{\bfD}^{\Kling} \in \calM_{G^{\hsuit}}^{\ord}(K_{G^{\hsuit}}^p, \II^{\ur}\lrbracket{\Gamma_{\calK}}) \otimes \Frac(\II^{\ur}).
\end{equation}
Moreover, $\bfE_{\bfD}^{\Kling} \in \calM_{G^{\hsuit}}^{1, \ord}(K_{G^{\hsuit}}^p, \II^{\ur}\lrbracket{\Gamma_{\calK}}) \otimes \Frac(\II^{\ur})$.

\begin{definition}
Let $\lambda_{\bfD}: \bfh_{\bfD} \rightarrow \Lambda_{\bfD}$ be the eigensystem of the Klingen Eisenstein family $\bfE^{\Kling}_{\bfD}$. We define $\calI_{\bfD}$ to be the ideal of $\bfh_{\bfD}$ generated by $t - \lambda_{\bfD}(t)$ for all $t \in \bfh_{\bfD}$. It is clear that the structure map $\Lambda_{\bfD} \rightarrow \bfh_{\bfD}/\calI_{\bfD}$ is surjective. We define the \emph{Klingen Eisenstein ideal} $\calE_{\bfD}$ attached to $\bfD$ as the kernel of this structure map:
\[
\calE_{\bfD} := \ker \left( \Lambda_{\bfD} \xrightarrow{\text{structure map}} \bfh_{\bfD}/\calI_{\bfD} \right).
\]
\end{definition}

\subsubsection{Constant terms of Klingen Eisenstein families}
Consider the Siegel operator $\Phi_{\frV}$ from Definition-Theorem \ref{thm:hidamain}(3) for cusp labels $\frV \in \frC_{\VV}/\Gamma$ of rank one. Then
\[
\Phi_{\frV}(\bfE^{\Kling}_{\bfD}) \in \calM_{G, \frV}^{0, \ord}(K_{G}^{\prime}, \II^{\ur}\lrbracket{\Gamma_{\calK}}) \otimes \Frac(\II^{\ur}),
\]
where $K_{G}^{\prime}$ is a compact open subgroup of $G(\AA)$. Evaluation at $g^{\prime} \in G(\AA_{\rmf}^{(p)})$ yields
\[
\Phi_{\frV}(\bfE^{\Kling}_{\bfD})(g^{\prime}) \in \II^{\ur}\lrbracket{\Gamma_{\calK}} \otimes \Frac(\II^{\ur}).
\]
Following the approach in \cite[Theorem 6.17]{skinner2014iwasawa}, we define $\calC_{\bfD}$ as \emph{the ideal of constant terms} of $\bfE^{\Kling}_{\bfD}$; specifically, the ideal generated in $\II^{\ur}\lrbracket{\Gamma_{\calK}} \otimes \Frac(\II^{\ur})$ by $\Phi_{\frV}(\bfE^{\Kling}_{\bfD})(g^{\prime})$ as $\frV$ ranges over all cusp labels in $\frC_{\VV}/\Gamma$ of rank one and $g^{\prime}$ ranges over $G(\AA_{\rmf}^{(p)})$.

\begin{theorem} \label{thm:constantterm}
With notation as above, $\calC_{\bfD} \subseteq (\calL_{\bfD})$, where
\[
\calL_{\bfD} := \calL^{\scS}_{\RibetSerre}((\chi_0^{\calF})^{-1} \epsilon_{\calK/\calF}^{m+n}, \calF) \cdot \calL^{\scS}_\Wan(\bff, \calK, \chi_0),
\]
and $\calL^{\scS}_{\RibetSerre}((\chi_0^{\calF})^{-1} \epsilon_{\calK/\calF}^{m+n}, \calF) \in \calO_{L}^{\ur}\lrbracket{\Gamma_{\calK}} \otimes_{\ZZ} \QQ$ is the $\scS$-imprimitive Deligne--Ribet $p$-adic $L$-function of the character $(\chi_0^{\calF})^{-1} \epsilon_{\calK/\calF}^{m+n}$ over the totally real field $\calF$.
\end{theorem}

\begin{proof}
This is proved in \cite[Theorem 1.1(iii)]{wan2015families}, with a correction provided in \cite[Theorem 6.8(iii)]{wan2019iwasawatheorymathrmursblochkato} (noting that since $\epsilon_{\calK/\calF}$ is a quadratic character, $\epsilon_{\calK/\calF}^{m+n} = \epsilon_{\calK/\calF}^{m-n}$).
\end{proof}

\subsection{$p$-primitive functionals}
We write $\calM_{\bfD} := \calM_{G^{\hsuit}}^{\ord}(K_{G^{\hsuit}}^p, \Lambda_{\bfD})$ for simplicity.

\begin{definition}[$p$-primitive functional]
\label{defn:primitivefunctional}
Let $\bfD$ be a $p$-adic Eisenstein datum.
\begin{enumerate}[label = \rm (\arabic*)]
    \item A $\Lambda_{\bfD}$-linear functional $\ell_{\bfD}: \calM_{\bfD} \rightarrow \Lambda_{\bfD}$ is called a \emph{$p$-primitive functional} for the datum $\bfD$ if, for any height-one prime $\frp$ of $\Lambda_{\bfD}$ dividing $\calL_{\bfD}$, we have $\ell_{\bfD}(\bfE^{\Kling}_{\bfD}) \notin \frp$.
    \item A $\Lambda_{\bfD}$-linear functional $\ell_{\bfD}: \calM_{\bfD} \otimes \II \rightarrow \Lambda_{\bfD} \otimes \II$ is called a \emph{weak $p$-primitive functional} for the datum $\bfD$ if, for any height-one prime $\frp$ of $\Lambda_{\bfD}$ such that $\frp \cap \II^{\ur} = (0)$ and $\frp$ divides $\calL_{\bfD}$, we have $\ell_{\bfD}(\bfE^{\Kling}_{\bfD}) \notin \frp$.
    \item The $p$-adic Eisenstein datum $\bfD$ is called \emph{$p$-primitive} (resp. \emph{weak $p$-primitive}) if there exists a $p$-primitive (resp. weak $p$-primitive) functional for $\bfD$.
\end{enumerate}
\end{definition}

The search for a $p$-primitive functional for $\bfD$ is the most technical aspect on the automorphic side of the Iwasawa main conjectures. We first show how the existence of a (weak) $p$-primitive functional leads to the inclusion \ref{item:1a} introduced in Section \ref{sec:eisensteincongruences}, following the method of \cite[Theorem 6.17]{skinner2014iwasawa}.

\begin{proposition} \label{prop:primitivity}
If $\bfD$ is weak $p$-primitive, then 
\[
\ord_{\frp}(\calE_{\bfD}) \geq \ord_{\frp}(\calL_{\bfD})
\]
for any height-one prime $\frp$ of $\Lambda_{\bfD}$ that is not a pullback of a height-one prime of $\II^{\ur}$. If, in addition, $\bfD$ is $p$-primitive, this inequality holds for all height-one prime ideals of $\Lambda_{\bfD}$.
\end{proposition}

\begin{proof}
The first step is to construct a cusp form $\bfH$ that is congruent to the Klingen Eisenstein family modulo $\frp$. It follows from Theorem \ref{thm:constantterm} that for any cusp label $\frV \in \frC_{\VV}/\Gamma$ of rank one, there exists $\bfF_{\frV} \in \calM^{0}_{G^{\hsuit}}(K_{G^{\hsuit}}, \Lambda_{\bfD})$ such that
\[
\Phi_{\frV} \bfE^{\Kling}_{\bfD} = \calL_{\bfD} \bfF_{\frV}.
\]
By the surjectivity of the Siegel operator $\Phi$ in the fundamental exact sequence of Eisenstein congruences (Corollary \ref{thm:fundaexact}), there exists $\bfF \in \calM^{1, \ord}_{G^{\hsuit}}(K_{G^{\hsuit}}^p, \Lambda_{\bfD})$ such that $\Phi_{\frV} \bfF = \bfF_{\frV}$ for every cusp label $\frV$ of rank one. Define
\[
\bfH := \bfE^{\Kling}_{\bfD} - \calL_{\bfD} \bfF \in \calM^{1, \ord}_{G^{\hsuit}}(K_{G^{\hsuit}}^{\prime}, \Lambda_{\bfD}) \otimes \Frac(\II^{\ur}).
\]
Direct verification shows that $\bfH \in \ker(\oplus \Phi_{\frV})$; thus, $\bfH \in \calM^{0, \ord}_{G^{\hsuit}}(K_{G^{\hsuit}}^{\prime}, \Lambda_{\bfD}) \otimes \Frac(\II^{\ur})$ by Corollary \ref{thm:fundaexact}.

Assume $\frp \mid \calL_{\bfD}$ (otherwise the result is trivial) and let $r := \ord_{\frp}(\calL_{\bfD})$. Using the primitive functional $\ell_{\bfD}$, we define a map
\[
\mu_{\bfD}: \bfh_{\bfD} \rightarrow \Lambda_{\bfD, \frp}/\frp^{r} \Lambda_{\bfD, \frp}, \quad h \mapsto \frac{\ell_{\bfD}(h \cdot \bfH)}{\ell_{\bfD}(\bfH)}.
\]
By the primitivity of $\ell_{\bfD}$ and the definition of $\bfH$,
$$
\ell_{\bfD}(\bfH) = \ell_{\bfD}(\bfE^{\Kling}_{\bfD} - \calL_{\bfD} \bfF) = \ell_{\bfD}(\bfE^{\Kling}_{\bfD}) - \ell_{\bfD}(\calL_{\bfD} \bfF) \equiv \ell_{\bfD}(\bfE^{\Kling}_{\bfD}) \not\equiv 0 \pmod{\frp},
$$
thus $\ell_{\bfD}(\bfH)$ is a unit in the localized ring $\Lambda_{\bfD, \frp}$, and hence a unit in the quotient ring $\Lambda_{\bfD, \frp}/\frp^{r} \Lambda_{\bfD, \frp}$, which ensures that the map $\mu_{\bfD}$ is well-defined. It is clear that $\mu_{\bfD}$ is surjective. Furthermore, 
\[
\ell_{\bfD}(h \cdot \bfH) \equiv \ell_{\bfD}(h \cdot \bfE^{\Kling}_{\bfD}) = \lambda_{\bfD}(h) \ell_{\bfD}(\bfE^{\Kling}_{\bfD}) \pmod{\frp^{r}},
\]
which implies
\[
\frac{\ell_{\bfD}(h \cdot \bfH)}{\ell_{\bfD}(\bfH)} \equiv \lambda_{\bfD}(h) \pmod{\frp^{r}}
\]
for any $h \in \bfh_{\bfD}$. It follows that $\calI_{\bfD} \subseteq \ker \mu_{\bfD}$. Thus, $\mu_{\bfD}$ induces a surjective $\Lambda_{\bfD, \frp}$-homomorphism
\[
\Lambda_{\bfD, \frp}/\calE_{\bfD} \Lambda_{\bfD, \frp} \xrightarrow{\sim} \bfh_{\bfD, \frp}/\calI_{\bfD} \bfh_{\bfD, \frp} \twoheadrightarrow \Lambda_{\bfD, \frp}/\frp^{r} \Lambda_{\bfD, \frp}.
\]
Since $\ord_{\frp}(\calE_{\bfD})$ is the length of $\Lambda_{\bfD, \frp}/\calE_{\bfD} \Lambda_{\bfD, \frp}$ over $\Lambda_{\bfD, \frp}$, the proposition is proved.
\end{proof}

%% file: AR03-Lattice.tex
\section{Lattice construction} \label{sec:lattice}
We introduce the lattice construction following the framework of \cite[Chapters 4, 7]{skinner2014iwasawa} and \cite[Section 9.2]{wan2020iwasawa}. For ease of reference, we summarize the setups and results appearing \textit{loc.\ cit.} We then prove several results regarding the Greenberg-type Iwasawa main theorem, assuming the (weak) $p$-primitivity of $\bfD$.

\subsection{Lattice construction}

\subsubsection{Axiomatic lattice construction}
Let $G$ be a group and $C$ be a ring. Let $r: G \rightarrow \Aut_{C}(V)$ be a representation of $G$ with $V \simeq C^{N}$. This extends to $r: C[G] \rightarrow \End_{C}(V)$. For any $x \in C[G]$, we define the \emph{characteristic polynomial}
\[
\Char(r, x, T) := \det(\id{} - r(x) T) \in C[T].
\]

Let $(r_1, V_1)$ and $(r_2, V_2)$ be two $C$-representations of $G$. Assume that both are defined over a local Henselian subring $B \subseteq C$. We say that $r_1$ and $r_2$ are \emph{residually disjoint} modulo the maximal ideal $\frm_{B}$ of $B$ if there exists $x \in B[G]$ such that $\Char(r_1, x, T) \pmod{\frm_{B}}$ and $\Char(r_2, x, T) \pmod{\frm_{B}}$ are coprime in $\kappa_{B}[T]$, where $\kappa_{B} := B/\frm_{B}$.

Let $H$ be a group with a semidirect product decomposition $H = G \rtimes \{1, \sfc\}$, where $\sfc \in H$ is an element of order two normalizing $G$. For any $C$-representation $(r, V)$ of $G$, we write $(r^{\sfc}, V)$ for the representation defined by $r^{\sfc}(g) := r(\sfc g \sfc)$ for any $g \in G$.

Let $\theta: G \rightarrow \GL_{L}(V)$ be a $G$-representation on an $L$-vector space, where $L$ is a field, and let $\psi: H \rightarrow L^{\times}$ be a character. We say that $\theta$ satisfies the \emph{$\psi$-polarization condition} if $\theta^{\sfc} \simeq \psi \otimes \theta^{\vee}$. By a \emph{$\psi$-polarization} of $\theta$, we mean an $L$-bilinear pairing $\Phi_{\theta}: V \times V \rightarrow L$ such that 
\[
\Phi_{\theta}(\theta(g)v, v^{\prime}) = \psi(g) \Phi_{\theta}(v, \theta^{\sfc}(g)^{-1} v^{\prime}).
\]
Let $\Phi_{\theta}^{\rmt}(v, v^{\prime}) := \Phi_{\theta}(v^{\prime}, v)$ be another $\psi$-polarization of $\theta$. We say that $\psi$ is \emph{compatible with the polarization} $\Phi_{\theta}$ (or simply with $\theta$ if the polarization is clear from the context) if $\Phi_{\theta}^{\rmt} = -\psi(\sfc)\Phi_{\theta}$.

\begin{definition}
An \emph{algebraic datum for the lattice construction} is a tuple 
\[
\pzL = (G, H, D; A_0, \frp; R_0, \frq; I_0, J_0),
\]
where:
\begin{itemize}
    \item $G$ and $H$ are profinite groups, and $D = \{D_1, \ldots, D_{g}\}$ is a set of subgroups of $G$.
    \item $A_0$ is a profinite $\ZZ_p$-algebra and a Krull domain.
    \item $\frp \subset A_0$ is a height-one prime ideal, and $A := \hatA_{0, \frp}$ is the completion of the localization of $A_{0}$ at $\frp$.
    \item $R_0$ is a local reduced finite $A_{0}$-algebra.
    \item $\frq \subset R_{0}$ is a prime ideal such that $\frq \cap A_{0} = \frp$, and $R := \hatR_{0, \frq}$ is the completion of the localization of $R_{0}$ at $\frq$.
    \item $J_{0} \subset A_{0}$ and $I_{0} \subset R_{0}$ are ideals such that $I_{0} \cap A_0 = J_0$, $A_0/J_0 = R_0/I_0$, and $J := J_0 A$, $I := I_0 A$ with $J_0 = J \cap A_0$ and $I_0 = I \cap R_0$.
\end{itemize}
A \emph{representation datum for the lattice construction} attached to $\pzL$ is a tuple $\pzR = (\nu, \chi, \rho, \sigma)$, where:
\begin{itemize}
    \item $\nu: H \rightarrow A_0^{\times}$ is a continuous character.
    \item $\chi: G \rightarrow A_{0}^{\times}$ is a continuous character such that $\barchi \neq \barnu \barchi^{-\sfc}$. Let $\chi^{\prime} := \nu \chi^{-\sfc}$.
    \item $\rho: G \rightarrow \Aut_{A}(V)$ with $V = A^{N}$ is a representation obtained by base change from a representation $\rho_0: G \rightarrow \Aut_{A_0}(V_0)$ such that:
    \begin{enumerate}
        \item $\rho^{\sfc} \simeq \rho^{\vee} \otimes \nu$;
        \item $\barrho$ is absolutely irreducible (hence $\rho$ is irreducible);
        \item $\rho$ is residually disjoint from $\chi$ and $\chi^{\prime}$.
    \end{enumerate}
    \item $\sigma: G \rightarrow \Aut_{R \otimes_{A} \Frac(A)}(M)$ with $M \simeq (R \otimes_{A} \Frac(A))^{N+2}$ is a representation defined in the image of $R_{0}$ in $R$ such that:
    \begin{enumerate}
        \item $\sigma^{\sfc} \simeq \sigma^{\vee} \otimes \nu$;
        \item $\tr(\sigma(g)) \in R$ for all $g \in G$;
        \item for any $v \in M$, $\sigma(R[G])v$ is a finitely generated $R$-module;
        \item $\tr \sigma(g) \equiv \chi^{\prime}(g) + \tr \rho(g) + \chi(g) \pmod{I}$ for all $g \in G$.
    \end{enumerate}
\end{itemize}
The datum $\pzR$ must satisfy the following conditions:
\begin{enumerate}[label = \rm (LC \arabic*)]
    \item $\nu$ is compatible with $\rho$. \label{item:LC1}
    \item (Local conditions for $\rho$) For each $i = 1, \ldots, g$, there exists a $D_i$-stable $A_0$-submodule $V_{0, i}^{+} \subseteq V_{0}$ such that $V_{0, i}^{+}$ and $V_{0, i}^{-} := V_{0}/V_{0, i}^{+}$ are free $A_0$-modules. \label{item:LC2}
    \item (Local conditions for $\sigma$) For each $i = 1, \ldots, g$, there exists a $D_i$-stable $R \otimes_{A} \Frac(A)$-submodule $M_{i}^{+} \subseteq M$ such that $M_{i}^{+}$ and $M_{i}^{-} := M/M_{i}^{+}$ are free $(R \otimes_{A} \Frac(A))$-modules. Moreover, $M_{i}^{+}$ and $M_{i}^{-}$ are disjoint modulo $I$. \label{item:LC3}
    \item (Compatibility with congruence condition) For each $i = 1, \ldots, g$, assume that for all $x \in R[D_i]$, the characteristic polynomial satisfies
    \[
    \Char(M_{i}^{+}, x, T) \equiv \Char(V_{i}^{+}, x, T)(1-T \chi(x)) \pmod{I}.
    \]
    It then follows that $\Char(M_{i}^{-}, x, T) \equiv \Char(V_{i}^{-}, x, T)(1-T \chi^{\prime}(x)) \pmod{I}$. \label{item:LC4}
    \item For each $\Frac(A)$-algebra homomorphism $\lambda: R \otimes_{A} \Frac(A) \rightarrow K$, where $K$ is a finite field extension of $\Frac(A)$, the representation $\sigma_{\lambda}: G \rightarrow \GL_{N+2}(M \otimes_{R \otimes_{A} \Frac(A), \lambda} K)$ obtained from $\sigma$ via $\lambda$ is either absolutely irreducible or contains an absolutely irreducible two-dimensional $K$-subrepresentation $\sigma_{\lambda}^{\prime}$ such that $\tr \sigma_{\lambda}^{\prime}(g) \equiv \chi(g)+\chi^{\prime}(g) \pmod{I}$ for any $g \in G$. \label{item:LC5}
\end{enumerate}
A \emph{datum for the lattice construction} is the tuple $(\pzL, \pzR)$.
\end{definition}

Given a datum for the lattice construction $(\pzL, \pzR)$, following \cite[Section 4.4.2]{skinner2014iwasawa}, we consider a character $\tilde{\chi}$ as the unique extension of the character $\chi^{-\sfc-1}$ to $H$ such that $\tilde{\chi}(\sfc) = 1$. Then $\tilde{\chi}\nu$ is a character of $H$ taking values in $A_{0}^{\times}$. We define the following Selmer groups:
\[
\Sel_{H}(\tilde{\chi}\nu) := \ker\left(\rmH^{1}(H, A_0^{\ast}(\tilde{\chi}\nu)) \rightarrow \prod_{i=1}^{g} \rmH^{1}(D_i, A_0^{\ast}(\tilde{\chi}\nu))\right)
\]
and
\[
\Sel_{G}(\rho_0 \otimes \chi^{-1}) := \ker\left(\rmH^{1}(G, V_0 \otimes_{A_0} A_0^{\ast}(\chi^{-1})) \rightarrow \prod_{i=1}^{g} \rmH^{1}(D_i, V_{0, i}^{-} \otimes_{A_0} A_0^{\ast}(\chi^{-1}))\right).
\]
Let $\Char_{H}(\tilde{\chi}\nu)$ and $\Char_{G}(\rho_0 \otimes \chi^{-1})$ be the characteristic ideals of the Pontryagin duals of $\Sel_{H}(\tilde{\chi}\nu)$ and $\Sel_{G}(\rho_0 \otimes \chi^{-1})$, respectively, viewed as $A_0$-modules.

We record the following theorems regarding the axiomatic lattice construction.

\begin{theorem}[{\cite[Corollary 4.16]{skinner2014iwasawa}}, cf.\ also {\cite[Proposition 9.2]{wan2020iwasawa}}] \label{thm:latticeconstruction}
Given a datum for the lattice construction $(\pzL, \pzR)$, if $\ord_{\frp}(\Char_{H}(\tilde{\chi}\nu)) = 0$, then
\[
\ord_{\frp}(\Char_{G}(\rho_0 \otimes \chi^{-1})) \geq \ord_{\frp}(J).
\]
\end{theorem}

\begin{theorem}[{\cite[Proposition 4.17]{skinner2014iwasawa}}] \label{thm:toblockkato}
Given a datum for the lattice construction $(\pzL, \pzR)$, if $\ord_{\frp}(J) > \ord_{\frp}(\Char_{H}(\tilde{\chi}\nu))$, then
\[
\ord_{\frp}(\Char_{G}(\rho_0 \otimes \chi^{-1})) \geq 1.
\]
\end{theorem}

\subsubsection{Application to Klingen Eisenstein ideals}
Assume that \ref{item:irredbff} and \ref{item:distbff} hold. Given a $p$-adic Eisenstein datum $\bfD$, we can define a semisimple representation\footnote{It might be simpler to twist $\rho_{\bfD}$ by $\sigma_{p}(\bpsi_0)^{-\sfc} \epsilon_{\cyc}^{m+n}$, which yields
\[
 \rho_{\bfD} \otimes \sigma_{p}(\bpsi_0)^{-\sfc} \epsilon_{\cyc}^{m+n} = \epsilon_{\cyc}^{-1} \oplus R_{p}(\bff) \otimes \sigma_{p}(\bxi_0)^{-\sfc} \oplus \sigma_{p}(\bxi_0^{\calF} \circ \Nm{}{})^{-1},
\]
which is more concise.}
\begin{multline} \label{eq:galoisklingen}
 \rho_{\bfD} := \sigma_{p}(\bpsi_0)^{\sfc} \cdot \epsilon_{\cyc}^{-(m+n+1)} \oplus R_{p}(\bff) \otimes \sigma_{p}(\bxi_0)^{-\sfc} \sigma_{p}(\bpsi_0)^{\sfc} \epsilon_{\cyc}^{-(m+n)} \\ \oplus \epsilon_{\cyc}^{-(m+n-1)} \det R_{p}(\bff) \sigma_{p}(\bxi_0^{\calF} \circ \Nm{}{})^{-1} \sigma_{p}(\bpsi_0)^{\sfc},   
\end{multline}
where $\sigma_p(-)$ denotes the reciprocity map from class field theory. This is the Galois representation associated with the Klingen Eisenstein family $\bfE_{\bfD}^{\Kling}$, following from a direct generalization of \cite[Proposition 9.14]{skinner2014iwasawa}.

We consider the following algebraic datum $\pzL_{\bfD}$ for the lattice construction:
\begin{itemize}
    \item $G := \Gal_{\calK, \Sigma}$ and $H := \Gal_{\calF, \Sigma}$, where $\sfc \in \Gal_{\calF}$ is the complex conjugation. 
    \item Let $v_1, \ldots, v_g$ be the $p$-adic CM type chosen at the outset, and let $D_i$ be the decomposition group at $v_{i}$.
    \item $A_{0} := \Lambda_{\bfD}$ and $A := \widehat{\Lambda}_{\bfD, \frp}$, where $\frp$ is a height-one prime ideal of $\Lambda_{\bfD}$.
    \item $R_{0} := \bfh_{\bfD}$ and $\frq$ is defined as the inverse image of $\frp \pmod{\calE_{\bfD}}$ under the map $\bfh_{\bfD} \twoheadrightarrow \bfh_{\bfD}/\calI_{\bfD} = \Lambda_{\bfD}/\calE_{\bfD}$.
    \item $J_{0} := \calE_{\bfD}$ and $I_{0} := \calI_{\bfD}$.
\end{itemize}

Guided by \eqref{eq:galoisklingen}, we consider the following representation datum $\pzR_{\bfD}$\footnote{Here we follow the normalization of \cite[Section 7.3.1]{skinner2014iwasawa} and \cite[Section 9.2.2]{wan2015hilbert} instead of \cite{wan2020iwasawa}. We note, however, that there appears to be a typographical error in the computation of $\nu$ in \cite[Section 7.3.1]{skinner2014iwasawa} and \cite[Section 9.2.2]{wan2015hilbert}. In the latter, the power of $\epsilon$ should be $-5$ to match the case $m = n = 1$, since $\det R_{p}(f) = \epsilon_{\cyc}^{-1}$.} for the lattice construction attached to $\pzL_{\bfD}$:
\begin{itemize}
    \item $\chi := \sigma_{p}(\bpsi_0)^{\sfc} \epsilon_{\cyc}^{-(m+n+1)}$ and $\chi^{\prime} := \epsilon_{\cyc}^{-(m+n-1)} \det R_{p}(\bff) \sigma_{p}(\bxi_0^{\calF} \circ \Nm{}{})^{-1} \sigma_{p}(\bpsi_0)^{\sfc}$.
    \item $\nu := \chi^{\sfc} \chi^{\prime} = \sigma_{p}(\bpsi_0^{\calF} \circ \Nm{}{}) \cdot \sigma_{p}(\bxi_0^{\calF} \circ \Nm{}{})^{-1} \cdot \epsilon_{\cyc}$.
    \item $\rho_0 := R_{p}(\bff) \otimes \sigma_{p}(\bxi_0)^{-\sfc} \sigma_{p}(\bpsi_0)^{\sfc} \epsilon_{\cyc}^{-(m+n)}$.
    \item $\sigma$ is the representation on $M$ obtained from the pseudorepresentation associated with $\bfh_{\bfD}$ as in \cite[Proposition 7.3]{skinner2014iwasawa}.
\end{itemize}
We take $V_{0, i}^{+}$ to be $T_{\bff, v_i}^{+} \otimes_{\II} A_0$, with $V_{0, i}^{-} := T_{\bff, v_i}^{-} \otimes_{\II} A_0$. Properties \ref{item:LC1}--\ref{item:LC4} follow from the descriptions in Section \ref{sec:Galoisreps}, and their verification is identical to the first paragraph of the proof of \cite[Theorem 7.6]{skinner2014iwasawa}. Property \ref{item:LC5} is proven in \cite[Proposition 7.4]{wan2019iwasawatheorymathrmursblochkato}: unlike the verification of \ref{item:LC1}--\ref{item:LC4}, this is no longer a direct generalization of the parallel conclusion in \cite{skinner2014iwasawa} (specifically, \cite[Proposition 7.5]{skinner2014iwasawa}), as we do not have modularity lifting results for general unitary groups. Nevertheless, $(\pzL_{\bfD}, \pzR_{\bfD})$ constitutes a valid datum for the lattice construction.

Applying Theorem \ref{thm:latticeconstruction}, we obtain the following proposition.
\begin{proposition} \label{prop:applyprop}
Suppose $\frp$ is a height-one prime ideal of $\II^{\ur}\lrbracket{\Gamma_{\calK}}$ such that $\frp \cap \II^{\ur}\lrbracket{\Gamma_{\calK}^{+}} = (0)$. Then
\begin{equation}\label{eq:latticefinalcyc}
\ord_{\frp}(X_{\Gr}^{\scS}(\bff/\calK_{\infty}, \xi_0^{-\sfc}\epsilon_{\cyc})) \geq \ord_{\frp}(\calE_{\bfD}).
\end{equation}
\end{proposition}
\begin{proof}
Direct calculation shows that
\[
\tilde{\chi}\nu = \sigma_{p}(\bxi_0^{\calF})^{-1} \epsilon_{\calK/\calF}^{m+n} \epsilon_{\cyc} \quad \text{and} \quad \rho_0 \otimes \chi^{-1} = R_{p}(\bff) \otimes \sigma_{p}(\bxi_0)^{-\sfc} \epsilon_{\cyc}.
\]
Recall that the Selmer groups in Theorem \ref{thm:latticeconstruction} are defined as $\Lambda_{\bfD}$-modules; thus, we must compare the various Selmer groups. One checks (see \cite[pp.\ 127--128]{skinner2014iwasawa}) that
\begin{equation} \label{eq:SelH}
\Sel_{H}(\tilde{\chi}\nu) = \Sel^{\scS}((\xi_0^{\calF})^{-1} \epsilon_{\calK/\calF}^{m+n} \epsilon_{\cyc}, \calF_{\infty}) \otimes_{\calO_{L}^{\ur}\lrbracket{\Gamma_{\calF}}} \Lambda_{\bfD},
\end{equation}
where $\Sel^{\scS}((\xi_0^{\calF})^{-1}\epsilon_{\calK/\calF}^{m+n} \epsilon_{\cyc}, \calF_{\infty}) := \Sel_{\Gr}^{\scS}(\calO_{L}((\xi_0^{\calF})^{-1}\epsilon_{\cyc}), \{0\}_{v \in \scV_{\calF}^{(p)}}, \calF, \calO_{L}^{\ur}\lrbracket{\Gamma_{\calF}})$ as in Definition \ref{defn:selmergroups}. Consequently, under the natural identification of $\Gamma_{\calK}^{+}$ with $\Gamma_{\calF}$, we see that $\Sel_{H}(\tilde{\chi}\nu)$ is only divisible by height-one prime ideals $\frp$ that are pullbacks of primes in $\II^{\ur}\lrbracket{\Gamma_{\calK}^{+}}$. Finally, note that $\Sel_{G}(\rho_0 \otimes \chi^{-1})$ is precisely the Selmer group $\Sel_{\Gr}^{\scS}(\bff/\calK_{\infty}, \xi_0^{-\sfc} \epsilon_{\cyc})$.
\end{proof}

\subsection{Application to Greenberg-type Iwasawa main conjectures}
Let 
\[
\II^{\ur}\lrbracket{\Gamma_{\calK}}^{(\cyc)} := \II^{\ur}\lrbracket{\Gamma_{\calK}} \otimes_{\II^{\ur}\lrbracket{\Gamma_{\calK}^{+}}} \Frac(\II^{\ur}\lrbracket{\Gamma_{\calK}^{+}})
\]
be the weight algebra of $\bfD$ with the cyclotomic direction deprived. We have the following partial result for the Iwasawa main conjecture.

\begin{theorem} \label{thm:imcdeprived}
Suppose $\bfD$ is weak $p$-primitive. Then:
\begin{enumerate}[label = \rm (\arabic*)]
    \item As ideals of $\II^{\ur}\lrbracket{\Gamma_{\calK}}^{(\cyc)}$, we have
\[
\Char_{\II^{\ur}\lrbracket{\Gamma_{\calK}}^{(\cyc)}}\left(X_{\Gr}^{\scS}(\bff/\calK_{\infty}, \xi_0^{-\sfc} \epsilon_{\cyc}) \otimes_{\II\lrbracket{\Gamma_{\calK}}} \II^{\ur}\lrbracket{\Gamma_{\calK}}^{(\cyc)}\right) \subseteq (\calL_{\Wan}^{\scS}(\bff, \calK, \chi_0)).
\]
    \item Furthermore, assume that $\calL_{\Wan}^{\scS}(\bff, \calK, \chi_0)$ is coprime to $\calL^{\scS}_{\RibetSerre}((\chi_0^{\calF})^{-1} \epsilon_{\calK/\calF}^{m+n}, \calF)$ in $\II^{\ur}\lrbracket{\Gamma_{\calK}^{+}}$.\footnote{This is satisfied, in particular, if $\calL_{\Wan}^{\scS}(\bff, \calK, \chi_0)$ is not divisible by any height-one ideals in $\II^{\ur}\lrbracket{\Gamma_{\calK}^{+}}$.} Then
    \begin{equation} \label{eq:finalimc}
    \Char_{\II^{\ur}\lrbracket{\Gamma_{\calK}}}\left(X_{\Gr}^{\scS}(\bff/\calK_{\infty}, \xi_0^{-\sfc} \epsilon_{\cyc})\right) \subseteq (\calL_{\Wan}^{\scS}(\bff, \calK, \chi_0))        
    \end{equation}
    as fractional ideals in $\II^{\ur}\lrbracket{\Gamma_{\calK}} \otimes \Frac(\II^{\ur})$. If $n = 0$, then \eqref{eq:finalimc} holds as ideals in $\II^{\ur}\lrbracket{\Gamma_{\calK}} \otimes \QQ$. If, in addition, $\bfD$ is $p$-primitive, \eqref{eq:finalimc} holds as ideals in $\II^{\ur}\lrbracket{\Gamma_{\calK}}$.
\end{enumerate}
\end{theorem}

\begin{proof}
Item (1) follows directly from Propositions \ref{prop:primitivity} and \ref{prop:applyprop}. To prove (2), it suffices to consider height-one primes $\frp$ of $\Lambda_{\bfD}$ such that $\frp \cap \II^{\ur}\lrbracket{\Gamma_{\calK}^{+}} \neq (0)$ and show that
\[
\ord_{\frp}(\calL_{\Wan}^{\scS}(\bff, \calK, \chi_0)) \leq \ord_{\frp}\Char_{\II^{\ur}\lrbracket{\Gamma_{\calK}}}\left(X_{\Gr}^{\scS}(\bff/\calK_{\infty}, \xi_0^{-\sfc} \epsilon_{\cyc})\right).
\]
Suppose $\frp$ divides $\calL_{\Wan}^{\scS}(\bff, \calK, \chi_0)$ (otherwise the claim is trivial). Then $\frp$ does not divide $\calL^{\scS}_{\RibetSerre}((\chi_0^{\calF})^{-1} \epsilon_{\calK/\calF}^{m+n}, \calF)$ by the coprimality assumption. The Iwasawa main conjecture for the character $\xi_{0}^{-1}$ over totally real fields was proved by Wiles \cite{wiles1990iwasawa}, which implies that $\frp$ does not divide the characteristic ideal of the Pontryagin dual of \eqref{eq:SelH}. The result follows from Theorem \ref{thm:latticeconstruction} by the same argument used in the proof of Proposition \ref{prop:applyprop}. For $n = 0$, everything is integral in $\II^{\ur}\lrbracket{\Gamma_{\calK}}$ by construction.
\end{proof}

%% file: TH1-Theta.tex
\section{Interlude: Preliminaries on theta theory} \label{sec:pretheta}

\subsection{Non-archimedean Langlands correspondences of unitary groups}
In this section, we briefly recall the basic notions and results concerning the non-archimedean Langlands correspondences for unitary groups. These have been established by the work of \cite{mok2015endoscopic, kaletha2014endoscopicclassificationrepresentationsinner}, and many others.

Let $F$ be a non-archimedean local field of characteristic zero, residue characteristic $p$, and residue cardinality $q_{F}$. Let $K$ be a quadratic field extension of $F$, and let $\epsilon_{K/F}$ be the quadratic character of $F^{\times}$ associated with $K/F$ via local class field theory.

\subsubsection{$\varepsilon$-Hermitian spaces and Witt towers} \label{sec:hermitianspace}
We fix $\varepsilon = \pm 1$. Let $V$ be a finite-dimensional vector space over $K$ equipped with a nondegenerate $\varepsilon$-Hermitian form
\[
\lrangle{-,-}_{V}: V \times V \rightarrow K.
\]
Denoting $n = \dim_{K} V$, we define
\[
\disc V := (-1)^{n(n-1)/2} \det V \in \begin{cases}
F^{\times}/\Norm_{K/F}(K^{\times}), &\quad \text{if } \varepsilon = +1, \\
\delta^{n} F^{\times}/\Norm_{K/F}(K^{\times}), &\quad \text{if } \varepsilon = -1,
\end{cases}
\]
and we define
\[
\epsilon(V) := \begin{cases}
\epsilon_{K/F}(\disc V), &\quad \text{if } \varepsilon = +1, \\
\epsilon_{K/F}(\delta^{-n} \disc(V)), &\quad \text{if } \varepsilon = -1.
\end{cases}
\]

Given a positive integer $n$, there are exactly two isometry classes of $n$-dimensional $\varepsilon$-Hermitian spaces $V$, distinguished by their signs $\epsilon(V)$. We denote an $\varepsilon$-Hermitian space by $V^{\epsilon}$ to emphasize its sign $\epsilon(V)$. Let $\rmU(V)$ be the unitary group of $V$, which is the connected reductive linear algebraic group over $F$ defined by
\[
\rmU(V)(R) := \{g \in \Aut_{K \otimes_{F} R}(V \otimes_{F} R): \lrangle{gv, gw}_{V} = \lrangle{v,w}_{V}, \text{ for all } v, w \in V \otimes_{F} R \}
\]
for any $F$-algebra $R$. When $\dim_{K} V = 0$, we regard $\rmU(V)$ as the trivial group $\{1\}$.

Let $V_{(0)}$ be an anisotropic space over $K$. For integers $r \geq 0$, let $V_{(r)} := V_{(0)} \oplus \scH^{r}$, where $\scH$ is the $\varepsilon$-Hermitian hyperbolic plane. The collection
\[
\calV := \{ V_{(r)} = V_{(0)} \oplus \scH^{r}: r \geq 0 \}
\]
is called a \emph{Witt tower} starting from $V_{(0)}$. We observe that the parity of $\dim_{K} V_{(r)}$ is independent of $r$, and $\disc(V_{r})$—or equivalently the sign $\epsilon(V_{(r)})$—is independent of $r$. Moreover, given any Witt tower $\calV$, there exists a companion Witt tower $\calV^{\prime} = (V_{(r)}^{\prime})$ with the opposite sign. We shall denote the two Witt towers by $\calV^{+} = (V_{(r)}^{+})$ and $\calV^{-} = (V_{(r)}^{-})$ accordingly. We note that any given $\varepsilon$-Hermitian space $V$ is a member of a unique Witt tower $\calV$.

\subsubsection{Representations of Weil--Deligne groups and local $L$-factors}
We denote by $\rmW_K$ and $\WD_{K} := \rmW_K \times \SL_2(\CC)$ the Weil group and the Weil--Deligne group of $K$, respectively. Let $I_{K}$ be the inertia group of $\rmW_{K}$. We fix a geometric Frobenius element $\Frob_{K}$ of $\rmW_{K}$. We regard $\rmW_{K}$ as a subgroup of $\rmW_F$ such that $\rmW_{F}/\rmW_{K} \cong \Gal(K/F)$, and we fix $\sfc \in \rmW_{F} \smallsetminus \rmW_{K}$.

Let $M$ be a finite-dimensional vector space over $\CC$. A homomorphism $\phi: \WD_{K} \rightarrow \GL(M)$ is a \emph{representation} of $\WD_{K}$ if:
\begin{enumerate}[label = \rm (\roman*)]
    \item $\phi(\Frob_{K})$ is semisimple;
    \item the restriction of $\phi$ to $\rmW_{K}$ is smooth;
    \item the restriction of $\phi$ to $\SL_2(\CC)$ is algebraic.
\end{enumerate}
Furthermore, we say a representation $\phi$ of $\WD_{K}$ is:
\begin{itemize}
    \item \emph{unramified} if it is trivial on $I_{K} \times \SL_2(\CC)$;
    \item \emph{tempered} if the image of $\rmW_{K}$ is bounded.
\end{itemize}

Let $\phi$ be a representation of $\WD_{K}$. The contragredient $\phi^{\vee}$ of $\phi$ is defined by $\phi^{\vee}(w)=\phi(w)^{-\rmt}$. We define a representation $\phi^{\sfc}$ of $\WD_{K}$ by $\phi^{\sfc}(w) = \phi(\sfc w \sfc^{-1})$. The equivalence class of $\phi^{\sfc}$ is independent of the choice of $\sfc$.

We fix $b \in \{\pm 1\}$. We say $\phi$ is \emph{conjugate self-dual} of sign $b$ if there exists a nondegenerate bilinear form $B: M \times M \rightarrow \CC$ such that
\[
B(\phi(w)x, \phi(\sfc w \sfc^{-1})y) = B(x,y), \quad \text{ and } B(y,x) = b \cdot B(x, \phi(\sfc^{2})y)
\]
for $x,y \in M$ and $w \in \WD_{K}$.

For each positive integer $k$, there exists a unique irreducible algebraic representation $\sfS_{k}$ of $\SL_2(\CC)$ of dimension $k$. A representation $\phi$ of $\WD_{K}$ can then be decomposed as
\begin{equation} \label{eq:decomposephi}
 \phi = \bigoplus_{n=1}^{\infty} \phi_n \boxtimes \sfS_{n},  
\end{equation}
where $(\phi_n, M_n)$ is a representation of $\rmW_{K}$ that is zero for all but finitely many positive integers $n$. Let $M_{n}^{I_{K}}$ be the subspace of $M_n$ consisting of $I_{K}$-fixed vectors. Note that $M_{n}^{I_{K}}$ is a subrepresentation of $M_n$, and $\phi_n(\Frob_{K}) \in \GL(M_{n}^{I_{K}})$ is independent of the choice of $\Frob_{K}$. We define the \emph{local $L$-factor} associated with $\phi$ by
\begin{equation} \label{eq:localLfactordefn}
L(s, \phi) := \prod_{n = 1}^{\infty} L \left( s + \dfrac{n-1}{2}, \phi_n \right) := \prod_{n = 1}^{\infty} \det \left( \bfone - q_{K}^{-s-(n-1)/2} \phi_n(\Frob_K) \mid M_{n}^{I_{K}} \right)^{-1}.    
\end{equation}

\subsubsection{Non-archimedean Langlands correspondence for unitary groups}
Let $V$ be an $\varepsilon$-Hermitian space of dimension $n$ over $K$. An \emph{$L$-parameter} for the unitary group $\rmU(V)$ is an $n$-dimensional conjugate self-dual representation of $\WD_{K}$
\[
\phi: \WD_{K} \rightarrow \GL_{n}(\CC)
\]
with sign $(-1)^{n-1}$. Let $\Phi(n)$ be the set of equivalence classes of $L$-parameters for unitary groups of $n$ variables. On the automorphic side, let $\Irr(\rmU(V^{\varepsilon}))$ be the set of equivalence classes of irreducible smooth representations of $\rmU(V^{\varepsilon})$.

The non-archimedean Langlands correspondence asserts a \emph{canonical} finite-to-one surjection
\[
\calL_n: \Irr(\rmU(V^{+})) \sqcup \Irr(\rmU(V^{-})) \rightarrow \Phi(n).
\]
For each $L$-parameter $\phi \in \Phi(n)$, the inverse image $\Pi_{\phi} := \calL_n^{-1}(\phi)$ is called the \emph{$L$-packet} associated with $\phi$. Among the many properties of $\calL_{n}$ (including their defining properties), we list those useful to us in the following proposition. The definition and properties of the local adjoint $L$-factors in this proposition will be introduced in the proof of Theorem \ref{thm:localLfactor}.

\begin{proposition} \label{assertion}
With notation as above:
\begin{enumerate}[label = \rm (\arabic*)]
    \item The map $\calL_{n}$ sends unramified representations to unramified $L$-parameters.
    \item The map $\calL_{n}$ preserves local $L$-factors and local adjoint $L$-factors for tempered representations, in the sense that
    \[
    L(s, \pi) = L(s, \calL_{n}(\pi)), \, L(s, \pi, \Ad) = L(s, \calL_{n}(\Ad(\pi))),
    \]
    where the left-hand side is the local $L$-factor of the representation $\pi$ defined via its base change to $\GL_n(K)$.
    \item Let $n_1, n_2$ be positive integers. Then $\calL_{n_1}$ and $\calL_{n_2}$ preserve local Rankin--Selberg $L$-factors for tempered representations $\pi_i$ of $\rmU(V_{n_i})$ for $i=1,2$, in the sense that
    \[
    L(s, \pi_1 \times \pi_2) = L(s, \calL_{n_1}(\pi_1) \otimes \calL_{n_2}(\pi_2)),
    \]
    where the left-hand side is the local $L$-factor of the Rankin--Selberg product defined in \cite{jacquet1983rankin} under the base change of $\pi_i$ to $\GL_{n_i}(K)$, and the right-hand side is the tensor product of $L$-parameters.
\end{enumerate}
\end{proposition}

\subsection{Non-archimedean theta correspondences} \label{sec:nonarchitheta}
Let $V$ be a Hermitian space and $W$ be a skew-Hermitian space (see Section \ref{sec:hermitianspace}) over $K$. We consider the theta correspondence for the reductive dual pair $(\rmU(V), \rmU(W))$, with the following additional data:
\begin{itemize}
    \item a nontrivial additive character $\psi: F \rightarrow \CC$,
    \item a pair of characters $\chi_V$ and $\chi_W$ of $K^{\times}$ such that
    \[
    \chi_V|_{F^{\times}} = \epsilon_{K/F}^{\dim_{K} V} \quad \text{and} \quad \chi_W|_{F^{\times}} = \epsilon_{K/F}^{\dim_{K} W}.
    \]
\end{itemize}
The tuple $\ulpsi := (\psi, \chi_V, \chi_W)$ is called a \emph{local datum for the theta correspondence}. We first equip the tensor product space $V \otimes_{K} W$ with a natural symplectic form defined by
\[
\lrangle{v_1 \otimes w_1, v_2 \otimes w_2} := \tr_{K/F}(\lrangle{v_1, v_2}_{V} \cdot \lrangle{w_1, w_2}_{W}).
\]

There is a natural map
\[
\rmU(V) \times \rmU(W) \rightarrow \Sp(V \otimes_{K} W).
\]
Let $\Mp(V \otimes_{K} W)$ be the metaplectic cover of $\Sp(V \otimes_{K} W)$. The additive character $\psi$, together with the symplectic form $\lrangle{-,-}$ on $V \otimes W$, determines a \emph{Weil representation} $\omega_{\psi}$ of $\Mp(V \otimes W)$. The datum $\ulpsi := (\psi, \chi_V, \chi_W)$ allows one to specify a splitting of the metaplectic cover, thereby defining the Weil representation $\omega_{\ulpsi, V, W}$ of $\rmU(V) \times \rmU(W)$.

Given an irreducible representation $\pi$ of $\rmU(W)$, the maximal $\pi$-isotypic quotient of $\omega_{\ulpsi, V, W}$ is of the form
\[
\Theta_{\ulpsi, V, W}(\pi) \boxtimes \pi
\]
for some smooth representation $\Theta_{\ulpsi, V, W}(\pi)$ of $\rmU(V)$ of finite length. This is called the \emph{big theta lift} of $\pi$ from $\rmU(W)$ to $\rmU(V)$ with respect to $\ulpsi$. By Howe duality, as proved in \cite{waldspurger1990demostration, gan2016howe, gan2016proof}, it has a unique irreducible quotient $\theta_{\ulpsi, V, W}(\pi)$, called the \emph{small theta lift}.

\subsubsection{The conservation law}
Let $\calV = (V_{(r)})$ be a Witt tower of Hermitian spaces over $K$ such that $V_{(r)} = V_{(0)} \oplus \scH^{r}$; then $\dim V_{(r)} = \dim V_{(0)} + 2r$. We consider the big theta lifts $\Theta_{\ulpsi, V_{(r)}, W}(\pi)$ along the Witt tower $\calV$ for an irreducible representation $\pi$ of $\rmU(W)$. Kudla \cite{kudla1986local} established the following results.

\begin{theorem} Let $\pi$ be an irreducible representation of $\rmU(W)$. Then:
\begin{enumerate}[label = \rm (\arabic*)]
    \item The number $m_{\calV}(\pi) := \min \{r : \Theta_{\ulpsi, V_{(r)}, W}(\pi) \neq 0 \}$ is finite. Moreover, $m_{\calV}(\pi) \leq \dim W$.
    \item For any $r > m_{\calV}(\pi)$, we have $\Theta_{\ulpsi, V_{(r)}, W}(\pi) \neq 0$.
    \item If $\pi$ is supercuspidal, then $\Theta_{\ulpsi, V_{(m_{\calV}(\pi))}, W}(\pi)$ is irreducible supercuspidal. For $r > m_{\calV}(\pi)$, $\Theta_{\ulpsi, V_{(r)}, W}(\pi)$ is irreducible but not supercuspidal.
\end{enumerate}
\end{theorem}

We call $m_{\calV}(\pi)$ the \emph{first occurrence index} of $\pi$ in the Witt tower $\calV$, relative to the datum $\ulpsi$. When $r \geq \dim W$, we say the theta lift is in the \emph{stable range}.

For any Witt tower $\calV$, there is a companion Witt tower $\calV^{\prime}$. We denote the two towers by $\calV^{+}$ and $\calV^{-}$ based on their signs and denote the first occurrence indices of $\pi$ by $m^{\pm}(\pi)$ accordingly. The following \emph{conservation relation}, proved by Kudla--Rallis \cite{kudla2005first} and Sun--Zhu \cite{sun2015conservation}, is fundamental.

\begin{theorem} \label{thm:conservation}
With notation as above, $\dim V_{(m^{+}(\pi))} + \dim V_{(m^{-}(\pi))} = 2 \dim W + 2$.
\end{theorem}

This conservation relation implies that the first occurrence on one tower determines the occurrence on the other, playing an essential role in the theory.

\subsubsection{Non-archimedean theta correspondence in the almost equal rank case}
By the "almost equal rank" case, we mean the case where $\dim V = n+1$ and $\dim W = n$, and we shall consider the theta correspondence for $(\rmU(W_{n}^{\epsilon^{\prime}}), \rmU(V_{n+1}^{\epsilon}))$. Here we add the subscripts to indicate the $K$-dimension of corresponding $\varepsilon$-Hermitian spaces. We summarize the results we need into the following theorem.

\begin{theorem} \label{thm:nonarchinonvanishing}
Let $\phi$ be an $L$-parameter for $\rmU(W_{n}^{\pm})$. Then we have:
\begin{enumerate}[label = \rm (\arabic*)] 
    \item Suppose that $\phi$ does not contain $\chi_V$, then for any $\pi \in \Pi_{\phi}^{\epsilon^{\prime}}$, the big theta lift $\Theta_{\ulpsi, V_{n+1}^{\epsilon}, W_{n}^{\epsilon^{\prime}}}(\pi)$ is nonzero, and the small theta lift $\theta_{\ulpsi, V_{n+1}^{\epsilon}, W_{n}^{\epsilon^{\prime}}}(\pi)$ has $L$-parameter
    \[
    \theta_{\ulpsi, V_{n+1}^{\epsilon}, W_{n}^{\epsilon^{\prime}}}(\phi) := (\phi \otimes \chi_V^{-1} \chi_W) \oplus \chi_W.
    \]
    \item Suppose that $\phi$ contains $\chi_V$. For any fixed $\pi \in \Pi_{\phi}^{\epsilon^{\prime}}$, exactly one of $\Theta_{\ulpsi, V_{n+1}^{+}, W_{n}^{\epsilon^{\prime}}}(\pi)$ or $\Theta_{\ulpsi, V_{n+1}^{-}, W_{n}^{\epsilon^{\prime}}}(\pi)$ is nonzero. Moreover, if $\Theta_{\ulpsi, V_{n+1}^{\epsilon}, W_{n}^{\epsilon^{\prime}}}(\pi)$ is nonzero, then the small theta lift $\theta_{\ulpsi, V_{n+1}^{\epsilon}, W_{n}^{\epsilon^{\prime}}}(\pi)$ has $L$-parameter
    \[
    \theta_{\ulpsi, V_{n+1}^{\epsilon}, W_{n}^{\epsilon^{\prime}}}(\phi) := (\phi \otimes \chi_V^{-1} \chi_W) \oplus \chi_W.
    \]
    \item For any tempered representation $\pi \in \Pi_{\phi}^{\epsilon^{\prime}}$ such that the big theta lift $\Theta_{\ulpsi, V_{n+1}^{\epsilon}, W_{n}^{\epsilon^{\prime}}}(\pi)$ is nonzero, then the big theta lift $\Theta_{\ulpsi, V_{n+1}^{\epsilon}, W_{n}^{\epsilon^{\prime}}}(\pi)$ is irreducible, and hence itself is the small theta lift. Moreover, all irreducible subquotients of $\Theta_{\ulpsi, V_{n+1}^{\epsilon}, W_{n}^{\epsilon^{\prime}}}(\pi)$ are tempered.
\end{enumerate}
\end{theorem}

\begin{proof}
Items (1), (2), and the first part of (3) are summarized in \cite[Theorem 4.4]{gan2016gross}, which is proved in \cite{gan2014formal}. The last part of (3) is proved in \cite[Proposition 5.5]{atobe2017local}.
\end{proof}

\subsubsection{Behaviours of local $L$-factors at finite unramified places}
In Theorem \ref{thm:nonarchinonvanishing}, we have the relation of local $L$-parameters under theta lifts. With Proposition \ref{assertion} on the non-archimedean Langlands correspondence of unitary groups, we can express the local $L$-factors of the theta lift $\Theta_{\ulpsi, V_{n+1}^{\epsilon}, W_{n}^{\epsilon^{\prime}}}(\pi)$, provided it is nonzero, with those of $\pi$, at least at the places where $\pi$, $\chi_V$ and $\chi_W$ are unramified.

\begin{theorem} \label{thm:localLfactor}
Suppose $\pi$ and $\sigma$ are unramified representations of $\rmU(W_{n})$, and $\chi_V$ and $\chi_W$ are unramified characters. Assume the theta lift $\Theta_{\ulpsi, V_{n+1}, W_{n}}(\pi)$ is nonzero, then
\begin{enumerate}[label = \rm (\arabic*)]
    \item $L(s, \Theta_{\ulpsi, V_{n+1}, W_{n}}(\pi)) = L(s, \pi, \chi_V^{-1}\chi_W) L(s, \chi_W)$,
    \item $L(s, \sigma \times \Theta_{\ulpsi, V_{n+1}, W_{n}}(\pi)) = L(s, \sigma \times \pi, \chi_V^{-1}\chi_W) L(s, \sigma, \chi_W)$,
    \item $L(s, \Theta_{\ulpsi, V_{n+1}, W_{n}}(\pi), \Ad) = L(s, \pi, \Ad) L(s, \BC(\pi) \otimes \chi_{V}^{-1}) L(s, \epsilon_{K/F})$.
\end{enumerate}
\end{theorem}
\begin{proof}
(1) and (2) are direct consequences of Theorem \ref{thm:nonarchinonvanishing} and Proposition \ref{assertion}(2)(3). Note that the unramified condition is essential because, according to Proposition \ref{assertion}(1), the corresponding $L$-parameters under $\calL_{1}$, $\calL_n$, and $\calL_{n+1}$ are trivial on $I_{K} \times \SL_2(\CC)$, which implies that no direct summands $\sfS_{n}$ with $n > 1$ appear in the decomposition \eqref{eq:decomposephi}. Now we prove (3). 

Previously, we defined the $L$-parameter of the unitary group $\rmU(V)$ as a $\WD_{K}$-representation $\phi: \WD_{K} \rightarrow \GL_{n}(\CC)$ satisfying the conjugate self-dual condition. To calculate the local adjoint $L$-factors, we first introduce an equivalent definition of $L$-parameters. Since we have assumed $\pi$ is an unramified representation, following the same reasoning as above, we will omit Deligne's $\SL_2(\CC)$ factor and consider only the representations of the Weil group.

Recall that $\rmW_{F}$ is the Weil group of the local field $F$, and the Weil group $\rmW_{K}$ of its quadratic extension $K$ is a subgroup of index $2$. Fix $w_{\sfc} \in \rmW_{F} \smallsetminus \rmW_{K}$ as an element projecting to the nontrivial Galois automorphism $\sfc \in \Gal(K/F)$. For an $n$-dimensional Hermitian space $W$, the $L$-group of its unitary group $\rmU(W)$ is defined as ${}^{L} \rmU(W) := \GL_{n}(\CC) \rtimes \rmW_{F}$, where the action of $\rmW_{F}$ on $\GL_{n}(\CC)$ is defined as follows: for $g \in \GL_{n}(\CC)$, $w \cdot g = g$ for $w \in \rmW_{K}$, and $w_{\sfc} \cdot g = \Phi_{n} g^{-\rmt} \Phi_n^{-1}$ for $w_{\sfc} \in \rmW_{F}$, where 
\[
\Phi_{n} := \begin{bmatrix}
0      & \cdots & 0      & 1 \\
0      & \cdots & -1     & 0 \\
\vdots & \iddots & \vdots & \vdots \\
(-1)^{n+1} & \cdots & 0 & 0
\end{bmatrix}.
\]
The $L$-parameter of the unitary group $\rmU(W)$ can be equivalently defined as a representation of the Weil group $\phi: \rmW_{F} \rightarrow {}^{L} \rmU(W)$.

For the $L$-group ${}^{L} \rmU(W)$ of the unitary group $\rmU(W)$, we can define its \emph{adjoint representation} on its Lie algebra $\Lie({}^{L} \rmU(W)) = \gl_{n}(\CC)$ \footnote{Refer to \cite[Example 1.13]{harris2011refined}.}: for an element $(g, w) \in {}^{L} \rmU(W)$, define
\[
\Ad(g,w) \cdot X = g (w \cdot X) g^{-1}, \quad X \in \gl_{n}(\CC),
\]
where for $w \in \rmW_{K}$, $w \cdot X := X$, and for $w_{\sfc} \in \rmW_{F}$, $w_{\sfc} \cdot X := - \Phi_{n} X^{\rmt} \Phi_{n}^{-1}$. For an $L$-parameter $\phi: \rmW_{F} \rightarrow {}^{L} \rmU(W)$, we can then define its \emph{adjoint $L$-parameter} as the composition
\[
\Ad(\phi) := \Ad \circ \phi: \rmW_{F} \rightarrow {}^{L} \rmU(W) \xrightarrow{\Ad} \Aut_{\CC}(\gl_{n}(\CC)).
\]
Given such a $\rmW_{F}$-representation, its \emph{local adjoint $L$-factor} is the local $L$-factor $L(s, \Ad(\phi))$ defined by \eqref{eq:localLfactordefn}, usually denoted by $L(s, \phi, \Ad)$.

Returning to the proof of (3), since Theorem \ref{thm:nonarchinonvanishing} only gives (equivalently) the values of the $L$-parameters on $\rmW_{K}$, we first calculate on $\rmW_{K}$. Theorem \ref{thm:nonarchinonvanishing} yields
\[
\theta(\phi) := \theta_{\ulpsi, V_{n+1}^{\epsilon}, W_{n}^{\epsilon^{\prime}}}(\phi) = (\phi \otimes \chi_V^{-1} \chi_W) \oplus \chi_W,
\]
where $\phi_{1} := \phi \otimes \chi_V^{-1} \chi_W$ is an $n$-dimensional representation and $\chi_{W}$ is a $1$-dimensional representation. Let $\sfW$ denote the representation space of the $L$-parameter $\phi_{1}$, and $\sfE$ denote the representation space of the $L$-parameter $\chi_{W}$. Following the above direct sum decomposition, the representation space $\sfV$ of the $L$-parameter $\theta(\phi)$ decomposes as $\sfV = \sfW \oplus \sfE$. We further decompose the Lie algebra $\Lie({}^{L} \rmU(V)) = \gl_{n+1}(\CC)$ of the $L$-group ${}^{L} \rmU(V)$ of the unitary group $\rmU(V)$ accordingly:
\[
\gl_{n+1}(\CC) = \End_{\CC}(\sfW) \oplus \Hom_{\CC}(\sfW, \sfE) \oplus \Hom_{\CC}(\sfE, \sfW) \oplus \End_{\CC}(\sfE).
\]
The action of $\rmW_{K}$ on these direct summands is given by the adjoint (i.e., conjugation) action:
\begin{itemize}
    \item On $\End_{\CC}(\sfW)$, the action of the Weil group $\rmW_{K}$ is $\Ad(\phi_{1}) = \Ad(\phi)$;
    \item On $\Hom_{\CC}(\sfW, \sfE)$, the action of the Weil group $\rmW_{K}$ is $\phi_1^{\vee} \otimes \chi_{W} = (\phi^{\vee} \otimes \chi_V \chi_W^{-1}) \otimes \chi_{W} = \phi^{\vee} \otimes \chi_V$;
    \item On $\Hom_{\CC}(\sfE, \sfW)$, the action of the Weil group $\rmW_{K}$ is $\phi_1 \otimes \chi_{W}^{-1} = (\phi \otimes \chi_V^{-1} \chi_W) \otimes \chi_{W}^{-1} = \phi \otimes \chi_V^{-1}$;
    \item On $\End_{\CC}(\sfE)$, the action of the Weil group $\rmW_{K}$ is the trivial action.
\end{itemize}
We further descend to the Weil group $\rmW_{F}$, noting that:
\begin{itemize}
    \item On $\End_{\CC}(\sfW)$, the action of the Weil group $\rmW_{F}$ is the previously defined $\Ad(\phi)$;
    \item The action of the Weil group element $w_{\sfc} \in \rmW_{F} \smallsetminus \rmW_{K}$ interchanges the direct summands $\Hom_{\CC}(\sfW, \sfE)$ and $\Hom_{\CC}(\sfE, \sfW)$, so the action of the Weil group $\rmW_{F}$ on $\sfX := \Hom_{\CC}(\sfW, \sfE) \oplus \Hom_{\CC}(\sfE, \sfW)$ is the induced representation $\Ind_{\rmW_{K}}^{\rmW_{F}}(\phi \otimes \chi_V^{-1})$;
    \item On $\End_{\CC}(\sfE)$, the action of the Weil group element $w_{\sfc} \in \rmW_{F} \smallsetminus \rmW_{K}$ is multiplication by $-1$. This corresponds (via class field theory) exactly to the quadratic character $\epsilon_{K/F}$ of the field extension $K/F$.
\end{itemize}
Thus, according to the direct sum property of local $L$-factors for $L$-parameters (cf. the textbook \cite[Eq. (12.8)]{getz2024}), we obtain
\[
L(s, \theta(\phi), \Ad) = L(s, \phi, \Ad) L(s, \Ind_{\rmW_{K}}^{\rmW_{F}}(\phi \otimes \chi_V^{-1})) L(s, \epsilon_{K/F}).
\]
Using the non-archimedean Langlands correspondence for unitary groups (Theorem \ref{assertion} \footnote{For the factor containing the induced representation, a more precise formulation is as follows: first note that, by the textbook \cite[Eq. (12.7)]{getz2024}, we have 
\[
L(s, \Ind_{\rmW_{K}}^{\rmW_{F}}(\phi|_{\rmW_{K}} \otimes \chi_V^{-1})) = L(s, \phi |_{\rmW_{K}} \otimes \chi_V^{-1}),
\]
and the right side of the equality, by the defining property of the base change representation $\BC(\pi)$, is $L(s, \BC(\pi) \otimes \chi_{V}^{-1})$.
}), we return to the local $L$-factors of automorphic representations, which translates to
\[
L(s, \Theta_{\ulpsi, V_{n+1}, W_{n}}(\pi), \Ad) = L(s, \pi, \Ad) L(s, \BC(\pi) \otimes \chi_{V}^{-1}) L(s, \epsilon_{K/F}),
\]
as desired.
\end{proof}

We remark that item (3) of the above theorem is a generalization of \cite[Proposition A.3]{harris2014refined}: in fact, the authors there consider the theta lifting from $\rmU(2,0)$ to $\rmU(3,0)$ and take the special case where the splitting characters $\chi_{W} = \gamma^{2}$ and $\chi_{V} = \gamma^{3}$ for a fixed character $\gamma$, noting that \footnote{To see \eqref{eq:conjugateselfdual} on the automorphic side, since the local representations and local field extensions have been assumed to be unramified, we can describe the local $L$-factors using the Satake parameters of the $\GL_{n}(K)$-representation $\BC(\pi) \otimes \chi_{V}$ (see \cite[Example 1.10]{harris2011refined}). At this point, note that the set of Satake parameters for the local representation $\BC(\pi)$ coming from the base change of the unitary group $\rmU(W)$ is invariant under inversion (i.e., the transition from $\pi$ to $\pi^{\vee}$). In fact, letting $n := \dim W$ be $2m$ or $2m+1$ accordingly, and denoting the uniformizers of $K$ and $F$ by $\varpi$, \cite[Theorem 4.1]{minguez2011unramified} proves that if the Satake parameters of $\pi$ are $(\chi_1(\varpi), \ldots, \chi_{m}(\varpi))$, then the Satake parameters of the base change $\BC(\pi)$ are as follows: if $n=2m$ is even, then
\[
(\chi_1(\varpi), \ldots, \chi_{m}(\varpi),  \chi_{m}^{-1}(\varpi), \ldots, \chi_{1}^{-1}(\varpi));
\]
if $n = 2m+1$ is odd, then
\[
(\chi_1(\varpi), \ldots, \chi_{m}(\varpi), 1,  \chi_{m}^{-1}(\varpi), \ldots, \chi_{1}^{-1}(\varpi)).
\]
Finally, the character $\chi_{V}$ as a splitting character satisfies $\chi_{V}|_{F^{\times}} = \varepsilon_{K/F}^{n+1}$, and thus $\chi_{V}(\varpi) = \chi_{V}^{-1}(\varpi) = -1$.
}
\begin{equation} \label{eq:conjugateselfdual}
 L(s, \BC(\pi^{\vee}) \otimes \chi_{V}^{-1}) = L(s, \BC(\pi) \otimes \chi_{V}).   
\end{equation}
\cite[Appendix A]{harris2014refined} also handles some low-rank cases which are not almost-equal rank. We will apply this theorem in Chapter \ref{sec:imcforcm} to the theta lift of a Hecke character $\eta: \calK^{\times} \bs \AA_{\calK}^{\times} \rightarrow \CC^{\times}$ from $\rmU(1,0)$ to $\rmU(2,0)$. In this case, $L(s, \BC(\eta) \otimes \chi_{V}^{-1}) = L(s, \eta\chi_{V}^{-1})$ is a Hecke $L$-function over the number field $\calK$.

\subsection{Archimedean theta correspondences} \label{sec:architheta}
In this section, we consider the theory of archimedean theta correspondences. Let $p \geq q$ and $r \geq s$ be nonnegative integers, and let $W_{p,q}$ and $V_{r,s}$ be a skew-Hermitian space and a Hermitian space over $\CC$ of signature $(p,q)$ and $(r,s)$, respectively. We denote $n = p+q$ and $m=r+s$. We consider the theta lifts for the reductive dual pair $(\rmU(p,q), \rmU(r,s))$. Similar to the non-archimedean case, we require the following additional data:
\begin{itemize}
    \item a nontrivial additive character $\psi: \RR \rightarrow \CC$, which we fix for certainty as $\psi(x) = \exp(-2 \pi \sqrt{-1} x)$,
    \item a pair of characters $\chi_V$ and $\chi_W$ of $\CC^{\times}$, which we fix as
    \[
    \chi_V(z) = \left( \dfrac{z}{\sqrt{z \bar{z}}} \right)^{m_0} \quad \text{and} \quad \chi_W(z) = \left( \dfrac{z}{\sqrt{z \bar{z}}} \right)^{n_0},
    \]
    with fixed integers $m_0, n_0$ such that $m_0 \equiv m \bmod 2$ and $n_0 \equiv n \bmod 2$.
\end{itemize}
Again, the tuple $\ulpsi := (\psi, \chi_V, \chi_W)$ is called a \emph{local datum for the theta correspondence}.

Parallel to the non-archimedean case, we obtain a symplectic space $V_{r,s} \otimes W_{p,q}$ and a natural map
\[
\rmU(r,s) \times \rmU(p,q) \rightarrow \Sp(V \otimes W).
\]
Let $\Mp(V \otimes W)$ be the metaplectic cover of $\Sp(V \otimes W)$. The additive character $\psi$, together with the symplectic form $\lrangle{-,-}$ on $V \otimes W$, determines a \emph{Weil representation} $\omega_{\psi}$ of $\Mp(V \otimes W)$. The datum $\ulpsi := (\psi, \chi_V, \chi_W)$ allows one to specify a splitting of the metaplectic cover, thereby defining a Weil representation $\omega_{\ulpsi, V, W}$ of $\rmU(r,s) \times \rmU(p,q)$.

Given an irreducible (unitary) representation $\pi$ of $\rmU(p,q)$, the maximal $\pi$-isotypic quotient of $\omega_{\ulpsi, V,W}$ is of the form
\[
\Theta_{\ulpsi, V_{r,s}, W_{p,q}}(\pi) \boxtimes \pi
\]
for some representation $\Theta_{\ulpsi, V_{r,s}, W_{p,q}}(\pi)$ of $\rmU(r,s)$. This is called the \emph{big theta lift} of $\pi$ from $\rmU(p,q)$ to $\rmU(r,s)$ with respect to the datum $\ulpsi$. By Howe duality \cite[Theorem 2.1]{howe1989transcending}, when $\Theta_{\ulpsi, V_{r,s}, W_{p,q}}(\pi)$ is nonzero, it has a unique irreducible quotient $\theta_{\ulpsi, V_{r,s}, W_{p,q}}(\pi)$, called the \emph{small theta lift}. When $\ulpsi$ and $\rmU(p,q)$ are clear from the context, we simply write $\Theta_{r,s}(\pi)$ (resp.\ $\theta_{r,s}(\pi)$) for the big (resp.\ small) theta lift.

In this article, we focus on discrete series representations of unitary groups at Archimedean places. These are automatically tempered by definition. We adopt the conventions from Section \ref{sec:iiarchi}.

\subsubsection{Archimedean theta correspondences}
We consider $\Theta_{\ulpsi, V_{r,s}, W_{p,q}}(\DS(\lambda))$, where $\DS(\lambda)$ is a discrete series representation of $\rmU(p,q)$ determined by the Harish-Chandra parameter $\ullam$. We fix $k_0 \in \{-1, 0\}$ such that $m \equiv n+k_0 \bmod 2$. We then write
\[
\ullam = \lambda_{\ulpsi} + \left( \dfrac{m_0}{2}, \ldots, \dfrac{m_0}{2} \right)
\]
with $\lambda_{\ulpsi} = (\lambda_{\ulpsi, 1}, \ldots, \lambda_{\ulpsi, n})$, such that
\begin{itemize}
    \item $\lambda_{\ulpsi, i} \in \ZZ + \frac{k_0 - 1}{2}$,
    \item $\lambda_{\ulpsi, i} \neq \lambda_{\ulpsi, j}$ if $i \neq j$,
    \item $\lambda_{\ulpsi, 1} > \cdots > \lambda_{\ulpsi, p}$ and $\lambda_{\ulpsi, p+1} > \cdots > \lambda_{\ulpsi, q}$.
\end{itemize}

The non-vanishing property of theta lifts of discrete series can be described in terms of Harish-Chandra parameters, which also yields an explicit description of the theta lifts themselves.

First, we introduce Jian-Shu Li's result \cite{li1990theta} on the non-vanishing of theta lifts for discrete series representations with sufficiently regular infinitesimal characters.

\begin{theorem} \label{thm:lioriginal}
Let $\pi:=\DS(\ullam)$ be an irreducible discrete series representation of $\rmU(p,q)$ with Harish-Chandra parameter $\ullam$. Then the theta lift $\Theta_{\ulpsi, V_{r,s}, W_{p,q}}(\pi)$ is nonzero if $m \geq n$ and the Harish-Chandra parameter $\lambda$ of $\pi$ can be written as
\begin{equation} \label{eq:writelambda}
\lambda = \lambda_{\ulpsi} + \left( \dfrac{m_0}{2}, \ldots, \dfrac{m_0}{2} \right), \quad \text{with } \lambda_{\ulpsi} = (\alpha_1, \ldots, \alpha_x, \beta_1, \ldots, \beta_y, \gamma_1, \ldots, \gamma_z, \delta_1, \ldots, \delta_w),
\end{equation}
satisfying
\begin{itemize}
    \item $x+y=p$ and $z+w=q$,
    \item $x+w \leq r$ and $z+y \leq s$,
    \item and $\alpha_x, -\beta_1, \gamma_z, -\delta_1 \geq \frac{m-n+1}{2}$.
\end{itemize}
\end{theorem}

\begin{remark} \label{rem:onjianshu}
We remark that Theorem \ref{thm:lioriginal} provides only a sufficient condition for the non-vanishing of $\Theta_{r,s}(\pi)$, and may not be necessary. There are many subsequent works on this aspect, as well as generalizations to wider classes of representations, such as tempered representations. We refer the reader to the article \cite{atobe2020non} for a comprehensive survey. For our current purposes, Li's result above suffices. We also note that we use the reinterpreted version of Li's theorem found in \cite[Remark 4.5]{ichino2020thetaliftingdiscreteseries}.
\end{remark}

Second, the theta lift $\Theta_{\ulpsi, V_{r,s}, W_{p,q}}(\pi)$ remains a discrete series representation. We have an explicit description of its Harish-Chandra parameters, and even the construction of the theta lift in the form of cohomologically induced representations. Due to the complexity of the combinatorial expressions involved, we refer interested readers directly to \cite[Section 4.2]{ichino2020thetaliftingdiscreteseries}; we shall introduce a simplified version for the reductive dual pairs $(\rmU(n,0), \rmU(n+1, 0))$ and $(\rmU(n,0), \rmU(n, 1))$ in the following subsection.

\subsubsection{Archimedean theta correspondence in the almost equal compact case}
By the ``almost equal rank'' case, we mean $m = n+1$, and by ``compact'', we mean theta lifts from compact unitary groups $\rmU(n,0)$. In this setting, two reductive dual pairs remain:
\begin{itemize}
    \item \textbf{Case (A)}: $(\rmU(n,0), \rmU(n+1, 0))$, and
    \item \textbf{Case (B)}: $(\rmU(n,0), \rmU(n, 1))$.
\end{itemize}
In both cases, $k_0 = -1$. For a discrete series representation $\DS(\ullam)$ of $\rmU(n,0)$, we can write
\[
\ullam = \lambda_{\ulpsi} + \left( \dfrac{m_0}{2}, \ldots, \dfrac{m_0}{2} \right) \text{, with } \lambda_{\ulpsi} = (\lambda_{\ulpsi, 1}, \ldots, \lambda_{\ulpsi, n})
\]
such that $\lambda_{\ulpsi, i} \in \ZZ$ and $\lambda_{\ulpsi, 1} > \cdots > \lambda_{\ulpsi, n}$.

In these compact almost equal rank cases, we have the following interpretations of Theorem \ref{thm:lioriginal}:

\begin{corollary} \label{coro:archinonvanishing}
Let $\pi := \DS(\ullam)$ be an irreducible discrete series representation of $\rmU(n,0)$ with Harish-Chandra parameter $\ullam$.
\begin{enumerate}[label = \rm (\arabic*)]
    \item In Case (A), the theta lift $\Theta_{\ulpsi, V_{n+1,0}, W_{n,0}}(\pi)$ is nonzero if $\lambda_{\ulpsi,n} \geq 1$.
    \item In Case (B), the theta lift $\Theta_{\ulpsi, V_{n,1}, W_{n,0}}(\pi)$ is nonzero if 
    \begin{enumerate}[label = \rm (\alph*)]
        \item \textbf{Case (B0)}: $\lambda_{\ulpsi,n} \geq 1$, or
        \item \textbf{Case (B1)}: $\lambda_{\ulpsi,n-1} \geq 1$ and $\lambda_{\ulpsi,n} \leq -1$.
    \end{enumerate}
\end{enumerate}
\end{corollary}

\begin{proof}
This is a direct application of Theorem \ref{thm:lioriginal}. We designate ``Case (B0)'' and ``Case (B1)'' because the value of $y$ in the notation of \eqref{eq:writelambda} is $0$ or $1$, respectively.
\end{proof}

Secondly, for the explicit expression of the Harish-Chandra parameters in this case, we have the following corollary of \cite[Theorem 4.1]{ichino2020thetaliftingdiscreteseries}.

\begin{corollary} \label{cor:hcgloballift}
Let $\pi := \DS(\ullam)$ be an irreducible discrete series representation of $\rmU(n,0)$ with Harish-Chandra parameter $\ullam$. Assume that the theta lift $\Theta_{\ulpsi, V, W}(\pi)$ is nonzero.
\begin{enumerate}[label = \rm (\alph*)]
    \item In Case (A), the theta lift $\Theta_{\ulpsi, V_{n+1,0}, W_{n,0}}(\pi)$ is a discrete series representation with Harish-Chandra parameter
    \[
    \Theta_{\psi, V_{n+1,0}, W_{n,0}}(\ullam) := (\lambda_{\ulpsi, 1}, \ldots, \lambda_{\ulpsi, n}, 0) + \left( \dfrac{n_0}{2}, \ldots, \dfrac{n_0}{2} \right).
    \]
    \item In Case (B), the theta lift $\Theta_{\ulpsi, V_{n,1}, W_{n,0}}(\pi)$ is a discrete series representation with Harish-Chandra parameter
    \[
    \Theta_{\psi, V_{n,1}, W_{n,0}}(\ullam) := (\lambda_{\ulpsi, 1}, \ldots, \lambda_{\ulpsi, n}, 0) + \left( \dfrac{n_0}{2}, \ldots, \dfrac{n_0}{2} \right)
    \]
    in Case (B0), and 
    \[
    \Theta_{\psi, V_{n,1}, W_{n,0}}(\ullam) := (\lambda_{\ulpsi, 1}, \ldots, \lambda_{\ulpsi, n-1}, 0, \lambda_{\ulpsi, n}) + \left( \dfrac{n_0}{2}, \ldots, \dfrac{n_0}{2} \right)
    \]
    in Case (B1).
\end{enumerate}
\end{corollary}

\begin{proof}
In \cite[Theorem 4.1]{ichino2020thetaliftingdiscreteseries}, it is shown that 
\[
\Theta_{\ulpsi, V, W}(\pi_{\lambda}) = A_{\frq}(\lambda^{\prime}),
\]
where $A_{\frq}(\lambda^{\prime})$ is the cohomologically induced representation given by parameters $\frq$ and $\lambda^{\prime}$, expressed in explicit terms (see \cite[Section 3.2]{ichino2020thetaliftingdiscreteseries} for notation). We note that the Harish-Chandra parameter of $A_{\frq}(\lambda^{\prime})$ is not $\lambda^{\prime}$, but rather $\lambda^{\prime} + \varrho$, where $\varrho$ is the half-sum of the positive roots $\Delta^{+}$ in $\Delta$\footnote{Evidence for this can be found in \cite[Proposition 6.1]{li1990theta}; note the relationship between infinitesimal characters and Harish-Chandra parameters. In this proof, the aforementioned root system $\Delta$ is that of either $\rmU(n+1, 0)$ or $\rmU(n,1)$, rather than $\rmU(n,0)$.}, which is 
\[
\varrho = \left( \dfrac{n}{2}, \dfrac{n-2}{2}, \ldots, \dfrac{-n}{2}  \right).
\]
The result then follows. We remark that this correspondence for discrete series representations is well-established; it appears, for instance, in \cite[Proposition 6.6]{adams2007theta} for the case where $m_0 = m$ and $n_0 = n$. The contribution of Ichino in \cite{ichino2020thetaliftingdiscreteseries} is the description of the theta lift specifically in the form of cohomological induction.
\end{proof}

\subsection{Global theta correspondences} \label{sec:globaltheta}
In this section, we consider global theta correspondences. As before, let $\calK/\calF$ be a CM extension of number fields and $d := [\calF:\QQ]$.

We denote by $\WW_{n}$ an $n$-dimensional skew-Hermitian space over $\calK$, and by $\VV_{m}$ an $m$-dimensional Hermitian space over $\calK$ with anisotropic kernel $\VV_{(0)}$, so that $\VV_{m} = \VV_{(0)} \oplus \HH^{r}$, where $m = \dim_{\calK} V_{(0)} + 2r$ and $\HH$ is the hyperbolic plane over $\calK$. The sequence of spaces $\calV := (\VV_{(r)} := \VV_{(0)} \oplus \HH^{r})_{r \geq 0}$ is a global version of a Witt tower.

Parallel to the local cases, the tensor product space $\WW_{n} \otimes_{\calK} \VV_{m}$ can be regarded as a $2mn$-dimensional symplectic space over $\calF$, and we have a natural map
\[
\iota_{\VV, \WW}: \rmU(\VV_{m})(\AA) \times \rmU(\WW_{n})(\AA) \rightarrow \Sp(\WW_{n} \otimes_{\calK} \VV_{m})(\AA).
\]
Let $\Mp(\WW_{n} \otimes_{\calK} \VV_{m})$ be the metaplectic cover of $\Sp(\WW_{n} \otimes_{\calK} \VV_{m})(\AA)$, and take a complete polarization
\[
\WW_{n} \otimes_{\calK} \VV_{m} = \XX + \YY.
\]
We fix a nontrivial additive character $\psi: \calF \bs \AA \rightarrow \CC^{\times}$ and obtain the Weil representation $\omega_{\psi}$ of $\Mp(\WW_{n} \otimes_{\calK} \VV_{m})$ associated with $\psi$, realized on the space $\scS(\XX(\AA))$ of Bruhat--Schwartz functions on $\XX(\AA)$. Consider the following functional:
\[
\Theta_{\psi, \VV_{m}, \WW_n}: \scS(\XX(\AA)) \rightarrow \CC, \quad \phi \longmapsto \Theta_{\psi, \VV_{m}, \WW_n}(\phi) := \sum_{x \in \XX(\calF)} \phi(x).
\]

We fix two Hecke characters $\chi_{\VV}$ and $\chi_{\WW}$ over $\calK$ such that 
\[
\chi_{\WW}|_{\AA^{\times}} = \epsilon_{\calK/\calF}^{n} \quad \text{and} \quad \chi_{\VV}|_{\AA^{\times}} = \epsilon_{\calK/\calF}^{m}.
\]
Similar to the local cases, we call the tuple $\ulpsi := (\psi, \chi_{\VV}, \chi_{\WW})$ a \emph{global datum for the theta correspondence}. From now on, we assume that any global datum $\ulpsi$ is a local datum for the theta correspondence at Archimedean places. There is an explicit splitting
\[
\tilde{\iota}_{\chi_{\VV}, \chi_{\WW}}: \rmU(\VV_{m})(\AA) \times \rmU(\WW_{n})(\AA) \rightarrow \Mp(\WW_{n} \otimes_{\calK} \VV_{m})
\]
determined by $\chi_{\VV}$ and $\chi_{\WW}$. We then have the pullback of the Weil representation $\omega_{\ulpsi, \VV, \WW} := \omega_{\psi} \circ \tilde{\iota}_{\chi_{\VV}, \chi_{\WW}}$ to the product group $\rmU(\VV_{m})(\AA) \times \rmU(\WW_{n})(\AA)$. For each $\phi \in \scS(\XX(\AA))$, we consider the theta function
\begin{equation} \label{eq:globalthetafunction}
 \Theta_{\ulpsi, \VV, \WW}(g,h; \phi) := \Theta_{\psi, \VV, \WW}(\omega_{\ulpsi, \VV, \WW}(g,h) \phi)   
\end{equation}
for $g \in \rmU(\VV_{m})(\AA)$ and $h \in \rmU(\WW_{n})(\AA)$, which is an automorphic function on $\rmU(\VV_{m})(\AA) \times \rmU(\WW_{n})(\AA)$. Let $(\pi, V_{\pi})$ be an irreducible cuspidal automorphic representation of $\rmU(\WW_{n})(\AA)$. The \emph{global theta lift} $\Theta_{\ulpsi, \VV, \WW}(\pi)$ is defined to be the automorphic representation of $\rmU(\VV_{m})$ generated by the functions
\begin{equation} \label{eq:globalthetalift}
  \Theta_{\ulpsi, \VV, \WW}(g; f, \phi) := \int_{\rmU(\WW)(\calF) \bs \rmU(\WW)(\AA)} \Theta_{\ulpsi, \VV, \WW}(g,h; \phi) \barf(h) \dif h,  
\end{equation}
where $f \in V_{\pi}$ and $\phi \in \scS(\XX(\AA))$, with $\dif h$ being the Tamagawa measure. When $\pi$ is a unitary representation realized in $L^{2}(\rmU(\WW)(\calF) \bs \rmU(\WW)(\AA))_{\cusp}$, the function $\Theta_{\ulpsi, \VV, \WW}(-; f, \phi)$ is an automorphic form on the group $\rmU(\VV)$, called the \emph{theta lift} of $f$ with respect to the datum $\ulpsi$ and the Bruhat--Schwartz function $\phi$.

Next, we introduce two essential properties of global theta correspondences: the tower property and non-vanishing results.

The first is the following tower property, which provides a criterion for the range of non-vanishing of global theta lifts and their cuspidality. This is a result of Rallis \cite{rallis1984howe}.
\begin{theorem} \label{thm:globaltower}
Let $\pi$ be a cuspidal automorphic representation of $\rmU(\WW_{n})$, and consider its global theta lifts $\Theta_{\ulpsi, \VV_{(r)}, \WW_{n}}(\pi)$ along the Witt tower $\calV := (\VV_{(r)})_{r \geq 0}$. Then:
\begin{enumerate}[label = \rm (\arabic*)]
    \item The number $m_{\calV}(\pi) := \min\{r: \Theta_{\ulpsi, \VV_{(r)}, \WW_{n}}(\pi) \neq 0\}$ is finite. Moreover, $m_{\calV}(\pi) \leq \dim \WW_{n}$.
    \item The global theta lift $\Theta_{\ulpsi, \VV_{(m_{\calV}(\pi))}, \WW_{n}}(\pi)$ is contained in the space of cuspidal automorphic forms; i.e., it is a cuspidal automorphic representation.
    \item For all $r > m_{\calV}(\pi)$, $\Theta_{\ulpsi, \VV_{(r)}, \WW_{n}}(\pi)$ is nonzero and non-cuspidal.
\end{enumerate}
\end{theorem}
As in the local cases, we call $m_{\calV}(\pi)$ the \emph{first occurrence index} of $\pi$ in the Witt tower $\calV$ relative to the fixed datum $\ulpsi$.

Secondly, the irreducibility and the local components of $\Theta_{\ulpsi, \VV_{(m_{\calV}(\pi))}, \WW_{n}}(\pi)$ were determined by Kudla--Rallis \cite[Corollary 7.1.3]{kudla1994regularized}.
\begin{theorem} \label{thm:localglobalcompa}
With notation as above, let $\pi = \otimes^{\prime}_{v \in \scV_{\calF}} \pi_v$. If $\Theta_{\ulpsi, \VV_{(r)}, \WW_{n}}(\pi)$ is nonzero and cuspidal, then it is irreducible and
\[
\Theta_{\ulpsi, \VV_{(r)}, \WW_{n}}(\pi) = \otimes_{v \in \scV_{\calF}}^{\prime} \theta_{\ulpsi_{v}, \VV_{(r), v}, \WW_{n,v}}(\pi_v),
\]
where $\theta_{\ulpsi_{v}, \VV_{(r), v}, \WW_{n,v}}(\pi_v)$ is the small theta lift of $\pi_v$.
\end{theorem}

Finally, we address the non-vanishing problem of global theta lifts. Rallis established a local-to-global criterion stating that the global theta lift is nonzero if and only if its local components are all nonzero and there are no global obstructions, which are reflected by certain global special $L$-values. Here we state Rallis' result for the almost equal rank case (namely $m=n+1$) using the formulation in \cite[Theorem 1.3]{gan2014regularized}.
\begin{theorem} \label{thm:globalnonvanishing}
Let $\pi$ be a cuspidal automorphic representation of $\rmU(\WW_{n})$ and let $\calV$ be a Witt tower. Then $\Theta_{\ulpsi, \VV_{n+1}, \WW_{n}}(\pi)$ is nonzero if and only if:
\begin{enumerate}[label = \rm (\alph*)]
    \item for all places $v \in \scV_{\calF}$, $\Theta_{\ulpsi_v, \VV_{(r), v}, \WW_{n,v}}(\pi_v) \neq 0$; and
    \item $L(1, \pi \times \chi_{\VV}) \neq 0$.
\end{enumerate}
\end{theorem}

We note that the almost equal rank case provides a significant advantage because, as discussed in \cite[bottom of page 745]{gan2016gross}, when $\pi$ is tempered, the non-vanishing of the value $L(1, \pi \times \chi_{\VV})$ is automatically satisfied according to \cite[Proposition 5]{lapid2005local}. Its holomorphicity follows from the fact that $s = 1$ is not the central critical value of $L(s, \pi \times \chi_{\VV})$. Therefore, there is no global obstruction in the almost equal rank case.

%% file: AZ01-IMCU1.tex
\section{Iwasawa main conjecture for CM fields} \label{sec:imcforcm}
In this section, we apply the machinery developed in Part \ref{part:automorphic} to the lowest rank case: the Iwasawa theory for $\rmU(1,0)$, or equivalently, the Iwasawa theory for CM fields $\calK/\calF$. This also applies to the Iwasawa theory of elliptic curves with complex multiplication (\textit{viz.} CM elliptic curves) over totally real fields $\calF$. It turns out that our new approach, which considers the Gan--Gross--Prasad type period integrals of Klingen Eisenstein series over $\rmU(2,1)$, enables us to remove some of the conditions on these Iwasawa main conjectures that were required in \cite{hsieh2014eisenstein}.

\subsection{Auxiliary CM family over $\rmU(2,0)$}
As in Section \ref{sec:introduction}, the computation of the Gan--Gross--Prasad type period integral of the Klingen Eisenstein family over $\rmU(2,1)$ involves a choice of an auxiliary cuspidal Hida family $\bfg$ over $\rmU(2,0)$. In this section, we shall consider the Hida family of theta lifts from $\rmU(1,0)$ to $\rmU(2,0)$ as a candidate, following the construction by Lee in \cite{lee2025hidafamilythetalift} as a pullback of the theta lift.

\subsubsection{Setup} We consider the following particular (skew-)Hermitian spaces $W$ and $V$, so that $\rmU(V) = \rmU(1,0)$ and $\rmU(W) = \rmU(2,0)$.
\begin{itemize}
    \item Let $(W, \lrangle{-,-}_{W})$ be a skew-Hermitian space over $\calK$ of dimension $2$. We fix a $\calK$-basis $\{\tte_1, \tte_{2}\}$ for which the skew-Hermitian matrix
    \[
    \bdelta := \left[\lrangle{\tte_i, \tte_j}_{W} \right] = \diag[\delta_1, \delta_{2}]
    \]
    is diagonal. We require $i_v(\sfi^{-1} \delta_i) > 0$, so that $W$ is an anisotropic space over $\calK$ of dimension $2$.
    \item Let $V := \calK$ be the Hermitian space defined by $(z,z^{\prime})_{V} = \barz z^{\prime}$. We identify $\GU(\AA_{\calF}) = \AA_{\calK}^{\times}$ and the similitude factor with the norm map $\Norm_{\calK/\calF}$, so that $\rmU(\AA_{\calF}) = \AA_{\calK}^{1}$. We adopt this identification for general unitary groups associated with any one-dimensional Hermitian or skew-Hermitian spaces.
\end{itemize}

We introduce a Hecke character $\lambda: \calK^{\times} \bs \AA_{\calK}^{\times} \rightarrow \CC^{\times}$ satisfying the condition $\lambda|_{\AA_{\calF}^{\times}} = \epsilon_{\calK/\calF}$. Following \cite[Section 3]{lee2025hidafamilythetalift}, we assume:
\begin{enumerate}[label = (Lee-\roman*)]
    \item The central $L$-value $L(\frac{1}{2}, \lambda)$ is nonzero; \label{item:leei}
    \item The conductor $\frc$ of $\lambda$ is coprime to $p$; \label{item:leeii}
    \item The Hecke character $\lambda \abs{-}_{\AA_{\calK}^{\times}}^{1/2}$ is algebraic and has infinity type $\Sigma$. \label{item:leeiii}
\end{enumerate}
After choosing such a character $\lambda$, which we refer to as the \emph{splitting character}, we choose $\delta = \delta_1 = \delta_2$ to satisfy conditions ($\delta 1$) and ($\delta 2$) in \cite[p.\ 12]{lee2025hidafamilythetalift}, which we shall not recall here. We also fix an additive character $\psi: \AA_{\calF}/\calF \rightarrow \CC$ such that $\psi_v(x) = \exp(-2 \pi \sfi x)$ for any archimedean place $v$ of $\calF$.

Consider the symplectic space $\WW$ over $\calF$ underlying $V \otimes_{\calK} W$ and fix any polarization $\WW = \XX + \YY$. 
Then the Weil representation $\omega_{\psi}$ of $\Mp(\WW)$ acts on the space of Bruhat--Schwartz functions $\calS(\XX(\AA_{\calF}))$. Referring back to Section \ref{sec:globaltheta} (but interchanging the roles of $V$ and $W$ there), we use $\chi_{V} := \lambda$ and $\chi_{W} := 1$ (the trivial character) as splitting characters to define the Weil representation $\omega_{\lambda} := \omega_{\psi, \lambda, 1, V, W}$ over $\rmU(W) \times \rmU(V)$ acting on $\calS(\XX(\AA_{\calF}))$. We shall consider the global theta lift
\[
\Theta_{\lambda}(g; \eta, \phi) := \int_{\rmU(V)(\calF) \bs \rmU(V)(\AA_{\calF})} \Theta_{\lambda}(g,h; \phi) \bareta(h) \dif h
\]
following \eqref{eq:globalthetalift} for the theta function $\Theta_{\lambda}(g,h; \phi)$ considered in \eqref{eq:globalthetafunction}, where $\eta: \AA_{\calK}^{\times} \rightarrow \CC^{\times}$ is an algebraic Hecke character over $\calK$.

\subsubsection{Interlude: Pullback theta lifts} \label{sec:pullbackthetalift}
Although it is not needed in our main results, we briefly recall the construction by Lee. Let $\bfV = W \oplus (-W)$ be the direct sum of the skew-Hermitian spaces. We fix a standard basis
    \[
    \ttx^{i} := (2\delta_{i})^{-1}(\tte_{i} - \tte_{i}^{-}), \quad \tty^{i} := \tte_{i}+\tte_{i}^{-}, \, i = 1, 2.
    \]
    The associated skew-Hermitian matrix is $w_1 := \begin{bmatrix}   & -1 \\ 1 & \end{bmatrix}$. There exists a natural map $\iota^{\dsuit}: \rmG(\rmU(W) \times \rmU(-W)) \rightarrow \GU(\bfV)$:
    \[
    (g_1, g_2) \mapsto \begin{bmatrix} (2 \bdelta)^{-1} & -(2\bdelta)^{-1} \\ 1 & 1 \end{bmatrix} \begin{bmatrix}
    g_1 & \\ & g_2
    \end{bmatrix} \begin{bmatrix} (2 \bdelta)^{-1} & -(2\bdelta)^{-1} \\ 1 & 1 \end{bmatrix}^{-1}.
    \]
Lee first constructed a $p$-adic family of theta lifts of certain characters $\eta$ over $\rmU(V)$ to the quasi-split unitary group $\rmU(\bfV)$, denoted by $\Theta_{\lambda}^{\square}(-;\eta, \phi^{\square})$, via $p$-adically interpolating its Fourier coefficients by carefully choosing the Schwartz function $\phi^{\square} \in \calS(\AA_{\calK}^{2})$. Consider the embedding
\[
\iota^{\square}: \rmU(W) \times (\rmU(W_1) \times \rmU(W_2)) \hookrightarrow \rmU(\bfV),
\]
where $\rmU(W_1) \times \rmU(W_2)$ is the diagonal torus of $\rmU(-W)$. Lee computed the \emph{pullback theta lift}
\begin{equation} \label{eq:pullbackthetalift}
\Theta^{\pb}_{\lambda}(g; \eta, \nu_1, \nu_2) := \int_{[\rmU(W_1)]} \int_{[\rmU(W_2)]} \Theta^{\square}_{\lambda}(\iota^{\square}(g, t_1, t_2); \eta, \phi^{\square}) \nu_{1}(t_1) \nu_{2}(t_{2}) \dif t_1 \dif t_{2},
\end{equation}
where $\nu_i$ are characters of $\rmU(W_i)$ for $i=1,2$. It turns out to be $\Theta_{\lambda}(-; \lambda^{2} \eta \nu_1 \nu_2; \phi)$ by \cite[Theorem 2.4]{lee2025hidafamilythetalift}, up to a constant which is nonzero under the assumptions \ref{item:leei}, ($\delta 1$) and ($\delta 2$) (see \cite[Remark 3.3]{lee2025hidafamilythetalift}). Taking $\nu_1 = \nu_2 = \lambda^{-1}$ gives the desired theta lift of $\eta$ from $\rmU(V)$ to $\rmU(W)$. The pullback integral in \eqref{eq:pullbackthetalift} can be $p$-adically interpolated using Proposition \ref{prop:padicpetersson}.

\begin{remark}
The settings in this paper involve slight modifications compared to Lee's article \cite{lee2025hidafamilythetalift}. To ensure consistency with the global theta lift defined in \eqref{eq:globalthetalift}, we have replaced $\eta$ in \cite{lee2025hidafamilythetalift} with $\bareta$. Correspondingly, we have replaced the CM type $\Sigma$ selected in \cite{lee2025hidafamilythetalift} with $\Sigma^{\sfc}$. Consequently, we have updated \cite[Section 3, Assumption $(\chi3)$]{lee2025hidafamilythetalift} and restated it as Assumption \ref{item:leeiii} in this paper. Furthermore, unlike \cite{lee2025hidafamilythetalift}, we use $\psi_v(x) = \exp(- 2 \pi \sqrt{-1} x)$ as the fixed additive character instead of $\psi_v(x) = \exp( 2 \pi \sqrt{-1} x)$ defined therein; this has no impact on the conclusions used in this paper.
\end{remark}

\begin{remark}
Although we closely follow Lee's construction in this article, the technique of pullback theta lifts was known prior to Lee's work. In \cite{wan2020iwasawa, castella2022iwasawa} (in particular \cite[Section 8B]{wan2020iwasawa}), the Hida family of theta lifts over $\rmU(2,0)$ from $\rmU(1,0)$ (referred to as a \emph{CM family} there) was already constructed in the same way, so as to prove the $p$-primitive property of the Klingen Eisenstein family over $\rmU(3,1)$ by computing its Fourier--Jacobi coefficients. In \cite{lee2025hidafamilythetalift}, this is generalized to the case $(\rmU(V), \rmU(W))$ for $W$ of arbitrary signature, but it remains unclear whether the resulting $p$-adic family over $\rmU(W)$ is $p$-ordinary except when $\rmU(W) = \rmU(2,0)$.
\end{remark}

\subsubsection{Hida family of theta lifts} \label{sec:constructhidatheta}
Let $\frs$ be an ideal of $\calO_{\calK}$ that is prime to $p\frc$ and is divisible only by split primes. We sometimes identify $\frs$ with an ideal of $\calO_{\calF}$. For each finite place $v$ of $\calF$ dividing $\frs$, we fix a choice of a finite place $w$ of $\calK$ above $v$. This amounts to fixing a decomposition $\frs = \frf + \bar{\frf}$ such that $\frf+\bar{\frf} = \calO_{\calK}$. Henceforth, if $v = w \barw$ is a split place and $v \mid p \frs$, we shall assume that either $w$ is a place above $p$ or $w \mid \frf$.

Let $\Cl(\fra)$ denote the $p$-part of the ray class group over $\calK$ of an ideal $\fra \subset \calO_{\calK}$. We define
\[
\frG_{\frs} := \varprojlim_{n} \Cl(p^n \frs),
\]
which is isomorphic to the Galois group over $\calK$ of the maximal abelian pro-$p$ extension over $\calK$ that is unramified outside $p \frs$, denoted by $\calK^{(p\frs)}$. We normalize the Artin reciprocity map $\Art: \AA_{\calK, \rmf} \rightarrow \frG_{\frs}$ so that a uniformizer $\varpi_{w}$ at $w \nmid p\frs$ is sent to the geometric Frobenius $\Fr_{w}$.

Since $p^{n}\frs = \overline{p^n \frs}$, the Galois group $\Gal(\calK/\calF)$ acts on $\frG_{\frs}$, and we define the subgroups
\[
\frG_{\frs}^{\pm} := \{ \gamma \in \frG_{\frs}: \gamma^{\sfc} = \gamma^{\pm 1} \},
\]
and the quotients
\[
\frG_{\frs}^{\cyc} := \frG_{\frs}/\frG_{\frs}^{-}, \quad \frG_{\frs}^{\ac} := \frG_{\frs}/\frG_{\frs}^{+},
\]
which we call the \emph{maximal cyclotomic quotient} and \emph{maximal anticyclotomic quotient} of $\frG_{\frs}$, respectively. Since $\Art(\AA_{\calF, \rmf}^{\times}) \subset \frG_{\frs}^{+}$, the group homomorphism $\AA_{\calK, \rmf}^{\times} \rightarrow \frG_{\frs} \rightarrow \frG_{\frs}^{\ac}$ factors through $\AA_{\calK, \rmf}^{1}$ and fits into the commutative diagram
\begin{equation} \label{tikz:cft}
\begin{tikzcd}
	{\AA_{\calK,\rmf}^{\times}} && {\AA_{\calK,\rmf}^{1}} && {\AA_{\calK,\rmf}^{\times}} \\
	{\frG_{\frs}} && {\frG_{\frs}^{\ac}} && {\frG_{\frs}}
	\arrow["{z \mapsto z/\barz}", from=1-1, to=1-3]
	\arrow["\Art"', from=1-1, to=2-1]
	\arrow[hook, from=1-3, to=1-5]
	\arrow["{\Art^{\ac}}", from=1-3, to=2-3]
	\arrow["\Art", from=1-5, to=2-5]
	\arrow[from=2-1, to=2-3]
	\arrow["{\gamma \mapsto \gamma/\gamma^{\sfc}}"', from=2-3, to=2-5]
\end{tikzcd}    
\end{equation}
which identifies the image of $\frG_{\frs}^{\ac}$ with the subgroup $\frG_{\frs}^{-}$.

Given an ideal $\frs$ as above, we say that an algebraic Hecke character $\eta: \calK^{\times} \bs \AA_{\calK}^{\times} \rightarrow \CC^{\times}$ is \emph{$\frs$-admissible} if the prime-to-$p$ part of the conductor of $\eta$ divides $\frs$, and $\eta$ has infinity type $k\Sigma$ for some integer $k \geq 0$. We say $\eta$ is \emph{admissible} if it is $\frs$-admissible for some ideal $\frs$ as above.

Suppose $\eta$ is an $\frs$-admissible Hecke character. Through the global Artin map $\Art: \calK^{\times}\backslash \AA_{\calK}^{\times} \longrightarrow \frG_{\frs}$, we regard $\eta$ as a character of $\frG_{\frs}$. We define the \emph{$\frs$-ramified weight space}
\[
\calY_{\frs} := \Hom_{\cts}(\frG_{\frs}^{-},\, \bar{\QQ}_p^{\times}),
\]
i.e., the space of continuous $\bar{\QQ}_p^{\times}$-valued characters of $\frG_{\frs}^{-}$. Let $\calY_{\frs}^{\adm} \subset \calY_{\frs}$ be the subset consisting of those characters $\eta_{\sfQ}: \frG_{\frs}^{-} \to \bar{\QQ}_p^{\times}$ that arise from an $\frs$-admissible Hecke character $\eta_{\sfQ}: \calK^{\times}\backslash \AA_{\calK}^{\times} \to \CC^{\times}$ via $\Art^{-1}$, and whose infinity type is $k\Sigma$ with $k \geq 0$. We refer to the elements of $\calY_{\frs}^{\adm}$ as \emph{anticyclotomic $\frs$-admissible characters}. The subset $\calY_{\frs}^{\adm}$ is Zariski dense in $\calY_{\frs}$.

Using the technique of pullback theta lifts reviewed in Section \ref{sec:pullbackthetalift}, Lee constructed the following $p$-ordinary family of theta lifts from $\rmU(V)$ to $\rmU(W)$ in \cite[Theorem 6.8]{lee2025hidafamilythetalift}.

\begin{theorem} \label{thm:lee2025hidafamilythetaliftcmfamily}
There exists a measure
\[
\bTheta_{\lambda, \frs} \in \Meas(\frG_{\frs}^{-}, V_{\rmU(W)}^{\ord}(K^{p}_{\rmU(W)}, \calO_{L}))
\]
where $L/\QQ_p$ is a sufficiently large finite extension, such that for any anticyclotomic $\frs$-admissible Hecke character $\eta$ of infinity type $k \Sigma$ with $k \geq 0$, we have
\[
\int_{\frG_{\frs}^{-}} \eta \dif \bTheta_{\lambda, \frs} = \Theta_{\lambda}(-, \eta, \phi_{\eta}),
\]
where $\phi_{\eta} \in \calS(\XX(\AA_{\calF}))$ is the Schwartz function chosen following \cite[Section 3.1]{lee2025hidafamilythetalift}, which depends on $\eta$.
\end{theorem}

\begin{remark}
The ordinary family $\bTheta_{\lambda, \frs}$ is often referred to as an \emph{ordinary CM family} because, given the Schwartz function $\phi^{\square}$ chosen in \cite[Section 3.1]{lee2025hidafamilythetalift}, its restriction to $\GL_2(\calF)$ coincides with the classical automorphic induction (i.e., the classical CM modular form) attached to the character $\eta^{-1}$. See \cite[Remark 3.8]{lee2025hidafamilythetalift} for further details.
\end{remark}

\begin{remark} \label{rem:weighttheta}
It is worth noting that the CM family $\bTheta_{\lambda, \frs}$ is \emph{not} a Hida family over $\mathrm{U}(2,0)$. This is because it is defined only over the one-dimensional weight space $\frG_{\frs}^{-}$, rather than the full weight space of $\mathrm{U}(2,0)$. More concretely, for an $\frs$-admissible Hecke character $\eta$ of infinity type $k\Sigma$ with $k \geq 0$, Corollary \ref{cor:hcgloballift} shows that the theta lift $\Theta_{\lambda}(-,\eta,\phi_{\eta})$ has highest weight $(k,0)$ at any real place of $\calF$. However, these weights do \emph{not} form a Zariski-dense subset of the weight algebra $\Lambda_{2,0}^{\circ}$; consequently, the resulting family cannot be a Hida family on $\mathrm{U}(2,0)$ (though it remains $p$-ordinary). Nevertheless, as will be seen in the proof of Theorem \ref{thm:hsiehprimitive}, this restricted family suffices for our purposes.
\end{remark}

By the global tower property of theta lifts (Theorem \ref{thm:globaltower}), the CM ordinary family $\bTheta_{\lambda, \frs}$ is either zero or a cuspidal ordinary family over $\rmU(W)$. Thus, it is applicable to our automorphic computations in Part \ref{part:automorphic} as the auxiliary family $\bfg$ over $\rmU(2,0)$.

\subsection{Non-vanishing of Hecke $L$-functions modulo $p$} \label{sec:nonvanishingrets}
In this section, we review existing results on the non-vanishing of Hecke $L$-functions modulo $p$, which will play an essential role in our proof of the Iwasawa main conjecture for CM fields.

Let $\chi: \AA_{\calK}^{\times} \rightarrow \CC^{\times}$ be an algebraic Hecke character over $\calK$ of infinity type $k \Sigma + \kappa (1-\sfc)$ with $\kappa = \sum_{\sigma \in \Sigma} \kappa_{\sigma} \sigma \in \ZZ_{\geq 0}[\Sigma]$ and $k \Sigma + \kappa \in \ZZ_{>0}[\Sigma]$; that is, 
\[
\chi_{\infty}(z) = \prod_{\sigma \in \Sigma} \sigma(z)^{-k-\kappa_{\sigma}} (\sigma^\sfc)(z)^{\kappa_{\sigma}}.
\]
We let $L(s, \chi)$ be the associated Hecke $L$-function as in Tate's thesis \cite{tate1967fourier}, with the center of the functional equation at $s=1$. Let $\cond(\chi)$ denote the conductor of $\chi$.

Following \cite{burungale2025modellnonvanishingselfdual}, we say a Hecke character $\chi$ over $K$ is \emph{self-dual} if
\[
\chi|_{\AA^{\times}} \cdot \abs{-}_{\AA^{\times}} = \epsilon_{\calK/\calF}.
\]
We say that $\chi$ is \emph{residually self-dual} if
\[
\hatchi|_{\calF^{\times}} \omega_{\calF} \equiv \epsilon_{\calK/\calF} \pmod{\frm_{p}},
\]
where $\hatchi$ is the $p$-adic avatar of $\chi$, $\omega_{\calF}$ is the Teichmüller character over $\calF$, and $\frm_{p}$ is the maximal ideal of $\barZZ_{(p)}$. 

\begin{remark} \label{rem:burungalehsieh}
Self-dual characters are \emph{conjugate symplectic self-dual} in the sense of \cite{gan2012symplectic}. Note that the definitions here differ from those of Hsieh in \cite{hsieh2012nonvanishing, hsieh2014mu}, where self-dual characters are defined to satisfy $\chi|_{\AA^{\times}} = \epsilon_{\calK/\calF} \cdot \abs{-}_{\AA^{\times}}$, and residually self-dual characters satisfy $\hatchi|_{\calF^{\times}} \equiv \omega_{\calF} \epsilon_{\calK/\calF} \pmod{\frm_{p}}$.
\end{remark}

For a character $\chi$ that is either self-dual or not residually self-dual, we define the invariant
\[
\mu_p(\chi) := \sum_{\substack{v \mid \Norm_{\calK/\calF}(\cond(\chi)) \\ v \text{ nonsplit}, v \nmid \ell}} \mu_p(\chi_v),
\]
where $\mu_p(\chi_v) := \inf_{x \in \calK_{v}^{\times}} v_{p}(\chi_v(x)-1)$ and "nonsplit" refers to primes in $\calF$ that do not split in $\calK$. This invariant often provides a lower bound for the $p$-adic valuation of $L^{\alg}(1, \chi \nu)$ for $\nu \in \Xi_{\frl}$. We record the observations from \cite[Remark 6.3]{burungale2025modellnonvanishingselfdual} in the following proposition.

\begin{proposition} \label{prop:mupzero}
Let $\chi$ be an algebraic Hecke character. Then we have the following:
\begin{enumerate}[label = \rm (\arabic*)]
    \item If $\chi$ is self-dual, then $\mu_p(\chi) = 0$ provided that $p$ is odd and, for each prime $\frq \nmid \frl$ of $\calF$ inert in $\calK$ such that $\frq \parallel \cond(\chi)$, we have $p \nmid (\Norm_{\calF/\QQ}(\frq) - 1)$.
    \item If $\chi$ is not residually self-dual, then $\mu_p(\chi) = 0$ provided that for each prime $\frq \parallel \cond(\chi)$ of $\calF$ nonsplit in $\calK$ and coprime to $\frl$, we have $p \nmid (\Norm_{\calF/\QQ}(\frq) - 1)$.
\end{enumerate}
\end{proposition}

We remark that for such a Hecke character $\chi$ of infinity type $k\Sigma + \kappa (1-\sfc)$ with $k \in \ZZ$ and $\kappa = \sum_{\sigma \in \Sigma} \kappa_{\sigma} \sigma \in \ZZ[\Sigma]$, it is self-dual if and only if $k = 1$, and it is a unitary character if and only if $k = 0$.

\subsubsection{Twisting by $\ell$-adic anticyclotomic characters}

Let $\frl$ be a fixed prime ideal of $\calF$ unramified in $\calK$, and let $\ell$ be the rational prime below $\frl$. We denote by $\calK_{\infty}^{(\frl)}/\calK$ the anticyclotomic $\ZZ_{\ell}^{[\calF_{\frl}:\QQ_{\ell}]}$-extension of $\calK$ unramified outside $\frl$. We let $\Xi_{\frl}$ be the set of $\barQQ^{\times}$-valued finite-order characters of $\Gamma_{\frl} := \Gal(\calK_{\infty}^{(\frl)}/\calK)$, and define
\[
\Xi_{\chi, \frl}^{+} := \{ \nu \in \Xi_{\frl}: \epsilon(\chi \nu) = 1 \},
\]
where $\epsilon(-)$ denotes the global epsilon factor.

The following is a main result of \cite{hida2022icts, hsieh2012nonvanishing, he2025stability}.
\begin{theorem} \label{thm:hsieh2012nonvanishing}
Let $\chi$ be an algebraic Hecke character over a $p$-ordinary CM quadratic extension $\calK/\calF$ which is either self-dual or not residually self-dual. Let $\frl \nmid p$ be a degree-one prime of $\calF$. Suppose that:
\begin{enumerate}[label = \rm (\roman*)]
    \item $p \nmid 2 D_{\calF}$;
    \item $\frl \nmid \cond(\chi|_{\AA^{\times}})$;
    \item if $\chi$ is self-dual and $\frl$ splits in $\calK$, then $\epsilon(\chi) = 1$.
\end{enumerate}
Then there exist infinitely many $\nu \in \Xi_{\frl}$ such that $v_p(L^{\alg}(1, \chi\nu)) = \mu_p(\chi)$.
\end{theorem}

We remark that this theorem was claimed in \cite{hida2004non, hida2007non, hsieh2012nonvanishing} for all but finitely many $\nu \in \Xi_{\frl}$, which is stronger than Theorem \ref{thm:hsieh2012nonvanishing}. However, in 2018, a gap was found in these strategies; in \cite{hida2022icts}, Hida recovered the results but proved them only for infinitely many characters $\nu$.

Combining Theorem \ref{thm:hsieh2012nonvanishing} with \cite[Theorem 1.4]{burungale2025modellnonvanishingselfdual}, we obtain the following theorem, which is a special case of \cite[Theorem 6.2]{burungale2025modellnonvanishingselfdual} that fits our application.

\begin{theorem} \label{thm:burungale2025modellnonvanishingselfdualmain}
Let $\chi_1$ and $\chi_2$ be algebraic Hecke characters over a $p$-ordinary CM quadratic extension $\calK/\calF$ such that $\chi_1$ is self-dual and $\chi_2$ is not residually self-dual. Let $\frl \nmid p$ be a degree-one prime of $\calF$ which splits in $\calK$. Suppose that:
\begin{enumerate}[label = \rm (\roman*)]
    \item $p \nmid 6 D_{\calF}$;
    \item $\frl \nmid \Norm_{\calK/\calF}(\cond(\chi_1)) \cond(\chi_2|_{\AA^{\times}})$;
    \item $\epsilon(\chi_1) = 1$.
\end{enumerate}
Then there exist infinitely many $\nu \in \Xi_{\frl}$ such that
\[
v_p(L^{\alg}(1, \chi_1\nu)) = \mu_p(\chi_1) \quad \text{and} \quad v_p(L^{\alg}(1, \chi_2\nu^{-1})) = \mu_p(\chi_2)
\]
hold simultaneously.
\end{theorem}

\begin{proof}
By \cite[Theorem 1.4]{burungale2025modellnonvanishingselfdual}, we have $v_p(L^{\alg}(1, \chi_1 \nu)) = \mu_p(\chi_1)$ for \emph{all except finitely many} $\nu \in \Xi_{\chi_1, \frl}^{+}$. On the other hand, by Theorem \ref{thm:hsieh2012nonvanishing}, we have $v_p(L^{\alg}(1, \chi_2 \nu^{-1})) = \mu_p(\chi_2)$ for infinitely many $\nu \in \Xi_{\frl}$. The theorem then follows from the observations implied by \cite[Lemma 4.11]{burungale2025modellnonvanishingselfdual} and the fact that $\epsilon(\chi_1) = 1$: the set $\Xi_{\chi_1, \frl}^{+}$ is infinite, and it differs from $\Xi_{\frl}$ by only finitely many elements since $\frl$ splits in $\calK$.
\end{proof}

\subsubsection{Twisting by $p$-adic anticyclotomic characters} \label{sec:twistbypanti}
Following Section \ref{sec:constructhidatheta}, we adopt a simplified setup where $\frs$ is divisible only by split primes. In this section, we recall results on the non-vanishing of Hecke $L$-functions modulo $p$ upon twists by anticyclotomic characters $\eta: \frG_{\frs} \rightarrow \CC^{\times}$. This is often formulated as the vanishing of anticyclotomic $\mu$-invariants.

In \cite{katz1978cm, hida1994anticyclotomic}, a $\barZZ_p$-valued measure on $\frG_{\frs}$, denoted by $\calL_{\Katz}(\frs, \Sigma)$, is constructed such that
\begin{equation}\label{eq:katzpadicl}
\int_{\frG_{\frs}} \widehat{\lambda} \dif \calL_{\Katz}(\frs, \Sigma) = L^{(p\frs)}(0, \lambda) \cdot \Eul_{p}(\lambda) \Eul_{\frs}(\lambda) \cdot \dfrac{\pi^{\kappa} \Gamma_{\Sigma}(k\Sigma+\kappa)[\calO_{\calK}^{\times}:\calO_{\calF}^{\times}]}{\sqrt{\abs{D_{\calF}}_{\RR}}} \left(\dfrac{\Omega_{p}}{\Omega_{\infty}}\right)^{k\Sigma+2\kappa},    
\end{equation}
where $\lambda$ is a Hecke character of $\calK$ of modulus $\frs p^{\infty}$ and infinity type $k \Sigma + \kappa (1-\sfc)$, with either $k \geq 1$ and $\kappa \in \ZZ_{\geq 0}[\Sigma]$ or $k \geq 1$ and $k\Sigma + \kappa \in \ZZ_{>0}[\Sigma]$. Here, $\widehat{\lambda}$ is its $p$-adic avatar, and $\Eul_{p}(\lambda), \Eul_{\frs}(\lambda)$ are certain modified Euler factors (see, for example, \cite[(4.16)]{hsieh2014mu}).

We fix an algebraic Hecke character $\chi$ over $\calK$ of infinity type $k \Sigma$ with $k \geq 1$, and let $\calL_{\Katz}(\chi, \Sigma)$ (resp. $\calL_{\Katz}^{\ac}(\chi, \Sigma)$) be the $p$-adic measure on $\Gamma_{\calK}$ (resp. $\Gamma_{\calK}^{\ac}$) obtained by the pullback of $\calL_{\Katz}(\frs, \Sigma)$ along $\chi$. In other words, for every locally constant function $\varphi$ on $\Gamma_{\calK}^{\bullet}$, we have
\[
\int_{\Gamma_{\calK}^{\bullet}} \varphi \dif \calL_{\Katz}^{\bullet}(\chi, \Sigma) = \int_{\frG_{\frs}} \varphi \widehat{\chi} \dif \calL_{\Katz}(\frs, \Sigma)
\]
for $\bullet \in \{\, , \ac\}$. The $p$-adic measure $\calL_{\Katz}(\chi, \Sigma)$ (resp. $\calL_{\Katz}^{\ac}(\chi, \Sigma)$) is called the \emph{Katz $p$-adic $L$-function} (resp. \emph{anticyclotomic Katz $p$-adic $L$-function}) attached to the branch character $\chi$. We denote by $\mu^{\bullet}(\chi, \Sigma)$ the $\mu$-invariant of $\calL_{\Katz}^{\bullet}(\chi, \Sigma)$ for $\bullet \in \{\, , \ac\}$, omitting the CM type $\Sigma$ if it is clear from the context. 

\begin{theorem} \label{thm:hsieh2014mu}
Let $\chi$ be an algebraic Hecke character over a $p$-ordinary CM quadratic extension $\calK/\calF$. Suppose that $p \nmid 2D_{\calF}$; then we have:
\begin{enumerate}[label = \rm (\arabic*)]
    \item $\mu(\chi, \Sigma) = 0$;
    \item if, moreover, $\chi$ is not residually self-dual in the sense of Hsieh (see Remark \ref{rem:burungalehsieh}), then $\mu^{\ac}(\chi, \Sigma) = 0$.
\end{enumerate}
\end{theorem}
\begin{proof}
The result (1) follows from \cite[Theorem A]{burungale2013vanishing}, and (2) follows from \cite[Theorem B]{hsieh2014mu}, noting that we have adopted a simplified setup where $\frs$ is divisible only by split primes, implying that $\frs^{-}$ is trivial.
\end{proof}

\begin{remark}
Note that (1) does \emph{not} imply (2). For example, consider $f(S,T) := ST+pS+pT \in \ZZ_p\lrbracket{S,T}$ as a two-variable power series over $\ZZ_p$, with $S$ regarded as the cyclotomic variable and $T$ as the anticyclotomic variable. Then the $\mu$-invariant of $f$ is zero, but its images under the ``modulo $S$'' map and the ``modulo $T$'' map both have $\mu$-invariant $1$.
\end{remark}

\subsection{Iwasawa main conjecture for CM fields}
With the preparations above, we are now ready to prove the Iwasawa main conjecture for CM fields. We consider the following objects:
\begin{itemize}
    \item Let $\calK_{\infty}/\calK$ be the compositum of the cyclotomic $\ZZ_p$-extension $\calK_{\infty}^{\cyc}/\calK$ and the anticyclotomic $\ZZ_p^{[\calF:\QQ]}$-extension $\calK_{\infty}^{\ac}$. We set $\Gamma_{\calK}^{\bullet} := \Gal(\calK_{\infty}^{\bullet}/\calK)$ for $\bullet \in \{ \, , \cyc, \ac \}$.
    \item Let $\calK^{\prime}$ be a finite extension of $\calK$ containing $\calK(\zeta_p)$ that is disjoint from $\calK_{\infty}$. Put $\Delta := \Gal(\calK^{\prime}/\calK)$ and $\calK_{\infty}^{\prime} = \calK_{\infty} \calK^{\prime}$.
    \item Let $\psi: \Delta \rightarrow \barZZ_{p}^{\times}$ be a finite-order Hecke character over $\calK$, which is automatically a unitary Hecke character. It is called a \emph{branch character}. We define $\calR := \barZZ_{p}^{\ur}[\psi]$ and $\Lambda := \calR\lrbracket{\Gamma_{\calK}}$.
\end{itemize}
Let $\calM_{\Sigma}$ be the maximal $p$-abelian $\Sigma_p$-ramified extension of $\calK_{\infty}^{\prime}$. The fundamental object is 
\[
X_{\Sigma} := \Gal(\calM_{\Sigma}/\calK_{\infty}^{\prime}) \otimes_{\ZZ_p[\Delta]\lrbracket{\Gamma_{\calK}}} \calR[\Delta]\lrbracket{\Gamma_{\calK}}.
\]
We denote by $X(\psi, \Sigma)$ the maximal $\psi$-isotypic quotient of $X_{\Sigma}$, which is a finitely generated torsion $\Lambda$-module, and let $F(\psi, \Sigma) \in \Lambda$ be its characteristic power series.

On the analytic side, we have the associated Katz $p$-adic $L$-function $\calL_{\Katz}(\psi, \Sigma) \in \Lambda$ as recalled in Section \ref{sec:twistbypanti}. We then have the following conjecture, known as the \emph{Iwasawa main conjecture for the CM field} $\calK/\calF$.

\begin{conjecture}[Iwasawa main conjecture for the CM field $\calK/\calF$]
With notation as above, the following equality holds as ideals of $\Lambda$:
\[
(F(\psi, \Sigma)) = (\calL_{\Katz}(\psi, \Sigma)).
\]
\end{conjecture}

To study this Iwasawa main conjecture, we consider the $p$-adic Eisenstein datum
\[
\bfD_{\Hs} := (L, \II := \barZZ_p[\psi], \psi^{-1}, 1, \bfone, 1)
\]
over $\rmU(1,0)$ in the sense of Definition \ref{defn:padiceisendatum}, following the construction of Klingen Eisenstein series over $\rmU(2,1)$ in \cite{hsieh2014eisenstein}\footnote{Note that in \cite{hsieh2014eisenstein}, the Siegel Eisenstein series is defined using Siegel Eisenstein sections in $I^{\Sieg}_{Q}(\chi^{-1}, s)$, while in \cite{wan2015families}, it is defined using Siegel Eisenstein sections in $I^{\Sieg}_{Q}(\chi, s)$. Therefore, to be consistent with the construction of the Klingen Eisenstein family in \cite{hsieh2014eisenstein}, we use $\psi^{-1}$ instead of $\psi$ in the $p$-adic Eisenstein datum $\bfD_{\Hs}$.}. Here, $L$ is a finite extension of $\QQ_p$, sufficiently large to contain all values of $\psi$, and $\bfone$ is the trivial Hida family over $\rmU(1,0)$ with the trivial character $1$ as its branching character. Then the weight algebra $\Lambda_{\bfD_{\Hs}} := \II^{\ur}\lrbracket{\Gamma_{\calK}}$ is precisely $\calR\lrbracket{\Gamma_{\calK}}$, and the characters $\xi_0 := \psi$ and $\bxi_0 := \bPsi = \psi \bvep_{\calK}: \Gal(\calK^{\Sigma}/\calK) \rightarrow \Lambda^{\times}$ are the $\Lambda$-valued Galois characters of $\calK$.

In this article, we prove the following result regarding the Iwasawa main conjecture for CM fields.

\begin{theorem}[Iwasawa main conjecture for the CM field $\calK/\calF$] \label{thm:myimccm}
Suppose that $p$ does not divide $6 D_{\calF} \abs{\Delta}$. Then we have the inclusion $(F(\psi, \Sigma)) \subseteq (\calL_{\Katz}(\psi, \Sigma))$ as ideals of $\Lambda$.
\end{theorem}

Compared with \cite[Theorem 2]{hsieh2014eisenstein}, we have been able to remove the following two conditions:
\begin{enumerate}[label = (Hs-\roman*)]
    \item $p \nmid h_{\calK}^{-} \abs{\Delta}$, where $h_{\calK}^{-} := h_{\calK}/h_{\calF}$ is the relative class number for $\calK/\calF$; \label{item:hsi}
    \item $\psi$ is unramified at $\Sigma_{p}^{\sfc}$, and $\psi \omega_{\calK}^{a}$ is unramified at $\Sigma_{p}$ for some integer $a \not\equiv 2 \pmod{p-1}$. \label{item:hsii}
\end{enumerate}
Conditions \ref{item:hsi} and \ref{item:hsii} in Hsieh's work were required to apply Hida's modulo $p$ non-vanishing theorem for Hecke $L$-functions. Our computation of Gan--Gross--Prasad type period integrals provides a new route to the $p$-primitive property of $\bfD_{\Hs}$; in this approach, Theorem \ref{thm:burungale2025modellnonvanishingselfdualmain} is applicable. Condition \ref{item:hsii} was also used in the lattice construction in \cite{hsieh2014eisenstein} (see in particular the proof of \cite[Theorem~8.14]{hsieh2014eisenstein}). In our argument, this hypothesis is removed by appealing to the lattice construction developed in \cite{wan2019iwasawatheorymathrmursblochkato}, specifically the improvement in \cite[Proposition~7.4]{wan2019iwasawatheorymathrmursblochkato}, as explained in Section \ref{sec:lattice}.

We devote the remainder of this section to the proof of Theorem \ref{thm:myimccm}.

\subsubsection{Step 1: Modified setup and the choice of the auxiliary Hida family}
We consider the $p$-adic Eisenstein datum $\bfD_{\Hs}$, to which is attached a Klingen Eisenstein family $\bfE_{\bfD_{\Hs}}^{\Kling}$ over $\rmU(2,1)$. Recall that in Part \ref{part:automorphic}, we analyzed the Gan--Gross--Prasad type period integral of $\bfE_{\bfD_{\Hs}}^{\Kling}$ against an auxiliary cuspidal Hida family $\bfg$ over $\rmU(2,0)$, denoted by $\bfP_\bfG^2$ as before. 

However, we noted in Remark \ref{rem:weighttheta} that the ordinary CM family $\bTheta_{\lambda,\frs}$ constructed in Theorem \ref{thm:lee2025hidafamilythetaliftcmfamily} is not a full Hida family over the weight space $\Lambda_{2,0}$, but rather a $p$-adic family of $p$-ordinary modular forms over $\calO_{L}\lrbracket{\frG_{\frs}^{-}}$. Nevertheless, we consider the restricted weight space $\Lambda_{\bfD_{\Hs}}\lrbracket{\frG_{\frs}^{-}}$, which has specializations $\sfP \times \sfQ \in \calX_{\bfD_{\Hs}} \times \calY_{\frs}$ with a Zariski-dense subset $\calX_{\bfD_{\Hs}}^{\gen} \times \calY_{\frs}^{\adm}$. Modifying the construction in Section \ref{sec:padicggp}, we can define
\[
\bfP^{2}_{\Hs} \in \Lambda_{\bfD_{\Hs}}\lrbracket{\frG_{\frs}^{-}}
\]
by setting $\bfg$ to be the CM ordinary family $\bTheta_{\lambda,\frs}$. Then $\bfP^2_{\Hs}$ satisfies the same interpolation formula as Theorem \ref{thm:interpoformula} for the specializations
\[
\sfP \times \sfQ \in \calX_{\bfD_{\Hs}}^{\gen} \times \calY_{\frs}^{\adm} \subseteq \calX_{\bfD_{\Hs}} \times \calY_{\frs}.
\]
Indeed, we note the following:
\begin{itemize}
    \item By Remark \ref{rem:weighttheta} (following Corollary \ref{cor:hcgloballift}), for any arithmetic point $\sfQ \in \calY_{\frs}^{\adm}$ corresponding to an $\frs$-admissible Hecke character $\eta_{\sfQ}$ of infinity type $k\Sigma$ with $k\geq 0$, the specialization $\bTheta_{\lambda,\frs}(\eta_{\sfQ})$ has highest weight $(k,0)$ at every real place of $\calF$. Consequently, for any pair $\sfP \times \sfQ \in \calX_{\bfD_{\Hs}}^{\gen} \times \calY_{\frs}^{\adm}$, the associated automorphic representations $\sigma_{\bff_{\sfP}}$ and $\pi_{\bfg_{\sfQ}}$ satisfy the weight-interlacing condition in Assumption \ref{ass:wtinterlacing} automatically, since $\bff$ is the trivial Hida family over $\rmU(1,0)$.
    \item The non-zero ideal $\frs$ of $\calO_{\calK}$ is chosen to satisfy the conditions of Section \ref{sec:constructhidatheta} and will be specified later. Since $\bff$ is the trivial Hida family in $\bfD_{\Hs}$, the automorphic assumption (Assumption \ref{ass:automorphic}) is automatically satisfied.
\end{itemize}
We emphasize that although not explicitly included in the notation, $\bfP_{\Hs}^{2}$ depends on the choice of the splitting character $\lambda$ and the ideal $\frs$, which will be selected later in the proof.

Therefore, applying Theorem \ref{thm:interpoformula} at all interpolation points 
\[
\sfP \times \sfQ \in \calX_{\bfD_{\Hs}}^{\gen} \times \calY_{\frs}^{\adm} \subsetneq \calX_{\bfD_{\Hs}} \times \calY_{\frs},
\]
we obtain
\begin{multline*}
(\sfP \times \sfQ)(\bfP^2_{\Hs}) = C \cdot \lrangle{1,1}_{\Pet} \lrangle{\bTheta_{\lambda, \frs}(\eta_{\sfQ}), \bTheta_{\lambda, \frs}(\eta_{\sfQ})^{\dagger}}_{\Pet} \\ \times 
\dfrac{L^{\alg}\left( \dfrac{1}{2}, \bTheta_{\lambda, \frs}(\eta_{\sfQ}) \right)}{L^{\alg}(1, \Ad) L^{\alg}(1, \bTheta_{\lambda, \frs}(\eta_{\sfQ}), \Ad)} L^{\alg}\left( s_{\sfP}+\dfrac{1}{2}, \bTheta_{\lambda^{-1}, \frs}(\eta_{\sfQ}^{-1}) \times \barchi_{\sfP}^{\sfc} \right)^2.
\end{multline*}
Here $C \in (L^{\ur})^{\times}$ denotes the extra factor appearing in the interpolation formula of Theorem \ref{thm:interpoformula}; it is a constant independent of the specialization point $\sfP \times \sfQ$ \footnote{It remains to be verified that the local integrals at primes above $p$ differ from the local Euler factors by at most a fixed $p$-adic unit. This was established satisfactorily through careful and elaborate calculations in \cite[Section 7]{castella2022iwasawa}.}. We then apply Theorem \ref{thm:localLfactor} to see that
\begin{multline*}
(\sfP \times \sfQ)(\bfP^2_{\Hs}) = C \cdot \lrangle{1,1}_{\Pet} \lrangle{\bTheta_{\lambda, \frs}(\eta_{\sfQ}), \bTheta_{\lambda, \frs}(\eta_{\sfQ})^{\dagger}}_{\Pet} \\
\times \dfrac{L^{\alg}\left( \dfrac{1}{2}, \lambda \eta_{\sfQ} \right) L^{\alg}\left( \dfrac{1}{2}, \lambda \right)}{L^{\alg}(1, \epsilon_{\calK/\calF})^3 L^{\alg}(1, \eta_{\sfQ})} L^{\alg}\left( s_{\sfP}+\dfrac{1}{2}, \barchi_{\sfP}^{\sfc} \lambda^{-1} \eta_{\sfQ}^{-1} \right)^2 L^{\alg}\left( s_{\sfP}+\dfrac{1}{2}, \barchi_{\sfP}^{\sfc}\lambda^{-1} \right)^2.
\end{multline*}
We write $\chi_{\sfP}^{\circ} := \barchi_{\sfP}^{\sfc} \abs{-}_{\AA_{\calK}^{\times}}^{\kappa_{\sfP}/2}$, which has infinity type $\kappa_{\sfP} \Sigma$ and is no longer a unitary character\footnote{Note that $\barchi_{\sfP}^{\sfc}$ has infinity type $\frac{\kappa_{\sfP}}{2} (\Sigma - \Sigma^{\sfc})$.}. We first fix a specialization $\sfP_{0} \in \calX_{\bfD_{\Hs}}^{\gen}$ with $\kappa_{\sfP} = 2$, so that
\begin{multline} \label{eq:p0qp}
(\sfP_0 \times \sfQ)(\bfP^2_{\Hs}) = C \cdot \lrangle{1,1}_{\Pet} \lrangle{\bTheta_{\lambda, \frs}(\eta_{\sfQ}), \bTheta_{\lambda, \frs}(\eta_{\sfQ})^{\dagger}}_{\Pet} \\
\times \dfrac{L^{\alg}\left(1, \lambda \abs{-}_{\AA^{\times}_{\calK}}^{1/2} \eta_{\sfQ} \right) L^{\alg}\left(1, \lambda \abs{-}_{\AA^{\times}_{\calK}}^{1/2}\right)}{L^{\alg}(1, \epsilon_{\calK/\calF})^3 L^{\alg}(1, \eta_{\sfQ})} L^{\alg}\left( 1, \chi_{\sfP_{0}}^{\circ} \eta_{\sfQ}^{-1} \lambda^{-1} \abs{-}_{\AA^{\times}_{\calK}}^{-1/2}\right)^2 L^{\alg}\left( 1, \chi_{\sfP_{0}}^{\circ} \lambda^{-1} \abs{-}_{\AA^{\times}_{\calK}}^{-1/2}\right)^2.
\end{multline}

\subsubsection{Step 2: weak $p$-primitivity of $\bfD_{\Hs}$} 
Building on the preparation in Step 1, we now prove Theorem \ref{thm:hsiehprimitive} by applying the nonvanishing results from Section \ref{sec:nonvanishingrets}. Throughout this second step, we assume that $p \nmid 6 D_{\calF}$.

Following the approach in \cite[Section 7.9, Section 7.11]{castella2022iwasawa}, we consider the following additional $p$-adic $L$-functions (whose existence is guaranteed by the existence of Katz $p$-adic $L$-functions):
\begin{itemize}
    \item $\calL_1 \in \calO_{L}^{\ur}\lrbracket{\frG_{\frs}^{-}}$: the Katz anticyclotomic $p$-adic $L$-function interpolating $L^{\alg}\left(1, \lambda \abs{-}_{\AA^{\times}_{\calK}}^{1/2} \eta_{\sfQ} \right)$ for admissible points $\sfQ \in \calY_{\frs}^{\adm}$.
    \item $\calL_2 \in L^{\ur}$: the value $L^{\alg}\left(1, \lambda \abs{-}_{\AA^{\times}_{\calK}}^{1/2}\right)$.
    \item $\calL_5 \in \calO_{L}^{\ur}\lrbracket{\frG_{\frs}^{-}}$: the Katz anticyclotomic $p$-adic $L$-function interpolating $L^{\alg}\left( 1, \chi_{\sfP_{0}}^{\circ} \eta_{\sfQ}^{-1} \lambda^{-1} \abs{-}_{\AA^{\times}_{\calK}}^{-1/2}\right)$ for admissible points $\sfQ \in \calY_{\frs}^{\adm}$.
    \item $\calL_6 \in L^{\ur}$: the value $L^{\alg}\left( 1, \chi_{\sfP_{0}}^{\circ} \lambda^{-1} \abs{-}_{\AA^{\times}_{\calK}}^{-1/2}\right)$.
    \item $\calL_{3} \in L^{\ur}$: the value $L^{\alg}(1, \epsilon_{\calK/\calF})$.
    \item $\calL_4 \in \calO_{L}^{\ur}\lrbracket{\frG_{\frs}^{-}}$: the Katz anticyclotomic $p$-adic $L$-function interpolating $L^{\alg}(1, \epsilon_{\calK/\calF})^2 L^{\alg}(1, \eta_{\sfQ})$ for admissible points $\sfQ \in \calY_{\frs}^{\adm}$.
\end{itemize}
According to \eqref{eq:p0qp} and \cite[Proposition 7.9.1]{castella2022iwasawa} (where we set the $\calL_4$ in the citation to be the $\calL_4$ above), we have
\[
(\sfP_{0} \times -) \bfP^2_{\Hs} = C^{\prime} \cdot \calL_1 \calL_2 \calL_5^2 \calL_6^{2} \in \calO_{L}^{\ur}\lrbracket{\frG_{\frs}^{-}} \otimes L
\]
where $C^{\prime} \in (L^{\ur})^{\times}$ is a constant independent of the specialization point $\sfQ \in \calY_{\frs}$ (arising from the constant $C \in (L^{\ur})^{\times}$ in Step 1, the comparison with the above interpolation formulas, and the cited \cite[Proposition 7.9.1]{castella2022iwasawa}). We set
\[
\calA(\lambda) := \calL_2 \calL_6^2 \in L^{\ur}, \quad \calB(\lambda) := \calL_1 \calL_5^2 \in \calO_{L}^{\ur}\lrbracket{\frG_{\frs}^{-}}.
\]
In a similar manner, we further have
\[
\widetilde{\calA}(\lambda) \in \calO_{L}^{\ur}\lrbracket{\Gamma_{\calK}} \otimes L^{\ur}, \quad \widetilde{\calB}(\lambda) \in \Lambda_{\bfD_{\Hs}}\lrbracket{\frG_{\frs}^{-}}, \quad \calC \in \Frac(\II^{\ur})
\]
satisfying 
\[
\sfP_{0}(\widetilde{\calA}(\lambda)) = \calA(\lambda), \quad \sfP_{0}(\widetilde{\calB}(\lambda)) = \calB(\lambda), \quad \sfP_{0}(\calC) = C^{\prime},
\]
where, for $\calC \in \Frac(\II^{\ur})$, one may refer to the explanation of the constant $\calC_{1}$ at \cite[the bottom of page 75 of Section 7.9]{castella2022iwasawa}. We emphasize that all the above terms depend on the choice of the splitting character $\lambda$ and the auxiliary ideal $\frs_{0}$.

Our strategy is as follows: first, choose a suitable $\lambda_{0}$ such that $\calA(\lambda_{0})$ is a $p$-adic unit; second, after selecting a specific ideal $\frs$ and a specific specialization point $\sfQ_{0} \in \calY_{\frs}^{\adm}$, prove that $\sfQ_{0}(\calB(\lambda_{0}))$ is also a $p$-adic unit. This will suffice to show that $\widetilde{\calA}(\lambda_{0}) \widetilde{\calB}(\lambda_{0})$ is a unit in $\Lambda_{\bfD_{\Hs}}\lrbracket{\frG_{\frs}^{-}}$.

For $\calA(\lambda)$, we follow \cite[Section 7.1]{hsieh2014eisenstein} and fix a Hecke character $\lambda_{\can}: \AA_{\calK}^{\times} \rightarrow \CC^{\times}$ in the set $\calX_{0}(\fra, m)$, with properties:
\begin{enumerate}[label = (\alph*)]
    \item $\lambda_{\can}|_{\AA_{\calF}^{\times}} = \epsilon_{\calK/\calF}$;
    \item $\lambda_{\can,\infty}(z_{\infty}) = \prod_{\sigma \in \Sigma} \abs{z_{\sigma}} z_{\sigma}^{-1}$;
    \item $\epsilon(\lambda_{\can} \abs{-}_{\AA^{\times}_{\calK}}^{-1/2}) = 1$.
\end{enumerate}
By (a), $\lambda_{\can}^{-1} \abs{-}_{\AA^{\times}_{\calK}}^{1/2}$ is self-dual and $\chi_{\sfP_{0}}^{\circ} \lambda_{\can} \abs{-}_{\AA^{\times}_{\calK}}^{-1/2}$ is not residually self-dual. By (a), (b), and \cite[Proposition 1]{rohrlich1982root}, $\lambda_{\can}^{-1} \abs{-}_{\AA^{\times}_{\calK}}^{1/2}$ is ramified only at primes where $\calK/\calF$ is ramified, and $\chi_{\sfP_{0}}^{\circ} \lambda_{\can} \abs{-}_{\AA^{\times}_{\calK}}^{-1/2}$ is ramified only at such primes and those above $p$. By Proposition \ref{prop:mupzero}, we have
\[
\mu_p \left(\lambda_{\can}^{-1} \abs{-}_{\AA^{\times}_{\calK}}^{1/2} \right) = \mu_p \left( \chi_{\sfP_{0}}^{\circ} \lambda_{\can} \abs{-}_{\AA^{\times}_{\calK}}^{-1/2} \right) = 0.
\]
Choosing an auxiliary degree-one prime $\frl$ of $\calF$ split in $\calK$ not lying above $p$, Theorem \ref{thm:burungale2025modellnonvanishingselfdualmain} implies there exist infinitely many $\nu \in \Xi_{\frl}$ such that
\[
v_p \left(L^{\alg}\left(1, \lambda_{\can}^{-1} \abs{-}_{\AA^{\times}_{\calK}}^{1/2} \nu \right) L^{\alg}\left(1, \chi_{\sfP_{0}}^{\circ} \lambda_{\can} \abs{-}_{\AA^{\times}_{\calK}}^{-1/2} \nu^{-1} \right) \right) = 0.
\]
We set $\lambda_{0} := \lambda_{\can}^{-1} \nu$ for such a $\nu$ satisfying $\epsilon(\lambda_{\can}^{-1} \abs{-}_{\AA^{\times}_{\calK}}^{1/2} \nu) = 1$. Therefore, $\calA(\lambda_0) \in L^{\ur}$ is a $p$-adic unit. In particular, noting that the denominator of $\widetilde{\calA}(\lambda)$ is introduced solely by the factor $\calL_{2}$, this implies that $\widetilde{\calA}(\lambda) \in \calO_{L}^{\ur}\lrbracket{\Gamma_{\calK}}$.

In fact, we need to verify that the chosen splitting character $\lambda_{0}$ satisfies Assumptions \ref{item:leei} through \ref{item:leeiii} required in the construction of the CM ordinary Hida family $\bTheta_{\lambda_0, \frs}$. First, by the choice $\lambda_{0} := \lambda_{\can}^{-1} \nu$, it is clear that \ref{item:leei} is automatically satisfied (indeed, the algebraic part of the $L$-value $L(\frac{1}{2}, \lambda_{0})$ is non-zero modulo $p$). Here we have applied Theorem \ref{thm:burungale2025modellnonvanishingselfdualmain} to include an $\ell$-adic anticyclotomic twist $\nu \in \Xi_{\frl}$ (for $\ell \neq p$, rather than the $p$-adic anticyclotomic twist provided by Theorem \ref{thm:hsieh2014mu}), so that the splitting character $\lambda_{0}$ satisfies condition \ref{item:leeii}. As for condition \ref{item:leeiii}, one can verify that $\lambda_{0} \abs{-}_{\AA_{\calK}^{\times}}^{1/2}$ has infinity type $\Sigma$ (noting that $\nu$ is a character of finite order). The above verification ensures that using $\lambda_{0}$ as the splitting character in Theorem \ref{thm:lee2025hidafamilythetaliftcmfamily} is feasible.

%\footnote{We remark that here we apply Theorem \ref{thm:burungale2025modellnonvanishingselfdualmain} to incorporate an $\ell$-adic anticyclotomic twist $\nu \in \Xi_{\frl}$ (for $\ell \neq p$), rather than the $p$-adic anticyclotomic twist provided by Theorem \ref{thm:hsieh2014mu}, to ensure that the splitting character $\lambda$ satisfies condition \ref{item:leeii}. With our choice of $\lambda := \lambda_{\can}^{-1} \nu$, one readily verifies that conditions \ref{item:leei} and \ref{item:leeiii} are also satisfied. Specifically, since $\nu$ is a character of finite order, $\lambda \abs{-}_{\AA_{\calK}^{\times}}^{1/2}$ possesses the algebraic infinity type $\Sigma$. To ensure compatibility with the construction in \cite{lee2025hidafamilythetalift} (particularly for the verification of \ref{item:leeiii}), we systematically transpose the roles of $\Sigma$ and $\Sigma^{\sfc}$ as they appear in Lee's original work; this modification is purely notational and does not affect the validity of our application. Consequently, $\lambda$ serves as a valid splitting character for the application of Theorem \ref{thm:lee2025hidafamilythetaliftcmfamily}.}

Next, for $\calB(\lambda_0)$, since $\nu \in \Xi_{\frl}$ is a finite-order anticyclotomic character, $\lambda_{\can}^{-1} \nu \abs{-}_{\AA^{\times}_{\calK}}^{1/2}$ remains self-dual and $\chi_{\sfP_{0}}^{\circ} \lambda_{\can} \nu^{-1} \abs{-}_{\AA^{\times}_{\calK}}^{-1/2}$ is not residually self-dual. As before,
\[
\mu_p \left(\lambda_{\can}^{-1} \nu \abs{-}_{\AA^{\times}_{\calK}}^{1/2} \right) = \mu_p \left(\chi_{\sfP_{0}}^{\circ} \lambda_{\can} \nu^{-1} \abs{-}_{\AA^{\times}_{\calK}}^{-1/2} \right) = 0.
\]
Since $\epsilon(\lambda_{\can}^{-1} \nu) = 1$, we choose an auxiliary degree-one prime $\frs_{0}$ of $\calF$ split in $\calK$ ($ \frs_0 \neq \frl, p$). Let $\frs := \frs_{0} \calO_{\calK}$. Theorem \ref{thm:burungale2025modellnonvanishingselfdualmain} implies that
\[
v_p \left(L^{\alg}\left(1, \eta_{\sfQ_0} \lambda_{\can}^{-1} \nu \abs{-}_{\AA^{\times}_{\calK}}^{1/2} \right) L^{\alg}\left( 1, \chi_{\sfP_{0}}^{\circ} \eta_{\sfQ_0}^{-1} \lambda_{\can} \nu^{-1} \abs{-}_{\AA^{\times}_{\calK}}^{-1/2} \right) \right) = 0
\]
for some $\eta_{\sfQ_0} \in \Xi_{\frs_0}$. This implies that $\sfQ_0(\calB(\lambda_0)) \in L^{\ur}$ is a $p$-adic unit.

Combining the two parts of the discussion above, we see that for the chosen splitting character $\lambda_{0}$ and ideal $\frs_{0}$, the specialization $(\sfP_{0} \times \sfQ_{0}) (\widetilde{\calA}(\lambda_0) \widetilde{\calB}(\lambda_{0}))$ is a $p$-adic unit. Crucially, we note that $\widetilde{\calA}(\lambda_0) \widetilde{\calB}(\lambda_{0}) \in \Lambda_{\bfD_{\Hs}}\lrbracket{\frG_{\frs}^{-1}}$, and therefore\footnote{We employ the following general statement: consider a multi-variable power series ring $\Lambda = \calO \lrbracket{T_1, \ldots, T_n}$ with coefficients in a discrete valuation ring $\calO$, whose unique maximal ideal is $\frm_{\Lambda} = (\varpi, T_1, \ldots, T_n)$. Here $\frm_{\calO}$ is the maximal ideal of $\calO$ and $\varpi$ is its uniformizer. Let $\phi: \Lambda \rightarrow \calO$ be a ring homomorphism (i.e., a specialization of $\Lambda$) such that $\alpha_i := \phi(T_i) \in \frm_{\calO}$ for all $1 \leq i \leq n$. For $f \in \Lambda$, we assert that if $f(\alpha_1, \ldots, \alpha_{n}) \in \calO^{\times}$, then $f \in \Lambda^{\times}$. Otherwise, $f \in \frm_{\Lambda}$, so it can be written as
\[
f = \varpi g + h_1 T_{1} + \cdots + h_n T_{n}, \quad g, h_1, \ldots, h_n \in \Lambda.
\]
Then,
\[
f(\alpha_1, \ldots, \alpha_n) = \varpi g(\alpha_1, \ldots, \alpha_n) + \sum_{i=1}^{n} \alpha_i h_i(\alpha_1, \ldots, \alpha_n).
\]
Since each $\alpha_i \in \frm_{\calO}$, this implies $f(\alpha_1, \ldots, \alpha_n) \in \frm_{\calO}$, which contradicts our assumption. The arithmetic specializations $\calX_{\bfD}^{\ari}$ in Definition \ref{defn:specializations} and the admissible specializations on the weight space $\calY_{\frs}$ in this chapter both satisfy the preconditions of the above argument (cf. also the explicit definition of arithmetic specialization points in \cite[Definition 5.3]{wan2015families}). The argument in this footnote is precisely the process intended for the reader in classical literature such as \cite[Proposition 7.11.3]{castella2022iwasawa}; the most critical point is that the element $\widetilde{\calA}(\lambda_0) \widetilde{\calB}(\lambda_{0})$ lies in $\Lambda_{\bfD_{\Hs}}\lrbracket{\frG_{\frs}^{-1}}$ rather than merely in $\Lambda_{\bfD_{\Hs}}\lrbracket{\frG_{\frs}^{-1}} \otimes \Frac(\II^{\ur})$.} $\widetilde{\calA}(\lambda_0) \widetilde{\calB}(\lambda_{0})$ is a unit in $\Lambda_{\bfD_{\Hs}}\lrbracket{\frG_{\frs}^{-1}}$. In particular, it is not divisible by any prime ideal $\frP$ of height 1 in $\Lambda_{\bfD_{\Hs}}\lrbracket{\frG_{\frs}^{-}}$. Returning to
\[
\bfP_{\Hs}^{2} = \calC \cdot \widetilde{\calA}(\lambda_0) \widetilde{\calB}(\lambda_{0}),
\]
due to the ambiguity introduced by $\calC \in \Frac(\II^{\ur})$, we can only conclude that $\bfP_{\Hs}^{2}$ is not divisible by any prime ideal $\frP$ of height 1 in $\Lambda_{\bfD_{\Hs}}\lrbracket{\frG_{\frs}^{-}}$ such that $\frP \cap \II^{\ur} = (0)$.

The above argument suffices to show that $\bfD_{\Hs}$ is weak $p$-primitive. More explicitly, we may enlarge the weight space $\Lambda_{\bfD_{\Hs}}$ in Definition \ref{defn:padiceisendatum} to $\Lambda_{\bfD_{\Hs}, \frs} := \Lambda_{\bfD_{\Hs}}\lrbracket{\frG_{\frs}^{-}}$, and then define the $\Lambda_{\bfD_{\Hs}, \frs}$-linear functional
    \begin{align*}
    \ell_{\bfD^{\Hs}}: \calM_{\bfD_{\Hs}} \otimes_{\Lambda_{\bfD_{\Hs}}} \Lambda_{\bfD_{\Hs}, \frs} \otimes \Frac(\II^{\ur}) &\rightarrow \Lambda_{\bfD_{\Hs}, \frs} \otimes \Frac(\II^{\ur}); \\  \bff &\mapsto \bfB_{K^{p}}(\bff, \bTheta_{\lambda_0,\frs_0})^{2},  
    \end{align*}
where the data $\lambda_{0}, \frs_{0}$, etc., are chosen according to the above proof, and a tame level group $K^{p} \subset G^{\hsuit}(\AA_{\calK, \rmf}^{p})$ is fixed throughout. Thus, $\ell_{\bfD^{\Hs}}(\bfE_{\bfD_{\Hs}}^{\Kling}) = \bfP_{\Hs}^{2}$. The preceding discussion proves that $\ell_{\bfD^{\Hs}}$ is a $\Lambda_{\bfD_{\Hs}, \frs}$-linear weak $p$-primitive functional in the sense of Definition \ref{defn:primitivefunctional}, provided we note that its $\Lambda_{\bfD_{\Hs}, \frs}$-linearity originates from the $\Lambda_{\bfD_{\Hs}}$-linearity of the pairing $\bfB_{K^{p}}(-,-)$ in the first component as stated in Proposition \ref{prop:padicpetersson}. Although we have enlarged the weight space, the same arguments in Section \ref{sec:genprimitivity} proceed by replacing $\Lambda_{\bfD_{\Hs}}$ with $\Lambda_{\bfD_{\Hs}, \frs}$. \footnote{This is consistent with the interpretation used in \cite{skinner2014iwasawa}, where the weight algebra is taken to be $\II\lrbracket{\Gamma_{\calK} \times \Gamma_{\calK}^{-}}$ enlarged from $\II\lrbracket{\Gamma_{\calK}}$ (see \cite[the definition on page 116]{skinner2014iwasawa}).}

This constitutes a key milestone of the present article, which we summarize in the following theorem.

\begin{theorem} \label{thm:hsiehprimitive}
Suppose that $p \nmid 6 D_{\calF}$. Then the $p$-adic Eisenstein datum $\bfD_{\Hs}$ is weak $p$-primitive (in the sense of Definition \ref{defn:primitivefunctional}).
\end{theorem}

\subsubsection{Step 3: The proof of Theorem \ref{thm:myimccm}}
With the weak $p$-primitive property of $\bfD_{\Hs}$ established, we are now ready to prove Theorem \ref{thm:myimccm} by applying Theorem \ref{thm:imcdeprived}. Let $R$ be a commutative Iwasawa algebra over $\calO_{\frp}$ and $X$ be a finitely generated $R$-module. For each height-one prime $\frP$ of $R$, we write
\[
\ell_{\frP}(X) := \length_{R_{\frP}}(X \otimes_{R} R_{\frP}).
\]
Motivated by Theorem \ref{thm:imcdeprived}, we first compare the formulation of Theorem \ref{thm:myimccm} with it, following the proof of \cite[Theorem 8.14]{hsieh2014eisenstein}.

On the algebraic side, we have the following comparison:
\begin{align*}
\ell_{\frP}(X_{\Gr}^{\scS}(\bfone/\calK_{\infty}, &\xi_0^{-\sfc} \epsilon_{\cyc}, \Sigma)) \\ &\leq \ell_{\frP}(X_{\Gr}(\bfone/\calK_{\infty}, \xi_0^{-\sfc} \epsilon_{\cyc}, \Sigma)) - \sum_{w \in \scS^{\calK}, w \nmid \cond(\xi_0)} \ell_{\frP}((1-\xi_0^{\sfc}(\Frob_{w}))) \\
 &= \ell_{\frP}(X_{\Gr}(\bfone/\calK_{\infty}, \xi_0^{\sfc}, \Sigma^{\sfc})) - \sum_{w \in \scS^{\calK}, w \nmid \cond(\xi_0)} \ell_{\frP}((1-\xi_0^{\sfc}(\Frob_{w}))),    
\end{align*}
by \cite[Proposition 8.12]{hsieh2014eisenstein}, where we have made explicit the chosen CM type $\Sigma$. Moreover, when $p \nmid \abs{\Delta}$, \cite[(8.1)]{hsieh2014eisenstein} gives
\[
X_{\Gr}(\bfone/\calK_{\infty}, \xi_0^{\sfc}, \Sigma^{\sfc}) \simeq X(\xi_0^{\sfc}, \Sigma^{\sfc}) \otimes_{\ZZ_p[\xi_0]} \calO_{\frp}.
\]
On the analytic side, the $\scS$-imprimitive Katz $p$-adic $L$-function is defined as
\[
\calL_{\Katz}^{\scS}(\xi_0^{\sfc}, \Sigma^{\sfc}) := \calL_{\Katz}(\xi_0^{\sfc}, \Sigma^{\sfc}) \cdot \prod_{w \in \scS^{\calK}, w \nmid \cond(\xi_0)} (1-\xi_0^{\sfc}(\Frob_{w})) \in \Lambda.
\]
By comparing the interpolation formulas \eqref{eq:wanpadicl} and \eqref{eq:katzpadicl}, we see that
\begin{equation} \label{eq:comparewankatz}
(\calL_{\Wan}^{\scS}(\bfone, \calK, \chi_0, \Sigma^{\sfc})) = (\calL_{\Katz}^{\scS}(\xi_0^{\sfc}, \Sigma^{\sfc}))    
\end{equation}
as ideals of $\Lambda \otimes L$ (recall that we have set $\chi_0 := \psi^{-1}$). Combined with Theorem \ref{thm:hsiehprimitive} and Theorem \ref{thm:imcdeprived}, we therefore have
\begin{equation} \label{eq:onesidesfc}
\ord_{\frP}(F(\xi_0^{\sfc}, \Sigma^{\sfc})) \geq \ord_{\frP} (\calL_{\Katz}(\xi_0^{\sfc}, \Sigma^{\sfc}))
\end{equation}
for all height-one prime ideals $\frP$ of $\Lambda$ except for the following two types:
\begin{itemize}
    \item Height-one prime ideals $\frP$ of $\Lambda$ such that
    \[
\ord_{\frP}(\calL_{\Katz}^{\scS}(\xi_0^{\sfc}, \Sigma^{\sfc})) > 0 \quad \text{and} \quad \ord_{\frP}(\calL^{\scS}_{\RibetSerre}((\chi_0^{\calF})^{-1} \epsilon_{\calK/\calF}, \calF)) > 0,
\]
which are called \emph{exceptional primes} following \cite[Definition 8.13]{hsieh2014eisenstein}.\footnote{Note that in \cite[Definition 8.13]{hsieh2014eisenstein}, the primitive Katz $p$-adic $L$-function appears instead of the $\scS$-imprimitive one.}
    \item Height-one prime ideals $\frP$ such that $\frP \cap \II^{\ur} \neq (0)$.\footnote{We exclude this type of height-one prime ideals since we only prove the weak $p$-primitive property of $\bfD_{\Hs}$ in Theorem \ref{thm:hsiehprimitive}, and the comparison \eqref{eq:comparewankatz} between Wan's $p$-adic $L$-function and Katz's $p$-adic $L$-function is only valid up to inverting $p$.}
\end{itemize}
We have thus recovered \cite[Theorem 8.14]{hsieh2014eisenstein} under our weaker assumptions. The arguments dealing with exceptional primes in the proof of \cite[Corollary 8.15]{hsieh2014eisenstein} still hold. This requires the extra assumption that $p$ does not divide $\abs{\Delta}$ to guarantee that $\mu_p(\psi) = 0$ when $\epsilon(\psi) = 1$, so as to apply \cite[Theorem A]{hsieh2014mu} in Hsieh's argument. This proves \eqref{eq:onesidesfc} for exceptional primes. Finally, by Theorem \ref{thm:hsieh2014mu}(1), we see that $\ord_{\frP} (\calL_{\Katz}(\xi_0^{\sfc}, \Sigma^{\sfc}))=0$ for $\frP$ such that $\frP \cap \II^{\ur} \neq (0)$; hence \eqref{eq:onesidesfc} holds trivially. 

Theorem \ref{thm:myimccm} follows by replacing $\xi_0^{\sfc}$ with $\xi_0$ and $\Sigma^{\sfc}$ with $\Sigma$ (note that the character $\xi_0$ attached to the $p$-adic Eisenstein datum $\bfD_{\Hs}$ is exactly the branch character $\psi$). \qed

\subsection{Two cases of Iwasawa main conjectures for CM fields}
Having established the one-sided divisibility $(F(\psi, \Sigma)) \subseteq (\calL_{\Katz}(\psi, \Sigma))$ of the Iwasawa main conjecture for CM fields in Theorem \ref{thm:myimccm}, by applying the work of Hida or Rubin on the other side of the divisibility, we obtain the equality $(F(\psi, \Sigma)) = (\calL_{\Katz}(\psi, \Sigma))$. The argument is parallel to \cite[Section 8.6.1]{hsieh2014eisenstein} with an identical proof, but replacing the corresponding theorem there with Theorem \ref{thm:myimccm} of this article. We summarize the results as follows.

\begin{theorem}[Iwasawa main conjecture for CM fields] \label{thm:imchida}
Suppose that $p$ does not divide $6 D_{\calF} \abs{\Delta}$, and that:
\begin{enumerate}[label = \rm (Hi-\roman*)]
    \item $\psi$ is an anticyclotomic character and has order prime to $p$; \label{item:hidai}
    \item the local character $\psi_w$ is nontrivial over $\calK_{w}^{\times}$ for all $w \in \Sigma_{p}$; \label{item:hidaii}
    \item $\psi|_{\Gal_{\calK(\sqrt{p^{\ast}})}}$ is nontrivial, with $p^{\ast} = (-1)^{(p-1)/2} p$. \label{item:hidaiii}
\end{enumerate}
Then as ideals of $\Lambda$, we have $(F(\psi, \Sigma)) = (\calL_{\Katz}(\psi, \Sigma))$.
\end{theorem}
The proof follows the same argument as \cite[Theorem 8.17]{hsieh2014eisenstein}, but replaces the application of \cite{hida2006anticyclotomic} with \cite[Theorem]{hida2009quadratic}, where assumptions \ref{item:hidai}, \ref{item:hidaii}, and \ref{item:hidaiii} are required. Comments are made in \cite[p.~914]{hida2009quadratic} on the role of these assumptions and how to possibly remove them.

Another approach to obtaining an equality is to apply Rubin's two-variable main conjecture for imaginary quadratic fields $\calM$ proved in \cite{rubin1991main} and utilize a base-change trick. This yields the following result, with the proof following \cite[Theorem 8.18]{hsieh2014eisenstein} using Theorem \ref{thm:myimccm}. Recall that $\psi: \Gal(\calK^{\prime}/\calK) \rightarrow \barZZ_{p}^{\times}$ is the branch character, and $\Delta := \Gal(\calK^{\prime}/\calK)$.

\begin{theorem}[Iwasawa main conjecture for CM fields] \label{thm:imcrubin}
Suppose that $p$ does not divide $6 D_{\calF} \abs{\Delta}$, and that:
\begin{enumerate}[label = \rm (Ru-\roman*)]
    \item $\calK = \calF \calM$, where $\calM$ is an imaginary quadratic field in which $p$ splits; \label{item:rubini}
    \item $\Sigma$ is the $p$-ordinary CM type of $\calK$ obtained by extending the inclusion $\iota_{\infty}: \calM \hookrightarrow \CC$; \label{item:rubinii}
    \item $\calK^{\prime}$ is abelian over $\calM$ and $p \nmid [\calK^{\prime}:\calM]$. \label{item:rubiniii}
\end{enumerate}
Then as ideals of $\Lambda$, we have $(F(\psi, \Sigma)) = (\calL_{\Katz}(\psi, \Sigma))$.
\end{theorem}
Unlike Theorem \ref{thm:imchida}, which only deals with anticyclotomic branch characters $\psi$, Theorem \ref{thm:imcrubin} deals with general branch characters but requires extra assumptions on the CM field $\calK/\calF$ and the branch field $\calK^{\prime}/\calK$.

\subsection{Iwasawa main conjecture for CM elliptic curves}
Next, we apply Theorem \ref{thm:myimccm} (and Theorem \ref{thm:imcrubin}) to study CM elliptic curves, following \cite[Section 8.6.2]{hsieh2014eisenstein}. Let $E$ be an elliptic curve defined over $\calF$ with complex multiplication by the ring of integers of an imaginary quadratic field $\calM$. Suppose that $E$ has good ordinary reduction at all places of $\calF$ above $p$, which implies that $p$ splits in $\calM$. We write $p \calO_{\calM} = \frp \overline{\frp}$, where $\frp$ is the prime induced by $\iota_p: \calM \hookrightarrow \CC_p$. Let $\calK = \calF \calM$, and let $\Sigma$ be the CM type of $\calK$ induced by $\iota_{\infty}: \calM \hookrightarrow \CC$.

Let $\calK^{\prime} := \calK(E[\frp])$ and let $\psi := \psi_E$ be the character obtained from the $\ZZ_p$-representation of $\Delta := \Gal(\calK^{\prime}/\calK)$ on $E[\frp]$. Let $\pzK/\calK$ be the unique $\ZZ_p^{2}$-extension of $\calK$ contained in $\calK(E[p^{\infty}])$, and let $\Lambda_{E} := \ZZ_p \lrbracket{\Gal(\pzK/\calK)}$ be the two-variable Iwasawa algebra over $\ZZ_p$. Note that $\calK_{\infty}^{\cyc} \subset \pzK \subset \calK_{\infty}$. We define the \emph{$\frp^{\infty}$-Selmer group} $\Sel_{\frp^{\infty}}(E/\pzK)$ attached to $E$ over $\pzK$ by
\[
\Sel_{\frp^{\infty}}(E/\pzK) := \ker \left( \rmH^{1}(\pzK, E[\frp^{\infty}]) \rightarrow \prod_{w} \rmH^{1}(\pzK_w, E) \right).
\]
Then $\Sel_{\frp^{\infty}}(E/\pzK)$ is a cofinitely generated discrete $\Lambda_{E}$-module; we let $X_{\frp^{\infty}}(E/\pzK)$ denote its Pontryagin dual. By \cite[Chapter II, Theorem 18]{perrinriou}, there is a $\Lambda_{E}$-module isomorphism
\[
\Sel_{\frp^{\infty}}(E/\pzK) \xrightarrow{\sim} \Hom_{\ZZ_p}(X(\psi, \Sigma)(\pzK), E[\frp^{\infty}]),
\]
where $X(\psi, \Sigma)(\pzK)$ is the $\psi$-isotypic quotient of
\[
X_{\Sigma}(\pzK) := \Gal(\calM_{\Sigma}/\pzK) \otimes_{\ZZ_{p}[\Delta]\lrbracket{\Gamma_{\calK}}} \calR[\Delta]\lrbracket{\Gamma_{\calK}}.
\]

Let $\alpha_{E, \frp}$ be the $\ZZ_p$-valued character given by the composition
\[
\alpha_{E, \frp}: \Gal(\pzK/\calK) \hookrightarrow \Gal(\calK(E[p^{\infty}])/\calK) \rightarrow \Aut(E[\frp^{\infty}]) = \ZZ_p^{\times},
\]
and let $\pi_{E, \frp}: \Lambda \rightarrow \Lambda_{E}$ be the morphism induced by the character $\widetilde{\alpha}: \Gal(\calK_{\infty}/\calK) \rightarrow \Lambda_{E}^{\times}$ defined by $g \mapsto \alpha_{E, \frp}(g|_{\pzK})g|_{\pzK}$. Let $\frf$ be the prime-to-$p$ conductor of $\psi$; it is well known that $\calK^{\prime} \subseteq \calK(\frf\frp)$ if $p \nmid 6$. Setting $\Lambda_{E, \calO_{\frp}} := \Lambda_{E} \otimes_{\ZZ_p} \calO_{\frp}$, we define the \emph{two-variable $p$-adic $L$-function} attached to $(E, \frp)$ as
\[
\calL_{\frp}(E) := \pi_{E, \frp}(\calL_{\Katz}(\psi_{E}^{\prime}, \Sigma)) \in \Lambda_{E, \calO_{\frp}},
\]
where $\calL_{\Katz}(\psi_{E}^{\prime}, \Sigma)$ is the Katz $p$-adic $L$-function attached to the branch character $\psi_{E}^{\prime}: \Delta_{\frf} \rightarrow \Delta \rightarrow \mu_{p-1} \hookrightarrow \ZZ_p^{\times}$.

\begin{theorem}[Two-variable Iwasawa main conjecture for CM elliptic curves] \label{thm:imccmtwo}
Suppose that $p \nmid 6 D_{\calF}$ and that $E$ has good ordinary reduction at all places above $p$. Then we have the following inclusion of ideals in $\Lambda_{E, \calO_{\frp}}$:
\[
\Char_{\Lambda_{E}}(X_{\frp^{\infty}}(E/\pzK)) \subseteq (\calL_{p}(E)).
\]
If we further assume that 
\begin{enumerate}[label = \rm (Ru)]
    \item $p$ does not divide $[\calF:\QQ]$ and $\calK(E[\frp])/\calM$ is abelian, \label{item:ru}
\end{enumerate}
then the equality $\Char_{\Lambda_{E}}(X_{\frp^{\infty}}(E/\pzK)) = (\calL_{p}(E))$ holds as ideals in $\Lambda_{E, \calO_{\frp}}$.
\end{theorem}
\begin{proof}
It follows by applying Theorem \ref{thm:myimccm} to $\psi := \psi_{E}$ and the argument of \cite[Corollary 8.19]{hsieh2014eisenstein} to obtain the inclusion (with the observation that $p \nmid [\calK(E[\frp]):\calK]$). For the equality, we apply Theorem \ref{thm:imcrubin}: the condition \ref{item:rubini} and \ref{item:rubinii} are satisfied in our setup, and since $[\calK^{\prime}:\calM] = [\calF: \QQ][\calM^{\prime}: \calK \cap \calM]$ with $p \nmid [\calM^{\prime}: \calK \cap \calM]$, the condition \ref{item:rubiniii} is satisfied if $p \nmid [\calF:\QQ]$ and $\calK(E[\frp])/\calM$ is abelian.
\end{proof}

\begin{remark}
Although not explicitly stated in \cite{hsieh2014eisenstein}, the equality in Theorem \ref{thm:imccmtwo} can also be proved using the Iwasawa main conjecture for CM fields from \cite[Theorem 8.18]{hsieh2014eisenstein} under the additional assumption that $p$ does not divide the relative class number $h_{\calK}^{-}$ of $\calK /\calF$. We note that when $\calF = \QQ$, the conditions $p \nmid 6 h_{\calK}^{-} D_{\calF}$ and \ref{item:ru} hold trivially for $p \geq 5$. In this case, our result coincides with \cite[Corollary 8.19]{hsieh2014eisenstein}.
\end{remark}

We can also consider the cyclotomic main conjecture for $E$ over $\calF$. Let $\Lambda_{\calF} := \ZZ_{p} \lrbracket{\Gal(\calF_{\infty}/\calF)}$ and let $\pi_{+}: \Lambda_{E} \rightarrow \Lambda_{+} \simeq \Lambda_{\calF}$ be the restriction map. Let $\Lambda_{\calF, \calO_{\frp}} := \Lambda_{\calF} \otimes_{\ZZ_p} \calO_{\frp}$ and let $\calL_{p}(E/\calF_{\infty}) := \pi_{+}(\calL_{\frp}(E)) \in \Lambda_{\calF, \calO_{\frp}}$ be the \emph{cyclotomic $p$-adic $L$-function} for $E$. On the analytic side, we note that the complex CM period $\Omega_{\infty}$ appearing in the interpolation formula for $\calL_{\Katz}(\psi_{E}^{\prime}, \Sigma)$ and the Néron period $\Omega_{E}$ of the elliptic curve $E$ differ only by a $p$-adic unit. On the algebraic side, by Shapiro's lemma, $\Sel_{p^{\infty}}(E/\calF_{\infty}) \simeq \Sel_{\frp^{\infty}}(E/\calK_{\infty}^{\cyc})$ as $\Lambda_{\calF}$-modules for odd $p$. Specializing Theorem \ref{thm:imccmtwo} to the cyclotomic line and using the descent of Selmer groups for CM elliptic curves \cite[II: Lemme 9, Proposition 12]{perrinriou}, we obtain the following result.

\begin{theorem}[Cyclotomic Iwasawa main conjecture for CM elliptic curves] \label{thm:imccmcyc}
Suppose that $p \nmid 6 D_{\calF}$ and that $E$ has good ordinary reduction at all places above $p$. Then we have the inclusion of ideals in $\Lambda_{\calF, \calO_{\frp}}$:
\[
\Char_{\Lambda_{\calF}}(\Sel_{\frp^{\infty}}(E/\calK_{\infty}^{\cyc})) \subseteq (\calL_{p}(E/\calF_{\infty})).
\]
If, in addition, \ref{item:ru} holds, then we have the equality of ideals in $\Lambda_{\calF, \calO_{\frp}}$:
\[
\Char_{\Lambda_{\calF}}(\Sel_{\frp^{\infty}}(E/\calK_{\infty}^{\cyc})) = (\calL_{p}(E/\calF_{\infty})).
\]
\end{theorem}

\subsection{Application to the Birch--Swinnerton-Dyer conjecture for CM elliptic curves}

Cyclotomic Iwasawa main conjectures for elliptic curves, including those with complex multiplication, have significant applications to rank-zero Birch--Swinnerton-Dyer conjectures. We highlight two such applications below.

\begin{theorem}[Rank-zero $p$-converse theorem for CM elliptic curves] \label{thm:pconverse}
Suppose that $p \nmid 6 D_{\calF}$ and that $E$ has good ordinary reduction at all places above $p$. Then we have the implication:
\[
\corank_{\ZZ_p} \Sel_{p^{\infty}}(E/\calF) = 0 \Rightarrow L(E/\calF, 1) \neq 0.
\]
If \ref{item:ru} holds, then the converse implication holds as well.
\end{theorem}

\begin{theorem}[Rank-zero refined $p$-part Birch--Swinnerton-Dyer formula for CM elliptic curves] \label{thm:pbsd}
Suppose that $p \nmid 6 D_{\calF}$ and that \ref{item:ru} holds. Assume further that $E$ has good ordinary reduction at all places above $p$ and that $L(E/\calF, 1) \neq 0$. Then
\[
v_p \left( \dfrac{L(E/\calF, 1)}{\Omega_{E}} \right) = v_p \left( \dfrac{\abs{\Sha(E/\calF)[p^{\infty}]}}{\abs{E(\calF)_{\tor}}^2} \prod_{v \not\in \scV_{\calF}^{(\infty)}} \Tam_{v}(E/\calF) \right).
\]
Here $\Tam_{v}(E/\calF)$ is the Tamagawa number at $v$ for a finite place $v$ of $\calF$, $\Sha(E/\calF)[p^{\infty}]$ is the $p$-primary part of the Tate--Shafarevich group of $E$ over $\calF$, and $E(\calF)_{\tor}$ is the finite group of torsion points of $E$ over $\calF$.
\end{theorem}

These two theorems follow directly from Theorem \ref{thm:imccmcyc} and the control theorems for $E$ over $\calF$ proved in \cite[Theorem 1.2, Proposition 4.1]{greenberg1999iwasawa}. For the convenience of the reader, we provide proofs here; these are valid for any elliptic curve over $\calF$ (not necessarily with CM) for which the cyclotomic Iwasawa main conjecture holds.

\begin{proof}[Proof of Theorem \ref{thm:pconverse}]
Suppose $L(E/\calF, 1) = 0$; we shall prove that $\Sel_{p^{\infty}}(E/\calF)$ has positive $\ZZ_p$-rank. For simplicity, let $X$ denote the Pontryagin dual of $\Sel_{p^{\infty}}(E/\calF_{\infty})$ and let $\gamma_{+} - 1$ be a topological generator of $\Lambda_{\calF}$. By \cite[Proposition 4.15]{greenberg1999iwasawa}, $X$ has no non-zero pseudo-null $\Lambda_{\calF}$-submodule. Hence we have the following commutative diagram with exact rows, where $X_{\std}$ is the standard $\Lambda_{\calF}$-module associated to $X$:
\[
\begin{tikzcd}
	0 & X & {X_{\std}} & {\coker \theta} & 0 \\
	0 & X & {X_{\std}} & {\coker \theta} & 0
	\arrow[from=1-1, to=1-2]
	\arrow["\theta", from=1-2, to=1-3]
	\arrow["\gamma_+-1", from=1-2, to=2-2]
	\arrow[from=1-3, to=1-4]
	\arrow["\gamma_+-1", from=1-3, to=2-3]
	\arrow[from=1-4, to=1-5]
	\arrow["\gamma_+-1", from=1-4, to=2-4]
	\arrow[from=2-1, to=2-2]
	\arrow["\theta", from=2-2, to=2-3]
	\arrow[from=2-3, to=2-4]
	\arrow[from=2-4, to=2-5]
\end{tikzcd}
\]
By the snake lemma, let $\delta$ be the connecting homomorphism; we have:
\[
0 \rightarrow \ker(\theta)_{\gamma_+} \rightarrow \im(\delta) \rightarrow X/(\gamma_{+} - 1)X \xrightarrow{\theta} X_{\std}/(\gamma_{+} - 1)X_{\std} \rightarrow \coker(\theta)/(\gamma_{+} - 1)\coker(\theta) \rightarrow 0.
\]
Since $\im(\delta)$ is finite by Mazur's control theorem (\cite[Theorem 1.2]{greenberg1999iwasawa}) and $\coker(\theta)/(\gamma_{+} - 1)\coker(\theta)$ is finite, it follows that
\begin{equation} \label{eq:infiniteequiv}
\abs{X_{\std}/(\gamma_{+} - 1)X_{\std}} = \infty \Longleftrightarrow \abs{X/(\gamma_{+} - 1)X} = \infty.   
\end{equation}

By assumption, $L(E/\calF, 1) = 0$, which implies that $\gamma_{+} - 1$ divides $\calL_{p}(E/\calF_{\infty})$ in $\Lambda_{\calF, \calO_{\frp}}$. By the inclusion of ideals in Theorem \ref{thm:imccmcyc}, $\gamma_{+} - 1$ must divide $\Char_{\Lambda_{\calF}}(X)$. Consequently, $X_{\std}/(\gamma_{+} - 1)X_{\std}$ is infinite. From \eqref{eq:infiniteequiv}, we conclude that $X/(\gamma_{+} - 1)X$ is also infinite. Applying Mazur's control theorem (\cite[Theorem 1.2]{greenberg1999iwasawa}) again, this is equivalent to $\Sel_{p^{\infty}}(E/\calF)$ having positive $\ZZ_p$-rank.

If the reverse divisibility of ideals (i.e., $\Char_{\Lambda_{\calF}}(\Sel_{\frp^{\infty}}(E/\calK_{\infty}^{\cyc})) \supseteq (\calL_{p}(E/\calF_{\infty}))$) holds in Theorem \ref{thm:imccmcyc}, the converse implication is straightforward. On the algebraic side, the structure theorem for finitely generated torsion $\Lambda_{\calF}$-modules implies:
\[
\rank_{\ZZ_p} X/(\gamma_{+} - 1)X \leq \ord_{\gamma_{+} - 1} \Char_{\Lambda_{\calF}}(X).
\]
By Mazur's control theorem, the left-hand side equals $\corank_{\ZZ_p}(\Sel_{p^{\infty}}(E/\calF))$. Assuming $L(E/\calF, 1) \neq 0$ (and \ref{item:ru}), Theorem \ref{thm:imccmcyc} implies the right-hand side is zero. Therefore, $\corank_{\ZZ_p}(\Sel_{p^{\infty}}(E/\calF)) = 0$.
\end{proof}

\begin{remark}
When $\calF = \QQ$, the rank-zero $p$-converse theorem was proved for all CM elliptic curves $E$ in \cite{burungale2025rankzeropconversetheorem} using a rational version of the Iwasawa main conjecture for $E$ formulated via Kato's zeta morphism in \cite{kato2004padic}. Due to the lack of Kato's zeta morphism for elliptic curves over more general totally real fields $\calF$, it is difficult to generalize their method to general $\calF$.
\end{remark}

\begin{proof}[Proof of Theorem \ref{thm:pbsd}]
We adopt the notation from the proof of Theorem \ref{thm:pconverse}. Since $L(E/\calF, 1) \neq 0$, it follows from Theorem \ref{thm:pconverse} that $\Sel_{p^{\infty}}(E/\calF)$ is finite. This allows us to apply the refined control theorem proposed by Greenberg in \cite[Theorem 4.1]{greenberg1999iwasawa} \footnote{Let $f_{E}(T) \in \ZZ_p\lrbracket{T}$ be the characteristic power series of $X = \Sel_{p^{\infty}}(E/\calF)^{\vee}$. In fact, \cite[Theorem 4.1]{greenberg1999iwasawa} proves the slightly stronger conclusion that
\begin{equation} \label{eq:fn1}
 f_{E}(0) \sim \prod_{v \mid p} \abs{\tilE(\FF_{v})[p^{\infty}]}^{2} \cdot \dfrac{\abs{\Sel_{p^{\infty}}(E/\calF)}}{\abs{E(\calF)[p^{\infty}]}^2} \prod_{v \neq \infty} \Tam_{v}(E/\calF)^{(p)}    
\end{equation}
(here, for $a, b \in \QQ_p^{\times}$, we write $a \sim b$ if $a$ and $b$ have the same $p$-adic valuation; for a positive integer $n \in \ZZ$, we write $n = n^{(p)} q$ where $q > 0$ and $(p,q) = 1$). In general (using the snake lemma, cf. \cite[Lemma 4]{skinner2018note}), let $\Lambda := \ZZ_p\lrbracket{T}$, and let $X$ be a finitely generated torsion $\Lambda$-module with no nonzero pseudo-null submodules; then for a nonzero element $\lambda \in \Lambda$, we have
\[
\abs{X/\lambda X} = \abs{\Lambda / (\Char(X), \lambda)}.
\]
In particular, if $f_{X} \in \Char(X)$ is a generator (i.e., the characteristic power series of $X$), then 
\begin{equation}\label{eq:fn2}
\abs{X/(\gamma-1)X} = \abs{\ZZ_p/f_{X}(0)} = \abs{f_{X}(0)}_{p}^{-1}.    
\end{equation}
This equality should be understood such that if one side is infinite, the other side is also infinite. Combining \eqref{eq:fn1} with \eqref{eq:fn2}, and noting again that by \cite[Proposition 4.15]{greenberg1999iwasawa}, $X$ has no nonzero pseudo-null $\Lambda_{\calF}$-submodule, we obtain the result cited in the text.}:
\[
\abs{X/(\gamma_{+} - 1)X} = \left\vert \dfrac{\abs{\Sel_{p^{\infty}}(E/\calF)}}{\abs{E(\calF)_{\tor}}^2} \prod_{v \mid p} \abs{\tilE(\FF_{v})}^{2} \cdot \prod_{v \neq \infty} \Tam_{v}(E/\calF) \right\vert_{p}^{-1}.
\]
Here, for each prime $v$ of $\calF$ lying above $p$, we denote by $\tilE$ the reduction of $E$ modulo $v$, which is defined over the residue field $\FF_{v}$. This is well-defined since $E$ has good ordinary reduction at $v$. On the other hand, through the equality of the cyclotomic Iwasawa main conjecture proven in Theorem \ref{thm:imccmcyc} and the interpolation formula for $\calL_p(E/\calF_{\infty})$, we have
\[
\abs{X/(\gamma_{+} - 1)X} = \left\vert \dfrac{L(E/\calF, 1)}{\Omega_{E}} \cdot  \prod_{v \mid p} \abs{\tilE(\FF_{v})}^{2} \right\vert_{p}^{-1}.
\]
Note that the term in the parentheses on the right is the Euler factor of $L(E/\calF, s)$ at primes $v \mid p$ and $s = 1$. The theorem then follows by noting that $\Sel_{p^{\infty}}(E/\calF) = \Sha(E/\calF)[p^{\infty}]$: indeed, this is implied by the short exact sequence
\[
0 \rightarrow E(\calF) \otimes \QQ_p/\ZZ_p \rightarrow \Sel_{p^{\infty}}(E/\calF) \rightarrow \Sha(E/\calF)[p^{\infty}] \rightarrow 0
\]
and the fact that $\corank_{\ZZ_p} \Sel_{p^{\infty}}(E/\calF) = 0$. 
\end{proof}

\begin{remark}
Without appealing to \cite[Theorem 4.1]{greenberg1999iwasawa}, one could obtain a weaker inequality without assuming \ref{item:ru}:
\[
\length_{\ZZ_p}(\Sel_{p^{\infty}}(E/\calF)) \geq v_p \left( \dfrac{L(E/\calF, 1)}{\Omega_{E}} \right),
\]
similar to \cite[Corollary 1]{hsieh2014eisenstein} but with the assumption $p \nmid h_{\calK}^{-}$ removed. However, the $\ZZ_p$-length of $\Sel_{p^{\infty}}(E/\calF)$ could potentially be infinite if the reverse divisibility in Theorem \ref{thm:imccmcyc} is not established.
\end{remark}

\begin{remark}
Our Theorem \ref{thm:pbsd} overlaps with existing results. The rank-zero $p$-part Birch--Swinnerton-Dyer formula was proved for CM elliptic curves $E$ over $\calF$ in \cite[Corollary 2]{burungale2024bsd} assuming $\calK(E_{\tor})/\calM$ is abelian, without constraints on the prime $p$ or the reduction type. In the specific case where $p \nmid 6D_{\calF} [\calF:\QQ]$ and $E$ has good ordinary reduction above $p$, our result is slightly stronger as we only assume $\calK(E[\frp])/\calM$ is abelian, which is a smaller extension than $\calK(E_{\tor})/\calM$.
\end{remark}

We conclude by noting that for rank-one CM elliptic curves, current approaches to the $p$-converse theorem and the $p$-part BSD formula (e.g., \cite{burungale2020pconverse, li2026bsd}) involve different Iwasawa main conjectures and are largely limited to $\calF = \QQ$ with good ordinary reduction at $p$.